\documentclass[sn-nature,Numbered]{sn-jnl}% Math and Physical Sciences Reference Style
\usepackage{graphicx}%
\usepackage{multirow}%
\usepackage{amsmath,amssymb,amsfonts}%
\usepackage{amsthm}%
\usepackage{mathrsfs}%
\usepackage[title]{appendix}%
\usepackage{xcolor}%
\usepackage{colortbl}
\usepackage{textcomp}%
\usepackage{manyfoot}%
\usepackage{booktabs}%
\usepackage{algorithm}%
\usepackage{algorithmicx}%
\usepackage{algpseudocode}%
\usepackage{listings}%
\usepackage{easybmat}
\usepackage{caption}
\definecolor{bluencs}{rgb}{0.0, 0.53, 0.74}
\definecolor{ceruleanblue}{rgb}{0.16, 0.32, 0.75}
\hypersetup{
	colorlinks=true,
	linkcolor=bluencs,
	filecolor=black,    
	citecolor=bluencs,
	urlcolor=bluencs
}
\usepackage{fourier}
\renewcommand\rmdefault{ptm}
\usepackage[bf]{subfigure}

\newcommand{\mbf}[1]{\mathbf{#1}}

\newcommand{\mbfM}{\mbf{M}}
\newcommand{\mbfC}{\mbf{C}}
\newcommand{\mbfK}{\mbf{K}}
\newcommand{\mbfF}{\mbf{F}}
\newcommand{\mbfu}{\mbf{U}}
\newcommand{\mbfv}{\dot{\mbf{U}}}
\newcommand{\mbfa}{\ddot{\mbf{U}}}
\newcommand{\dt}{\varDelta\! t}

\newcommand{\rhoinf}{\rho_\infty}

\renewcommand{\Omega}{\varOmega}
\usepackage{colortbl}
\definecolor{mygray}{gray}{.9}
\definecolor{mygraye}{gray}{.85}
\definecolor{dcolor}{rgb}{0.5,0.5,0.8}
\allowdisplaybreaks[4]
\usepackage{algorithm,algpseudocode,float}
\makeatletter
\newenvironment{breakablealgorithm}
{% \begin{breakablealgorithm}
		\begin{center}
			\refstepcounter{algorithm}% New algorithm
			\hrule height.8pt depth0pt \kern2pt% \@fs@pre for \@fs@ruled
			\renewcommand{\caption}[2][\relax]{% Make a new \caption
				{\raggedright\textbf{\ALG@name~\thealgorithm} ##2\par}%
				\ifx\relax##1\relax % #1 is \relax
				\addcontentsline{loa}{algorithm}{\protect\numberline{\thealgorithm}##2}%
				\else % #1 is not \relax
				\addcontentsline{loa}{algorithm}{\protect\numberline{\thealgorithm}##1}%
				\fi
				\kern2pt\hrule\kern2pt
			}
		}{\kern2pt\hrule\relax% \@fs@post for \@fs@ruled
		\end{center}
	}
	\makeatother
\theoremstyle{thmstyleone}%
\newtheorem{theorem}{Theorem}%  meant for continuous numbers
\newtheorem{proposition}[theorem]{Proposition}% 

\theoremstyle{thmstyletwo}%
\newtheorem{remark}{Remark}%

\theoremstyle{thmstylethree}%

\begin{document}

\title[High-order accurate multi-sub-step implicit integration algorithms with dissipation control]{High-order accurate multi-sub-step implicit integration algorithms with dissipation control for second-order hyperbolic problems}


\author[1]{\fnm{Jinze} \sur{Li}}\email{pinkie.ljz@gmail.com}

%\author[1]{\fnm{Naigang} \sur{Cui}}\email{cui\_naigang@163.com}

\author[2]{\fnm{Hua} \sur{Li}}\email{lihua@ntu.edu.sg}

%\author[1]{\fnm{Yiwei} \sur{Lian}}\email{lianyiwei\_hit@qq.com}

\author*[1]{\fnm{Kaiping} \sur{Yu}}\email{kaipingyu1968@gmail.com}

\author[1]{\fnm{Rui} \sur{Zhao}}\email{ruizhao@hit.edu.cn}
%\author*[1,2]{\fnm{Third} \sur{Author}}\email{iiiauthor@gmail.com}

\affil[1]{\orgdiv{School of Astronautics}, \orgname{Harbin Institute of Technology}, \orgaddress{\street{No.~92 West Dazhi Street}, \city{Harbin}, \postcode{150001}, \country{China}}}

\affil[2]{\orgdiv{School of Mechanical and Aerospace Engineering}, \orgname{Nanyang Technological University}, \orgaddress{\street{50 Nanyang Avenue}, \postcode{639798}, \country{Singapore}}}
%
%\affil[3]{\orgdiv{Department}, \orgname{Organization}, \orgaddress{\street{Street}, \city{City}, \postcode{610101}, \state{State}, \country{Country}}}


\abstract{This paper proposes an implicit family of sub-step integration algorithms grounded in the explicit singly diagonally implicit Runge-Kutta (ESDIRK) method. The proposed methods achieve third-order consistency per sub-step and thus the trapezoidal rule is always employed in the first sub-step. This paper demonstrates for the first time that the proposed $ s $-sub-step implicit method with $ s\le6 $ can reach $ s $th-order accuracy when achieving dissipation control and unconditional stability simultaneously. Hence, this paper develops, analyzes, and compares four cost-optimal high-order implicit algorithms within the present $ s $-sub-step method using three, four, five, and six sub-steps. Each high-order implicit algorithm shares identical effective stiffness matrices to achieve optimal spectral properties. Unlike the published algorithms, the proposed high-order methods do not suffer from the order reduction for solving forced vibrations. Moreover, the novel methods overcome the defect that the authors' previous algorithms require an additional solution to obtain accurate accelerations. Linear and nonlinear examples are solved to confirm the numerical performance and superiority of four novel high-order algorithms.}




\maketitle

\section{Introduction}
The numerical simulation of dynamic structures is one of the central problems in computational dynamics. In the past decades, a great number of numerical techniques have been proposed to solve various dynamical problems. When considering linear elastic structures and after using spatial discretizations \cite{hughesFiniteElement2000}, the following second-order differential equations of motion can be obtained as
\begin{equation}\label{eq:mck}
	\mbf{M}\ddot{\mbf{U}}(t)+\mbf{C}\dot{\mbf{U}}(t)+\mbf{K}\mbf{U}(t)=\mbf{F}(t)
\end{equation}
with the appropriate initial conditions $ \mbf{U}_0=\mbf{U}(t_0) $ and $ \dot{\mbf{U}}_0=\dot{\mbf{U}}(t_0) $. In Eq.~(\ref{eq:mck}), $ \mbf{U}(t)$ collects unknown nodal displacements and a dot denotes differentiation with respect to time $ t $; $ \mbf{M},~\mbf{C} $, and $ \mbf{K} $ are the global mass, damping, and stiffness matrices, respectively, and $ \mbf{F}(t) $ stands for the load vector, as a known function of time $ t $. Among all numerical techniques to solve the system (\ref{eq:mck}), the preferred one is of the single-step type consisting of updating displacement, velocity, and acceleration vectors at current time $ t_{n} $ to the next instant $ t_{n+1}=t_n+\dt $, where $ \dt $ denotes the time increment. This type is often called direct time integration algorithms \cite{hughesFiniteElement2000}, or step-by-step time marching schemes. In general, according to the computational efficiency, time integration algorithms can be further divided into two parts: explicit and implicit methods, and they possess own advantages and disadvantages. 

In explicit methods, numerical solutions at discrete instants depend only on responses in preceding time step(s) \cite{rezaiee-pajandMoreAccurate2015}, and they can significantly save the computational cost but achieve only conditional stability, so the used integration steps are strongly limited by their stability limits. The single-step explicit scheme \cite{liIdenticalSecondorder2021} and the composite two-sub-step explicit algorithms \cite{nohExplicitTime2013,liTwoThirdorder2022} are recommended in this paper. These explicit methods all achieve, at least, second-order of accuracy and flexible dissipation control at the bifurcation point, and they are superior to the earlier explicit algorithms. On the other hand, implicit methods compute numerical solutions by (iteratively) solving an equation involving both known and unknown states of the structure \cite{rezaiee-pajandNewExplicit2016} and thus need the direct/iterative solver to solve the resulting linearized equations. Importantly, they can achieve unconditional stability at the expense of computational costs, so the used integration steps are often chosen based on accuracy instead of stability. The main attention of this paper is paid to developing implicit algorithms.

Many of implicit integration schemes have been developed by using several desgin ideas.
%, such as the time finite elements \cite{fungExtrapolatedGalerkin1996,fungUnconditionallyStable1996,liStructuralDynamic1996}, collocation schemes \cite{hilberCollocationDissipation1978,golleyTimeSteppingProcedure1996,fungUnconditionallyStable2001}, finite difference schemes \cite{parkImprovedStiffly1975,houboltRecurrenceMatrix1950,yuNewFamily2008}, and weighted residual methods \cite{golleyWeightedResidual1998,rezaiee-pajandEfficientWeighted2021}. 
Some well-known implicit algorithms are the Newmark \cite{newmarkMethodComputation1959}, HHT-$\alpha$ \cite{hilberImprovedNumerical1977}, WBZ-$\alpha$ \cite{woodAlphaModification1980}, and TPO/G-$\alpha$ \cite{shaoThreeParameters1988,chungTimeIntegration1993} algorithms. These mentioned algorithms are single-step, so they often achieve second-order accuracy, unconditional stability, the solution of linear systems once within each time step, and self-starting features. It should be emphasized that these algorithms all suffer from unexpected overshoots in displacement and/or velocity when imposing numerical dissipation in the high-frequency range, and they are reduced to either the trapezoidal rule or the midpoint rule in the non-dissipative case. Among all singe-step implicit schemes without subsidiary variables, only the trapezoidal rule and the midpoint rule are highly recommended since they possess desired numerical properties.% except for the dissipation control. %In the dissipative case, the GS4-2 algorithms \cite{shimadaNovelDesign2015,zhouDesignAnalysis2004,maxamReevaluationOvershooting2022} are also acceptable, although some algorithm members with the high-frequency dissipation suffer from overshoots in either displacement or velocity and algorithm members without overshoots impose numerical dissipation only in the middle-frequency range.

\subsection{Second-order sub-step schemes}
Recently, some novel implicit integration methods, namely the so-called composite sub-step algorithms, have attracted the attention of some researchers. A composite two-sub-step algorithm was analyzed by Bathe \cite{batheCompositeImplicit2005} to successfully solve the strongly nonlinear problems that the non-dissipative trapezoidal rule fails to solve. The Bathe algorithm uses the composite two-sub-step technique. Firstly, the current integration interval $ t\in[t_n,~t_{n+1}] $ is split into two sub-intervals $ [t_n,~t_n+\gamma\dt] $ and $ [t_n+\gamma\dt,~t_{n+1}] $, where $ \gamma $ denotes the splitting ratio of sub-step size. Then, the non-dissipative trapezoidal rule and the three-point backward difference formula are used in the first and second sub-steps, respectively. Finally, the resulting scheme achieves L-stability, second-order accuracy, non-overshooting, and self-starting features. The unique free-parameter $ \gamma $ adjusts numerical dissipation imposed in the low-frequency range. It has been shown in \cite{LiFurtherAssessment2021} that the Bathe algorithm \cite{batheCompositeImplicit2005} corresponds to the non-dissipative trapezoidal rule in either $ \gamma=0 $ or $ 1 $, while the first-order backward Euler scheme is covered in $ \gamma=2 $. These facts also explain the reason why $\gamma$ can adjust numerical dissipation in the low-frequency range. When the parameter $\gamma$ gets close to either $ 0 $ or $ 1 $, the numerical behavior of the Bathe algorithm \cite{batheCompositeImplicit2005} naturally tends to that of the non-dissipative trapezoidal rule, so controlling dissipation in the low-frequency range. More importantly, the Bathe method \cite{batheCompositeImplicit2005,batheConservingEnergy2007} corresponds essentially to the two-stage diagonally implicit Runge-Kutta (RK) method proposed by Bank et al.~\cite{bankTransientSimulation1985}. Hence, the composite sub-step methods, featuring consistent displacement and velocity updating formulas, essentially form part of the RK family. Over the past decades, some composite sub-step implicit methods developed in computational mechanics have predominantly belonged to the implicit RK schemes and are briefly viewed as follows. 

The composite multi-sub-step implicit algorithms using the non-dissipative trapezoidal rule and more general multi-point backward difference scheme were constructed and analyzed in the literature \cite{liThreeOpitmal2023,nohBatheTime2019,liNovelFamily2019,rezaiee-pajandMixedMultistep2010}. Besides, other integration algorithms without adopting the trapezoidal rule in the first sub-step were successfully developed in \cite{liAlternativeBathe2019,liTrulySelfstarting2020,liSimpleTruly2020} and they can predict more accurate solutions for solving some dynamical problems. Two and three sub-steps have been widely used to construct composite sub-step algorithms since these sub-step schemes can be readily derived and analyzed. For example, another two Bathe algorithms, named as $ \rhoinf $-Bathe \cite{nohBatheTime2019} and $ \beta_1/\beta_2 $-Bathe \cite{malakiyehBatheTime2019}, are composite two-sub-step schemes, and they actually share the completely same algorithm structure \cite{LiFurtherAssessment2021}. The theoretical analysis \cite{LiFurtherAssessment2021} has shown that the second-order $ \beta_1/\beta_2 $-Bathe algorithm is algebraically identical to the $ \rhoinf $-Bathe algorithm, and the first-order $ \beta_1/\beta_2 $-Bathe algorithm is not competitive with other higher-order methods. The researches have clearly shown that the composite multi-sub-step implicit algorithms attain optimal spectral properties, namely the maximum numerical dissipation but the minimum relative period errors, when achieving identical effective stiffness matrices within each sub-step. Hence, the optimal composite two-sub-step implicit algorithm achieving dissipation control, second-order accuracy, and no overshoots has been proposed in \cite{liNovelFamily2020}. Other composite two-sub-step schemes, such as $ \rhoinf $-Bathe \cite{nohBatheTime2019}, can cover it when achieving optimal spectral properties. Similarly, an optimal three-sub-step implicit method has been also considered in \cite{liNovelFamily2019}, where identical effective stiffness matrices and controllable algorithmic dissipation in the high-frequency range are achieved. %It should be pointed out that another composite three-sub-step implicit algorithm \cite{jiOptimizedThreesubstep2020} is algebraically identical to the earlier published one \cite{liNovelFamily2019}. 
All composite multi-sub-step implicit algorithms mentioned above are only second-order accurate.

\subsection{High-order sub-step schemes}
High-order accurate integration algorithms often require much more computational cost than the common second-order ones. So far, there have been a great number of studies \cite{fungExtrapolatedGalerkin1996,tarnowHowRender1994,kimEffectiveHigherOrder2017,fungUnconditionallyStable1997,fanComprehensiveUnified1997a,fungComplextimestepNewmark1998,mancusoEfficientIntegration2003,krenkConservativeFourthorder2015,zhangOptimizationNsubstep2020,rezaiee-pajandImplicitHigherorder2008,rezaiee-pajandHighlyAccurate2018,liDirectlySelfstarting2022,kimNewFamily2017,rezaiee-pajandNumericalTime2008} reporting high-order accurate integration algorithms. A preferred way to develop high-order schemes is to operate the first-order differential systems by converting the original second-order equations of motion. The studies \cite{fungExtrapolatedGalerkin1996,fanComprehensiveUnified1997a,krenkConservativeFourthorder2015,mancusoEfficientIntegration2003} used the time finite elements and the weighted residual approach to address the resulting first-order systems. These high-order schemes suffer mainly from two drawbacks. One is that the resulting integration algorithms cannot output acceleration responses, although they are directly self-starting. Moreover, for solving structures subjected to the external loads, the lack of external load analysis could cause the order reduction in displacement and velocity. For example, the third/fourth-order accurate schemes \cite{fungExtrapolatedGalerkin1996} from the extrapolated Galerkin time finite elements are only of second-order for solving forced vibrations. The other is how to compute external loads effectively and accurately. Because the load terms are often expressed as the integral form in the current time interval, how to compute these loads in the integral form is crucial to maintain the expected order of accuracy. It is a known fact \cite{kimNewFamily2017,kimEffectiveHigherOrder2017,wangOverviewHighOrder2021} that the $ (p+1) $th-order approximation for the second-order dynamical problems imposes $ p $ unknown variables, and thus the order of accuracy can be optimized up to $ (2p-1) $ and $ 2p $. In general, the developed high-order algorithms can be designed to have controllable dissipation capability. The high-order methods are generally $ (2p-1) $th- and $ 2p $th-order accurate in the dissipative and non-dissipative cases, respectively. It should be emphasized that the high-order methods derived from the transformed first-order systems have to solve the equation system of dimension $ pd\times pd $ where $ p $ and $ d $ denote the number of the involved variables and the degree-of-freedom of the spatially discretized model, respectively. 
%For example, Kim and Reddy \cite{kimNewFamily2017,kimEffectiveHigherOrder2017} have recently proposed two new families of high-order integration algorithms based on a modified time-weighted residual method. 
In the work \cite{kimNewFamily2017}, the $ p $th-order Lagrange polynomial and Gauss-Lobatto integration points are used to develop $ (2p-1) $th- and $ 2p $th-order accurate algorithms, and a family \cite{kimEffectiveHigherOrder2017} of high-order algorithms considers $ p $th-order Hermite interpolation polynomials and higher-order time derivatives of the equilibrium equations. 
%In his thesis \cite{kimImprovedTime2016a}, another high-order implicit method that adopts the $ p $th-order Lagrange polynomial and equally spaced integration points is also described. 
Note that these families \cite{kimNewFamily2017,kimEffectiveHigherOrder2017} of high-order algorithms belonging to the fully implicit RK family must solve the equation systems with higher dimensions per time step, so increasing computational costs. It should be pointed out that the effective stiffness matrices of these two methods \cite{kimNewFamily2017,kimEffectiveHigherOrder2017} are neither symmetric nor sparse. Some existing high-order implicit methods \cite{fungWeightingParameters1999a,idesmanNewHighorder2007,mancusoEfficientIntegration2003,fanComprehensiveUnified1997a,argyrisDynamicResponse1973} suffer from these issues. Hence, some researchers \cite{liStructuralDynamic1996,mancusoEfficientIntegration2003,kimNewFamily2017} proposed various iterative solvers to solve the resulting higher dimensional equation systems. A recent overview of these high-order implicit algorithms can refer to the literature \cite{wangOverviewHighOrder2021}.  

On the other hand, high-order integration algorithms \cite{tarnowHowRender1994,zhangOptimizationNsubstep2020,liDirectlySelfstarting2022} directly operating the second-order differential systems show some advantages over other high-order ones, but there are some shortcomings. One of the main advantages of these high-order implicit methods is avoiding solving higher dimensional equation systems per time step, so significantly reducing computational costs. For instance, the Tarnow and Simo technique \cite{tarnowHowRender1994} is often applied to the single-step non-dissipative schemes, such as the trapezoidal rule and the midpoint rule. When some single-step dissipative schemes are used in the Tarnow and Simo technique, the resulting methods suffer from the order reduction. Hence, the high-order schemes from the Tarnow and Simo technique are generally non-dissipative. If the Tarnow and Simo technique is applied to the midpoint rule, the resulting scheme also covers a directly self-starting, fourth-order accurate, and non-dissipative method \cite{defrutosEasilyImplementable1992}. A single-step high-order implicit method based on Pad$ \acute{\text{e}} $ expansion was developed in \cite{songHighorderImplicit2022}. The developed method \cite{songHighorderImplicit2022} is always non-dissipative due to the use of diagonal Pad$ \acute{\text{e}} $ approximation, and it involves the complicated calculations of external loads and complex operations. Soares developed a directly self-starting fourth-order integration scheme \cite{soaresStraightforwardHighorder2020} for solving undamped models. It is reduced to be third- and second-order accurate for numerically and physically damped models, respectively. However, it is only conditionally stable, so its effective applications may be limited. Zhang et al. \cite{zhangOptimizationNsubstep2020} adopted the composite sub-step technique to develop an implicit family of high-order algorithms with controllable high-frequency dissipation and acceptable computational cost, but these algorithms also present lower-order of accuracy for solving forced vibrations \cite{liDirectlySelfstarting2022}. The second-order accurate $ \rhoinf $-Bathe algorithm \cite{nohBatheTime2019} is reanalyzed to achieve third-order accuracy \cite{kwonSelectingLoad2021} via considering optimal load selections. The resulting third-order scheme \cite{kwonSelectingLoad2021,choiTimeSplitting2022} cannot achieve a full range of numerical dissipation and the load information in the previous time step has to be used in the current step. The $ \rhoinf $-Bathe method has been revisited \cite{choiTimeSplitting2022} to analyze the influence of the time splitting ratio $ \gamma $ on accuracy. When a complex-valued $ \gamma $ is adopted, the $ \rhoinf $-Bathe method can achieve a full range of dissipation control and the fourth-order accuracy is obtained in the non-dissipative case. 

\subsection{Objectives and outlines}
Adopting the composite sub-step technique, the authors have already developed an implicit family \cite{liDirectlySelfstarting2022} of high-order algorithms which correspond essentially to the singly diagonally implicit Runge-Kutta methods (SDIRK). The developed algorithms are directly self-starting, avoiding both calling any starting procedures and computing the initial acceleration vector. In general, the directly self-starting $ s $-sub-step implicit methods \cite{liDirectlySelfstarting2022} can achieve $ s $th-order accuracy, unconditional stability, no overshoots, and controllable high-frequency dissipation. As stressed by Li et al.~\cite{liDirectlySelfstarting2022}, the directly self-starting high-order algorithms require using the second-order equations of motion once additionally to guarantee identical $ s $th-order accuracy, which is computationally expensive. Otherwise, the methods \cite{liDirectlySelfstarting2022} cannot output the same acceleration accuracy as those in displacement and velocity. In addition to the considerable computational effort required for the output of acceleration responses, the directly self-starting high-order methods also show poorer robustness for solving complicated dynamical problems due to losing the acceleration $\mbfa_n$ in all sub-steps. Almost all directly self-starting algorithms are faced with these computational issues, and further studies on analyzing the pros and cons of directly self-starting algorithms are needed in the near future. 

To overcome these issues, the authors will construct and develop a novel implicit family of high-order algorithms, corresponding to the explicit singly diagonally implicit Runge-Kutta methods (ESDIRK), which achieves the following desirable numerical characteristics. %The developed algorithms are expected to achieve identical high-order accuracy without additional solutions and with better robustness for solving complicated dynamical problems than the published schemes.
%Notice that the reviewed high-order implicit methods \cite{tarnowHowRender1994,defrutosEasilyImplementable1992,songHighorderImplicit2022,soaresStraightforwardHighorder2020,zhangOptimizationNsubstep2020,kwonSelectingLoad2021,liDirectlySelfstarting2022} only require solving the equation system of the same dimension as the original equilibrium equations per time step, so the resulting effective stiffness matrix inherits all desired properties of the global mass, damping, and stiffness matrices, such as sparsity, symmetry, and positive definiteness. Therefore, this paper focuses mainly on developing this type of high-order implicit methods. 
%To this end, the manuscript will design an implicit family of composite multi-sub-step algorithms to achieve the following characteristics. 
\begin{itemize}%[(a)]
	\item The self-starting property, avoiding the imposition of an additional starting procedure, has been emphasized in the work of Hilber and Hughes \cite{hilberCollocationDissipation1978}. They underscored that a non-self-starting algorithm often necessitates more in-depth analysis and consideration and the incorporation of an additional starting procedure may result in heightened coding and computational costs. Thus, the algorithm should be crafted with the intention of being (directly) self-starting. However, as mentioned earlier, the directly self-starting algorithms come with certain computational disadvantages due to the loss of $\mbfa_n$ in all sub-steps. Consequently, the novel algorithms are intentionally designed to be self-starting but not directly self-starting.
%	Self-starting property, avoiding imposing an additional starting procedure. Hilber and Hughes \cite{hilberCollocationDissipation1978} stressed that a non-self-starting algorithm often requires more analysis and consideration, and the additional starting procedure may engender more coding and computational costs. The algorithm should be designed to be (directly) self-starting. As mentioned above, the directly self-starting algorithms have some computational disadvantages due to losing $\mbfa_n$ in all sub-steps, so the novel algorithms are designed to be self-starting, not directly self-starting.
	
	\item Identical effective stiffness matrices across all sub-steps, incorporating optimal spectral properties, play a crucial role. The identity of effective stiffness matrices not only decreases computational expenses in solving linear elastic problems within multi-sub-step implicit algorithms but also integrates optimal spectral properties. When traditional implicit algorithms are considered as composite single-sub-step schemes, this property is automatically realized in all traditional algorithms. Remarkably, this property has been conceptually extended to the advancement of multi-sub-step explicit methods \cite{liSecondorderSsubstep2023,liSuiteSecondorder2023}.
%	Identical effective stiffness matrices within all sub-steps, embedding optimal spectral properties. The identical effective stiffness matrices can make multi-sub-step implicit algorithms not only reduce computational costs for solving linear elastic problems but also embed optimal spectral properties. If the traditional implicit algorithms are viewed as composite single-sub-step schemes, all traditional algorithms automatically achieve this property. Interestingly, this property has been extended conceptually to the development of multi-sub-step explicit methods \cite{liSecondorderSsubstep2023,liSuiteSecondorder2023}.
	
	\item The ability to control numerical dissipation in the high-frequency range, aiding in the damping out of spurious modes and stabilizing iterative procedures, is recognized as a desirable property in time integration algorithms. To accommodate complicated solutions, the novel methods should be meticulously designed to achieve comprehensive dissipation control, ranging from the non-dissipative scenario to the asymptotically annihilating case.
%	Controllable numerical dissipation in the high-frequency range, suppressing spurious high-frequency components. The capability of dissipation control is considered one of the desirable properties of time integration algorithms, such as damping out spurious modes and stabilizing iterative procedures. To accommodate complicated solutions, the novel methods should be well-designed to achieve a full range of dissipation control from the non-dissipative case to the asymptotically annihilating case.
	
	\item Higher-order accuracy in displacement, velocity, and acceleration is crucial for solving general structures. In this paper, novel methods are formulated with real-valued parameters to attain identical high-order accuracy. The sustained high-order accuracy achieved by these novel methods is demonstrated in the resolution of a comprehensive benchmark problem.
%	Higher-order accuracy in displacement, velocity, and acceleration for solving general structures. In this paper, the novel methods are designed with real-valued parameters to achieve identical high-order accuracy. The high-order accuracy achieved in the novel methods persists in solving the general benchmark problem.
	
	\item Unconditional stability, rendering the integration step constrained by accuracy rather than stability. Implicit algorithms must unquestionably be devised to possess unconditional stability, as lacking this feature would diminish their competitive edge in practical applications. Note that, achieving identical effective stiffness matrices, controllable numerical dissipation, and high-order accuracy, the proposed methods lack additional parameter design flexibility for further enhancement of BN-stability \cite{jiUnconditionallyStable2021}.
	
%	Unconditional stability, making the used integration step limited by accuracy instead of stability. Undoubtedly, implicit algorithms should be designed to be unconditionally stable, or they have no competitive advantages in practical applications.
	
	\item Achieving a favorable equilibrium between computational cost and high-order accuracy is essential. Similar to the existing high-order implicit methods \cite{tarnowHowRender1994,defrutosEasilyImplementable1992,songHighorderImplicit2022,soaresStraightforwardHighorder2020,zhangOptimizationNsubstep2020,kwonSelectingLoad2021,liDirectlySelfstarting2022}, the novel methods necessitate solving linearized equations with the same dimension as the second-order equations of motion per (sub-)step. Consequently, the resulting effective stiffness matrix inherits all desired properties of the global mass, damping, and stiffness matrices, including sparsity, symmetry, and positive definiteness.
%	The good balance between computational cost and high-order accuracy. Like the reviewed high-order implicit methods \cite{tarnowHowRender1994,defrutosEasilyImplementable1992,songHighorderImplicit2022,soaresStraightforwardHighorder2020,zhangOptimizationNsubstep2020,kwonSelectingLoad2021,liDirectlySelfstarting2022}, the novel methods require solving the linearized equations with the same dimension as the second-order equations of motion per (sub-)step, so the resulting effective stiffness matrix inherits all desired properties of the global mass, damping, and stiffness matrices, such as the sparsity, symmetry, and positive definiteness.
\end{itemize}

The remaining sections of this work are structured as follows. Section \ref{sec:development} introduces the construction of a novel composite $ s $-sub-step implicit scheme, where the conditions for achieving $ p $th-order accuracy and controllable numerical dissipation are derived. In accordance with these conditions, Section \ref{sec:alg} details the development of novel high-order implicit algorithms. Section \ref{sec:sp} provides spectral analysis and comparisons with existing high-order algorithms, confirming the well-established unconditional stability and dissipation control of the proposed high-order algorithms. Section \ref{sec:example} presents numerical examples to validate the numerical performance and superiority of the four novel high-order schemes. Lastly, Section \ref{sec:conclusion} summarizes essential conclusions drawn from this study.

%The remainder of this work is organized as follows. A novel composite $ s $-sub-step implicit scheme is constructed in Section \ref{sec:development}, where the conditions for achieving $ p $th-order accuracy and controllable numerical dissipation are derived. According to these conditions, four novel high-order implicit algorithms are developed in Section \ref{sec:alg}. Spectral analysis and comparisons with the existing high-order algorithms are given in Section \ref{sec:sp}, where the unconditional stability and dissipation control are confirmed well for the proposed high-order algorithms. Numerical examples are solved in Section \ref{sec:example} to validate the numerical performance and superiority of four novel high-order schemes. Finally, some important conclusions are summarized in Section \ref{sec:conclusion}.

\section{Development and analysis}\label{sec:development}
This section will firstly present a novel $ s $-sub-step implicit algorithm based on the ESDIRK method to solve the system (\ref{eq:mck}). Then, the conditions for achieving the designed order of accuracy are derived by analyzing a damped single-degree-of-freedom (SDOF) system subjected to the external load. The section ends up with the controllable dissipation analysis in the high-frequency range.

\subsection{The novel s-sub-step implicit algorithm}
For the sake of simplify and clarity, the constructed $ s $-sub-step implicit algorithm is described to solve the linear second-order equations of motion (\ref{eq:mck}). Considering the current integration internal $ t\in\left[t_n,~t_{n}+\dt\right] $, the composite sub-step technique firstly divides the whole time increment into $ s $ sub-intervals $ \cup_{i=1}^s[t_n+\gamma_{i-1}\dt,~t_n+\gamma_{i}\dt]$, where $ \gamma_i $ denotes the splitting ratio of the $ i $th sub-step. In this paper, $ \gamma_0=0\ \text{and}\ \gamma_s=1 $ are always used to provide acceleration responses. Then, an integration scheme in the $i$th sub-step $ t\in[t_n+\gamma_{i-1}\dt,~t_n+\gamma_{i}\dt]~(i=1,~\dots,~s) $ is constructed as
\begin{subequations}\label{eq:nsubstep}
	\begin{align}
		\mbf{M}\mbf{\ddot{U}}_{n+\gamma_i}&+\mbf{C}\mbf{\dot{U}}_{n+\gamma_i}+\mbf{K}\mbf{U}_{n+\gamma_i}  =\mbf{F}(t_n+\gamma_i\dt)                                                    \\
		\mbf{U}_{n+\gamma_i}                                                                             & =\mbf{U}_n+\dt\sum_{j=0}^{i}\alpha_{ij}\mbf{\dot{U}}_{n+\gamma_j}        \\
		\mbf{\dot{U}}_{n+\gamma_i}                                                                       & =\mbf{\dot{U}}_n+\dt\sum_{j=0}^{i}\alpha_{ij}\mbf{\ddot{U}}_{n+\gamma_j}.
	\end{align}
\end{subequations}
Notice that numerical solutions ($ \mbf{U}_{n+1}, ~\mbf{\dot{U}}_{n+1}$, and  $\mbf{\ddot{U}}_{n+1} $) at time $ t_{n+1} $ are obtained in the last sub-step $ t\in[t_{n}+\gamma_{s-1}\dt,~t_{n}+\gamma_s\dt] $ due to $ \gamma_s=1 $. That is,
\begin{subequations}\label{eq:202103041}
	\begin{align}
		\mbf{M}\mbf{\ddot{U}}_{n+1}&+\mbf{C}\mbf{\dot{U}}_{n+1}+\mbf{K}\mbf{U}_{n+1}  =\mbf{F}(t_{n+1})                                                         \\
		\mbf{U}_{n+1}                                                               & =\mbf{U}_n+\dt\sum_{j=0}^{s}\alpha_{sj}\mbf{\dot{U}}_{n+\gamma_j}        \\
		\mbf{\dot{U}}_{n+1}                                                         & =\mbf{\dot{U}}_n+\dt\sum_{j=0}^{s}\alpha_{sj}\mbf{\ddot{U}}_{n+\gamma_j}.
	\end{align}
\end{subequations}
As stressed previously, the composite sub-step methods can also be viewed as a special case of the Runge-Kutta family, and thus the proposed multi-sub-step method \eqref{eq:nsubstep} can be described in the Butcher tableau \cite{butcherNumericalMethods2016} as 
\begin{equation}\label{eq:suci_tab}
	\begin{BMAT}(b){c|c}{c|c}
		\mbf{c} & \mbf{A}\\
		& \mbf{b}
	\end{BMAT}=\begin{BMAT}(@,30pt,15pt){c|cccccc}{cccccc|c}
		0 & 0 &&&&& \\
		\gamma_1 & \alpha_{10} & \alpha_{11}&&&&\\
		\gamma_2 & \alpha_{20} & \alpha_{21} & \alpha_{22}&&&\\
		\vdots & \vdots & \vdots & \vdots &  \ddots&&\\
		\gamma_{s-1} & \alpha_{(s-1)0} & \alpha_{(s-1)1} & \alpha_{(s-1)2} & \cdots &\alpha_{(s-1)(s-1)}&\\
		\gamma_s & \alpha_{s0} & \alpha_{s1} & \alpha_{s2} & \cdots & \alpha_{s(s-1)} &\alpha_{ss}\\
		& \alpha_{s0} & \alpha_{s1} & \alpha_{s2} & \cdots & \alpha_{s(s-1)} &\alpha_{ss}
	\end{BMAT}
\end{equation}
This paper endeavors to identify algorithmic parameters conducive to achieving identical high-order accuracy and controllable numerical dissipation within the framework of computational structural dynamics. As a result, the developed high-order algorithms given by Eq.~\eqref{eq:suci_tab} are also applicable to the first-order transient problems.

In the initial time $ t_0 $, the initial acceleration vector $ \mbf{\ddot{U}}_0 $ is given by solving the equilibrium equation at $ t_0 $ with given initial conditions $ \mbf{U}_0 $ and $ \mbf{\dot{U}}_0 $.
\begin{equation}\label{eq:a0}
	\mbf{\ddot{U}}_0=\mbf{M}^{-1}\left\{\mbf{F}(t_0)-\mbf{K}\mbf{U}_0-\mbf{C}\mbf{\dot{U}}_0 \right\}
\end{equation}
Hence, the proposed $ s $-sub-step implicit algorithm (\ref{eq:nsubstep}) is obviously self-starting, like the known WBZ-$\alpha$ \cite{woodAlphaModification1980}, HHT-$\alpha$ \cite{hilberImprovedNumerical1977}, and TPO/G-$\alpha$ \cite{shaoThreeParameters1988,chungTimeIntegration1993} algorithms.

\subsubsection{The previous work}
Before the novel implicit method (\ref{eq:nsubstep}) is optimized to determine unknown algorithmic parameters, it is necessary to demonstrate differences between the present method (\ref{eq:nsubstep}) and the authors' previous work \cite{liDirectlySelfstarting2022}. A directly self-starting $ s $-sub-step implicit method was constructed and developed in \cite{liDirectlySelfstarting2022} to achieve unconditional stability, $ s $th-order accuracy, no overshoots, and controllable algorithmic dissipation in the high-frequency range. The published directly self-starting $ s $-sub-step method \cite{liDirectlySelfstarting2022} is written in the $ i $th sub-step $ t\in\left[t_n+\gamma_{i-1}\dt,~t_n+\gamma_i\dt\right] $ as
\begin{subequations}\label{eq:dsuci}
	\begin{align}
		\mbf{M}\widetilde{\ddot{\mbf{U}}}_{n+\gamma_i}&+\mbf{C}\widetilde{\dot{\mbf{U}}}_{n+\gamma_i}+\mbf{K}\widetilde{\mbf{U}}_{n+\gamma_i}=\mbf{F}(t_n+\gamma_i\dt)\\
		\widetilde{\mbf{U}}_{n+\gamma_i}&=\mbf{U}_n+\dt\sum_{j=1}^i\alpha_{ij}\widetilde{\dot{\mbf{U}}}_{n+\gamma_j}\\
		\widetilde{\dot{\mbf{U}}}_{n+\gamma_i}&=\mbf{\dot{U}}_n+\dt\sum_{j=1}^i\alpha_{ij}\widetilde{\ddot{\mbf{U}}}_{n+\gamma_j}
	\end{align}
	with the displacement and velocity updates at $t_{n+1}$:
	\begin{align}
		\mbf{U}_{n+1}&=\mbf{U}_n+\dt\sum_{i=1}^{s}\beta_i\widetilde{\dot{\mbf{U}}}_{n+\gamma_i}\\
		\mbf{\dot{U}}_{n+1}&=\mbf{\dot{U}}_n+\dt\sum_{i=1}^s\beta_i\widetilde{\ddot{\mbf{U}}}_{n+\gamma_i}.
	\end{align}
\end{subequations}
It should be emphasized that the directly self-starting $s$-sub-step implicit method \eqref{eq:dsuci} can also be viewed as the Runge-Kutta method with the following Butcher tableau:
\begin{equation}\label{eq:dsuci_tab}
	\begin{BMAT}(b){c|c}{c|c}
		\mbf{c} & \mbf{A}\\
		& \mbf{b}
	\end{BMAT}=\begin{BMAT}(@,30pt,15pt){c|ccccc}{ccccc|c}
		\gamma_1  & \alpha_{11} & & & & \\
		\gamma_2  & \alpha_{21} & \alpha_{22} &&&\\
		\vdots  & \vdots & \vdots &  \ddots&&\\
		\gamma_{s-1} & \alpha_{(s-1)1} & \alpha_{(s-1)2} & \cdots &\alpha_{(s-1)(s-1)}&\\
		\gamma_s  & \alpha_{s1} & \alpha_{s2} & \cdots & \alpha_{s(s-1)} &\alpha_{ss}\\
		&  \beta_{1} & \beta_{2} & \cdots & \beta_{s-1} &\beta_{s}\\
	\end{BMAT}
\end{equation}
Obviously, the present method differs essentially from the previous method at the algorithm level. In addition, there are other differences between them.  
\begin{itemize}%[(a)]
	\item The previous method (\ref{eq:dsuci}) is directly self-starting since the acceleration vector $ \mbf{\ddot{U}}_n $ does not participate in calculating displacement and velocity updates, whereas the present method (\ref{eq:nsubstep}) is only described to be self-starting since it involves the acceleration vector $ \mbf{\ddot{U}}_n $. As a result, the previous method \eqref{eq:dsuci} can directly start the transient analysis without computing the initial acceleration $\mbfa_0$, whereas the present method \eqref{eq:nsubstep} must accurately compute $\mbfa_0$ by solving Eq.~\eqref{eq:a0} to avoid the accuracy reduction \cite{hulbertErrorAnalysis1987}.
	
	\item For the previous method (\ref{eq:dsuci}), numerical solutions at $ t_{n+1} $ do not generally satisfy the second-order equations of motion (\ref{eq:mck}), whereas numerical solutions of the present method (\ref{eq:nsubstep}) exactly satisfy the system (\ref{eq:mck}) at the instant $ t_{n+1} $.  
	
	\item Although both of the methods are required to achieve $ s $th-order accuracy for the output responses ($ \mbf{U}_{n+1} $, $ \mbf{\dot{U}}_{n+1} $, and $ \mbf{\ddot{U}}_{n+1} $), the previous method (\ref{eq:dsuci}) only requires second-order consistency in each sub-step whereas the present method (\ref{eq:nsubstep}) enhances third-order consistency in each sub-step; see Eq.~(\ref{eq:2orderi}). 
	
	\item In terms of the acceleration accuracy, the previous method (\ref{eq:dsuci}) requires an additional procedure to achieve the same acceleration accuracy as those in displacement and velocity (see Section 4 of the work \cite{liDirectlySelfstarting2022}), whereas the present method (\ref{eq:nsubstep}) automatically provides identical high-order accuracy in displacement, velocity, and acceleration. 
\end{itemize}
In what follows, this paper will further develop high-order implicit members within the $s$-sub-step method \eqref{eq:nsubstep} or the Runge-Kutta method \eqref{eq:suci_tab} by employing the same analysis techniques as the authors' previous work \cite{liDirectlySelfstarting2022}.%the present method \eqref{eq:nsubstep} is a completely new construction and development using the composite sub-step technique. 

\subsection{Accuracy analysis}
It is cumbersome and difficult to analyze numerical properties of an integration algorithm by directly manipulating the multi-degree-of-freedom system (\ref{eq:mck}). However, the modal decomposition technique \cite{hughesFiniteElement2000} can be employed to uncouple the system (\ref{eq:mck}) by using the orthogonality properties of the free-vibration mode shapes of the undamped system. As a result, the modal damping is often used. It is therefore convenient to analyze numerical properties of an integration method by solving the standard SDOF system \cite{rezaiee-pajandImprovingStability2011,liTwoThirdorder2022,rezaiee-pajandNovelTime2017}, which is written as
\begin{equation}\label{eq:sdof}
	\ddot{u}(t)+2\xi\omega\dot{u}(t)+\omega^2u(t)=f(t)
\end{equation}
where $ \xi $, $ \omega $, and $ f(t) $ are the damping ratio, the undamped natural frequency of the system, and the modal force, respectively.

For the subsequent use, considering the above SDOF system with given initial displacement $ u_n $ and velocity $ \dot{u}_n $ and external load $ f(t)=\exp(t) $, its exact displacement and velocity are analytically obtained as
\begin{subequations}
	\begin{align}
		&\begin{aligned}
			u(t)= & \exp(-\xi\omega \tau)(\cos(\omega_d\tau)+\frac{\xi\omega}{\omega_d}\sin(\omega_d\tau))u_n+\frac{1}{\omega_d}\exp(-\xi\omega \tau)\sin(\omega_d\tau)\dot{u}_n+\frac{\exp(\tau)}{(\omega^2+2\xi\omega+1)} \\
			&\quad -\exp(-\xi\omega \tau)\frac{\omega_d\cos(\omega_d\tau)+(1+\xi\omega)\sin(\omega_d\tau)}{\omega_d(\omega^2+2\xi\omega+1)}
		\end{aligned}\\
		&\begin{aligned}
			\dot{u}(t)= & -\frac{\omega^2}{\omega_d}\exp(-\xi\omega \tau)\sin(\omega_d\tau)u_n+\frac{1}{\omega_d}\exp(-\xi\omega \tau)(-\xi\omega\sin(\omega_d\tau)+\omega_d\cos(\omega_d\tau))\dot{u}_n+\frac{\exp(\tau)}{(\omega^2+2\xi\omega+1)} \\
			&\quad +\exp(-\xi\omega \tau)\frac{(\omega^2+\xi\omega)\sin(\omega_d\tau)-\omega_d\cos(\omega_d \tau)}{\omega_d(\omega^2+2\xi\omega+1)}
		\end{aligned}
	\end{align}
\end{subequations}
where $ \tau=t-t_n $ and $ \omega_d=\sqrt{1-\xi^2}\cdot\omega $. Then, the exact solutions at time $ t_{n+1} $ can be collected as
\begin{equation}\label{eq:uvexa}
	\begin{bmatrix}
		u(t_{n+1}) \\ \dot{u}(t_{n+1})
	\end{bmatrix}=\mbf{D}_{exa}\begin{bmatrix}
		u_n \\ \dot{u}_n
	\end{bmatrix}+\mbf{L}_{exa}
\end{equation}
where $ \mbf{D}_{exa} $ is the exact amplification matrix \cite{liDirectlySelfstarting2022,rezaiee-pajandEfficientWeighted2021}, which is
\begin{subequations}
	\begin{equation}\label{key}
		\mbf{D}_{exa}=\exp(-\xi\omega\dt)\begin{bmatrix}
			\displaystyle\cos(\omega_d\dt)+\frac{\xi\omega}{\omega_d}\sin(\omega_d\dt) & \displaystyle\frac{1}{\omega_d}\sin(\omega_d\dt)                         \\[2mm]
			\displaystyle-\frac{\omega^2}{\omega_d}\sin(\omega_d\dt)                     & \displaystyle -\frac{\xi\omega}{\omega_d}\sin(\omega_d\dt)+\cos(\omega_d\dt)
		\end{bmatrix},
	\end{equation}
	and $ \mbf{L}_{exa} $ is the exact load operator associated with the external load $ \exp(t) $.
	\begin{equation}\label{key}
		\mbf{L}_{exa}=\begin{bmatrix}
			\displaystyle\frac{\exp(\dt)}{(\omega^2+2\xi\omega+1)}-\exp(-\xi\omega\dt)\frac{\omega_d\cos(\omega_d\dt)+(1+\xi\omega)\sin(\omega_d\dt)}{\omega_d(\omega^2+2\xi\omega+1)} \\[2mm]
			\displaystyle\frac{\exp(\dt)}{(\omega^2+2\xi\omega+1)}+\exp(-\xi\omega \dt)\frac{(\omega^2+\xi\omega)\sin(\omega_d\dt)-\omega_d\cos(\omega_d \dt)}{\omega_d(\omega^2+2\xi\omega+1)}
		\end{bmatrix}
	\end{equation}
\end{subequations}
Note that the external load $f(t)=\exp(t)$ is not the only choice in accuracy analysis, and other appropriate external loads, such as $f(t)=\sin(t)+\cos(t)$, can also be used to derive the same conditions as $f(t)=\exp(t)$. Appendix A provides the theoretical explanation for the rationality of using $f(t)=\exp(t)$.

On the other hand, applying the proposed implicit algorithm (\ref{eq:nsubstep}) to solve the SDOF system (\ref{eq:sdof}) can yield 
\begin{subequations}\label{eq:ljz}
	\begin{align}
		\ddot{u}_{n+\gamma_i}&+2\xi\omega\dot{u}_{n+\gamma_i}+\omega^2u_{n+\gamma_i}  =f(t_n+\gamma_i\dt)                                              \\
		{u}_{n+\gamma_i}                                                            & ={u}_n+\dt\sum_{j=0}^{i}\alpha_{ij}\dot{{u}}_{n+\gamma_j}        \\
		\dot{{u}}_{n+\gamma_i}                                                      & =\dot{{u}}_n+\dt\sum_{j=0}^{i}\alpha_{ij}\ddot{{u}}_{n+\gamma_j}
	\end{align}
	with varying the index $ i\in\{1,~2,~\dots,~s\} $.
\end{subequations}
Inserting Eq.~(\ref{eq:ljz}b) into Eq.~(\ref{eq:ljz}a) to eliminate $ u_{n+\gamma_i} $ gives
\begin{equation}\label{eq:ddY}
	\begin{aligned}
		f(t_n+\gamma_i\dt) & =\ddot{u}_{n+\gamma_i}+2\xi\omega\dot{{u}}_{n+\gamma_{i}}+\omega^2\left({u}_n+\dt\sum_{j=0}^{i}\alpha_{ij}\dot{{u}}_{n+\gamma_j} \right)         =\left(2\xi\omega\dot{{u}}_{n+\gamma_{i}}+\omega^2\dt\sum_{j=0}^i\alpha_{ij}\dot{{u}}_{n+\gamma_{j}}\right)+\ddot{{u}}_{n+\gamma_{i}}+\omega^2u_n.
	\end{aligned}
\end{equation}
Define the notations $ \dot{\mbf{Y}} $ and $ \ddot{\mbf{Y}} $ as $
\dot{\mbf{Y}}=\begin{bmatrix}
	\dot{u}_{n+\gamma_0} &
	\dot{u}_{n+\gamma_1} &
	\dots                &
	\dot{u}_{n+\gamma_{s}}
\end{bmatrix}^\mathsf{T}$ and $ \ddot{\mbf{Y}}=\begin{bmatrix}
	\ddot{u}_{n+\gamma_0} &
	\ddot{u}_{n+\gamma_1} &
	\dots                 &
	\ddot{u}_{n+\gamma_{s}}
\end{bmatrix}^\mathsf{T}
$, respectively, and also define $ \mbf{f} $ as
$
\mbf{f}=\begin{bmatrix}
	f(t_n+\gamma_0\dt) & f(t_n+\gamma_1\dt) & \dots & f(t_n+\gamma_s\dt)
\end{bmatrix}^\mathsf{T}
$.
Then, Eq.~(\ref{eq:ddY}) with varying $ i=\{0,~1,~\cdots,~s\} $ can be written as
\begin{equation}\label{eq:dotY}
	\left(2\xi\omega\mbf{I}+\omega^2\dt\mbf{A}\right)\dot{\mbf{Y}}+\ddot{\mbf{Y}}=\mbf{f}-\omega^2u_n\mbf{1}
\end{equation}
where $ \mbf{A} $ is defined in Eq.~\eqref{eq:suci_tab}; $ \mbf{1} $ is the column vector whose all elements are unity, and $ \mbf{I} $ is the unity matrix with dimension $ s+1 $. Similarly, Eq.~(\ref{eq:ljz}c) can be rewritten as
\begin{equation}\label{eq:dotY1}
	\dot{\mbf{Y}}-\dt\mbf{A}\ddot{\mbf{Y}}=\dot{{u}}_n\mbf{1}.
\end{equation}
Solving Eqs.~(\ref{eq:dotY}) and (\ref{eq:dotY1}) yields
\begin{equation}\label{eq:YY}
	\begin{bmatrix}
		\dot{\mbf{Y}} \\ \ddot{\mbf{Y}}
	\end{bmatrix}=\begin{bmatrix}
		2\xi\omega\mbf{I}+\omega^2\dt\mbf{A} & \mbf{I}     \\
		\mbf{I}                              & -\dt\mbf{A}
	\end{bmatrix}^{-1}\left\{\begin{bmatrix}
		\mbf{f} \\ \mbf{0}
	\end{bmatrix}+\begin{bmatrix}
		-\omega^2\mbf{1} & \mbf{0} \\ \mbf{0} & \mbf{1}
	\end{bmatrix}\begin{bmatrix}
		u_n \\ \dot{u}_n
	\end{bmatrix}\right\}.
\end{equation}
where $ \mbf{0} $ denotes a column vector whose all elements are zero.

Finally, numerical displacement $ u_{n+1} $ and velocity $ \dot{{u}}_{n+1} $ are updated in the last sub-step due to $ \gamma_s=1 $, that is
\begin{subequations}
	\begin{align}
		u_{n+1}&=u_n+\dt\sum_{j=0}^s\alpha_{sj}\dot{{u}}_{n+\gamma_j}\\
		\dot{u}_{n+1}&=\dot{u}_n+\dt\sum_{j=0}^s\alpha_{sj}\ddot{{u}}_{n+\gamma_j}.
	\end{align}
\end{subequations}
The above equations can be further written as
\begin{equation}\label{eq:20210304}
	\begin{bmatrix}
		u_{n+1} \\ \dot{u}_{n+1}
	\end{bmatrix}=\begin{bmatrix}
		1 & 0 \\0 & 1
	\end{bmatrix}\begin{bmatrix}
		u_n \\ \dot{u}_n
	\end{bmatrix}+\begin{bmatrix}
		\dt\mbf{b}^\mathsf{T} & \mbf{0}^\mathsf{T} \\ \mbf{0}^\mathsf{T} & \dt\mbf{b}^\mathsf{T}
	\end{bmatrix}\begin{bmatrix}
		\dot{\mbf{Y}} \\ \ddot{\mbf{Y}}
	\end{bmatrix}
\end{equation}
where the column vector $ \mbf{b} $ is defined in Eq.~\eqref{eq:suci_tab}. Substituting Eq.~(\ref{eq:YY}) into Eq.~(\ref{eq:20210304}) yields
\begin{equation}\label{eq:uvnum1}
	\begin{bmatrix}
		u_{n+1} \\ \dot{u}_{n+1}
	\end{bmatrix}=\mbf{D}_{num}\begin{bmatrix}
		u_n \\ \dot{u}_n
	\end{bmatrix}+\mbf{L}_{num}
\end{equation}
where $ \mbf{D}_{num} $ and $ \mbf{L}_{num} $ are called the numerical amplification matrix and load operator, respectively.
\begin{subequations}\label{eq:DL}
	\begin{align}
		\mbf{D}_{num} & =\begin{bmatrix}
			1 &\quad 0 \\0 &\quad 1
		\end{bmatrix}+\begin{bmatrix}
			\dt\mbf{b}^\mathsf{T} & \mbf{0}^\mathsf{T} \\ \mbf{0}^\mathsf{T} & \dt\mbf{b}^\mathsf{T}
		\end{bmatrix}\begin{bmatrix}
			2\xi\omega\mbf{I}+\omega^2\dt\mbf{A} & \mbf{I}     \\
			\mbf{I}                              & -\dt\mbf{A}
		\end{bmatrix}^{-1}\begin{bmatrix}
			-\omega^2\mbf{1} & \mbf{0} \\ \mbf{0} & \mbf{1}
		\end{bmatrix} \\
		\mbf{L}_{num} & =\begin{bmatrix}
			\dt\mbf{b}^\mathsf{T} & \mbf{0}^\mathsf{T} \\ \mbf{0}^\mathsf{T} & \dt\mbf{b}^\mathsf{T}
		\end{bmatrix}\begin{bmatrix}
			2\xi\omega\mbf{I}+\omega^2\dt\mbf{A} & \mbf{I}     \\
			\mbf{I}                              & -\dt\mbf{A}
		\end{bmatrix}^{-1}\begin{bmatrix}
			\mbf{f} \\ \mbf{0}
		\end{bmatrix}
	\end{align}
\end{subequations}

With the exact and numerical iterative relationships in hand, the conditions for achieving the novel method' $ p $th-order accuracy can be defined as follows.
\begin{proposition}\label{pos:accuracy}\emph{
		The proposed $ s $-sub-step method (\ref{eq:nsubstep}) achieves $ p $th-order accuracy for the standard SDOF system \eqref{eq:sdof}} if and only if both $ u_{n+1}-u(t_{n+1})=\mathcal{O}(\dt^{p+1}) $ and $ \dot{u}_{n+1}-\dot{u}(t_{n+1})=\mathcal{O}(\dt^{p+1}) $ are satisfied, namely
		\begin{subequations}\label{eq:accuracy}
			\begin{align}
				\mbf{D}_{num}-\mbf{D}_{exa} & =\mbf{O}(\dt^{p+1}) \\
				\mbf{L}_{num}-\mbf{L}_{exa} & =\mbf{O}(\dt^{p+1}).
			\end{align}
	\end{subequations}
\end{proposition}
Notice that Proposition \ref{pos:accuracy} actually requires the numerical scheme in the last sub-step to be $ p $th-order accurate due to $ \gamma_s=1 $.
To eliminate the redundant algorithmic parameters, an additional requirement is imposed in this paper.
In the $ i $th sub-step $ t\in\left[t_n+\gamma_{i-1}\dt,~t_n+\gamma_i\dt\right] $, the displacement predicted by Eq.~(\ref{eq:ljz}b) is defined to be $ q $th-order consistency if $ u(t_n+\gamma_i\dt)-u_{n+\gamma_i}=\mathcal{O}(\dt^{q}) $ is fulfilled, namely,
\begin{equation}\label{eq:ljz1}
	u(t_n+\gamma_i\dt)-\left[u(t_n)+\dt\sum_{j=0}^i\alpha_{ij}\dot{{u}}(t_n+\gamma_j\dt)\right]=\mathcal{O}(\dt^{q}).
\end{equation}
Calculating Taylor series expansions of $ u(t_n+\gamma_i\dt) $ and $ \dot{u}(t_n+\gamma_j\dt) $ at time $ t_n $ gives, respectively,
	\begin{subequations}\label{eq:20230407}
		\begin{equation}\label{key}
			u(t_n+\gamma_i\dt)=u(t_n)+\gamma_i\dt\dot{u}(t_n)+\dfrac{\gamma_i^2}{2}\dt^2\ddot{u}(t_n)+O(\dt^3)
		\end{equation}
		and 
		\begin{equation}\label{key}
			\dot{u}(t_n+\gamma_j\dt)=\dot{u}(t_n)+\gamma_j\dt\ddot{u}(t_n)+O(\dt^2).
		\end{equation}
	\end{subequations}

Substituting Eq.~\eqref{eq:20230407} into Eq.~\eqref{eq:ljz1} and then simplifying the result yield
\begin{equation}\label{eq:ljz2}
	u(t_n+\gamma_i\dt)-u_{n+\gamma_i}=\left(\gamma_i-\sum_{j=0}^i\alpha_{ij}\right)\dt\dot{{u}}(t_n)+\left(\frac{\gamma_i^2}{2}-\sum_{j=0}^i\alpha_{ij}\gamma_j\right)\dt^2\ddot{{u}}(t_n)+\mathcal{O}(\dt^3).
\end{equation}
Then, the following conditions are satisfied to eliminate low-order terms, reaching a third-order consistency \cite{liTrulySelfstarting2020}.
\begin{equation}\label{eq:2orderi}
	\sum_{j=0}^i\alpha_{ij}=\gamma_i,\quad \sum_{j=0}^i\alpha_{ij}\gamma_j=\frac{\gamma_{i}^2}{2},\qquad \forall i\in1,~\cdots,~s
\end{equation}
Similarly, the velocity predicted by Eq.~(\ref{eq:ljz}c) is defined to be $q$th-order consistency if $ \dot{u}(t_n+\gamma_i\dt)-\dot{u}_{n+\gamma_i}=\mathcal{O}(\dt^{q}) $ is fulfilled. Adopting the same analysis as Eq.~\eqref{eq:ljz2}, one can easily get
\begin{equation}\label{eq:20230404}
	\dot{u}(t_n+\gamma_i\dt)-\dot{u}_{n+\gamma_i}=\left(\gamma_i-\sum_{j=0}^i\alpha_{ij}\right)\dt\ddot{{u}}(t_n)+\left(\frac{\gamma_i^2}{2}-\sum_{j=0}^i\alpha_{ij}\gamma_j\right)\dt^2\dddot{{u}}(t_n)+\mathcal{O}(\dt^3).
\end{equation}
Eq.~\eqref{eq:20230404} demonstrates that the conditions (\ref{eq:2orderi}) also ensure a third-order consistency in velocity (\ref{eq:ljz}c). As stated previously, the published directly self-starting implicit algorithms \cite{liDirectlySelfstarting2022} only require a second-order consistency in each sub-step.
\begin{remark}
In the case of $ i=1 $, namely in the first sub-step, the following relations given by Eq.~\eqref{eq:2orderi} hold:
\begin{equation}\label{key}
	\alpha_{10}+\alpha_{11}=\gamma_1\quad\text{and}\quad \alpha_{10}\gamma_0+\alpha_{11}\gamma_1=\frac{\gamma_1^2}{2}.
\end{equation}
Due to $ \gamma_0=0 $, the above equations give $ \alpha_{10}=\alpha_{11}=\gamma_1/2 $, implying that the non-dissipative trapezoidal rule \cite{hughesFiniteElement2000} should be used in the first sub-step $ t\in\left[t_n,~t_n+\gamma_1\dt\right] $. 
\end{remark}

\subsection{Dissipation control}
The previous subsection has derived the numerical amplification matrix $ \mbf{D}_{num} $, thus it is very easy to design controllable numerical dissipation in the high-frequency range. Firstly, the characteristic polynomial of $ \mbf{D}_{num} $ is computed as
\begin{equation}\label{eq:cp}
	\det(\mbf{I}_2-\zeta\mbf{D}_{num})=\zeta^2-2A_1\zeta+A_2=0
\end{equation}
where $ \mbf{I}_2 $ denotes the unity matrix with dimension $ 2 $, and coefficients $ A_1 $ and $ A_2 $ are two principal invariants of $ \mbf{D}_{num} $. When considering the high-frequency limit ($ \omega\to\infty $), the characteristic polynomial (\ref{eq:cp}) can be expressed as
\begin{equation}\label{eq:cp1}
	\zeta_{\infty}^2-2A_{1\infty}\zeta_{\infty}+A_{2\infty}=0
\end{equation}
where $ A_{1\infty}=\lim\limits_{\omega\to\infty} A_1 $ and $ A_{2\infty}=\lim\limits_{\omega\to\infty} A_2 $.
On the other hand, the optimal high-frequency dissipation \cite{chungTimeIntegration1993} is achieved if and only if eigenvalues of the amplification matrix $ \mbf{D}_{num} $ become the same real roots, denoted by $ \rhoinf $, in the high-frequency limit. Then,
\begin{equation}\label{eq:cp2}
	(\zeta_{\infty}-\rhoinf)^2=\zeta_{\infty}^2-2\rhoinf\zeta_{\infty}+\rhoinf^2=0.
\end{equation}

Comparing Eqs.~(\ref{eq:cp1}) and (\ref{eq:cp2}) yields the conditions to achieve controllable high-frequency dissipation:
\begin{equation}\label{eq:optimaldissipation}
	A_{1\infty}  =\rhoinf\quad\text{and}\quad
	A_{2\infty}  =\rhoinf^2.
\end{equation}

It is necessary to point out that two invariants $ A_1 $ and $ A_2 $ given by Eq.~(\ref{eq:cp}) are often employed to analyze the stability of an integrator. When the following conditions \cite{hughesFiniteElement2000} are fulfilled, the proposed algorithms are said to be unconditionally stable.
\begin{subequations}\label{eq:uc}
	\begin{align}
		|2A_1|\le A_2+1,\quad & -1\le A_2<1 \\
		|A_1|<1,\quad         & A_2=1
	\end{align}
\end{subequations}
After embedding high-order accuracy and controllable numerical dissipation, the developed implicit algorithm (\ref{eq:nsubstep}) will achieve unconditional stability via selecting proper $ \gamma_1 $. It should be emphasized that, after achieving identical high-order accuracy and controllable numerical dissipation, the developed $s$-sub-step method can only achieve fundamental spectral stability instead of the BN-stability \cite{jiUnconditionallyStable2021}. %In Section \ref{sec:sp}, these numerical properties will be validated well.

\section{High-order accurate schemes}\label{sec:alg}
%Via the conditions derived in the previous section, the proposed $ s $-sub-step implicit method (\ref{eq:nsubstep}) will be analyzed to achieve controllable numerical dissipation, high-order accuracy, and unconditional stability. As a result, the composite $ s $-sub-step method (\ref{eq:nsubstep}) with $ s\le6 $ can generally achieve $ s $th-order accuracy after simultaneously embedding controllable high-frequency dissipation and unconditional stability. 
Apart from the conditions derived in the previous section, the identity of effective stiffness matrices is also required for the present $ s $-sub-step implicit method (\ref{eq:nsubstep}) to attain optimal spectral properties \cite{liNovelFamily2019} and to reduce computational costs for solving linear elastic problems, which requires
\begin{equation}\label{eq:iesm}
	\alpha_{11}=\alpha_{22}=\cdots=\alpha_{ss}.
\end{equation}
Using the conditions \eqref{eq:iesm}, the proposed $s$-sub-step method given by Eq.~\eqref{eq:suci_tab} reduces to the ESDIRK method in mathematics. The integration algorithms developed in this section are named as SUCI$n$, where the first four letters denote the sub-step, unconditional stability, controllable numerical dissipation, and identical effective stiffness matrices, while the final letter $ n $ represents either the order of accuracy or the number of sub-steps (\textit{these two quantities are the same in this paper}). % Since the composite algorithms developed in this paper correspond to the Runge-Kutta family, the algorithm structure  

Obviously, the proposed $ s $-sub-step method (\ref{eq:nsubstep}) corresponds simply to the non-dissipative trapezoidal rule in the case of $ s=1 $. For the two-sub-step case ($ s=2 $), the resulting scheme covers the published two-sub-step implicit algorithm \cite{liNovelFamily2020}, also named as SUCI2 in this paper, rewritten herein in the Butcher tableau as
\begin{equation}\label{eq:liyu}
		\begin{BMAT}(b){c|c}{c|c}
		\mbf{c} & \mbf{A}\\ 
		& \mbf{b}
	\end{BMAT}=\begin{BMAT}{c|ccc}{ccc|c}
		0 & 0 & & \\
		\gamma_1 & \dfrac{\gamma_1}{2} & \dfrac{\gamma_1}{2}&\\
		1 & \dfrac{-\gamma_1^2+3\gamma_1-1}{2\gamma_1} & \dfrac{1-\gamma_1}{2\gamma_1} & \dfrac{\gamma_1}{2}\\
		&\dfrac{-\gamma_1^2+3\gamma_1-1}{2\gamma_1} & \dfrac{1-\gamma_1}{2\gamma_1} & \dfrac{\gamma_1}{2}\\
	\end{BMAT}
\end{equation}
where the parameter $ \gamma_1 $ is correlated with the ultimate spectral radius ($\rhoinf$) through $\gamma_1=\left( 2-\sqrt{2(1+\rhoinf)} \right)/(1-\rhoinf)$. % controls numerical dissipation in the high-frequency range.
%\begin{equation}\label{key}
%	\gamma_1=\frac{2-\sqrt{2(1+\rhoinf)}}{1-\rhoinf}
%\end{equation}
%where $ \rhoinf\in[-1,~1] $. 
In the case of $ \rhoinf=1 $, the parameter $ \gamma_1 $ given above should be set as $ 1/2 $, covering the composite equal sub-step trapezoidal rule \cite{LiFurtherAssessment2021}.

SUCI2 is optimal with respect to spectral properties since some existing two-sub-step methods, such as the $ \rhoinf $-Bathe \cite{nohBatheTime2019} algorithm, correspond to it when embedding identical effective stiffness matrices within two sub-steps. When achieving a full range of controllable numerical dissipation, the second-order accuracy is only obtained for SUCI2. For the proposed $ s $-sub-step method (\ref{eq:nsubstep}), more sub-steps are generally necessary to achieve higher-order accuracy, dissipation control, and unconditional stability. The $ s $-sub-step schemes with $ 2\le s\le 6 $ will be developed in the following. It should be pointed out that the sub-step splitting ratios $ \gamma_i~(i=2,~\cdots,~s-1) $ are not pre-assumed to satisfy the relations $ \gamma_i=i\cdot\gamma_1 $, so the derivations below are general so that the expressions of $ a_{ij} $ are a little complicated for the five- and six-sub-step implicit members.

\subsection{Three-sub-step third-order scheme: SUCI3}
For the three-sub-step scheme, the detailed analysis will be carried out to clearly demonstrate how to determine the unknown algorithmic parameters $ \alpha_{ij} $. The proposed $ s $-sub-step method \eqref{eq:suci_tab} in the case of $ s=3 $ can be rewritten as ($\gamma_3=1$)

\begin{equation}\label{eq:3substep}
		\begin{BMAT}(b){c|c}{c|c}
		\mbf{c} & \mbf{A}\\ %\hline
		& \mbf{b}
	\end{BMAT}=\begin{BMAT}(@,25pt,10pt){c|cccc}{cccc|c}
		0 & 0 &&&\\
\gamma_1 & \dfrac{\gamma_1}{2} & \dfrac{\gamma_1}{2}&&\\
\gamma_2 & \alpha_{20} & \alpha_{21} & \dfrac{\gamma_1}{2}&\\
1 & \alpha_{30} & \alpha_{31} & \alpha_{32} & \dfrac{\gamma_1}{2}\\%		\hline
& \alpha_{30} & \alpha_{31} & \alpha_{32} & \dfrac{\gamma_1}{2}\\
	\end{BMAT}
\end{equation}
Notice that the above scheme (\ref{eq:3substep}) has used the conditions (\ref{eq:iesm}) to achieve identical effective stiffness matrices within three sub-steps, i.e.,~$ \alpha_{33}=\alpha_{22}=\alpha_{11}=\gamma_1/2 $. 
%In this case, the resulting scheme can save computational costs for solving linear elastic problems since the triangular factorization of the effective stiffness matrix is required once during the whole simulation. More importantly, identical effective stiffness matrices within each sub-step guarantee optimal spectral properties, that is the maximum numerical dissipation but minimum period elongation error. 
When the three-sub-step implicit scheme \eqref{eq:3substep} is developed to achieve second-order accuracy, the design idea covers the authors' previous work \cite{liThreeOpitmal2023}. This paper will use the three-sub-step implicit scheme \eqref{eq:3substep} to achieve third-order accuracy, controllable numerical dissipation, zero-order overshoots, and unconditional stability. 

\subsubsection{Third-order accuracy}
In what follows, the accuracy analysis is carried out to determine algorithmic parameters $ \alpha_{ij},0\le j<i\le 3 $. Firstly, the conditions (\ref{eq:2orderi}) with $ i=2$ and $3 $ are explicitly expressed as
\begin{subequations}
	\begin{align}
		\alpha_{20}+\alpha_{21}+\frac{\gamma_1}{2}                                             & =\gamma_2& \alpha_{30}+\alpha_{31}+\alpha_{32}+\frac{\gamma_1}{2}                                  &=\gamma_3  &&&& &&          \\
		\alpha_{20}\gamma_0+\alpha_{21}\gamma_1+\frac{\gamma_1}{2}\gamma_2                     & =\frac{\gamma_2^2}{2}&
		\alpha_{30}\gamma_0+\alpha_{31}\gamma_1+\alpha_{32}\gamma_2+\frac{\gamma_1}{2}\gamma_3  &=\frac{\gamma_3^2}{2}.&&&&
	\end{align}
\end{subequations}
Solving the above conditions with $ \gamma_0=0$ and $\gamma_3=1 $ gives
\begin{equation}\label{eq:alpha3sub1}
	\alpha_{20}=\frac{-\gamma_1^2+3\gamma_1\gamma_2-\gamma_2^2}{2\gamma_1}    \ \                 \alpha_{21}=\frac{\gamma_2(\gamma_2-\gamma_1)}{2\gamma_1}      \ \
	\alpha_{30}=\frac{-\gamma_1^2+(3-2\alpha_{32})\gamma_1+2\alpha_{32}\gamma_2-1}{2\gamma_1} \ \ \alpha_{31}=\frac{-2\alpha_{32}\gamma_2-\gamma_1+1}{2\gamma_1}.
	%	\begin{aligned}
		%		\alpha_{20}&=\frac{-\gamma_1^2+3\gamma_1\gamma_2-\gamma_2^2}{2\gamma_1}                    & \alpha_{21}&=\frac{\gamma_2(\gamma_2-\gamma_1)}{2\gamma_1}      \\
		%		\alpha_{30}&=\frac{-\gamma_1^2+(3-2\alpha_{32})\gamma_1+2\alpha_{32}\gamma_2-1}{2\gamma_1} & \alpha_{31}&=\frac{-2\alpha_{32}\gamma_2-\gamma_1+1}{2\gamma_1}
		%	\end{aligned}
\end{equation}
With the known conditions (\ref{eq:alpha3sub1}) in hand, the numerical amplification matrix and load operator can be further simplified by using Eq.~(\ref{eq:DL}), in which $ \mbf{A} $, $ \mbf{b} $ and $ \mbf{f} $ are given, respectively, as
\begin{equation}\label{eq:ljz4}
	\mbf{A}=\left[\begin{BMAT}(@,15pt,15pt){cccc}{cccc}
		0                  & 0                  & 0                  & 0                  \\
		\dfrac{\gamma_1}{2} & \dfrac{\gamma_1}{2} & 0                  & 0                  \\
		\alpha_{20}        & \alpha_{21}        & \dfrac{\gamma_1}{2} & 0                  \\
		\alpha_{30}        & \alpha_{31}        & \alpha_{32}        & \dfrac{\gamma_1}{2}
	\end{BMAT}\right],\quad \mbf{b}=\left[\begin{BMAT}(@,20pt,20pt){c}{cccc}
		\alpha_{30} \\ \alpha_{31}\\ \alpha_{32} \\ \dfrac{\gamma_1}{2}
	\end{BMAT}\right],\quad\text{and}\quad \mbf{f}=\left[\begin{BMAT}(@,20pt,20pt){c}{cccc}
		1 \\ \exp(\gamma_1\dt) \\ \exp(\gamma_2\dt) \\ \exp(\dt)
	\end{BMAT}\right].
\end{equation}

Using Proposition \ref{pos:accuracy} and Eq.~(\ref{eq:alpha3sub1}), numerical errors in the amplification matrix and load operator are calculated as
\begin{subequations}\label{eq:ljz3}
	\begin{align}
		\mbf{D}_{num}-\mbf{D}_{exa} & =\begin{bmatrix}
			\displaystyle-\frac{\xi\omega^3c_{30}}{6}\dt^3+\mathcal{O}(\dt^4)        & \displaystyle\frac{(1-4\xi^2)\omega^2c_{30}}{12}\dt^3+\mathcal{O}(\dt^4)   \\[2mm]
			\displaystyle\frac{(4\xi^2-1)\omega^4c_{30}}{12}\dt^3+\mathcal{O}(\dt^4) & \displaystyle\frac{\xi(2\xi^2-1)\omega^3c_{30}}{3}\dt^3+\mathcal{O}(\dt^4)
		\end{bmatrix} \\
		\mbf{L}_{num}-\mbf{L}_{exa} & =\begin{bmatrix}
			\displaystyle\frac{(2\xi\omega-1)c_{30}}{12}\dt^3+\mathcal{O}(\dt^4) \\[2mm] 
			\displaystyle\frac{(2\xi\omega-1+\omega^2(1-4\xi^2))c_{30}}{12}\dt^3+\mathcal{O}(\dt^4)
		\end{bmatrix}
	\end{align}
\end{subequations}
where  $
c_{30}=6(\gamma_1-\gamma_2)\gamma_2\alpha_{32}+3\gamma_1^2-6\gamma_1+2
$. %{\fr The detailed derivations about Eq.~\eqref{eq:ljz3} are given in Appendix A.} 
Obviously, $ c_{30}=0 $ is solved to achieve third-order accuracy, namely
\begin{equation}\label{key}
	\alpha_{32}=\frac{3\gamma_1^2-6\gamma_1+2}{6\gamma_2(\gamma_2-\gamma_1)}.
\end{equation}

\subsubsection{Unconditional stability}
The unconditional stability of SUCI3 is analyzed herein as a demonstration. After achieving third-order accuracy, the characteristic polynomial (\ref{eq:cp}) in the case of $ s=3 $ can be simplified and its principal invariants ($ A_1 $ and $ A_2 $) are given in the absence of $ \xi $ as
\begin{subequations}
	\begin{align}
		A_1&=\frac{\gamma_1^3(3\gamma_1^3-18\gamma_1^2+18\gamma_1-4)\Omega^6+\gamma_1(36\gamma_1^3+24\gamma_1^2-144\gamma_1+48)\Omega^4+(144\gamma_1^2-96)\Omega^2+192}{3(\gamma_1^2\Omega^2+4)^3}\\
		A_2&=\frac{3(3\gamma_1^3-18\gamma_1^2+18\gamma_1-4)\Omega^6+12(9\gamma_1^4+12\gamma_1^3-36\gamma_1^2+24\gamma_1-4)\Omega^4+432\gamma_1^2\Omega^2+576}{9(\gamma_1^2\Omega^2+4)^3}
	\end{align}
\end{subequations}
where $ \Omega=\omega\cdot\dt $. Due to the quantity $ A_1^2-A_2=-576\left((\gamma_1-1)\gamma_1^2(\gamma_1-1/3)\Omega^4+(2\gamma_1^2-4/9)\Omega^2+8/3\right)^2\Omega^2/(\gamma_1^2\Omega^2+4)^6\le 0 $, two eigenvalues of the amplification matrix can be expressed as
\begin{equation}\label{key}
	\zeta_{1,2}=A_1\pm\sqrt{A_2-A_1^2}\cdot\text{Im}
\end{equation}
where $ \text{Im} $ denotes the imaginary unit being defined as $ \text{Im}=\sqrt{-1} $. According to the definition \cite{hughesFiniteElement2000} of the spectral radius ($ \rho $), SUCI3's spectral radius ($ \rho $) is simply calculated as
\begin{equation}\label{key}
	\rho=\max\left(|\zeta_{1,2}|\right)=\sqrt{A_2}.
\end{equation}
Note that $ A_2\ge0 $ is indicated by $ A_1^2-A_2\le0 $. In this case, the proposed SUCI3 algorithm is said to be unconditionally stable if and only if $ \rho\le1 $ for all $\Omega\in[0,~\infty)$, equivalently,
\begin{equation}\label{eq:20231204}
	A_2-1=-\frac{24\Omega^4\aleph}{(\gamma_1^2\Omega^2+4)^6}\cdot\left[\left(\gamma_1-\frac13\right)\left(\gamma_1-\frac23\right)\left(\gamma_1^3-3\gamma_1^2+3\gamma_1-\frac23\right)\Omega^2-\frac43\gamma_1^3+4\gamma_1^2-\frac83\gamma_1+\frac49\right]\le 0
\end{equation}
where $ \aleph $ denotes a complicated factor and it is always greater than zero. Eq.~\eqref{eq:20231204} illustrates that SUCI3 is unconditionally stable if and only if the parameter $\gamma_1$ satisfies 
\begin{subequations}
	\begin{align}
		\left(\gamma_1-\frac13\right)\left(\gamma_1-\frac23\right)\left(\gamma_1^3-3\gamma_1^2+3\gamma_1-\frac23\right)&\ge 0\\
		-\frac43\gamma_1^3+4\gamma_1^2-\frac83\gamma_1+\frac49&\ge 0.
	\end{align}
\end{subequations}
Solving the inequalities above yields\begin{equation}\label{eq:suci3_uc_g1}
	\gamma_1\in\left[\dfrac23,~2.137158043\right].
\end{equation}
Hence, SUCI3 achieves unconditional stability if the parameter $\gamma_1$ satisfies Eq.~\eqref{eq:suci3_uc_g1}.

\subsubsection{Dissipation control}
\begin{figure}[htbp]
	\centering
	\subfigtopskip=2pt 
	\subfigbottomskip=-4pt
	\subfigcapskip=-5pt 
	\subfigure[SUCI3 ]{
		\includegraphics[scale=0.32]{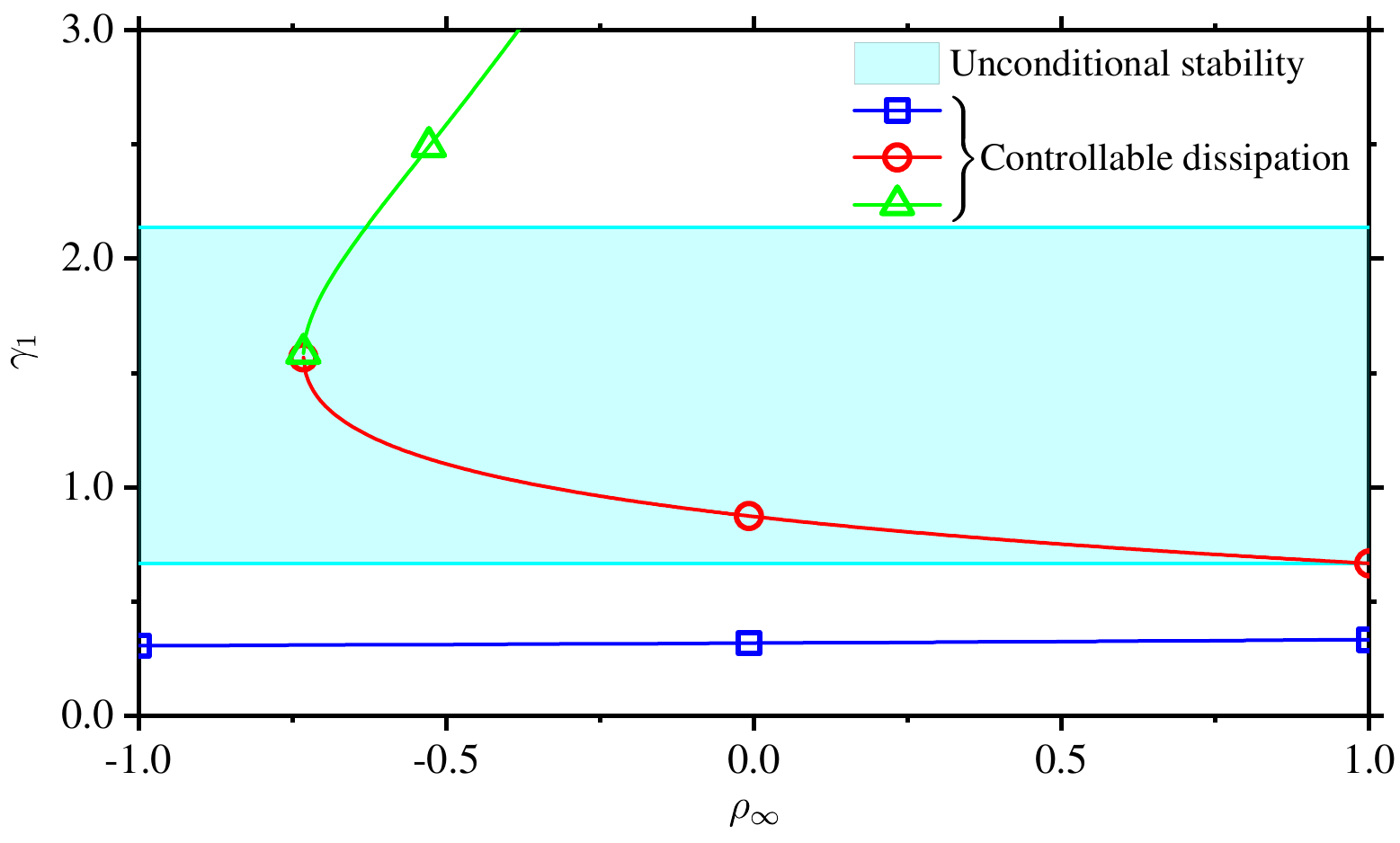}}
	\subfigure[SUCI4 ]{
		\includegraphics[scale=0.32]{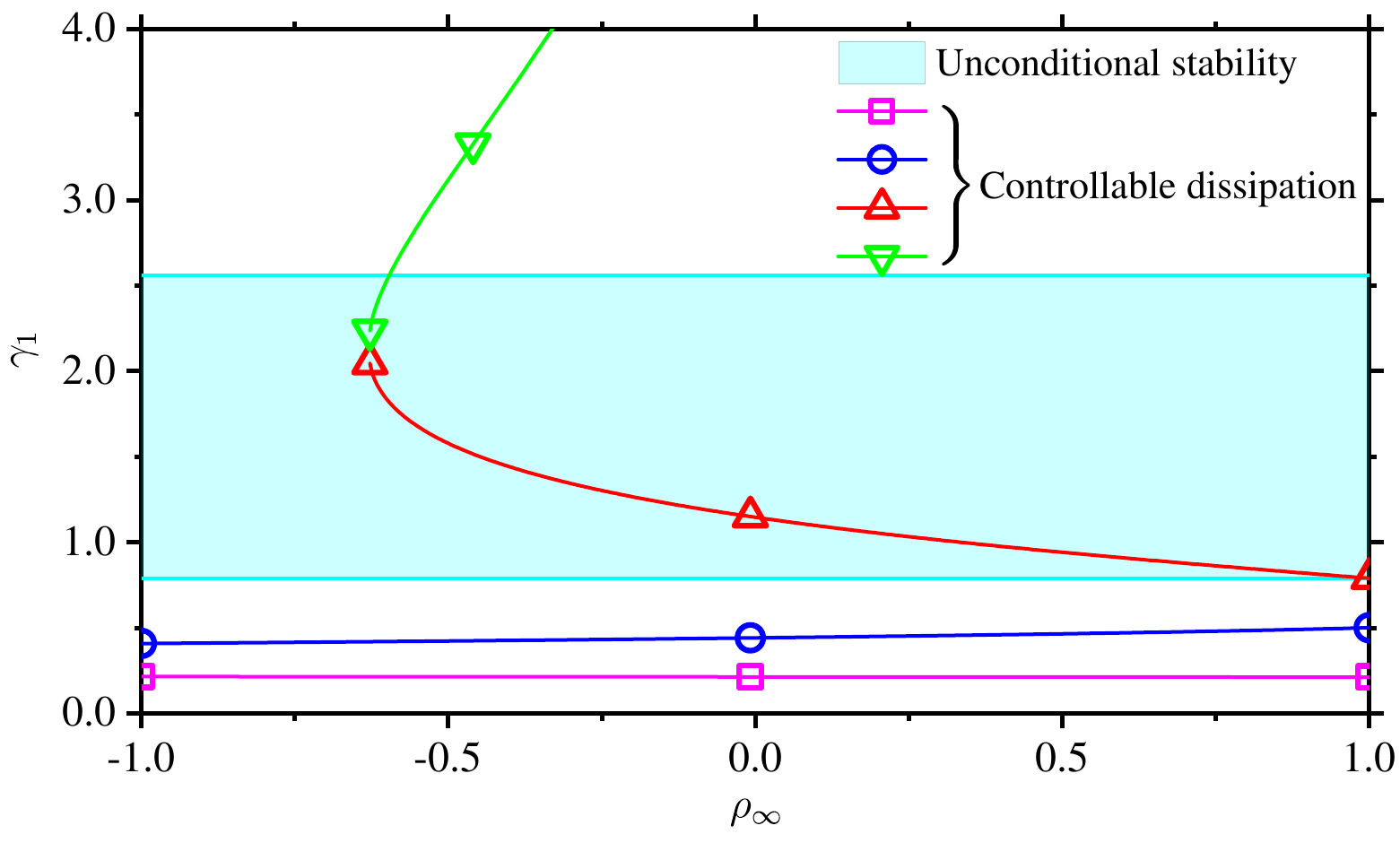}}
	\subfigure[SUCI5 ]{
		\includegraphics[scale=0.32]{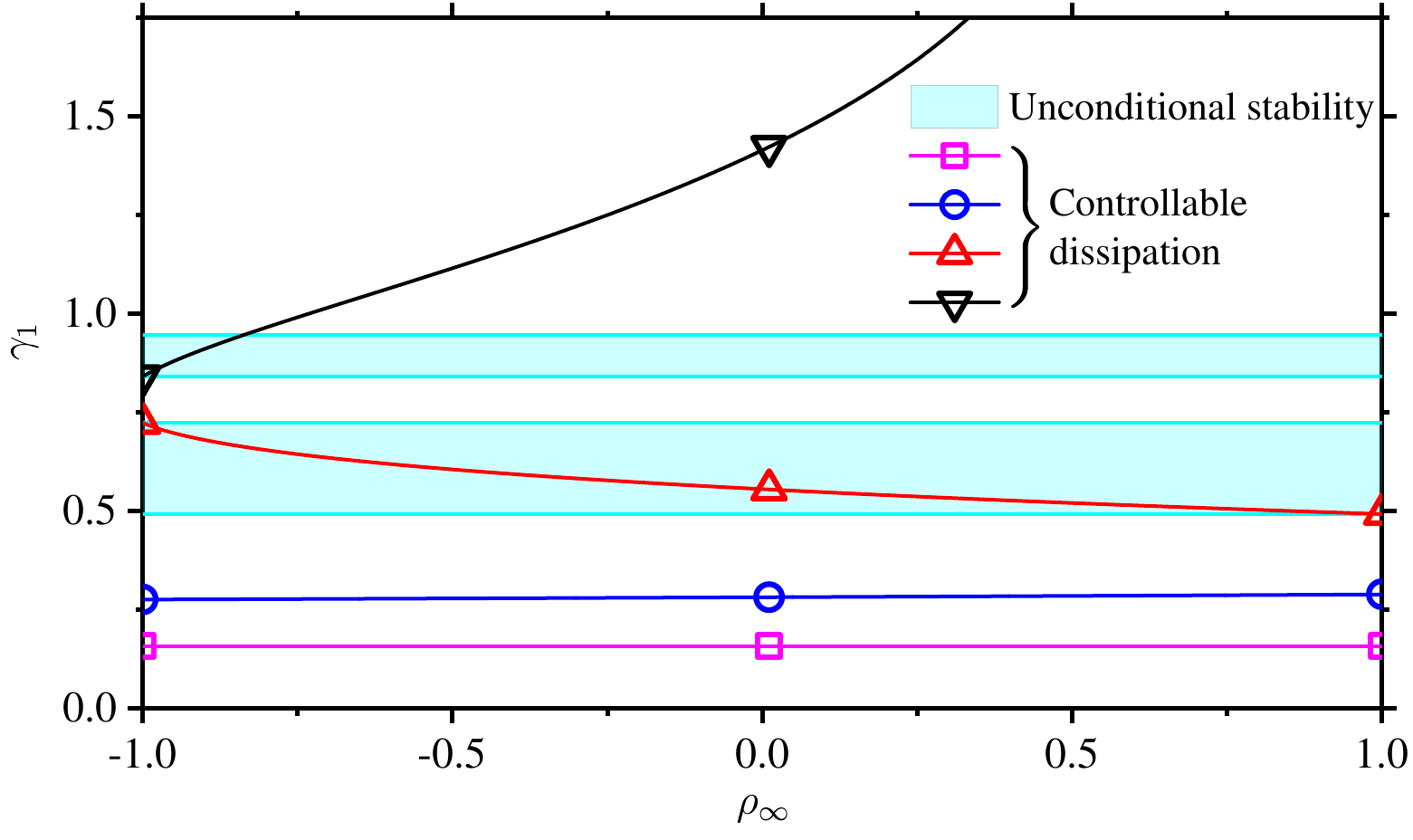}}
	\subfigure[SUCI6 ]{
			\includegraphics[scale=0.32]{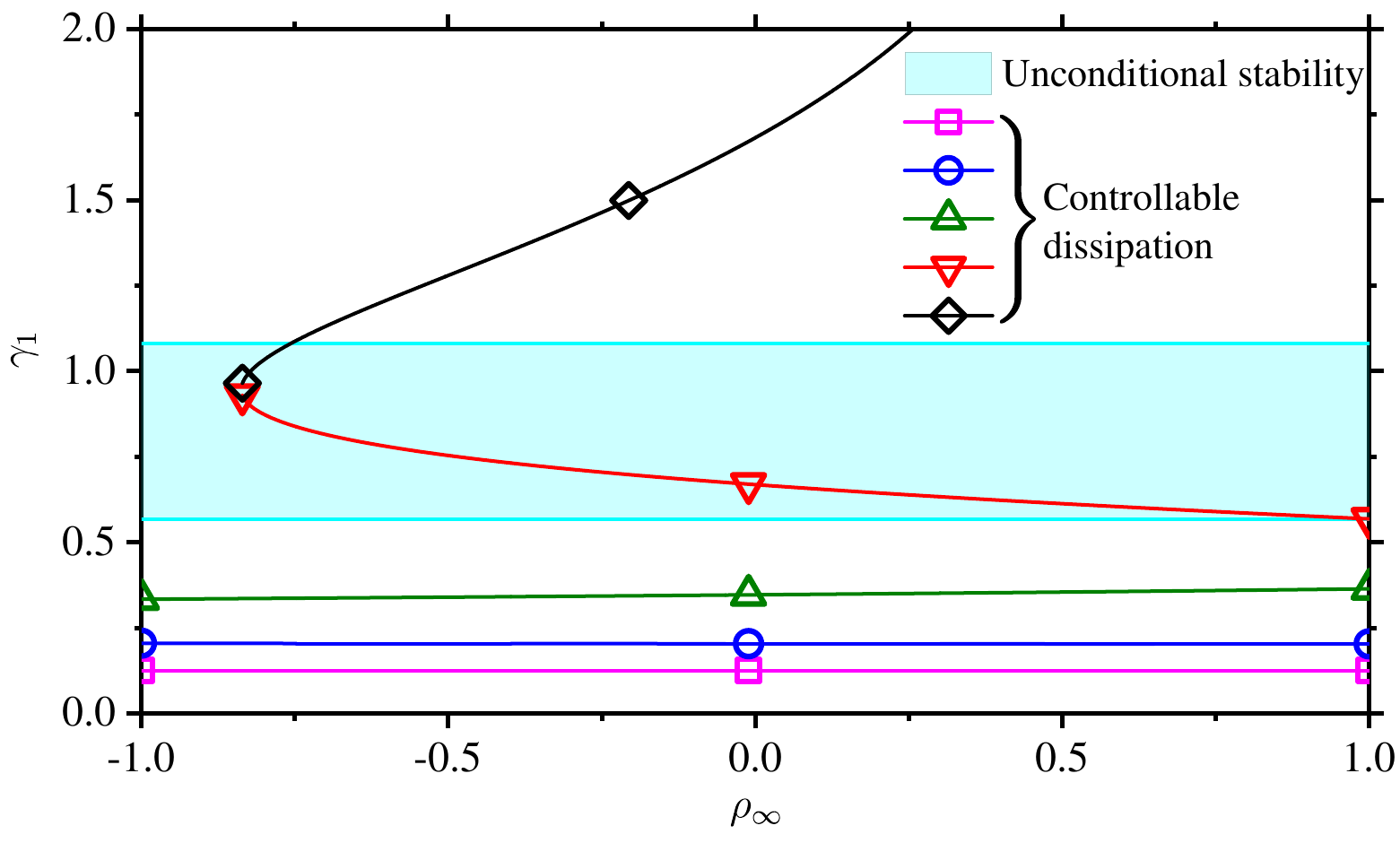}}
	\caption{Sections of $\gamma_1$ for SUCI$n$ based on unconditional stability and dissipation analyses.}
	\label{fig:sucin_g1}
\end{figure}

%Now, the present three-sub-step scheme (\ref{eq:3substep}) leaves $ \gamma_1 $ and $ \gamma_2 $ to be determined. The controllable numerical dissipation is achieved to determine $ \gamma_1 $. 
After complicated computations, the characteristic polynomial (\ref{eq:cp1}) of the numerical amplification matrix $ \mbf{D}_{num} $ can be simplified in the high-frequency limit ($ \omega\to\infty $) as
\begin{equation}\label{key}
	\left(\zeta_{\infty}-\frac{3\gamma_1^3-18\gamma_1^2+18\gamma_1-4}{3\gamma_1^3}\right)^2=0.
\end{equation}
Then, the conditions (\ref{eq:optimaldissipation}) can give
\begin{equation}\label{eq:3subgamm1}
	\frac{3\gamma_1^3-18\gamma_1^2+18\gamma_1-4}{3\gamma_1^3}=\rhoinf.
\end{equation}
It is difficult to obtain the analytical expression of $ \gamma_1 $ in terms of $ \rhoinf $ by solving Eq.~(\ref{eq:3subgamm1}). For each given $\rhoinf\in\left[-1,~1\right]$, numerically solving Eq.~\eqref{eq:3subgamm1} yields three roots of $\gamma_1$ which are plotted in Fig.~\ref{fig:sucin_g1}(a). It is reasonable that the values of $\gamma_1$ should be selected to achieve unconditional stability and controllable numerical dissipation. Hence, the circled red curve in Fig.~\ref{fig:sucin_g1}(a) should be used and the corresponding values of $\gamma_1$ in $\rhoinf\in[0,~1]$ are recorded in Table \ref{tab:gamma1}.

\begin{remark}
	The parameter $ \gamma_2 $ is free except for $ \gamma_2\neq\gamma_1 $ and it has no influence on linear numerical properties. In this paper, $ \gamma_2=\left(3+\sqrt{3}\right)\gamma_1/3 $ is used to achieve fourth-order consistency in time $ t_n+\gamma_2\dt $. It is necessary to stress that the proposed SUCI3 method is more generalized than MSSTH3 \cite{zhangOptimizationNsubstep2020} since the former does not preinstall $ \gamma_2=2\cdot\gamma_1 $ to make the first two sub-steps identical. For example, the $ \gamma_2=(1+\gamma_1)/2 $ can also be set in SUCI3 to make the last two sub-steps identical.
\end{remark}

%When the parameter $ \gamma_1 $ is determined to control high-frequency dissipation (see Table \ref{tab:gamma1}), the quantity $ A_2-1 $ is always less than zero for all $ \Omega\in[0,~\infty) $, thus validating unconditional stability of SUCI3. Note that $ A_2 -1$ cannot reach zero for all $ \Omega\in[0,~\infty) $ so that SUCI3 with $ \rhoinf=1 (\gamma_1=2/3) $ is not completely non-dissipative. As shown later, SUCI3 with $ \rhoinf=1 $ inevitably imposes some dissipation in the middle-frequency range. This is the case for the published high-order implicit methods \cite{liDirectlySelfstarting2022,zhangOptimizationNsubstep2020}. 

\begin{table}[htpb]
	\centering
	\caption{The parameter $ \gamma_1 $ of the high-order accurate $ s $-sub-step algorithms for each given $ \rhoinf\in[0,~1] $.}\vskip 1mm
	\small\begin{tabular}{ccccc}
		\toprule
		$ \boldsymbol{\rhoinf} $ & \textbf{SUCI3}      & \textbf{SUCI4}     & \textbf{SUCI5}     & \textbf{SUCI6}      \\
		\midrule
		
		0.0         & 0.8717330430 & 1.1456321252  & 0.5561076823 & 0.6682847341 \\ 
		
		0.1         & 0.8429736308 & 1.0967332903  & 0.5482826121 & 0.6557502542 \\ 
		
		0.2         & 0.8170015790 & 1.0527729141  & 0.5409197735 & 0.6440471963 \\ 
		
		0.3         & 0.7932944182 & 1.0126602385  & 0.5339560879 & 0.6330349995 \\ 
		
		0.4         & 0.7714620009 & 0.9755949496 & 0.5273404634 & 0.6226034838 \\ 
		
		0.5         & 0.7512044500 & 0.9409611552 & 0.5210308332 & 0.6126639724 \\ 
		
		0.6         & 0.7322856202 & 0.9082615701 & 0.5149920597 & 0.6031433531 \\ 
		
		0.7         & 0.7145156239 & 0.8770723798 & 0.5091944163 & 0.5939799400 \\ 
		
		0.8         & 0.6977389062 & 0.8470075321 & 0.5036124624 & 0.5851204729 \\ 
		
		0.9         & 0.6818258455 & 0.8176837322 & 0.4982241931 & 0.5765178426 \\ 
		
		1.0         & 0.6666666666 & 0.7886751346 & 0.4930103863 & 0.5681292760 \\ 
		\bottomrule
	\end{tabular}\label{tab:gamma1}
\end{table}

	\subsubsection{Zero-order overshoots}
	Although the overshoot analysis does not provide additional conditions for SUCI$n$, it is still necessary to confirm zero-order overshoots of SUCI$n$. When using SUCI3 to solve the SDOF system \eqref{eq:sdof} with $f(t)=0$, one can get numerical solutions at the first step as
	\begin{equation}\label{eq:re1:uv1}
		\begin{bmatrix}
			u_1\\ \dot{u}_1
		\end{bmatrix}=\mbf{D}_{num}\cdot\begin{bmatrix}
			u_0\\ \dot{u}_0
		\end{bmatrix}\quad\Longrightarrow\quad \begin{cases}
			u_1=d_{11}u_0+d_{12}\dot{u}_0\\
			\dot{u}_1=d_{21}u_0+d_{22}\dot{u}_0
		\end{cases}
	\end{equation} 
	where four coefficients $d_{ij}$ are not explicitly written herein for brevity. Further, Eq.~\eqref{eq:re1:uv1} is calculated in the limit of $\dt\to\infty$ as
	\begin{equation}\label{eq:re1:uv1inf}
		\begin{cases}
			u_1^\infty=\dfrac{3\gamma_1^3-18\gamma_1^2+18\gamma_1-4}{3\gamma_1^3}u_0=\rhoinf u_0\\[3mm]
			\dot{u}_1^\infty=\dfrac{3\gamma_1^3-18\gamma_1^2+18\gamma_1-4}{3\gamma_1^3}\dot{u}_0=\rhoinf\dot{u}_0.
		\end{cases}
	\end{equation}	
	Note that the equation above has used the condition \eqref{eq:3subgamm1} to eliminate the parameter $\gamma_1$. It follows from Eq.~\eqref{eq:re1:uv1inf} that numerical solutions at the first step do not tend to infinity in the limit of $\dt\to\infty$. In other words, SUCI3 does not suffer from overshoots at the first step. Similarly, numerical solutions at subsequent time steps do not also tend to infinity in the limit of $\dt\to\infty$. Hence, it is concluded that SUCI3 achieves zero-order overshoots in displacement and velocity. 
	
	\begin{remark}
		As with SUCI3, other SUCI$n$ algorithms developed in this paper do not suffer from overshoots in either displacement or velocity, and one can easily confirm these facts by using the same analysis as above. Hence, the overshooting analyses for other SUCI$n$ algorithms are omitted for brevity. 
	\end{remark}

\subsection{Four-sub-step fourth-order scheme: SUCI4}
Following the same development path as SUCI3, the four-sub-step implicit member can be developed to achieve fourth-order accuracy, controllable numerical dissipation, and unconditional stability. However, the detailed analysis is omitted to save the length of this paper, but some important results are given herein.

The four-sub-step scheme is explicitly described in the Butcher tableau as 
\begin{equation}\label{eq:4substep}
	\begin{BMAT}(b){c|c}{c|c}
		\mbf{c} & \mbf{A}\\ %\hline
		& \mbf{b}
	\end{BMAT}=\begin{BMAT}(@,25pt,10pt){c|ccccc}{ccccc|c}
		0 & 0 &&&&\\
		\gamma_1 & \dfrac{\gamma_1}{2} & \dfrac{\gamma_1}{2}&&&\\
		\gamma_2 & \alpha_{20} & \alpha_{21} & \dfrac{\gamma_1}{2}&&\\
		\gamma_3 & \alpha_{30} & \alpha_{31} & \alpha_{32} & \dfrac{\gamma_1}{2}&\\
		1 & \alpha_{40} & \alpha_{41} & \alpha_{42} & \alpha_{43} & \dfrac{\gamma_1}{2}\\
		& \alpha_{40} & \alpha_{41} & \alpha_{42} & \alpha_{43} & \dfrac{\gamma_1}{2}\\
	\end{BMAT}
\end{equation}
Obviously, the above scheme has shared identical effective stiffness matrices within each sub-step. There are twelve algorithmic parameters to be determined. The accuracy conditions given by Eqs.~(\ref{eq:2orderi}) and (\ref{eq:accuracy}) can make all $ \alpha_{ij} $ as functions of $ \gamma_i $.
\begin{subequations}\label{eq:alpha4sub1}
	\begin{align}
		\alpha_{20} &=\frac{-\gamma_1^2+3\gamma_1\gamma_2-\gamma_2^2}{2\gamma_1}&                                                     \alpha_{21}&=\frac{\gamma_2(\gamma_2-\gamma_1)}{2\gamma_1}                           \\
		\alpha_{30} &=\frac{-\gamma_1^2+(3\gamma_3-2\alpha_{32})\gamma_1+2\alpha_{32}\gamma_2-\gamma_3^2}{2\gamma_1}&                  \alpha_{31}&=\frac{-2\alpha_{32}\gamma_2-\gamma_1\gamma_3+\gamma_3^2}{2\gamma_1}     \\
		\alpha_{40} &=\frac{-\gamma_1^2+(3-2\alpha_{42}-2\alpha_{43})\gamma_1+2\alpha_{42}\gamma_2+2\alpha_{43}\gamma_3-1}{2\gamma_1}& \alpha_{41}&=\frac{-2\alpha_{42}\gamma_2-2\alpha_{43}\gamma_3-\gamma_1+1}{2\gamma_1}\\
		\alpha_{32}&=\frac{-3\gamma_1^3+9\gamma_1^2-6\gamma_1+1}{12\alpha_{43}\gamma_2(\gamma_2-\gamma_1)}& \alpha_{42}&=\frac{6\alpha_{43}\gamma_1\gamma_3-6\alpha_{43}\gamma_3^2+3\gamma_1^2-6\gamma_1+2}{6\gamma_2(\gamma_2-\gamma_1)}\\
		\alpha_{43}&=\frac{6(1-\gamma_2)\gamma_1^2+12\gamma_1\gamma_2-10\gamma_1-4\gamma_2+3}{12\gamma_3(\gamma_3-\gamma_2)(\gamma_3-\gamma_1)}
	\end{align}
\end{subequations}
%Like the previous SUCI3 scheme, the fourth-order SUCI4 leaves three splitting ratios of sub-step size, namely $ \gamma_1,~\gamma_2 $ and $ \gamma_3 $, to be determined. 
Next, the controllable numerical dissipation in the high-frequency range is imposed to determine $ \gamma_1 $. In the high-frequency limit ($ \omega\to\infty $), the resulting characteristic polynomial (\ref{eq:cp1}) is further simplified as
\begin{equation}\label{key}
	\left(\zeta_{\infty}-\frac{3\gamma_1^4-24\gamma_1^3+36\gamma_1^2-16\gamma_1+2}{3\gamma_1^4}\right)^2=0.
\end{equation}
Then, the conditions (\ref{eq:optimaldissipation}) for achieving controllable numerical dissipation give
\begin{equation}\label{eq:suci4_d}
	\frac{3\gamma_1^4-24\gamma_1^3+36\gamma_1^2-16\gamma_1+2}{3\gamma_1^4}=\rhoinf.
\end{equation}
On the other hand, the condition for achieving unconditional stability is given in the stability analysis as 
\begin{equation}
	\gamma_1\in\left[0.7886751346,~2.561159523\right].
\end{equation}
Numerically solving Eq.~\eqref{eq:suci4_d} gives four values of $ \gamma_1 $ for each given $ \rhoinf\in[0,~1] $, as shown in Fig.~\ref{fig:sucin_g1}(b). Fig.~\ref{fig:sucin_g1}(b) reveals that, among the four values of $\gamma_1$, precisely one falls within the unconditionally stable region and it is recorded in Table \ref{tab:gamma1} for each $\rhoinf$.

\begin{remark}
	The sub-step splitting ratios $\gamma_2$ and $\gamma_3$ are free parameters except for $\gamma_3\neq\gamma_2\neq\gamma_1$. In this work, $ \gamma_2=2\cdot\gamma_1 $ and $ \gamma_3=3\cdot\gamma_1 $ are adopted.
\end{remark}
%which are recorded in Table \ref{tab:gamma1}. Besides, 

\subsection{Five-sub-step fifth-order scheme: SUCI5}
Analogously to the previous SUCI3 and SUCI4 schemes, the five-sub-step member is developed by using the theoretical results in Section \ref{sec:development} to achieve fifth-order accuracy, dissipation control, and unconditional stability. The developed five-sub-step scheme are explicitly given in the Butcher tableau as
\begin{equation}\label{eq:5substep}
	\begin{BMAT}(b){c|c}{c|c}
		\mbf{c} & \mbf{A}\\ %\hline
		& \mbf{b}
	\end{BMAT}=\begin{BMAT}(@,25pt,10pt){c|cccccc}{cccccc|c}
		0 & 0 &&&&&\\
		\gamma_1 & \dfrac{\gamma_1}{2} & \dfrac{\gamma_1}{2}&&&&\\
		\gamma_2 & \alpha_{20} & \alpha_{21} & \dfrac{\gamma_1}{2}&&&\\
		\gamma_3 & \alpha_{30} & \alpha_{31} & \alpha_{32} & \dfrac{\gamma_1}{2}&&\\
		\gamma_4 & \alpha_{40} & \alpha_{41} & \alpha_{42} & \alpha_{43} & \dfrac{\gamma_1}{2}&\\
		1 & \alpha_{50} & \alpha_{51} & \alpha_{52} & \alpha_{53} & \alpha_{54}& \dfrac{\gamma_1}{2}\\
		&  \alpha_{50} & \alpha_{51} & \alpha_{52} & \alpha_{53} & \alpha_{54}& \dfrac{\gamma_1}{2}
	\end{BMAT}
\end{equation}
%\begin{subequations}\label{eq:5substep}
%	\begin{align}
%		\mbf{M}\ddot{\mbf{U}}_{n+\gamma_4}&+\mbf{C}\dot{\mbf{U}}_{n+\gamma_4}+\mbf{K}\mbf{U}_{n+\gamma_4}  =\mbf{F}(t_n+\gamma_4\dt)   & \mbf{M}\ddot{\mbf{U}}_{n+1}&+\mbf{C}\dot{\mbf{U}}_{n+1}+\mbf{K}\mbf{U}_{n+1}                       =\mbf{F}(t_{n+1})                                                                                                                                                                                                                             \\
%		\mbf{U}_{n+\gamma_4}                                                                             & =\mbf{U}_n+\dt\left(\sum_{j=0}^3\alpha_{4j}\dot{\mbf{U}}_{n+\gamma_j}+\frac{\gamma_1}{2}\dot{\mbf{U}}_{n+\gamma_4}\right)       & \mbf{U}_{n+1}                                                                                    & =\mbf{U}_n+\dt\left(\sum_{j=0}^4\alpha_{5j}\dot{\mbf{U}}_{n+\gamma_j}+\frac{\gamma_1}{2}\dot{\mbf{U}}_{n+1}\right)                  \\
%		\dot{\mbf{U}}_{n+\gamma_4}                                                                       & =\dot{\mbf{U}}_n+\dt\left(\sum_{j=0}^3\alpha_{4j}\ddot{\mbf{U}}_{n+\gamma_j}+\frac{\gamma_1}{2}\ddot{\mbf{U}}_{n+\gamma_4}\right) & \dot{\mbf{U}}_{n+1}                                                                              & =\dot{\mbf{U}}_n+\dt\left(\sum_{j=0}^4\alpha_{5j}\ddot{\mbf{U}}_{n+\gamma_j}+\frac{\gamma_1}{2}\ddot{\mbf{U}}_{n+1}\right).
%	\end{align}
%\end{subequations}
Again, the identity of effective stiffness matrices within each sub-step has been achieved by requiring $ \alpha_{ii}=\gamma_1/2~(i=1,~2,~\cdots,~5) $.

There are nineteen algorithmic parameters to be determined. The conditions given by Eqs.~(\ref{eq:accuracy}) and (\ref{eq:2orderi}) for achieving fifth-order accuracy are used to determine $ \alpha_{ij} $, which are
\begin{subequations}\label{eq:alpha5sub1}
	%	\begin{equation}
		\begin{align}
			\alpha_{20}&=\frac{-\gamma_1^2+3\gamma_1\gamma_2-\gamma_2^2}{2\gamma_1}                          \hspace{5.65cm} \alpha_{21}=\frac{\gamma_2(\gamma_2-\gamma_1)}{2\gamma_1}                                            \\
			\alpha_{30}&=\frac{-\gamma_1^2+(3\gamma_3-2\alpha_{32})\gamma_1+2\alpha_{32}\gamma_2-\gamma_3^2}{2\gamma_1}                                    \hspace{3.1cm} \alpha_{31}=\frac{-2\alpha_{32}\gamma_2-\gamma_1\gamma_3+\gamma_3^2}{2\gamma_1}                      \\
			\alpha_{40}&=\frac{-\gamma_1^2+(3\gamma_4-2\alpha_{42}-2\alpha_{43})\gamma_1+2\alpha_{42}\gamma_2+2\alpha_{43}\gamma_3-\gamma_4^2}{2\gamma_1}  \qquad\ \  \alpha_{41}=\frac{-2\alpha_{42}\gamma_2-2\alpha_{43}\gamma_3-\gamma_1\gamma_4+\gamma_4^2}{2\gamma_1}\\
			\alpha_{50} & =\frac{-\gamma_1^2+(3-2\alpha_{52}-2\alpha_{53}-2\alpha_{54})\gamma_1+2\alpha_{52}\gamma_2+2\alpha_{53}\gamma_3+2\alpha_{54}\gamma_4-1}{2\gamma_1}     \\
			\alpha_{51} & =\frac{-2\alpha_{52}\gamma_2-2\alpha_{53}\gamma_3-2\alpha_{54}\gamma_4-\gamma_1+1}{2\gamma_1}  \hspace{3cm}
			\alpha_{32}  =\frac{15\gamma_1^4-60\gamma_1^3+60\gamma_1^2-20\gamma_1+2}{120\alpha_{43}\alpha_{54}\gamma_2(\gamma_2-\gamma_1)}                                      \\ \alpha_{42}&=\frac{\left\{\begin{aligned}
					15&\alpha_{53}\gamma_1^4 +30(\alpha_{43}\alpha_{54}-2\alpha_{53})\gamma_1^3+30(2\alpha_{53}-3\alpha_{43}\alpha_{54})\gamma_1^2\\
					& +20(-6\alpha_{43}^2\alpha_{54}^2\gamma_3+3\alpha_{43}\alpha_{54}-\alpha_{53})\gamma_1+120\alpha_{43}^2\alpha_{54}^2\gamma_3^2-10\alpha_{43}\alpha_{54}+2\alpha_{53}
				\end{aligned}\right\}}{120\alpha_{43}\alpha_{54}^2\gamma_2(\gamma_1-\gamma_2)}\\
			\alpha_{52} & =\frac{3\gamma_1^2+6(\alpha_{53}\gamma_3+\alpha_{54}\gamma_4-1)\gamma_1-6\alpha_{53}\gamma_3^2-6\alpha_{54}\gamma_4^2+2}{6\gamma_2(\gamma_2-\gamma_1)}\\
			\alpha_{43} & =\frac{\gamma_4(\gamma_4-\gamma_1)(\gamma_4-\gamma_2)(\gamma_4-\gamma_3)(15\gamma_1^3\gamma_2-15\gamma_1^3-45\gamma_1^2\gamma_2+35\gamma_1^2+30\gamma_1\gamma_2-20\gamma_1-5\gamma_2+3)}
			{\gamma_3(\gamma_3-\gamma_1)(\gamma_3-\gamma_2)G(\gamma_3)} \\
			\alpha_{53} & =\frac{G(\gamma_4)}{60\gamma_3(\gamma_3-\gamma_4)(\gamma_2-\gamma_3)(\gamma_1-\gamma_3)}
			\hspace{3.5cm}
			\alpha_{54} =\frac{G(\gamma_3)}
			{60\gamma_4(\gamma_4-\gamma_3)(\gamma_4-\gamma_2)(\gamma_4-\gamma_1)}
		\end{align}
		%	\end{equation}
	where $ G(x)=30(1-x)(1-\gamma_2)\gamma_1^2+(50x-45+10(5-6x)\gamma_2)\gamma_1+5(4x-3)\gamma_2-15x+12 $.
\end{subequations}

%The parameter $ \gamma_1 $ is given to control numerical dissipation in the high-frequency range. 
In the high-frequency limit ($ \omega\to\infty $), the characteristic polynomial (\ref{eq:cp1}) is simplified using the known parameters above as
\begin{equation}
	\left(\zeta_{\infty}-\frac{15\gamma_1^5-150\gamma_1^4+300\gamma_1^3-200\gamma_1^2+50\gamma_1-4}{15\gamma_1^5}\right)^2=0.
\end{equation}
The controllable numerical dissipation is achieved by solving
\begin{equation}\label{eq:suci5_d}
	\frac{15\gamma_1^5-150\gamma_1^4+300\gamma_1^3-200\gamma_1^2+50\gamma_1-4}{15\gamma_1^5}=\rhoinf.
\end{equation}
On the other hand, the unconditional stability requires $\gamma_1$ to satisfy 
\begin{equation}
	\gamma_1\in\left[0.4930103863,~0.7236067977\right]\cup\left[0.8415650255,~0.9465367825\right].
\end{equation}
For each given $ \rhoinf\in[-1,~1] $, numerically solving Eq.~\eqref{eq:suci5_d} yields five values of $\gamma_1$ and one of them happens to fall precisely within the unconditionally stable region, as depicted in Fig.~\ref{fig:sucin_g1}(c). Table \ref{tab:gamma1} lists numerical values of $\gamma_1$ for SUCI5 with $\rhoinf\in\left[0,~1\right]$. Note that Fig.~\ref{fig:sucin_g1}(c) illustrates that the values of $\gamma_1$ corresponding to $\rhoinf\in[-1,~0)$ also makes SUCI5 unconditionally stable and controllably dissipative, but they often provide worse spectral properties, such as large period errors, than those corresponding to $\rhoinf\in\left[0,~1\right]$.
\begin{remark}
	The sub-step splitting ratios $\gamma_2$, $\gamma_3$, and $\gamma_4$ are free parameters except for $\gamma_4\neq\gamma_3\neq\gamma_2\neq\gamma_1$. In this work, $ \gamma_2=2\cdot\gamma_1 $, $ \gamma_3=3\cdot\gamma_1 $, and $\gamma_4=4\cdot\gamma_1$ are adopted.
\end{remark}

%The remaining $ \gamma_2,~\gamma_3 $, and $ \gamma_4 $ can be taken by default as $ \gamma_2=2\cdot\gamma_1,~\gamma_3=3\cdot\gamma_1 $, and $ \gamma_4=4\cdot\gamma_1 $, respectively, since they have no influence on linear numerical properties.

%\begin{remark}
Following the same development path as the three-, four-, and five-sub-step schemes, one can propose the six-sub-step scheme with sixth-order accuracy; see Appendix B for more details. Note that this paper does not pre-assume that the first ($ s-1 $) sub-steps in the general $ s $-sub-step implicit method \eqref{eq:nsubstep} share the same length, that is, $ \gamma_j=j\cdot\gamma_1~(j=2,~\cdots,~s-1)$, so the developed high-order algorithms are general so that the parameters $ a_{ij} $ are a little complicated. As a result, users can make the high-order members possess the same length either in the first $ (s-1) $ sub-steps by forcing $ \gamma_j=j\cdot\gamma_1~(j=2,~\cdots,~s-1)$ or in the last $ (s-1) $ sub-steps by forcing $ 1-\gamma_{s-1}=\gamma_{s-1}-\gamma_{s-2}=\cdots=\gamma_2-\gamma_1$. This paper adopts the former to make the first $ (s-1) $ sub-steps identical.
%\end{remark}

In summary, this section, as well as Appendix B, has developed four novel high-order implicit algorithms. For $ 2\le s\le6 $, the proposed $ s $-sub-step method (\ref{eq:nsubstep}) can achieve simultaneously $ s $th-order accuracy, unconditional stability, and a full range of dissipation control. However, this case is not true for $ s\ge 7 $. For instance, the seven-sub-step implicit scheme within the framework of (\ref{eq:nsubstep}) cannot achieve seventh-order accuracy and controllable numerical dissipation since the resulting scheme is unconditionally unstable. Hence, the high-order accurate algorithms proposed in this work end up with $ s=6 $.% In particular, the proposed four high-order accurate algorithms actually remain some splitting ratios of sub-step size as free parameters, such as $ \gamma_2$ and $\gamma_3 $ of SUCI4. In general, the composite $ s $-sub-step algorithms leave $ (s-2) $ splitting ratios as free parameters, which are $ \gamma_i~(i=2,~\cdots,~s-1) $. As mentioned previously, these free parameters do not influence linear numerical properties, and they are taken by default as $ \gamma_i=i\cdot\gamma_1 $.

\section{Spectral properties}\label{sec:sp}
In this section, spectral properties of SUCI$n$ will be analyzed and compared. The percentage amplitude decay and relative period error are employed to measure numerical dissipation and dispersion of an integrator, respectively, and their mathematical derivations refer to the literature \cite{hughesFiniteElement2000}. 
%In addition, the unconditional stability will be validated by plotting spectral radii in the undamped and damped cases. 
The published second-order two-sub-step implicit scheme (SUCI2) \cite{liNovelFamily2020} and the well-known TPO/G-$\alpha$ method \cite{shaoThreeParameters1988,chungTimeIntegration1993} are used herein for reference purposes.

\begin{figure}[htbp]
	\centering
	\includegraphics[scale=0.7]{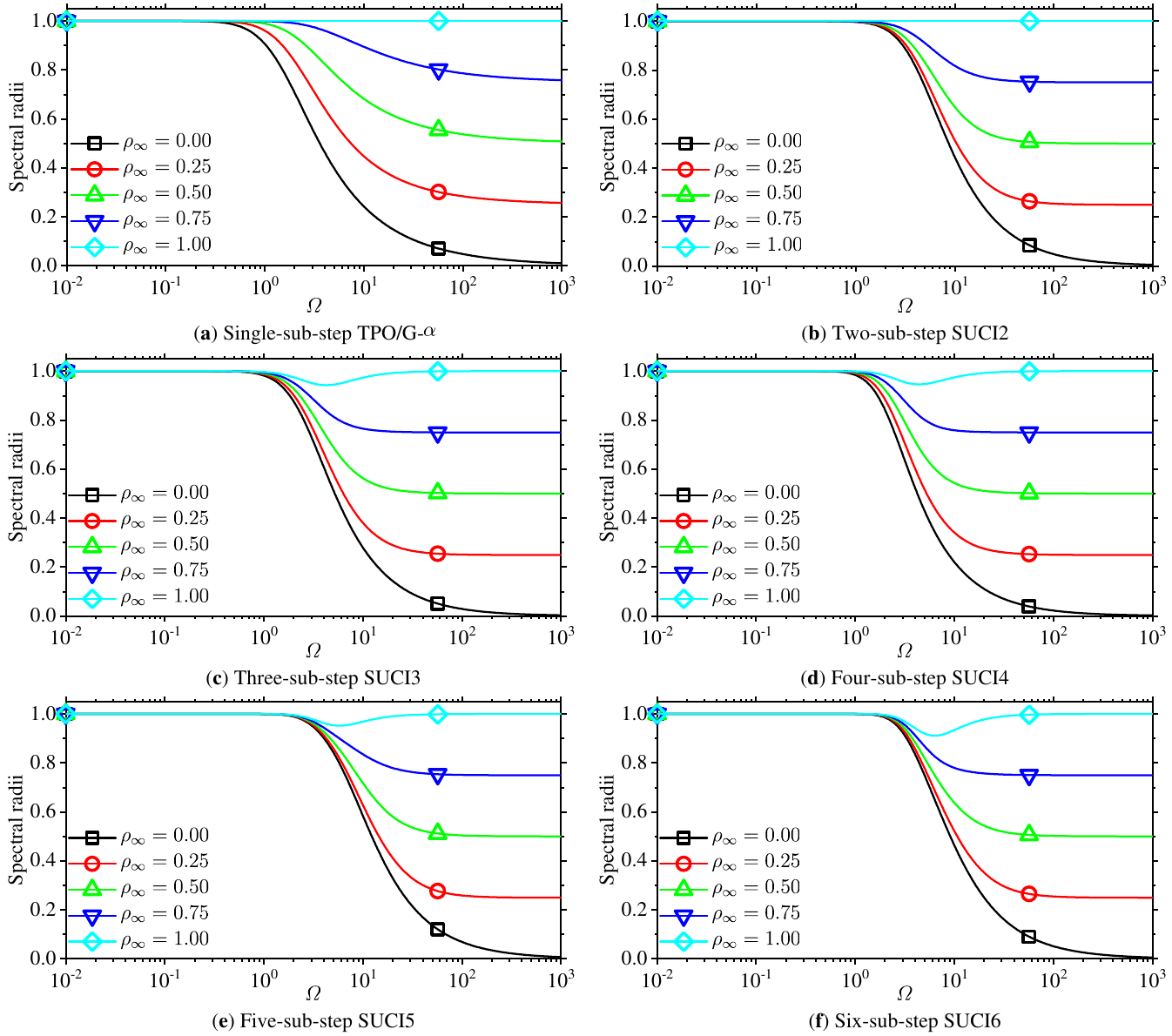}%\vskip -2mm
	\caption{Spectral radii of the TPO/G-$\alpha$ \cite{shaoThreeParameters1988,chungTimeIntegration1993}, two-sub-step SUCI2 \cite{liNovelFamily2020}, and multi-sub-step SUCI$n$ algorithms in the absence of $ \xi $.}
	\label{fig:spSUCI2}
\end{figure}

\begin{figure}[htbp]
	\centering
	\includegraphics[scale=0.7]{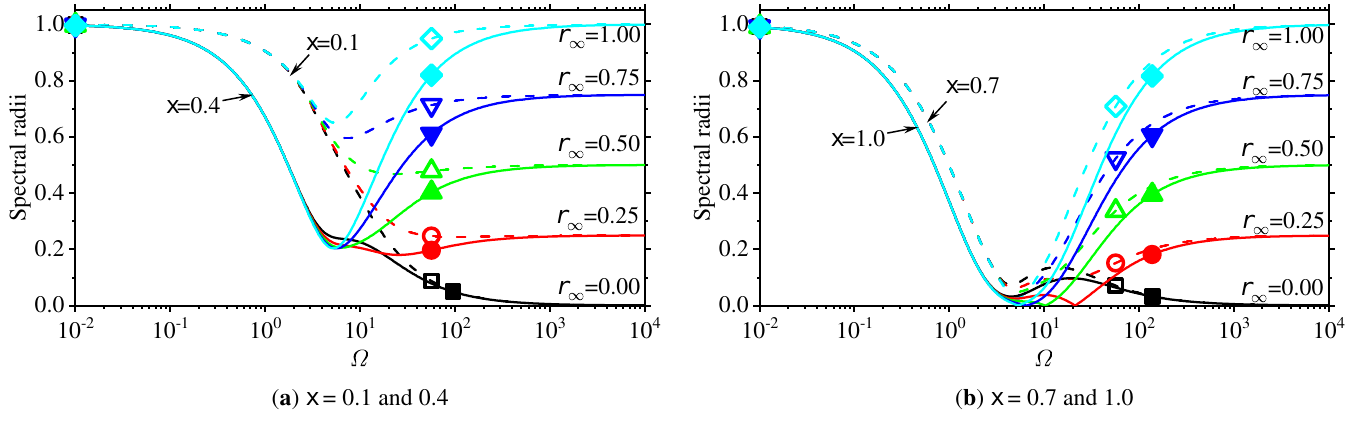}\vskip -2mm
	\caption{Spectral radii of six-sub-step SUCI6 when considering four cases of $ \xi $: (\textbf{a}) $ \xi=0.1 $ and $ 0.5 $ and (\textbf{b}) $ \xi=0.7 $ and $ 1.0 $.}
	\label{fig:spxi}
\end{figure}

The subplots (a) and (b) in Fig.~\ref{fig:spSUCI2} first plot spectral radii of TPO/G-$\alpha$ \cite{shaoThreeParameters1988,chungTimeIntegration1993} and SUCI2 \cite{liNovelFamily2020}, respectively, in the absence $ \xi $. It can be seen that the parameter $ \rhoinf $ denoting the spectral radius in the high-frequency limit exactly controls spectral radii in the high-frequency range and the two-sub-step SUCI2 scheme imposes less dissipation in the low-frequency range. Next, spectral radii of four high-order integration algorithms in the case of $ \xi=0 $ are plotted in the subplots (c-f) of Fig.~\ref{fig:spSUCI2}, where four high-order algorithms achieve controllable spectral radii in the high-frequency range by adjusting $ \rhoinf $, thus controlling numerical high-frequency dissipation. Note that in the non-dissipative case ($ \rhoinf=1 $), four high-order algorithms impose some dissipation in the middle-frequency range (about $ \Omega\in[2,10] $) since their spectral radii are mildly less than unity. This phenomenon has also been observed for other high-order algorithms \cite{liNovelFamily2019,zhangOptimizationNsubstep2020,liDirectlySelfstarting2022}. The subplots (c-f) in Fig.~\ref{fig:spSUCI2} clearly demonstrate that the four novel high-order algorithms achieve unconditional stability in the undamped case. When considering damped cases, the novel algorithms are still unconditionally stable, as shown in Fig.~\ref{fig:spxi} where spectral radii of SUCI6 are presented in four damped cases. Fig.~\ref{fig:spxi} also illustrates that the viscous damping ratio $ \xi $ has an influence on spectral radii only in the middle-frequency range, and with the increase of $ \xi $, more dissipation is imposed in the middle-frequency range instead of the high-frequency range.

\begin{figure}[H]
	\centering
	\includegraphics[scale=0.7]{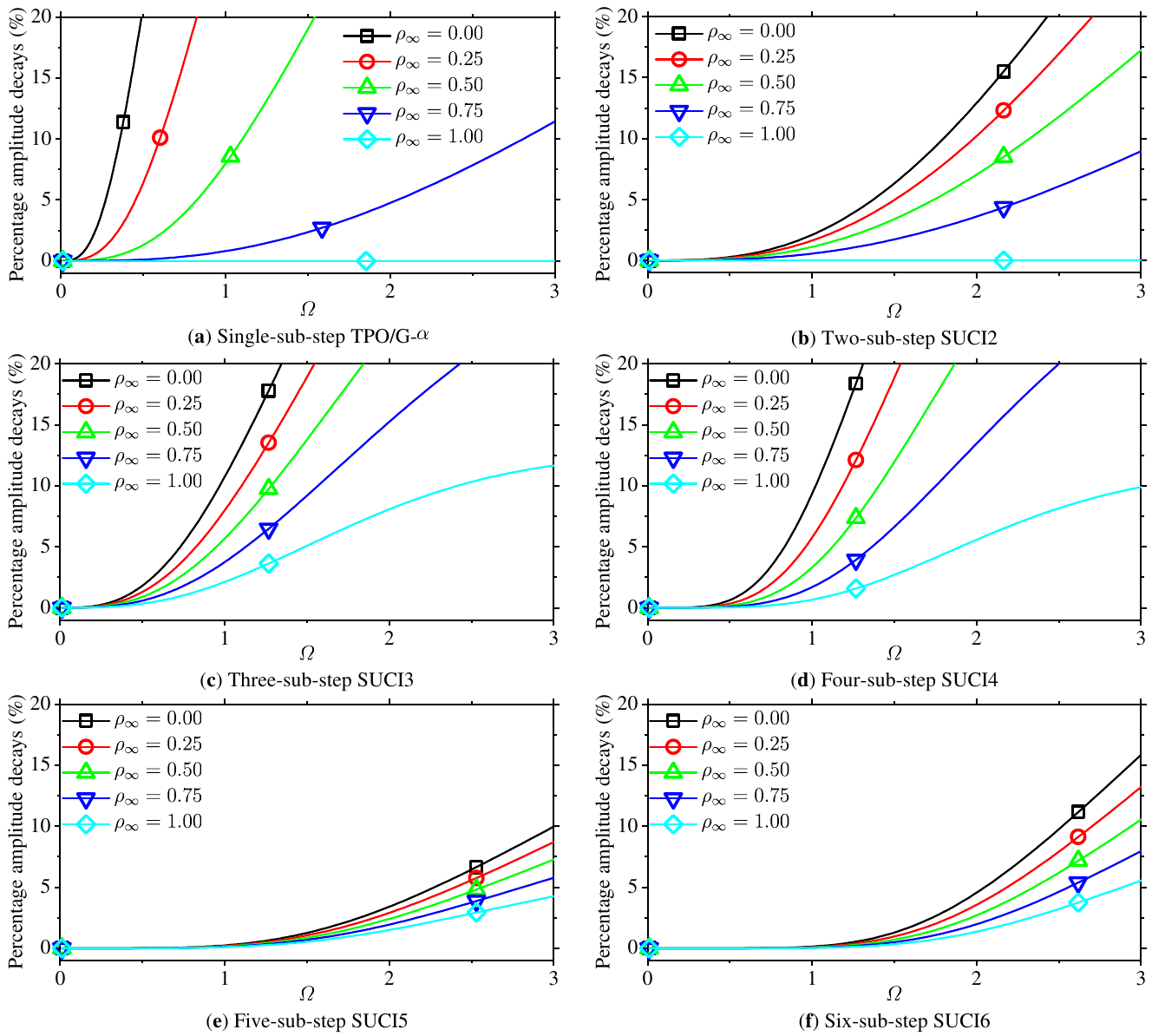}\vskip -2mm
	\caption{Percentage amplitude decays of the TPO/G-$\alpha$ \cite{shaoThreeParameters1988,chungTimeIntegration1993}, two-sub-step SUCI2 \cite{liNovelFamily2020}, and multi-sub-step SUCI$n$ algorithms in the absence of $ \xi $.}
	\label{fig:adSUCI2}
\end{figure}

The percentage amplitude decays of various algorithms are plotted in Fig.~\ref{fig:adSUCI2}, where some important observations are found and collected as follows.
\begin{itemize}%[(a)]
	\item Among the algorithms shown herein, the single-step TPO/G-$\alpha$ method \cite{shaoThreeParameters1988,chungTimeIntegration1993} produces the largest amplitude decays in the low-frequency range while the fifth-order accurate SUCI5 scheme provides the smallest amplitude decays. The previous study \cite{liNovelFamily2019} has already reported that the second-order multi-sub-step implicit schemes generally produce fewer amplitude decays in the low-frequency range as the number of sub-steps increases. However, this fact is not seemly true for the higher-order sub-step algorithms since the higher-order SUCI3 and SUCI4 schemes obviously provide larger amplitude decays than the second-order SUCI2 scheme. In addition, SUCI6 also produces larger amplitude decays than SUCI5.
	
	\item Unlike the second-order TPO/G-$\alpha$ \cite{shaoThreeParameters1988,chungTimeIntegration1993} and SUCI2 \cite{liNovelFamily2020} algorithms, the novel high-order methods impose some amplitude decays in the non-dissipative ($ \rhoinf=1 $) case since their spectral radii in the middle-frequency range $ \Omega\in[2,~10] $ are less than unity, as shown in the subplots (c-f) of Fig.~\ref{fig:spSUCI2}. However, amplitude decays in $ \rhoinf=1 $ finally decrease to zero with the increase of $ \Omega $, although they are not shown herein for brevity.
	
	\item Amplitude decays of each integration algorithm increase with the decrease of $ \rhoinf $, thus controlling numerical dissipation by changing $ \rhoinf $.
\end{itemize}

\begin{figure}[htbp]
	\centering
	\includegraphics[scale=0.7]{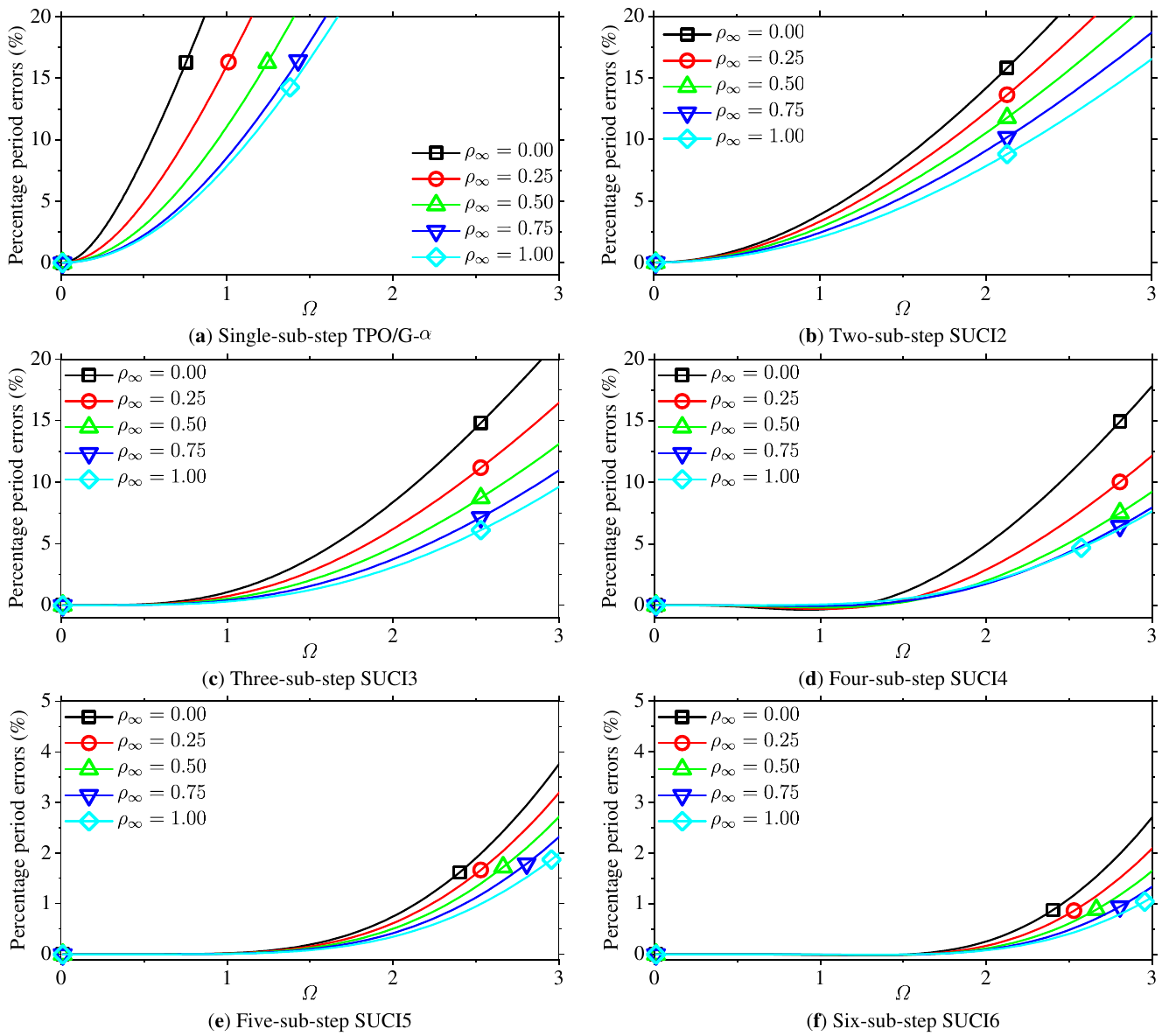}\vskip -2mm
	\caption{Percentage period errors of the TPO/G-$\alpha$ \cite{shaoThreeParameters1988,chungTimeIntegration1993}, two-sub-step SUCI2 \cite{liNovelFamily2020}, and multi-sub-step SUCI$n$ algorithms in the absence of $ \xi $.}
	\label{fig:peSUCI2}
\end{figure}

The period error actually characterizes the precision of an integrator, so a higher-order accurate algorithm should produce fewer period errors in the low-frequency range than the lower-order scheme. This fact is confirmed well for SUCI$n$. Fig.~\ref{fig:peSUCI2} shows percentage period errors of the TPO/G-$\alpha$ and SUCI$ n$ algorithms. It is pronounced that the higher-order accurate algorithm produces fewer period errors and the non-dissipative ($\rhoinf=1$) scheme provides the smallest period errors for each integration algorithm. 

%\begin{figure}[htbp]
%	\centering
%	\includegraphics[scale=0.7]{peSUCIn}
%	\caption{Percentage period errors of four novel high-order integration algorithms in the absence of $ \xi $.}
%	\label{fig:peSUCIn}
%\end{figure}

\subsection{Comparisons}
Spectral properties of SUCI$ n$ have been analyzed and compared with the common second-order schemes such as TPO/G-$\alpha$ \cite{shaoThreeParameters1988,chungTimeIntegration1993} and SUCI2 \cite{liNovelFamily2020}. Herein, the third-order dissipative EG3 \cite{fungExtrapolatedGalerkin1996} using the extrapolation technique, the third-order $ \rhoinf $-Bathe \cite{kwonSelectingLoad2021,choiTimeSplitting2022} using the four-point load selections, the fourth-order non-dissipative trapezoidal rule (TR-TS) using the Tarnow and Simo technique \cite{tarnowHowRender1994}, and the third-order trapezoidal rule (TR-CS) using the complex sub-step strategy \cite{fungComplextimestepNewmark1998,fungUnconditionallyStable1997} are compared with SUCI$ n$. Notice that EG3 can impose slight numerical dissipation via $ \beta_2=\left(\rhoinf-2-\sqrt{\rhoinf^2+2\rhoinf-2}\right)$ $/(4\rhoinf-4) $ where $ \rhoinf\in\left(\sqrt{3}-1,~1\right) $ and $ \rhoinf=3/4 $ is highly recommended in \cite{fungExtrapolatedGalerkin1996} reaching the smallest period errors, while TR-CS employs a complex-valued parameter and achieves controllable numerical dissipation. As with EG3, the third-order $ \rhoinf $-Bathe method with real-valued parameters cannot present a full range of dissipation due to $ \rhoinf\in\left(-1,~1-\sqrt{3}\right]$ and it attains the minimum period errors at $ \rhoinf=1-\sqrt{3}\approx-0.732 $. In this paper, $ \rhoinf=-3/4 $ is used for the third-order $ \rhoinf $-Bathe algorithm since it provides almost the same spectral properties as $ \rhoinf=1-\sqrt{3} $.

\begin{figure}[H]
	\centering
	\includegraphics[scale=0.7]{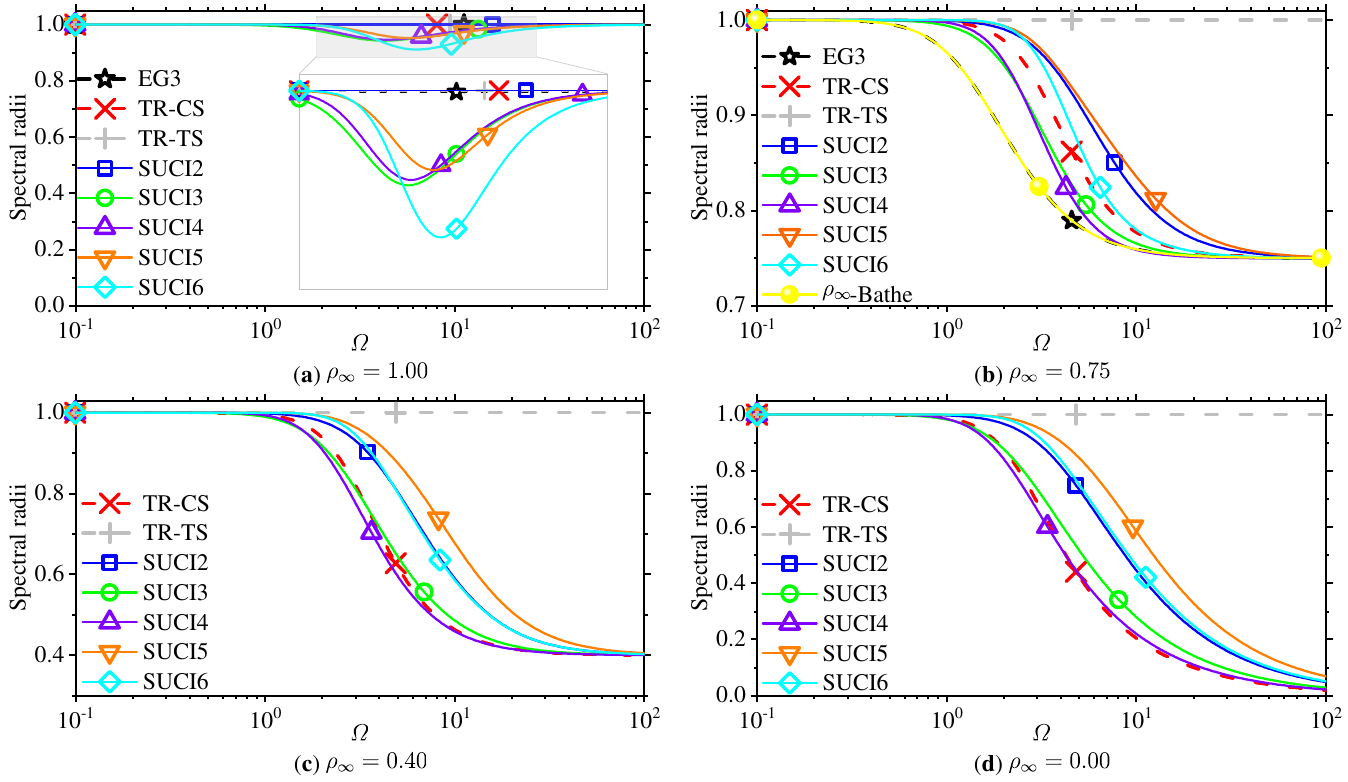}\vskip -2mm
	\caption{Comparisons of spectral radii among various algorithms when considering different $ \rhoinf $ in the absence of $ \xi $. Notice that EG3 \cite{fungExtrapolatedGalerkin1996} uses $ \rhoinf=0.999 $ in the subplot (\textbf{a}) and $ \rhoinf $-Bathe \cite{kwonSelectingLoad2021} should use the negative value of $ \rhoinf $ in the subplot (\textbf{b}).}
	\label{fig:spCom}
\end{figure}
\begin{figure}[hbtp]
	\centering
	\includegraphics[scale=0.7]{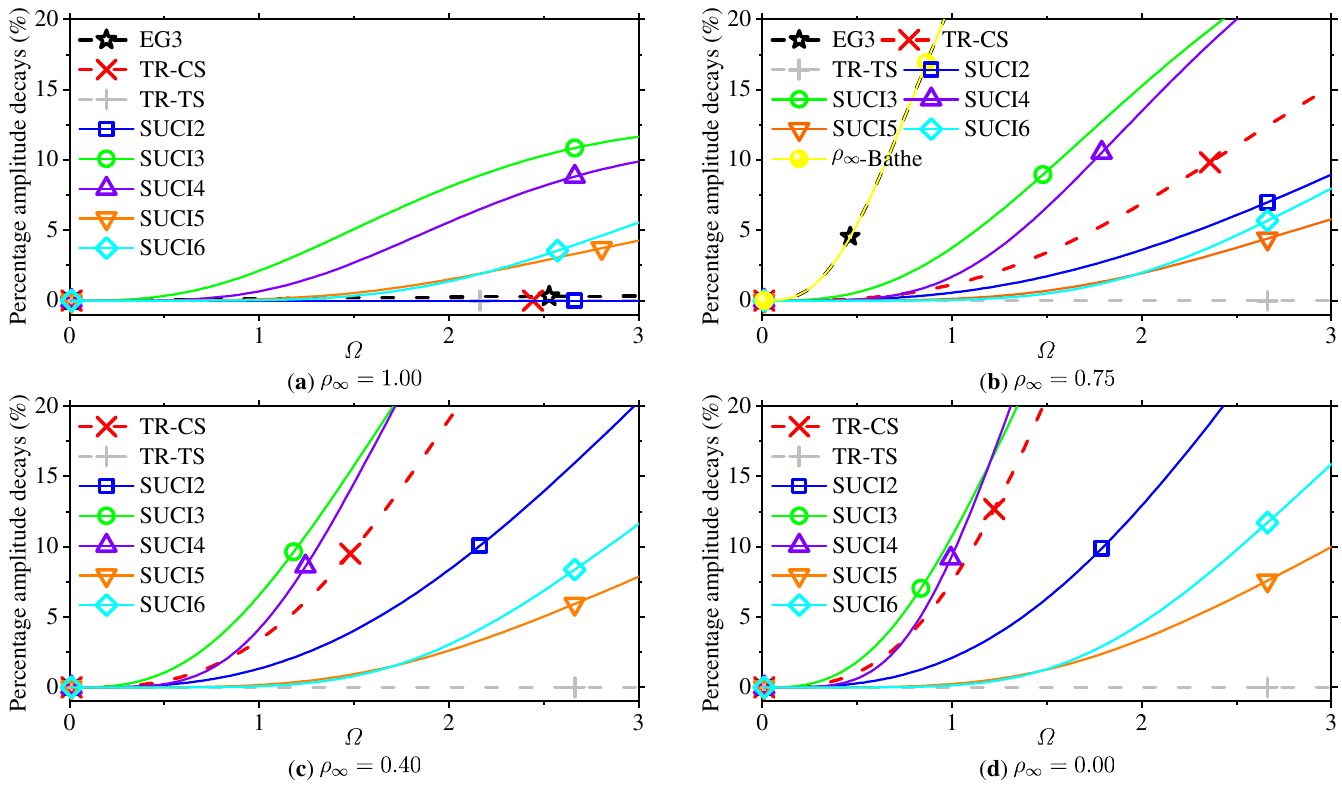}\vskip -2mm
	\caption{Comparisons of amplitude decays among various algorithms when considering different $ \rhoinf $ in the absence of $ \xi $. Notice that EG3 \cite{fungExtrapolatedGalerkin1996} uses $ \rhoinf=0.999 $ in the subplot (\textbf{a}) and $ \rhoinf $-Bathe \cite{kwonSelectingLoad2021} should use the negative value of $ \rhoinf $ in the subplot (\textbf{b}).}
	\label{fig:adCom}
\end{figure}

\begin{figure}[htbp]
	\centering
	\includegraphics[scale=0.7]{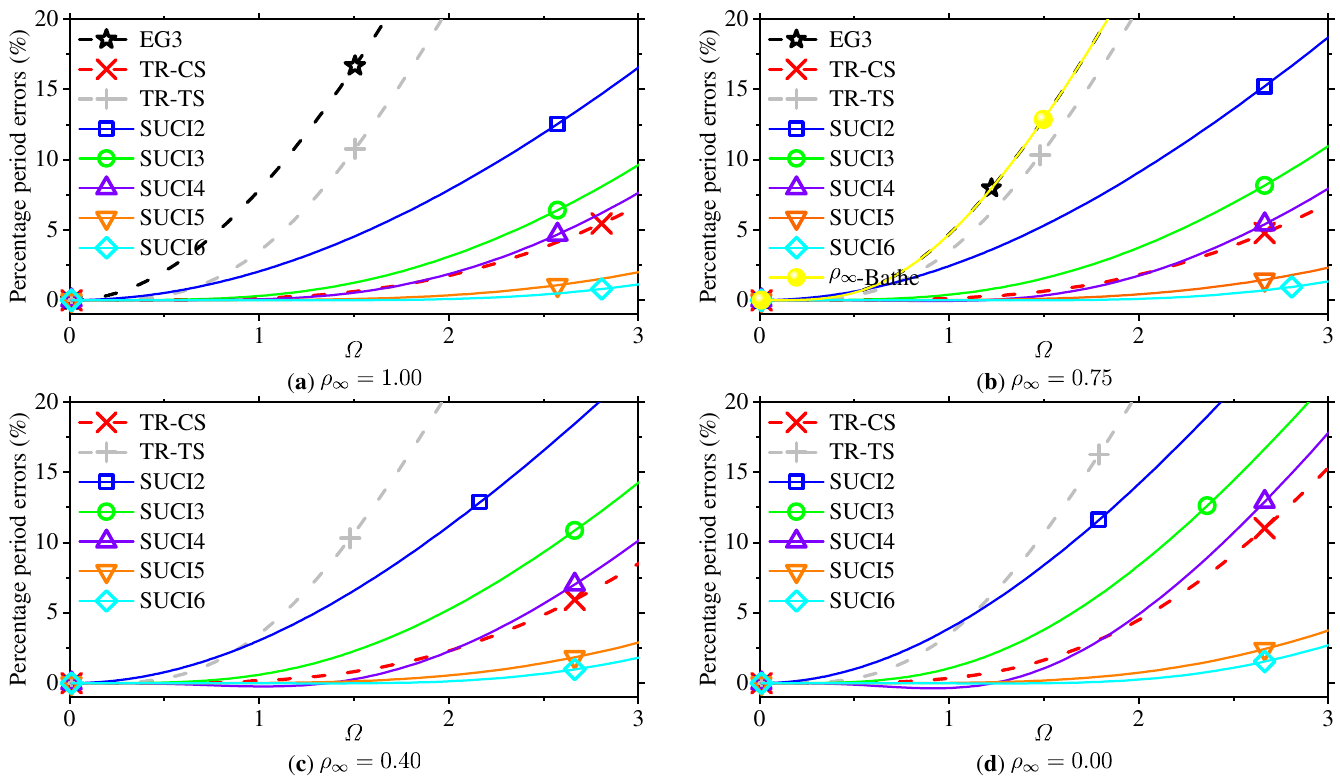}\vskip -2mm
	\caption{Comparisons of period errors among various algorithms when considering different $ \rhoinf $ in the absence of $ \xi $. Notice that EG3 \cite{fungExtrapolatedGalerkin1996} uses $ \rhoinf=0.999 $ in the subplot (\textbf{a}) and $ \rhoinf $-Bathe \cite{kwonSelectingLoad2021} should use the negative value of $ \rhoinf $ in the subplot (\textbf{b}).}
	\label{fig:peCom}
\end{figure}

Spectral radii, percentage amplitude decays, and period errors of various integration algorithms with  four different values of $ \rhoinf $ are plotted in Figs.~\ref{fig:spCom}-\ref{fig:peCom}. %Some important observations from these plots are collected as follows.
%\begin{itemize}%[(a)]
%	\item 
	In the non-dissipative ($ \rhoinf=1 $) case, the novel high-order algorithms produce some amplitude decays in the middle-frequency range, as shown in the subplots (\textbf{a}) of Figs.~\ref{fig:spCom} and \ref{fig:adCom}. In addition, the third-order EG3 scheme \cite{fungExtrapolatedGalerkin1996} only imposes slight dissipation in the high-frequency range due to $ \rhoinf\in\left(\sqrt{3}-1,~1\right) $, thus it is not included in the cases of $ \rhoinf=0.4 $ and $ 0.0 $. Similarly, the third-order $ \rhoinf $-Bathe algorithm \cite{kwonSelectingLoad2021} is compared only in the case of $ |\rhoinf|=0.75 $.
%	\item 
	Among the third-order algorithms, EG3 \cite{fungExtrapolatedGalerkin1996} and $ \rhoinf $-Bathe \cite{kwonSelectingLoad2021} provide the worst spectral behavior, such as the largest amplitude decays and period errors in the low-frequency range, while TR-CS \cite{fungComplextimestepNewmark1998,fungUnconditionallyStable1997} produces better spectral accuracy than SUCI3. However, it will be shown later that TR-CS requires more considerable computational costs than SUCI3 due to involving complex computations. In other words, TR-CS possesses better spectral accuracy at the expense of computational costs.
%	\item 
	Among the fourth-order schemes, SUCI4 shows greater advantages than TR-TS \cite{tarnowHowRender1994} since SUCI4 not only achieves controllable numerical dissipation but also produces significantly fewer period errors.
%	\item 
	Fig.~\ref{fig:peCom} illustrates that within the framework of Eq.~(\ref{eq:nsubstep}), the novel high-order schemes generally provide fewer period errors, so predicting more accurate numerical solutions.
%\end{itemize}

\begin{figure}[htbp]
	\centering
	\includegraphics[scale=0.8]{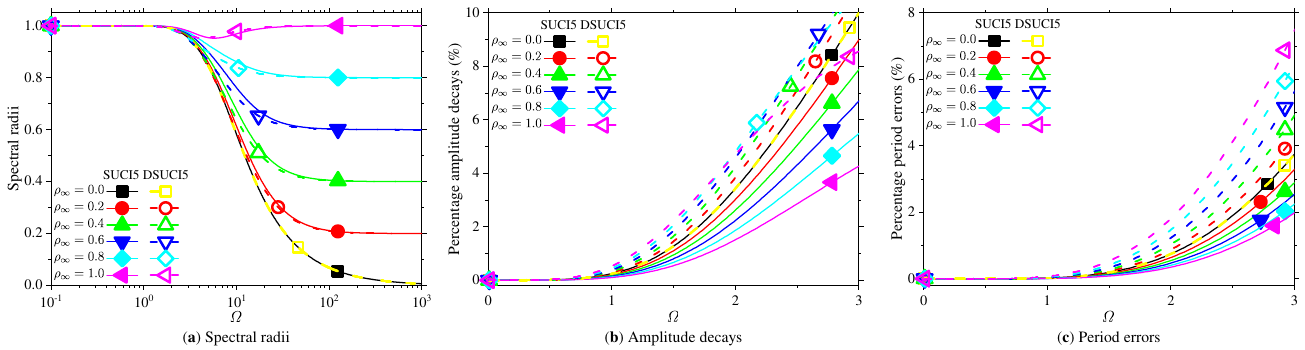}%\vskip -2mm
	\caption{Comparisons of spectral properties between DSUCI5 \cite{liDirectlySelfstarting2022} and SUCI5 in the absence of $ \xi $.}
	\label{fig:spd5}
\end{figure}

The directly self-starting high-order implicit methods \cite{liDirectlySelfstarting2022} (DSUCI$n$) share similar spectral properties to SUCI$n$ developed in this paper. In detail, except for the five-sub-step implicit schemes, other sub-step algorithms between DSUCI$n$ and SUCI$n$ present the same spectral properties and their comparisons are not thus given herein. For the five-sub-step members, SUCI5 and DSUCI5 \cite{liDirectlySelfstarting2022} provide different spectral properties and their comparisons are plotted in Fig.~\ref{fig:spd5}. As one can observe, the novel SUCI5 algorithm imposes less numerical dissipation and fewer period elongation errors than DSUCI5. Hence, the novel SUCI5 outperforms DSUCI5 in terms of spectral properties.

\section{Numerical examples}\label{sec:example}
This section will solve some linear and nonlinear examples to demonstrate that
\begin{itemize}%[(a)]
	\item SUCI$n$ achieves identical designed order of accuracy for solving general structures and does not suffer from the order reduction. For example, SUCI6 provides identical sixth-order accuracy in displacement, velocity, and acceleration for solving the damped forced vibration;
	\item SUCI$n$ embedding controllable numerical dissipation can not only effectively filter out spurious high-frequency components but also accurately integrate important low-frequency modes;
	\item SUCI$n$ possesses significant advantages over some published methods for solving nonlinear dynamics, and
	\item the computational cost of SUCI$n$ is acceptable compared with the published high-order schemes. For instance, the third-order SUCI3 algorithm requires less computational time than the third-order TR-CS scheme \cite{fungComplextimestepNewmark1998,fungUnconditionallyStable1997}.
\end{itemize}

The step-by-step solution procedure to simulate linear problem \eqref{eq:mck} using the novel $s$-sub-step implicit method \eqref{eq:nsubstep} is summarized in Algorithm \ref{code:suci_linear}. Note that the novel high-order methods only requires the calculation and decomposition of the effective stiffness matrix $\widetilde{\mbfK}$ once for solving linear problems due to achieving identical effective stiffness matrices. 
	\begin{breakablealgorithm}%[h]
		\caption{The novel $s$-sub-step implicit method \eqref{eq:nsubstep} for solving linear problems \eqref{eq:mck}}
		\begin{algorithmic}[1]
			\State Select the number of sub-steps $s\in\{1,~2,~3,~4,~5,~6\}$; set $\rhoinf\in[0,~1]$; calculate $\gamma_i$ and $\alpha_{ij}$.
			\State Solve: $\mbfa_0$ by $\mbfM\mbfa_0=\mbfF(t_0)-\mbfC\mbfv_0-\mbfK\mbfu_0$.
			\State Compute: $\widetilde{\mbfK}=\mbfM+\alpha_{ii}\dt\mbfC+{\alpha_{ii}}^2\dt^2\mbfK:=\mbf{L}\mbf{U}^\mathsf{T}$. {\color{dcolor}\qquad//$\alpha_{ii}=\gamma_1/2$ due to identical effective stiffness matrices.} 
			\For{$n=0$ to $n=N$} {\color{dcolor}\qquad\qquad\qquad\ \ \   //Loops for integration steps.}
			\For{$i=1$ to $i=s$} {\color{dcolor}\qquad\qquad\qquad\  //Loops for $s$ sub-steps.}
			\State Predict: $\widetilde{\dot{\mbf{U}}}_{n+\gamma_i}=\mbfv_n+\dt\sum_{j=0}^{i-1}\alpha_{ij}\mbfa_{n+\gamma_j}$.
			\State Predict: $\widetilde{\mbfu}_{n+\gamma_i}=\mbfu_n+\dt\sum_{j=0}^{i-1}\alpha_{ij}\mbfv_{n+\gamma_j}+\alpha_{ii}\dt\widetilde{\dot{\mbf{U}}}_{n+\gamma_i}$. 
			\State Solve: $\mbfa_{n+\gamma_i}$ by $\mbf{L}\mbf{U}^\mathsf{T}\mbfa_{n+\gamma_i}=\mbfF(t_n+\gamma_i\dt)-\mbfK\widetilde{\mbfu}_{n+\gamma_i}-\mbfC\widetilde{\dot{\mbf{U}}}_{n+\gamma_i}$.
			\State Compute: $\mbfv_{n+\gamma_i}=\widetilde{\dot{\mbf{U}}}_{n+\gamma_i}+\alpha_{ii}\dt\mbfa_{n+\gamma_i}$.
			\State Compute: $\mbfu_{n+\gamma_i}=\widetilde{\mbfu}_{n+\gamma_i}+{\alpha_{ii}}^2\dt^2\mbfa_{n+\gamma_i}$.
			\EndFor {\color{dcolor}\qquad\ \  // Sub-step schemes.}
			\State {\color{dcolor}// Solutions at the discrete instants $t_n$ are provided by the $s$th sub-step scheme due to $\gamma_s=1$.}
			\EndFor
		\end{algorithmic}\label{code:suci_linear}
	\end{breakablealgorithm}

\subsection{A damped SDOF problem}
The damped SDOF system subjected to the nonzero external load is solved to test the numerical accuracy for various algorithms. It is rigorous for confirming the numerical accuracy to solve the damped forced system. For example, EG3 \cite{fungExtrapolatedGalerkin1996} that achieves third-order accuracy for solving free vibrations is only second-order accurate for simulating forced vibrations. The damped SDOF system is described as
\begin{equation}\label{eq:sdof1}
	\ddot{u}(t)+4\dot{u}(t)+5u(t)=\sin(2t)
\end{equation}
with initial conditions $ u(0)=57/65$ and $\dot{u}(0)=2/65 $. The exact displacement is calculated as $ u(t)=\exp(-2t)(\cos(t)+2\sin(t))-(8\cos(2t)-\sin(2t))/65$. The global error is computed by using
\begin{equation}\label{eq:globalError}
	\text{Error} = \left[\sum_{j=1}^{N}\left(x(t_j)-x_j \right)^2\bigg/\sum_{j=1}^{N}\left(x(t_j)\right)^2 \right]^{1/2}
\end{equation}
where $ N $ denotes the total number of time steps in the analysis; $ x(t_j) $ and $ x_j $ stand for the exact and numerical solutions at time $ t_j $, respectively. In this test, the total simulation time is assumed to be $ t\approx5.62 $s.

\begin{figure}[htbp]
	\centering
	\includegraphics[scale=1.]{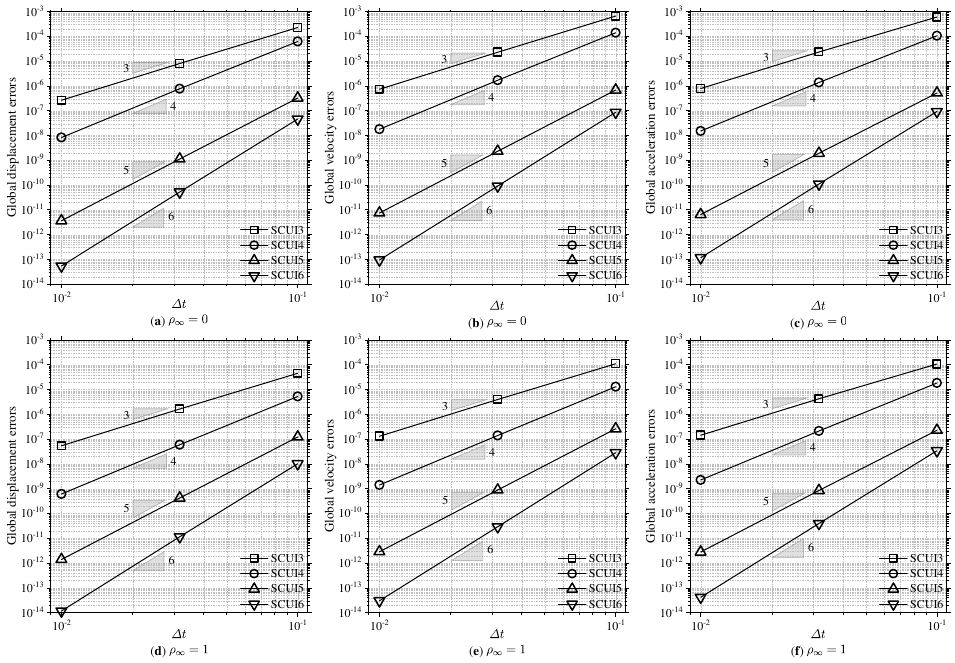}
	\caption{Convergence ratios of four novel high-order SUCI$n$ algorithms: (\textbf{a-c}):$ \rhoinf=0 $ and (\textbf{d-f}): $ \rhoinf=1 $.}
	\label{fig:sdof}
\end{figure}
\begin{figure}[htbp]
	\centering
	\includegraphics[scale=1.0]{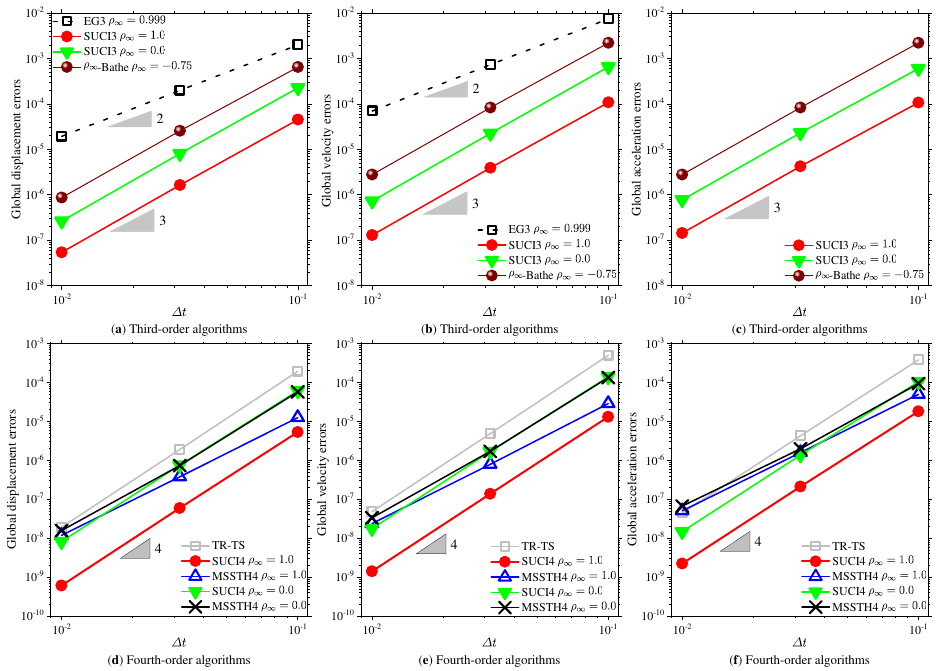}\vskip -1mm
	\caption{Comparisons of global errors for solving the damped forced vibration: (\textbf{a-c}) third-order algorithms and (\textbf{d-f}) fourth-order algorithms.}
	\label{fig:sdof1}
\end{figure}

Fig.~\ref{fig:sdof} first plots global errors in displacement, velocity, and acceleration versus the integration step using SUCI$n$. It is apparent that the novel high-order algorithms achieve the designed order of accuracy in displacement, velocity, and acceleration. For comparison, global errors of EG3 \cite{fungExtrapolatedGalerkin1996}, %TR-CS \cite{fungComplextimestepNewmark1998,fungUnconditionallyStable1997}, 
$ \rhoinf $-Bathe \cite{kwonSelectingLoad2021}, and SUCI3 are plotted in subplots (a-c) of Fig.~\ref{fig:sdof1}, where SUCI3 is obviously superior to EG3 since EG3 presents second-order accuracy for solving forced vibrations. 
%TR-CS produces fewer global errors than SUCI3 due to requiring much more computational costs, while 
The third-order $ \rhoinf $-Bathe method presents the largest global errors among all third-order algorithms. The subplots (d-f) of Fig.~\ref{fig:sdof1} plot global errors of the SDOF system (\ref{eq:sdof1}) using the fourth-order algorithms. It is obvious that SUCI4 is superior to TR-TS \cite{tarnowHowRender1994} and MSSTH4 \cite{zhangOptimizationNsubstep2020}. Particularly, MSSTH4 suffers from the order reduction for solving forced vibrations in the non-dissipative case, showing third-order accuracy. This case also holds for MSSTH5 \cite{zhangOptimizationNsubstep2020}. Hence, the proposed SUCI5 algorithm is obviously superior to MSSTH5.

%\begin{figure}[htbp]
%	\centering
%	\includegraphics[scale=1.0]{sdof3}
%	\caption{Convergence ratios of MSSTH5 \cite{zhangOptimizationNsubstep2020} for solving the damped forced vibration.}
%	\label{fig:sdof3}
%\end{figure}

\begin{figure}[htbp]
	\centering
	\subfigtopskip=2pt %?????????????????
	\subfigbottomskip=-4pt %??????????????????????????????
	\subfigcapskip=-5pt %?????????????
	\subfigure[DSUCI3 ]{
		\includegraphics[scale=1.0]{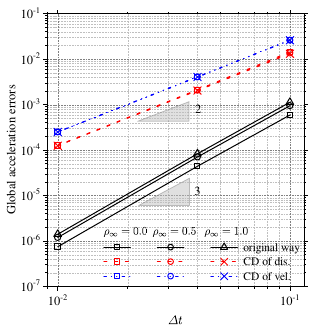}}
	\subfigure[DSUCI4 ]{
		\includegraphics[scale=1.0]{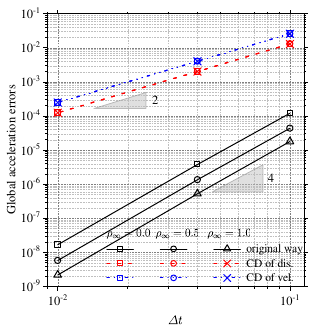}}
	\subfigure[DSUCI5 ]{
		\includegraphics[scale=1.0]{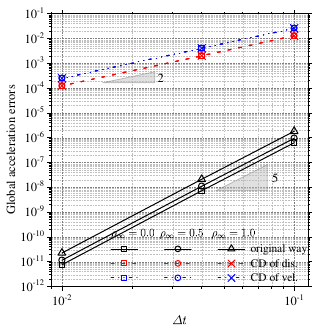}}
	\caption{Convergence ratios of accelerations predicted by DSUCI$n$ with different post-processing ways. For the original way, DSUCI3 adopts the weight combination of sub-step accelerations whereas other DSUCI$n$ algorithms use the balance equation at the discretized instants.}
	\label{fig:con_cd}
\end{figure}

Since the previous DSUCI$n$ algorithms \cite{liDirectlySelfstarting2022} are designed to be directly self-starting, some post-processing techniques are needed to output the acceleration response whenever necessary. The original study \cite{liDirectlySelfstarting2022} provides two approaches for users to output accelerations. In addition, after obtaining accurate displacement and velocity responses, users can employ the central difference (CD) of displacement or velocity to output accelerations. As shown in Fig.~\ref{fig:con_cd}, using the central difference of displacement or velocity, DSUCI$n$ always predicts second-order accurate accelerations even if the displacement and velocity responses are high-order accurate. It is also observed that  accelerations produced by the central difference of displacement are slightly more accurate than those from velocity. By default, DSUCI$n$ uses the original techniques given in \cite{liDirectlySelfstarting2022} to output accelerations, and the central difference technique could show its advantage for some solutions, such as the stiff-flexible dynamic problem in Section \ref{sec:2dofs}. 

This damped SDOF system subjected to the load clearly shows that the proposed high-order algorithms (SUCI$n$) achieve the identical order of accuracy in displacement, velocity, and acceleration. Importantly, the novel methods do not suffer from the order reduction for the benchmark test. These conclusions are naturally valid for simulating other simpler cases, such as the undamped and/or free vibrations.

%\begin{figure}[htbp]
%	\centering
%	\includegraphics[scale=1.0]{sdof2}
%	\caption{Comparisons of global errors among the fourth-order accurate algorithms for solving the damped forced vibration.}
%	\label{fig:sdof2}
%\end{figure}

\begin{figure}[htpb]
	\centering
	\includegraphics[scale=0.4]{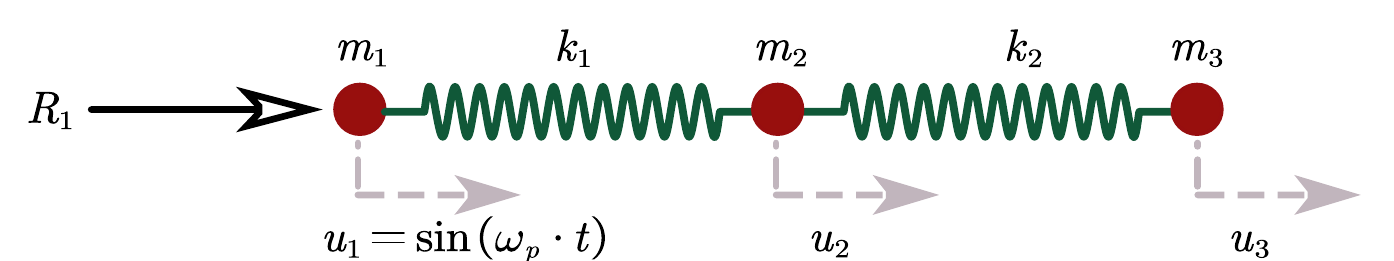}
	\caption{A linear mass-spring system where $k_1=10^7$N/m, $k_2=1$N/m, $m_1=0$kg, $m_2=m_3=1$kg, and $ \omega_p=1.2$rad/s.}
	\label{fig:modalproblem}
\end{figure}
\subsection{A double-degree-of-freedom system with spurious high-frequency component}\label{sec:2dofs}
A standard modal problem \cite{choiTimeSplitting2022,batheInsightImplicit2012} shown in Fig.~\ref{fig:modalproblem} has been solved by various dissipative algorithms to show their capabilities of eliminating spurious high-frequency modes. This mass-spring system represents the complex structural dynamic problems, consisting of stiff and flexible parts. The initial governing equation is expressed as
\begin{equation}\label{eq:modeproblem}
	\begin{bmatrix}
		m_1 & 0 & 0 \\ 0 & m_2 & 0\\ 0 & 0 & m_3\\
	\end{bmatrix}\begin{bmatrix}
		\ddot{u}_1 \\ \ddot{u}_2 \\ \ddot{u}_3\\
	\end{bmatrix}+\begin{bmatrix}
		k_1 & -k_1 & 0 \\ -k_1 & k_1+k_2 & -k_2\\ 0 & -k_2 & k_2\\
	\end{bmatrix}\begin{bmatrix}
		u_1 \\ u_2\\ u_3\\
	\end{bmatrix}=\begin{bmatrix}
		R_1 \\ 0 \\ 0\\
	\end{bmatrix}
\end{equation}
where the displacement at node 1 is assumed to be $u_1=\sin(\omega_pt)=\sin(1.2t)$m and $R_1$ is the reaction force at node 1. All initial conditions are set to be zero. By using the known $u_1$ at node 1, Eq.~(\ref{eq:modeproblem}) is further simplified as
\begin{equation}
	\begin{bmatrix}
		m_2 & 0 \\ 0 & m_3\\
	\end{bmatrix}\begin{bmatrix}
		\ddot{u}_2 \\ \ddot{u}_3\\
	\end{bmatrix}+\begin{bmatrix}
		k_1+k_2 & -k_2 \\ -k_2 & k_2\\
	\end{bmatrix}\begin{bmatrix}
		u_2 \\ u_3\\
	\end{bmatrix}=\begin{bmatrix}
		k_1u_1 \\ 0\\
	\end{bmatrix}.\label{eq:standardproblem}
\end{equation}
Note that the high-order methods compared in this paper adopt the different number of sub-steps and their algorithmic parameters are either real-valued or complex-valued, so the compared methods should set different time steps to approximate the same computational cost. Typically, assuming that the single-sub-step algorithms with real-valued parameters use the time step as $\dt$, the time steps of $n$-sub-step algorithms with real- and complex-valued parameters should be set as $n\times \dt$ and $4n\times \dt$, respectively. This strategy has been used by Choi et al. \cite{choiTimeSplitting2022}, who developed and analyzed the third-order $\rhoinf$-Bathe algorithm with complex-valued parameters.

\begin{figure}[hbtp]
	\centering
	\subfigtopskip=2pt %?????????????????
	\subfigbottomskip=-4pt %??????????????????????????????
	\subfigcapskip=-5pt %?????????????
	\includegraphics[scale=0.37]{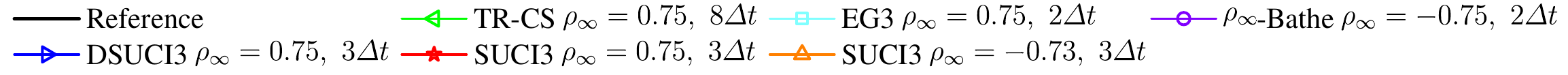}
	\subfigure[ ]{
		\includegraphics[scale=0.37]{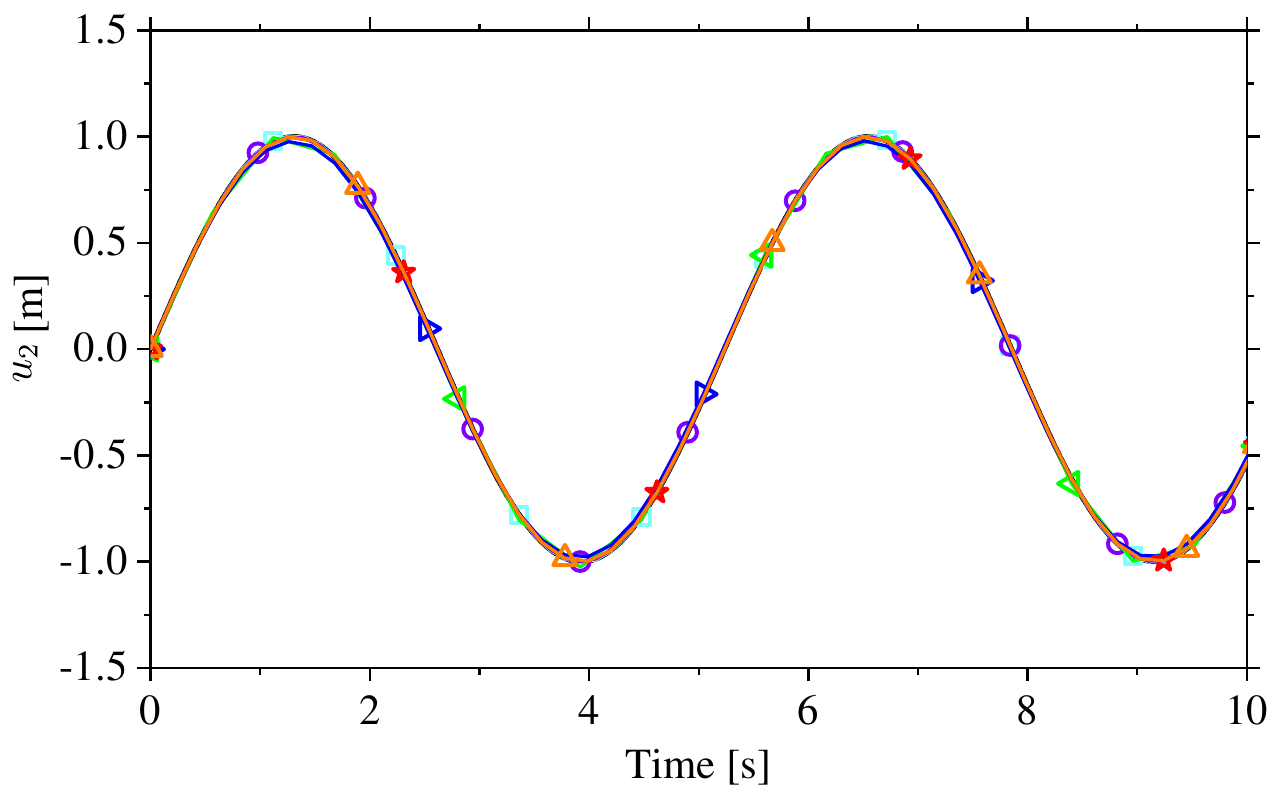}}
	\subfigure[ ]{
		\includegraphics[scale=0.37]{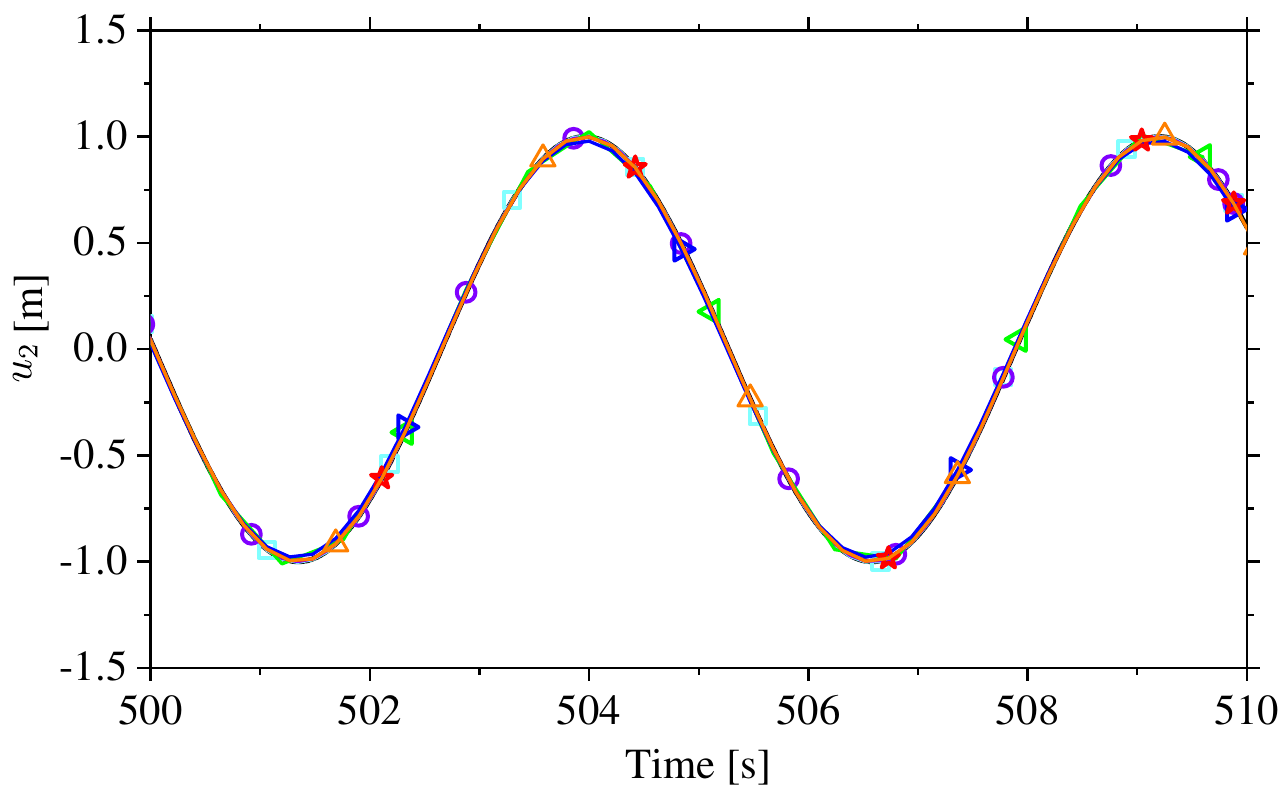}}
	\subfigure[ ]{
		\includegraphics[scale=0.37]{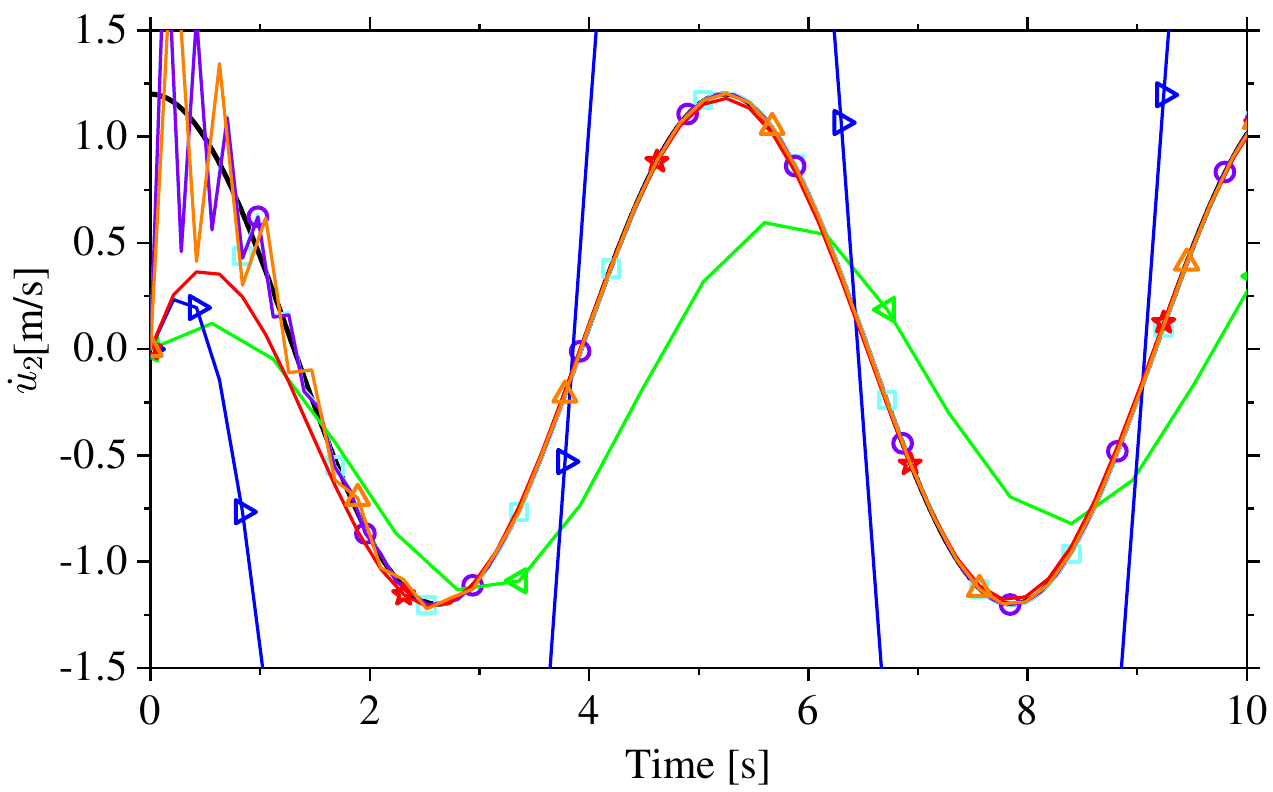}}
	\subfigure[ ]{
		\includegraphics[scale=0.37]{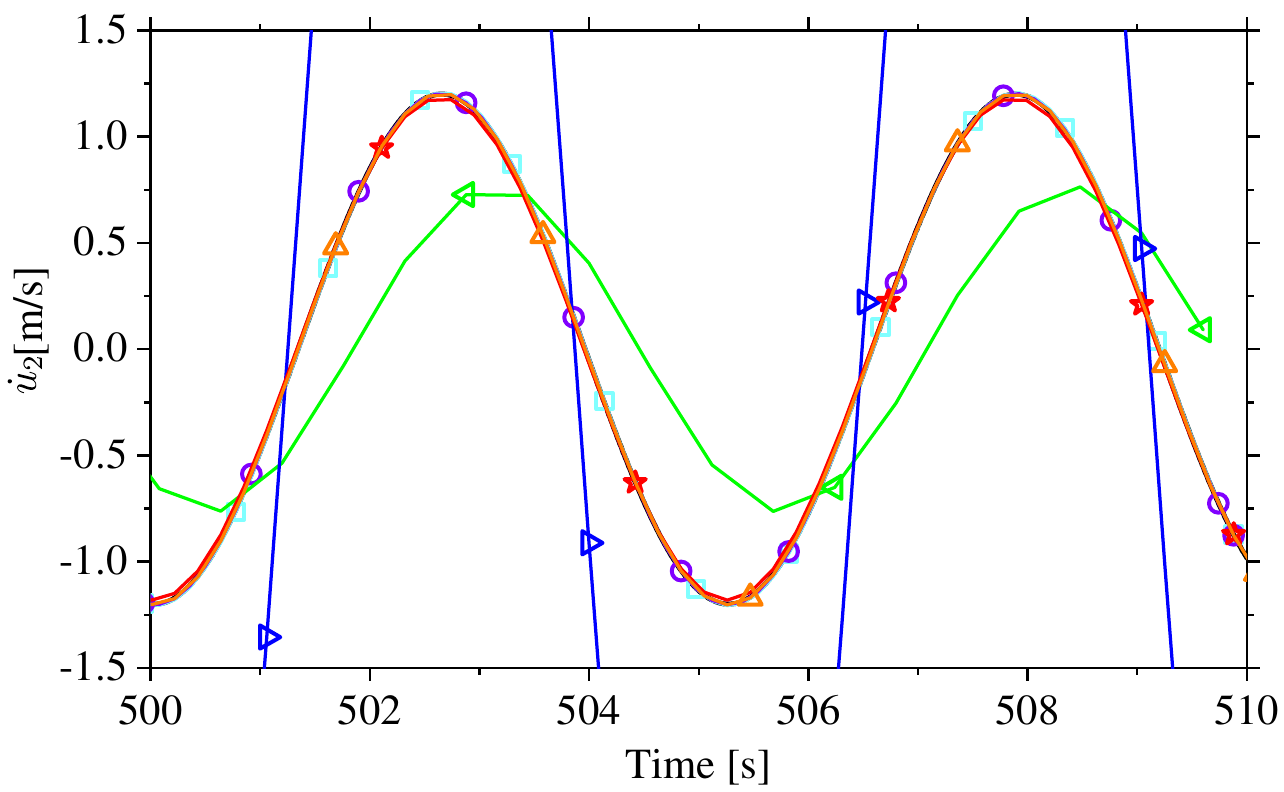}}
	\subfigure[ ]{
		\includegraphics[scale=0.37]{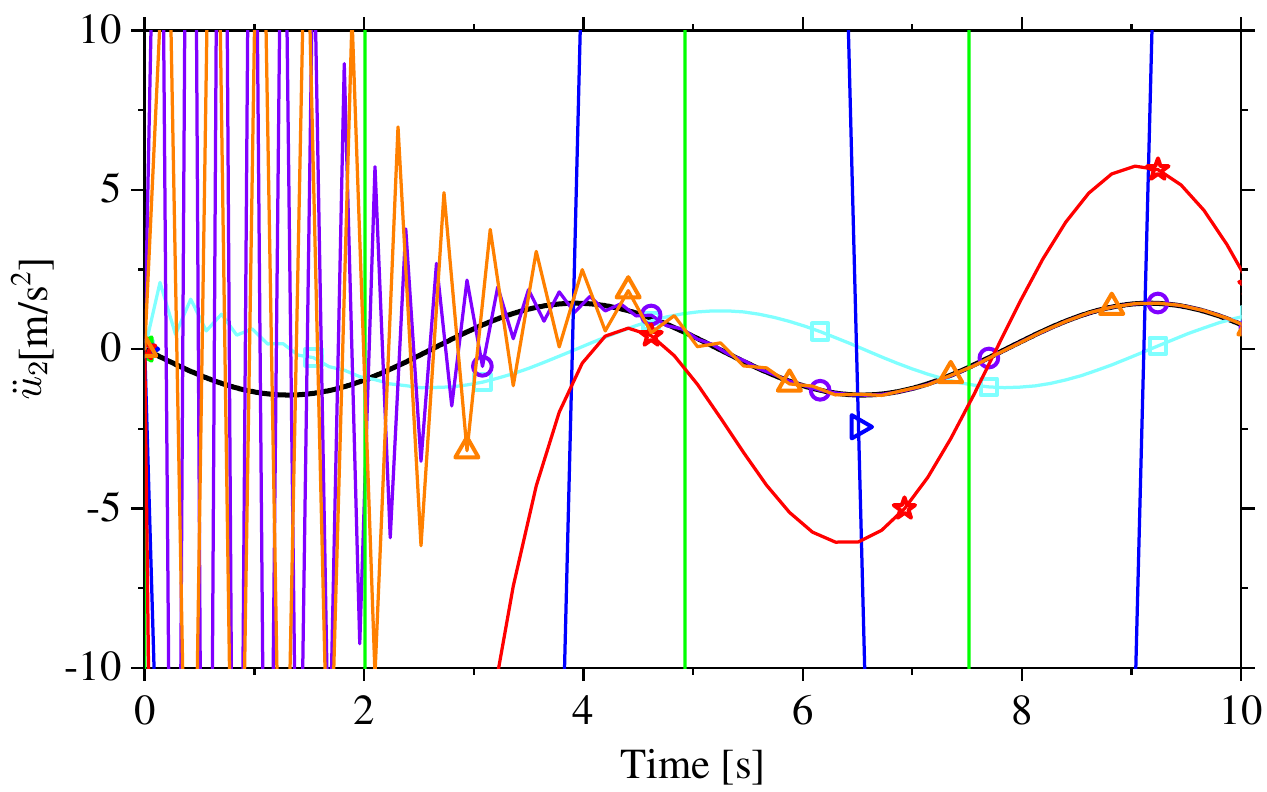}}
	\subfigure[ ]{
		\includegraphics[scale=0.37]{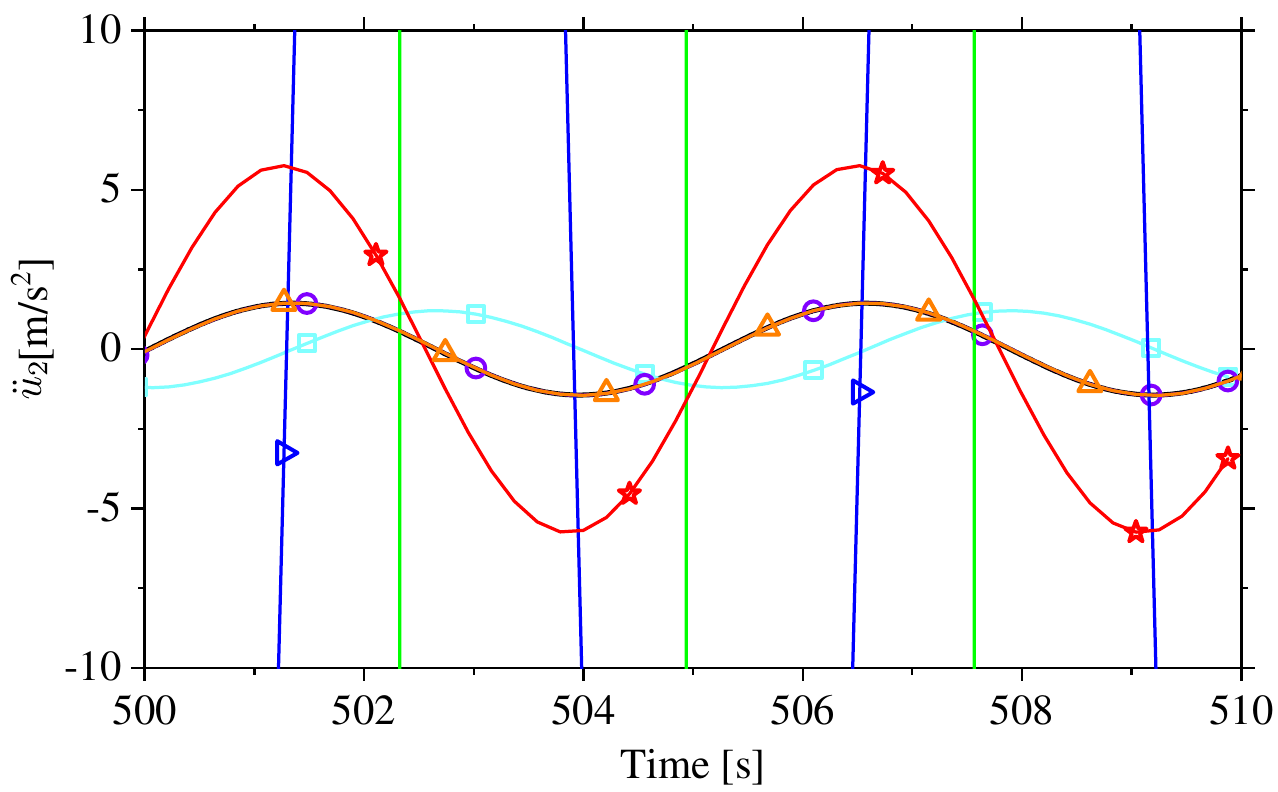}}
	\caption{Numerical displacements, velocities, and accelerations at node 2,  predicted by various third-order implicit algorithms with $\dt=0.07$s.}
	\label{fig:2dof_3rd_u2}
\end{figure}
\begin{figure}[htbp]
	\centering
	\subfigtopskip=2pt %?????????????????
	\subfigbottomskip=-4pt %??????????????????????????????
	\subfigcapskip=-5pt %?????????????
	{
		\includegraphics[scale=0.37]{2dof_3rd_leg}}
	\subfigure[ ]{
		\includegraphics[scale=0.37]{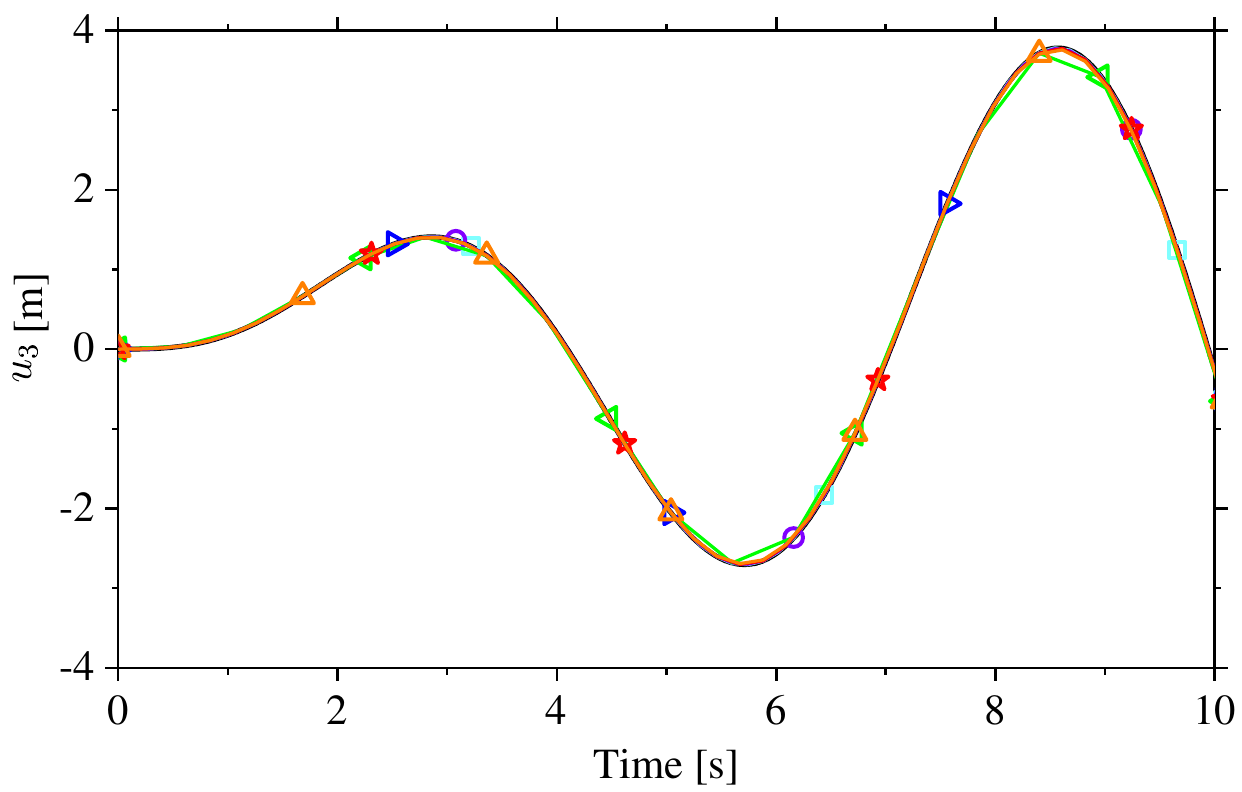}}
	\subfigure[ ]{
		\includegraphics[scale=0.37]{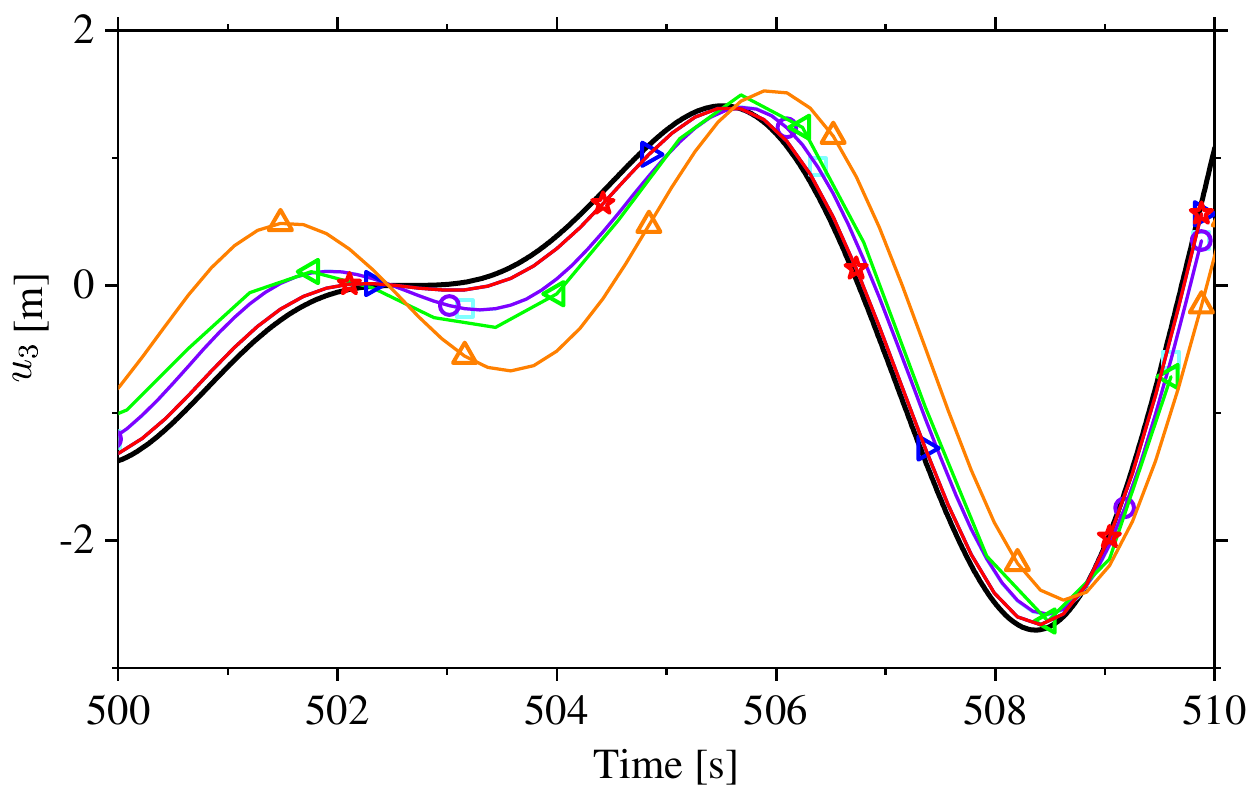}}
	\subfigure[ ]{
		\includegraphics[scale=0.37]{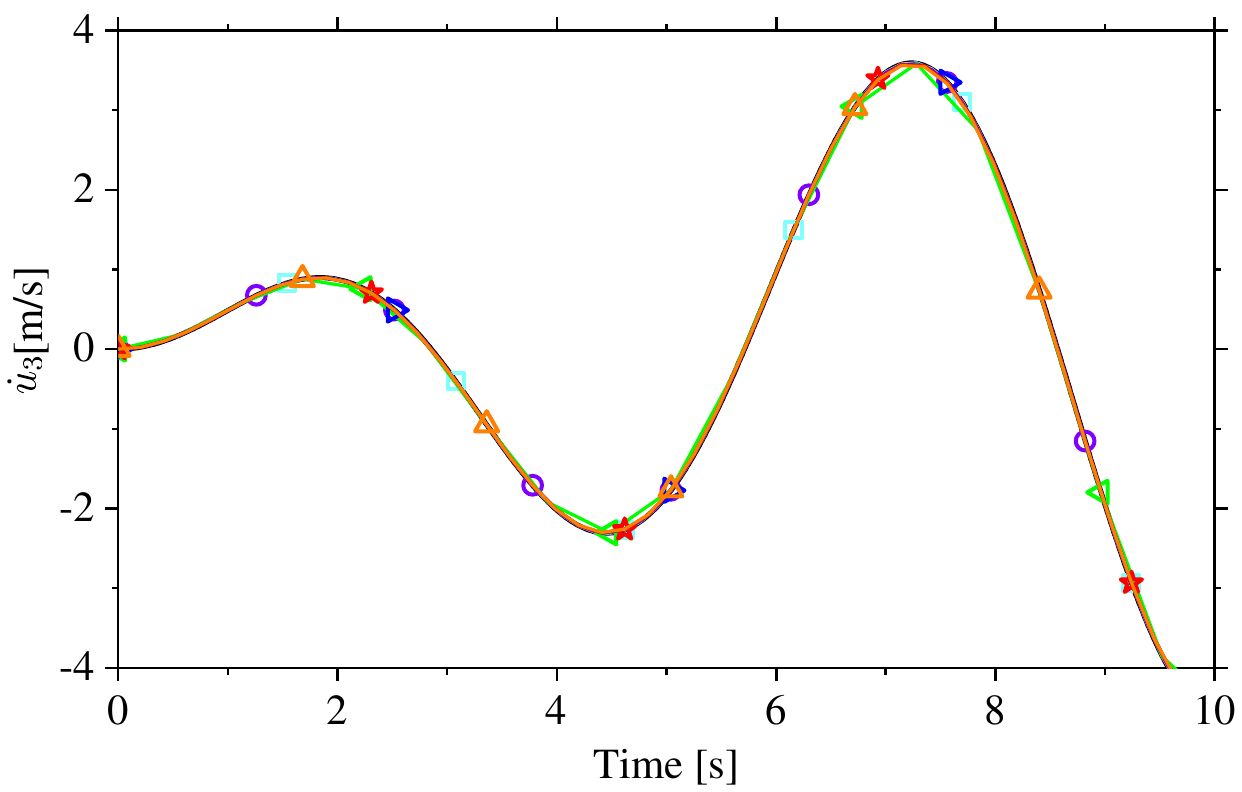}}
	\subfigure[ ]{
		\includegraphics[scale=0.37]{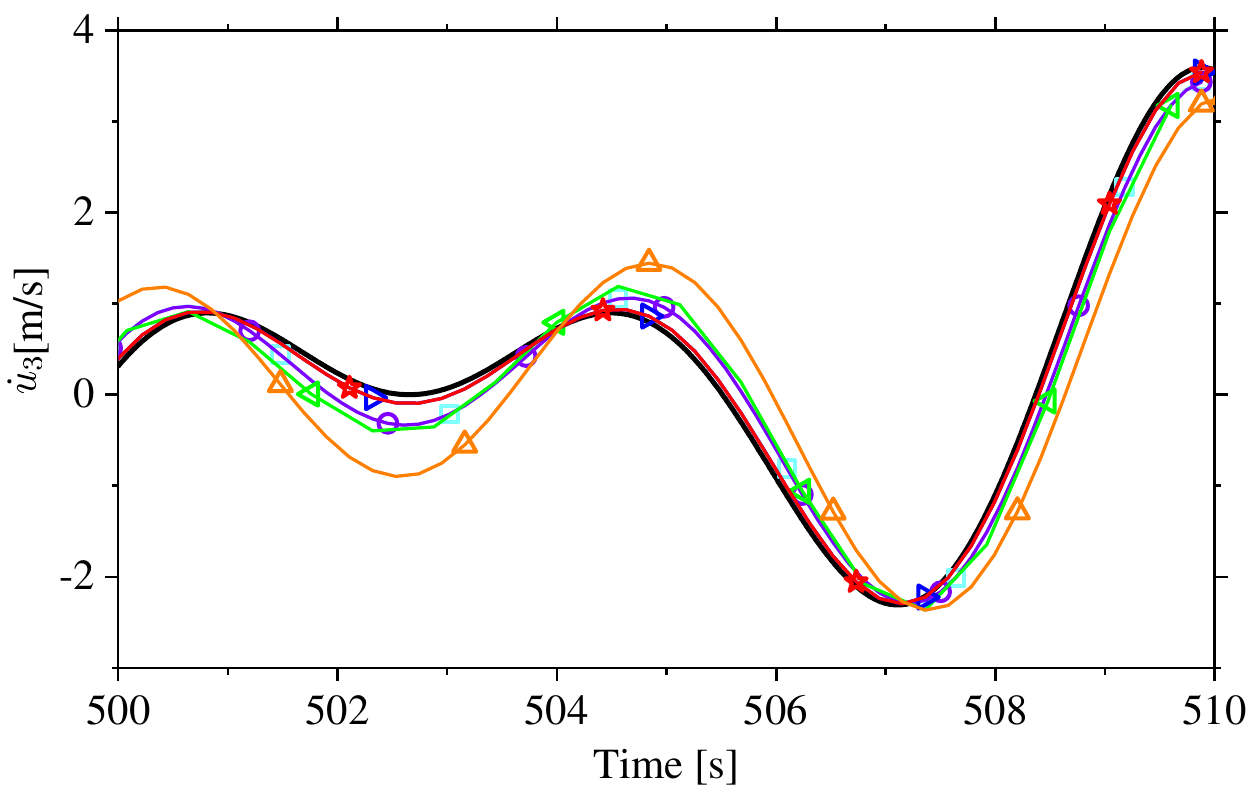}}
	\subfigure[ ]{
		\includegraphics[scale=0.37]{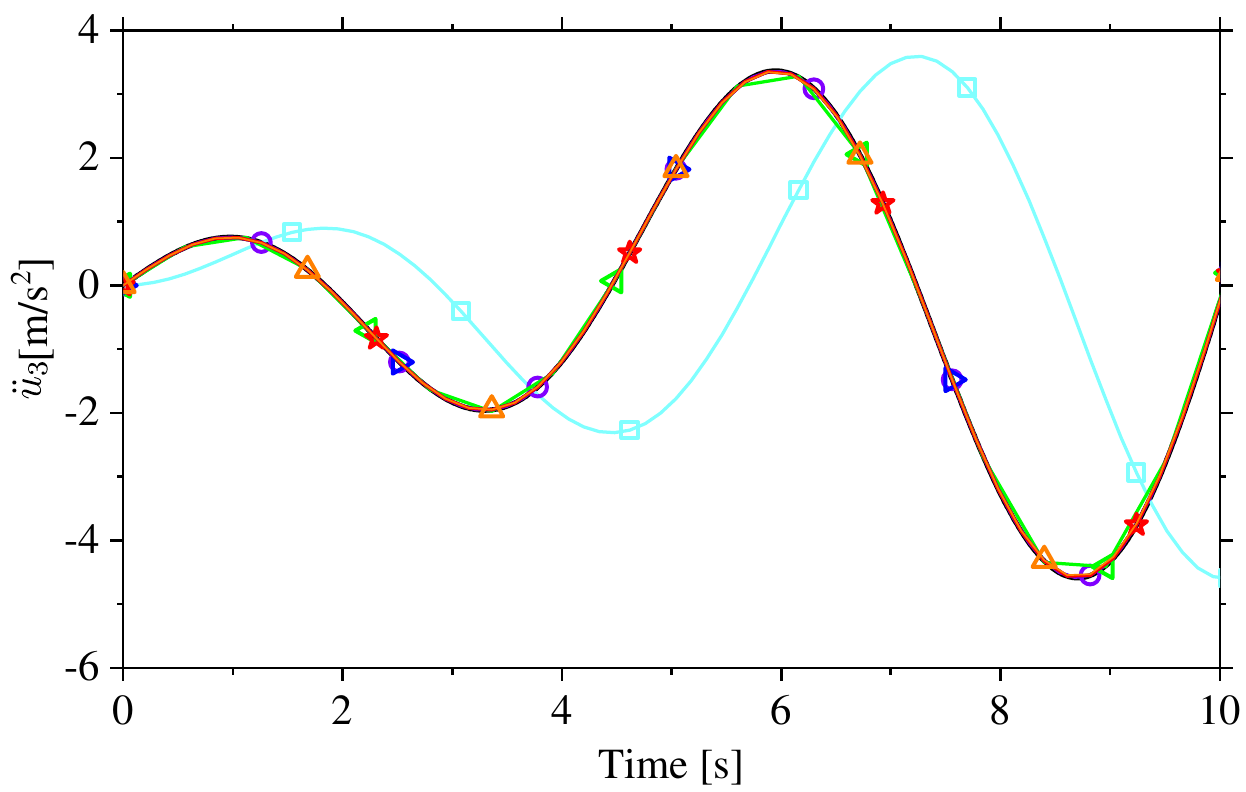}}
	\subfigure[ ]{
		\includegraphics[scale=0.37]{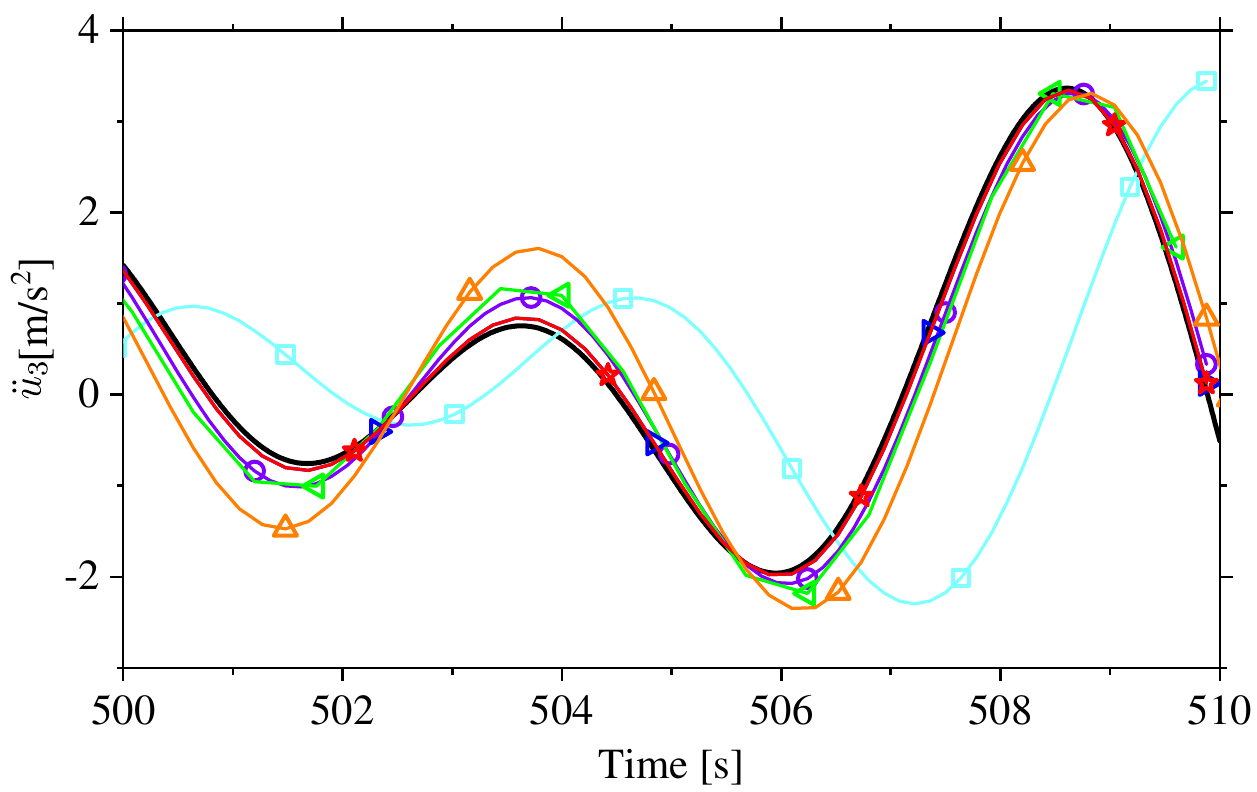}}
	\subfigure[ ]{
		\includegraphics[scale=0.37]{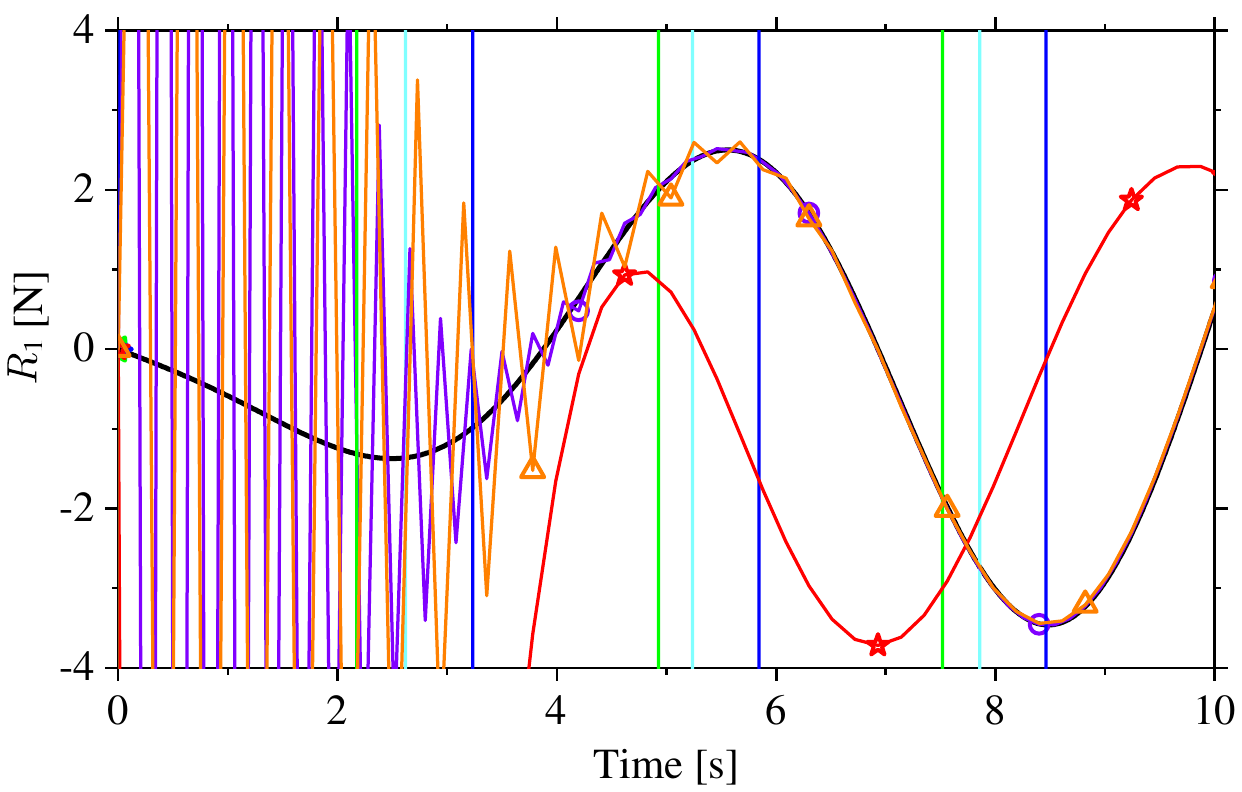}}
	\subfigure[ ]{
		\includegraphics[scale=0.37]{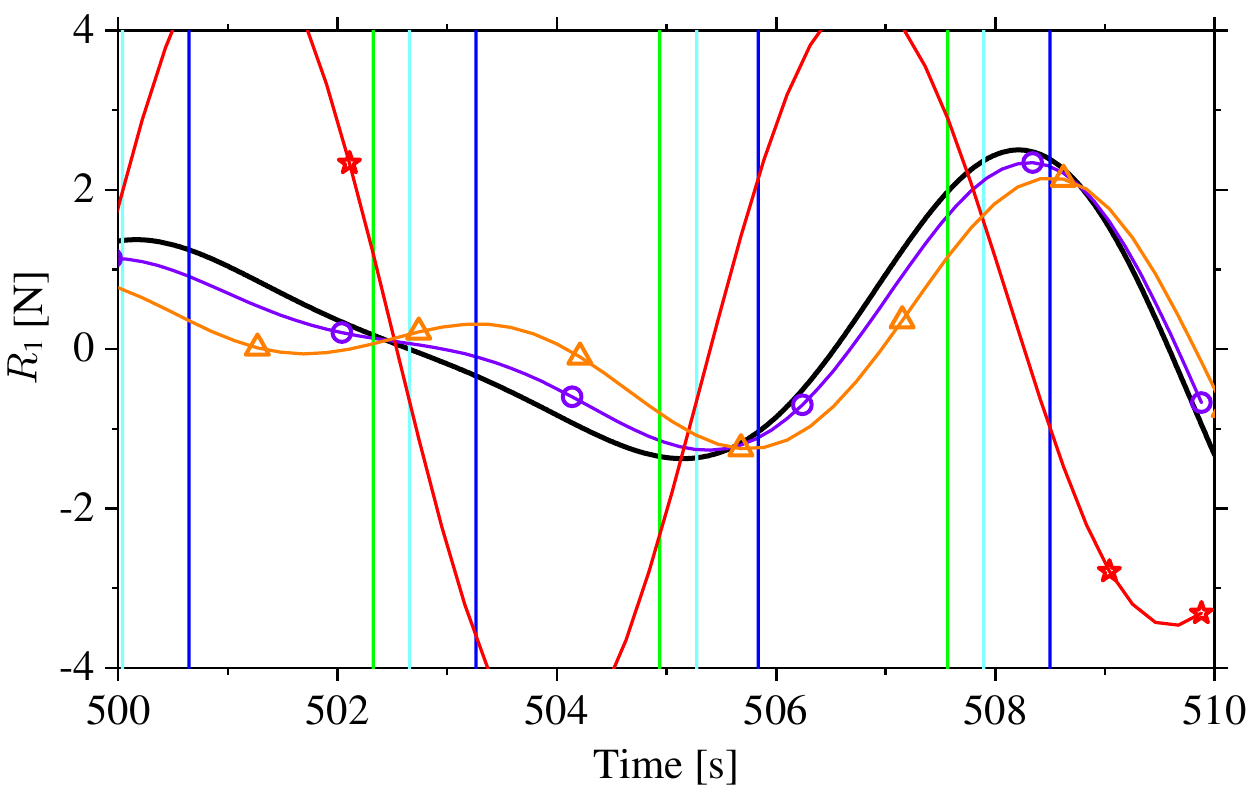}}
	\caption{Numerical displacements, velocities, and accelerations at node 3 and reaction force $R_1$,  predicted by various third-order implicit algorithms with $\dt=0.07$s.}
	\label{fig:2dof_3rd_u3}
\end{figure}
Figs.~\ref{fig:2dof_3rd_u2} and \ref{fig:2dof_3rd_u3} present a comparison of the numerical responses of nodes 2 and 3, as well as the reaction force $R_1$ at node 1, among the third-order implicit algorithms with moderate dissipation. It is observed that the TR-CS \cite{fungComplextimestepNewmark1998,fungUnconditionallyStable1997} and DSUCI3 \cite{liDirectlySelfstarting2022} algorithms fail to accurately predict velocities and accelerations of node 2, as well as the reaction force at node 1. Additionally, the EG3 algorithm \cite{fungExtrapolatedGalerkin1996} does not provide accurate predictions for accelerations and the reaction force. 
By utilizing the recommended positive value of $\rhoinf$, the SUCI3 algorithm introduces noticeable errors in the acceleration of node 2 and the reaction force. However, these errors can be effectively reduced by employing negative values of $\rhoinf$, as exemplified by the use of $\rhoinf=-0.73$. As shown in Fig.~\ref{fig:adCom}, the $\rhoinf$-Bathe method \cite{kwonSelectingLoad2021} with $\rhoinf=-0.75$ exhibits the largest amplitude decays among the compared algorithms, enabling it to accurately predict the responses of this model. This fact also highlights why SUCI3 with $\rhoinf=-0.73$ is capable of providing more accurate accelerations and reaction forces.
Overall, when moderate dissipation is considered, the $\rhoinf$-Bathe algorithm yields the best performance, followed by SUCI3, in solving this stiff-flexible dynamic problem. However, the superiority of SUCI3 over $\rhoinf$-Bathe becomes evident when $\rhoinf=0$ is utilized.
%Figs.~\ref{fig:2dof_3rd_u2} and \ref{fig:2dof_3rd_u3} compare the numerical responses of nodes 2 and 3, and the reaction force $R_1$ at node 1, among the third-order implicit algorithms with moderate dissipation. It is seen that TR-CS \cite{fungComplextimestepNewmark1998,fungUnconditionallyStable1997} and DSUCI3 \cite{liDirectlySelfstarting2022} do not predict reasonable velocities and accelerations of node 2, and the reaction force at node 1, and EG3 \cite{fungExtrapolatedGalerkin1996} does not provide accurate predictions for accelerations and the reaction force. Using the recommended positive value of $\rhoinf$, SUCI3 produces noticeable errors in the acceleration of node 2 and the reaction force, but these errors can be effectively reduced by using the negative values of $\rhoinf$, such as the used $\rhoinf=-0.73$. As analyzed in Fig.~\ref{fig:adCom}, the $\rhoinf$-Bathe method \cite{kwonSelectingLoad2021} with $\rhoinf=-0.75$ imposes the largest amplitude decays among the compared algorithms, so it is able to predict more accurate responses for solving this model. This is the reasons why SUCI3 with $\rhoinf=-0.73$ can predict more accurate accelerations and the reaction force. Overall, with moderate dissipation, the $\rhoinf$-Bathe algorithm performs best, followed by SUCI3, for the present stiff-flexible dynamics. However, the superiority of SUCI3 over $\rhoinf$-Bathe will be highlighted in the case of $\rhoinf=0$. 

\begin{figure}[htbp]
	\centering
	\subfigtopskip=2pt %?????????????????
	\subfigbottomskip=-4pt %??????????????????????????????
	\subfigcapskip=-5pt %?????????????
	\includegraphics[scale=0.37]{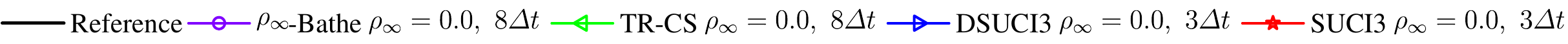}
	\subfigure[ ]{
		\includegraphics[scale=0.37]{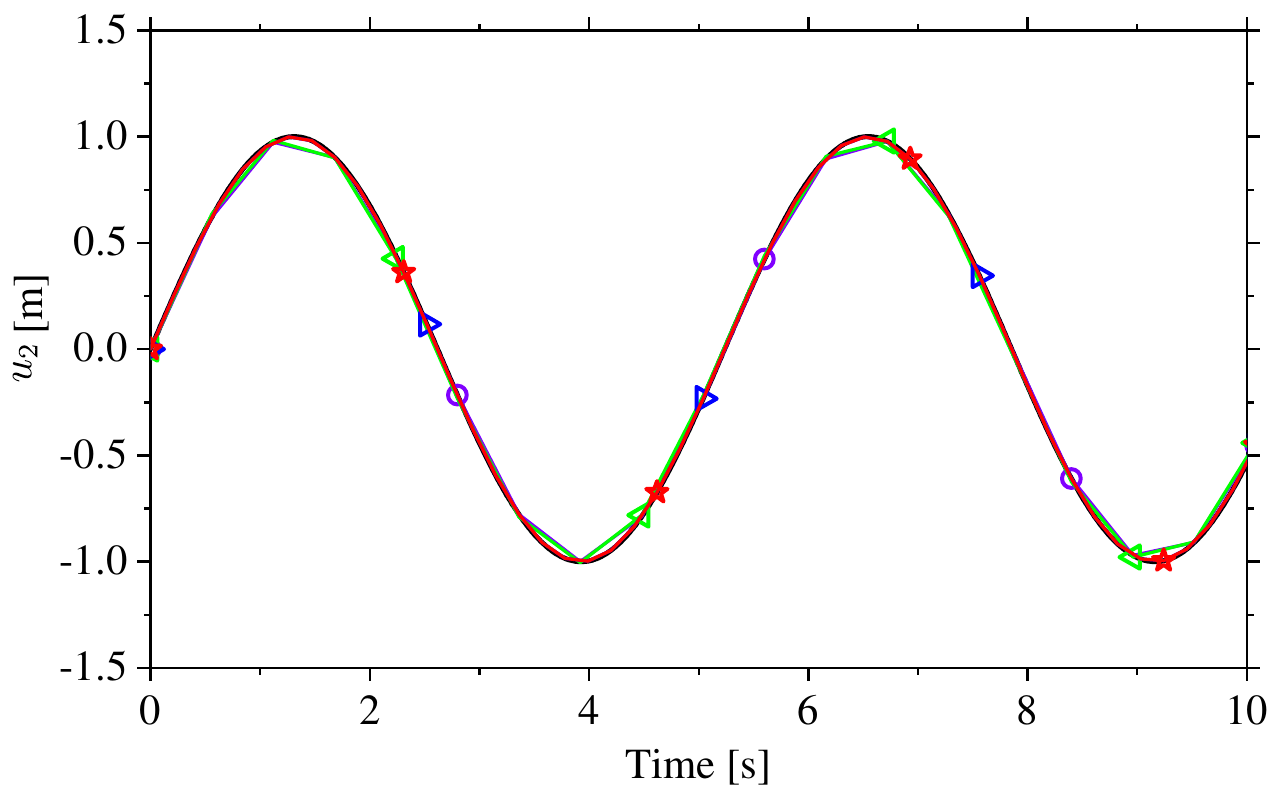}}
	\subfigure[ ]{
		\includegraphics[scale=0.37]{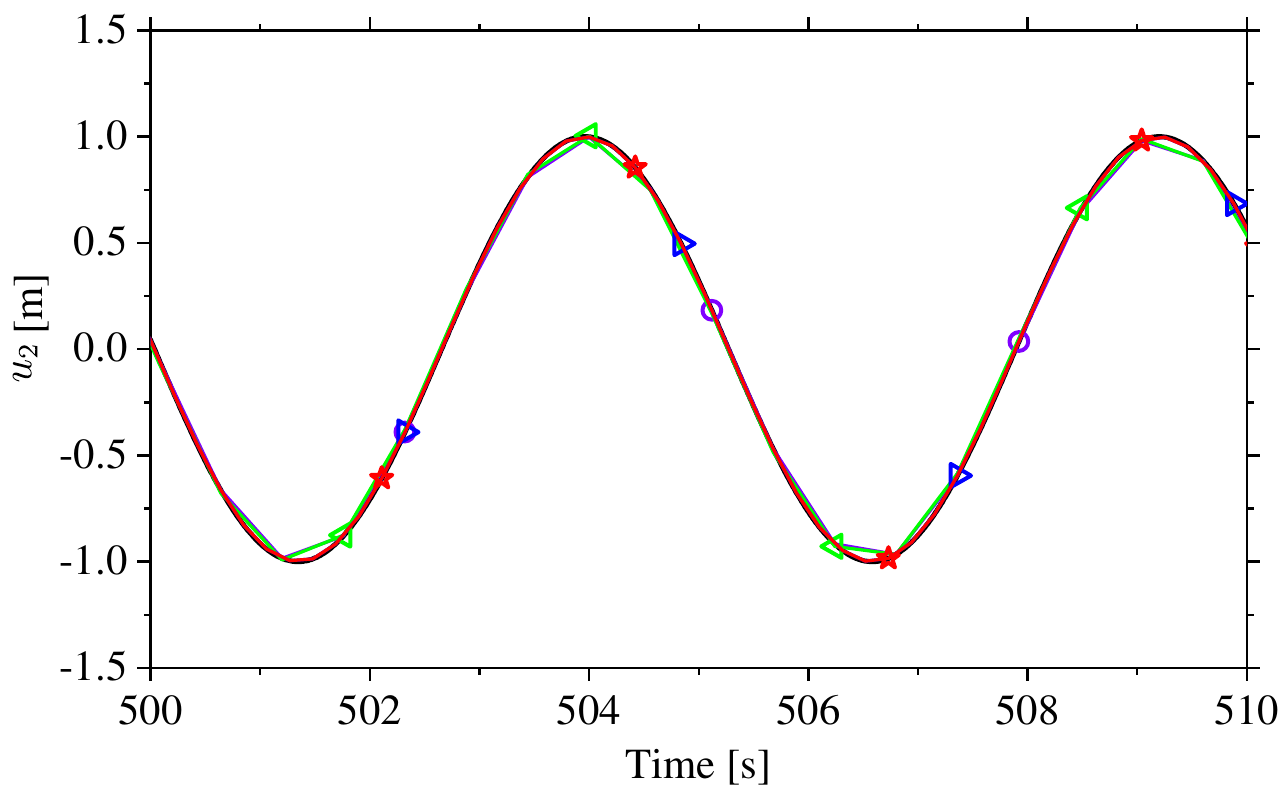}}
	\subfigure[ ]{
		\includegraphics[scale=0.37]{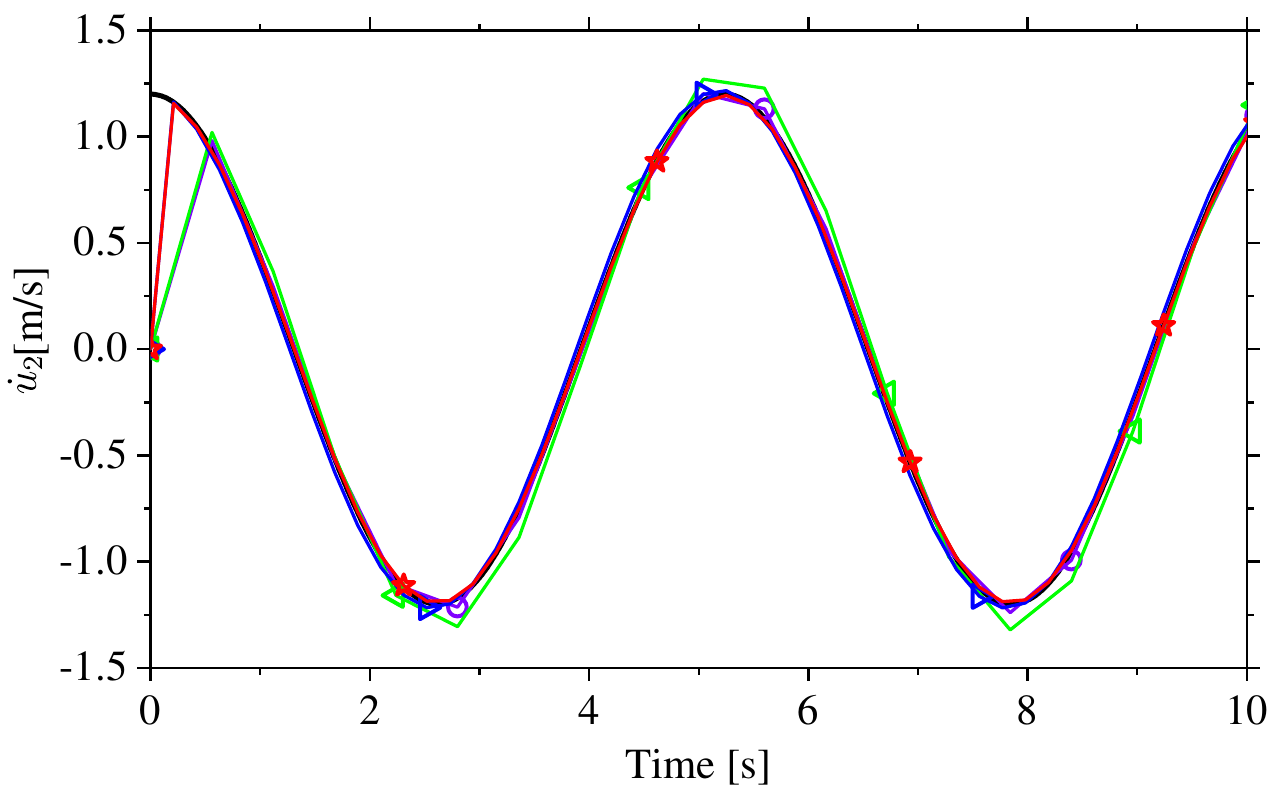}}
	\subfigure[ ]{
		\includegraphics[scale=0.37]{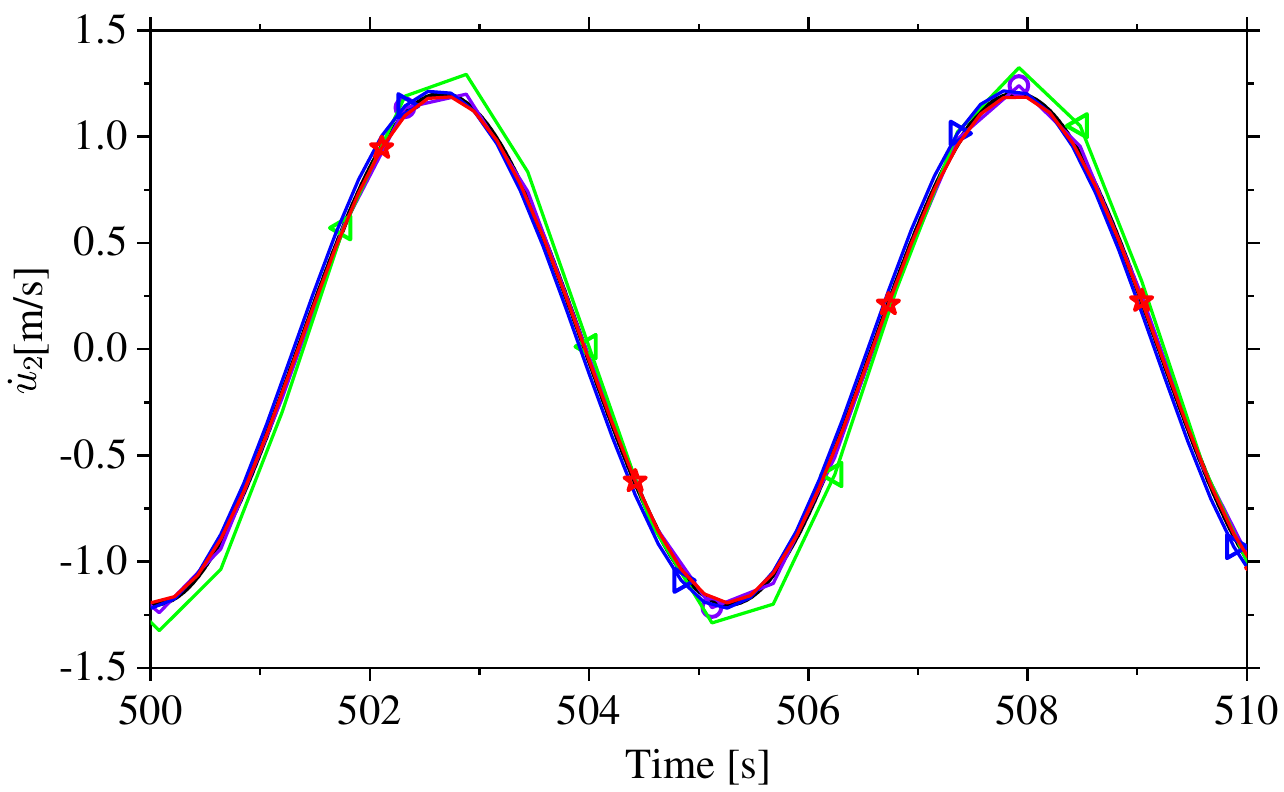}}
	\subfigure[ ]{
		\includegraphics[scale=0.37]{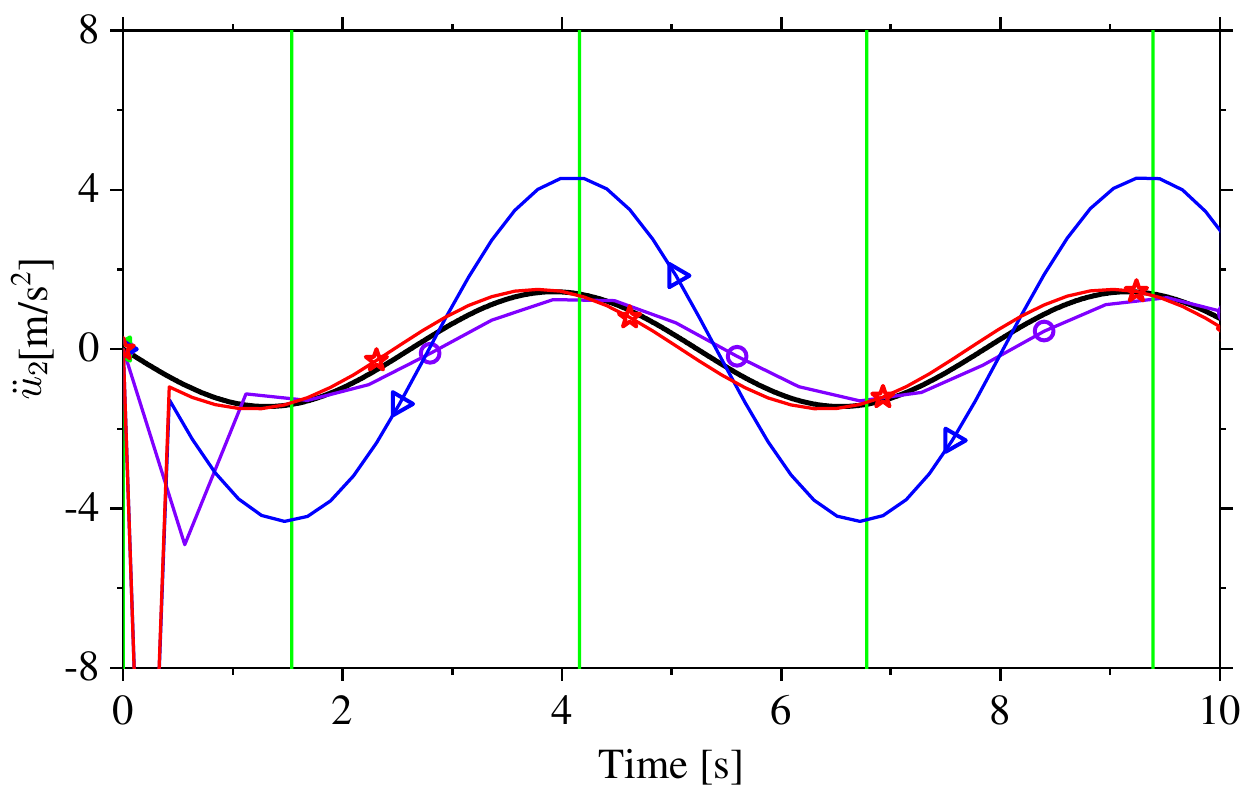}}
	\subfigure[ ]{
		\includegraphics[scale=0.37]{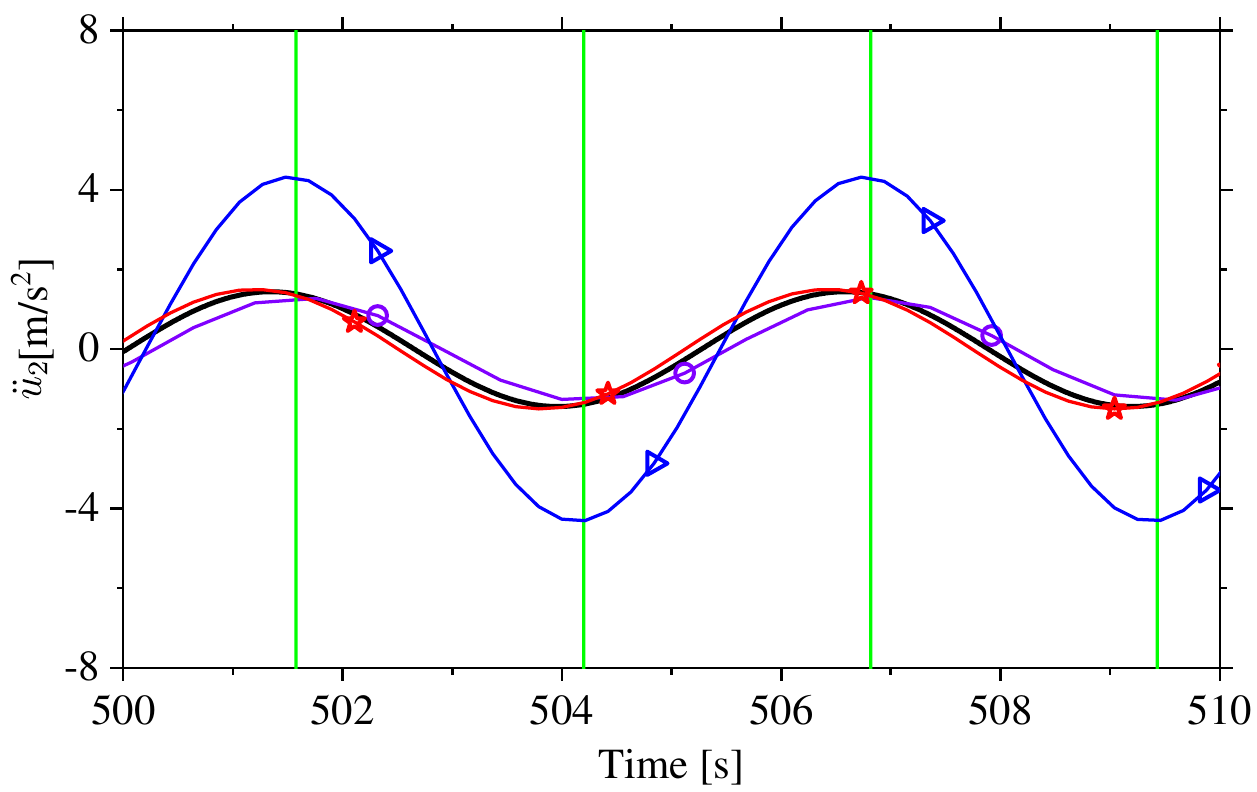}}
	\caption{Numerical displacements, velocities, and accelerations at node 2,  predicted by various third-order implicit algorithms with $\rhoinf=0$ and $\dt=0.07$s.}
	\label{fig:2dof_3rd_u2_1}
\end{figure}

\begin{figure}[htbp]
	\centering
	\subfigtopskip=2pt %?????????????????
	\subfigbottomskip=-4pt %??????????????????????????????
	\subfigcapskip=-5pt %?????????????
	\includegraphics[scale=0.37]{2dof_3rd_leg_1}
	\subfigure[ ]{
		\includegraphics[scale=0.37]{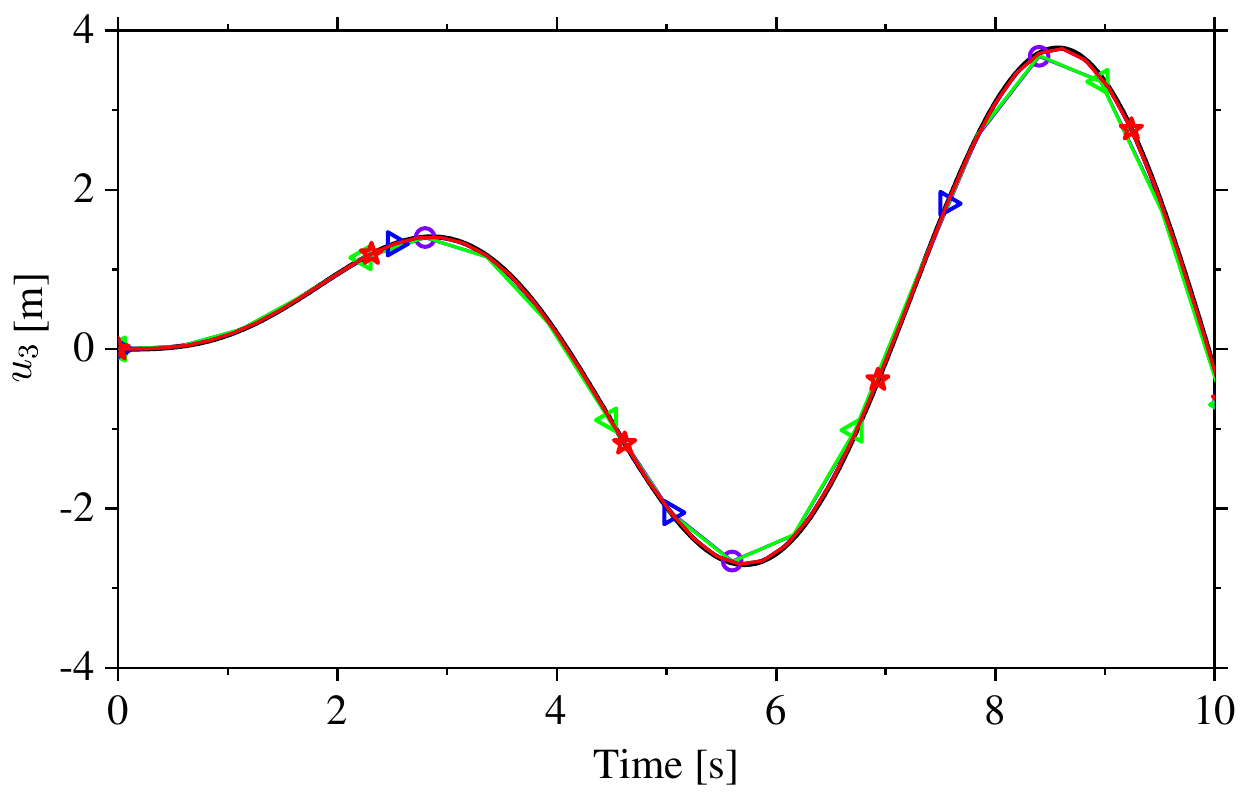}}
	\subfigure[ ]{
		\includegraphics[scale=0.37]{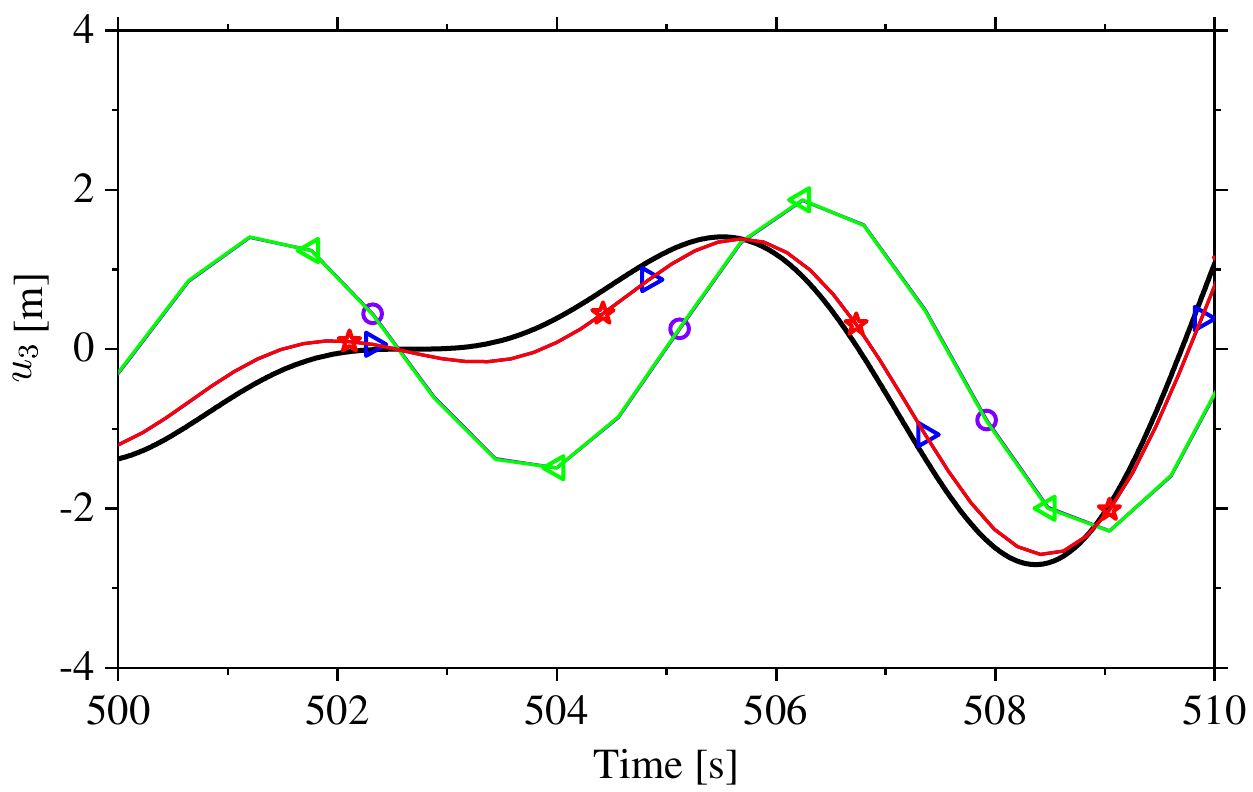}}
	\subfigure[ ]{
		\includegraphics[scale=0.37]{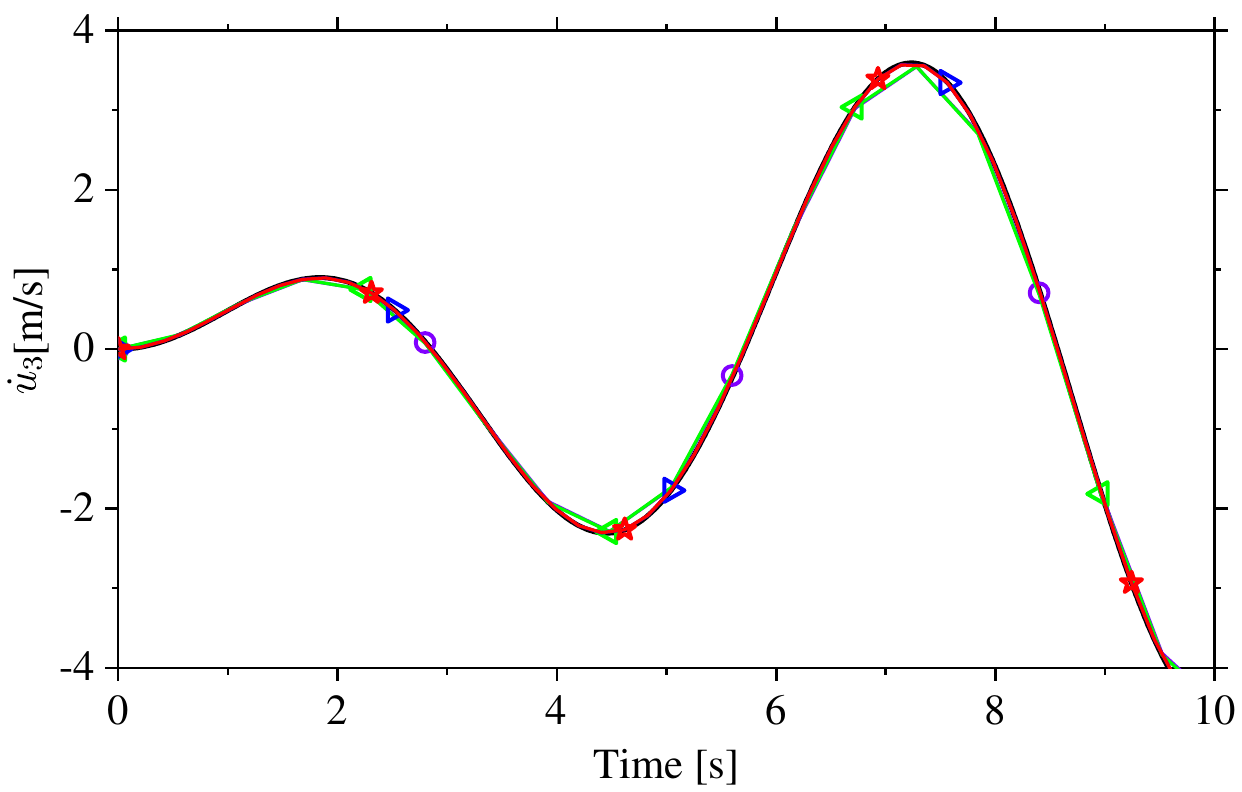}}
	\subfigure[ ]{
		\includegraphics[scale=0.37]{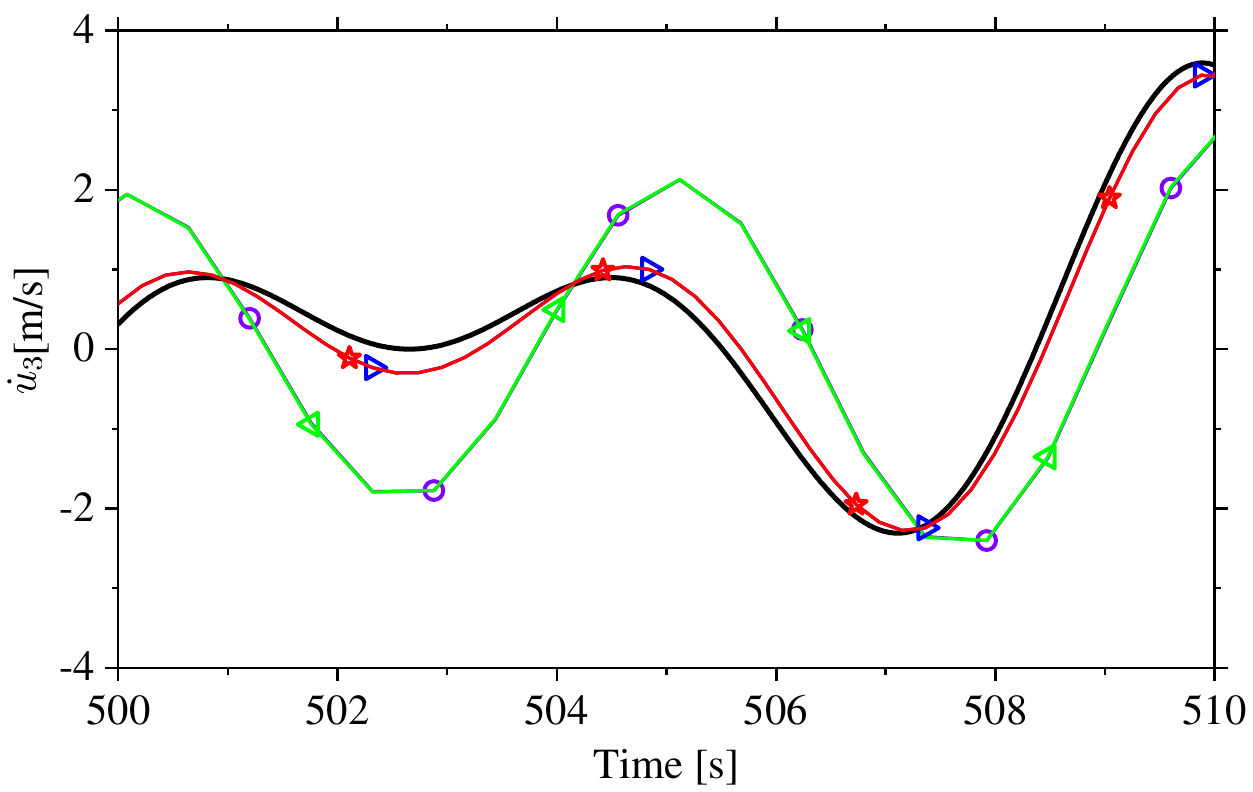}}
	\subfigure[ ]{
		\includegraphics[scale=0.37]{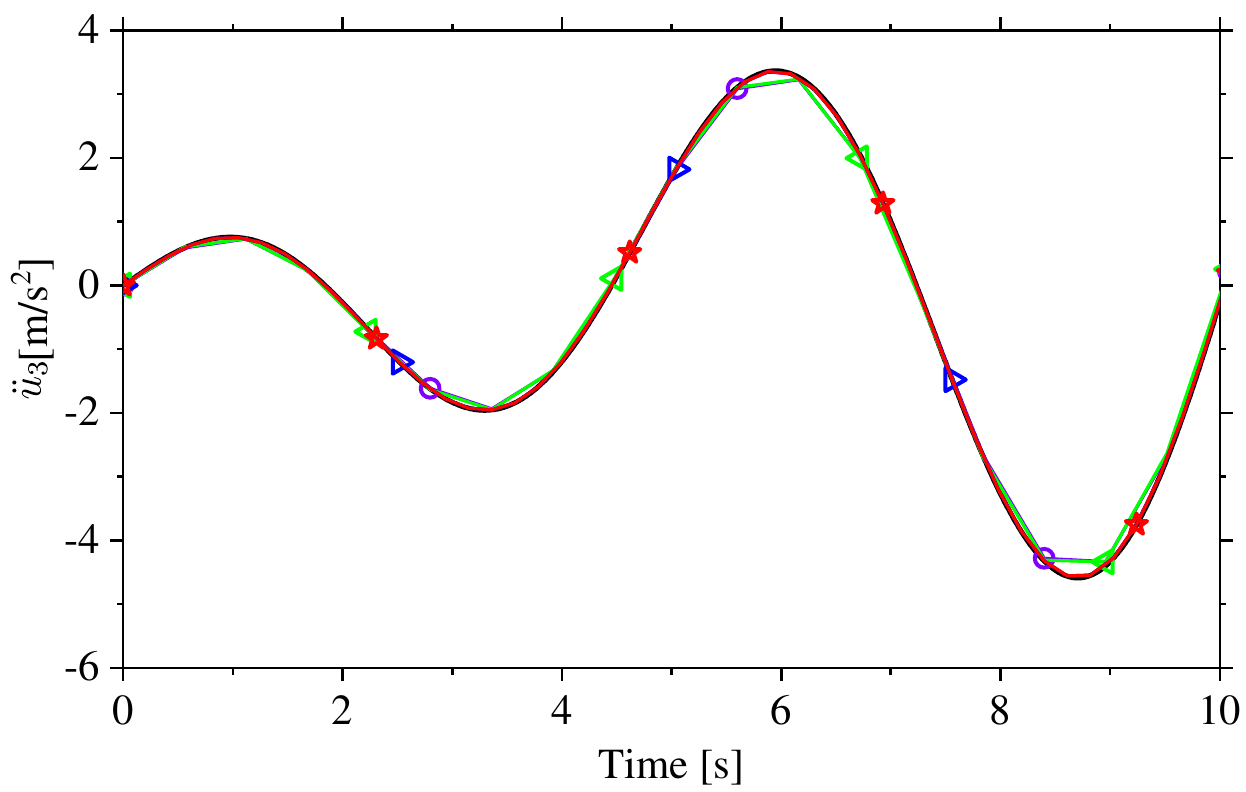}}
	\subfigure[ ]{
		\includegraphics[scale=0.37]{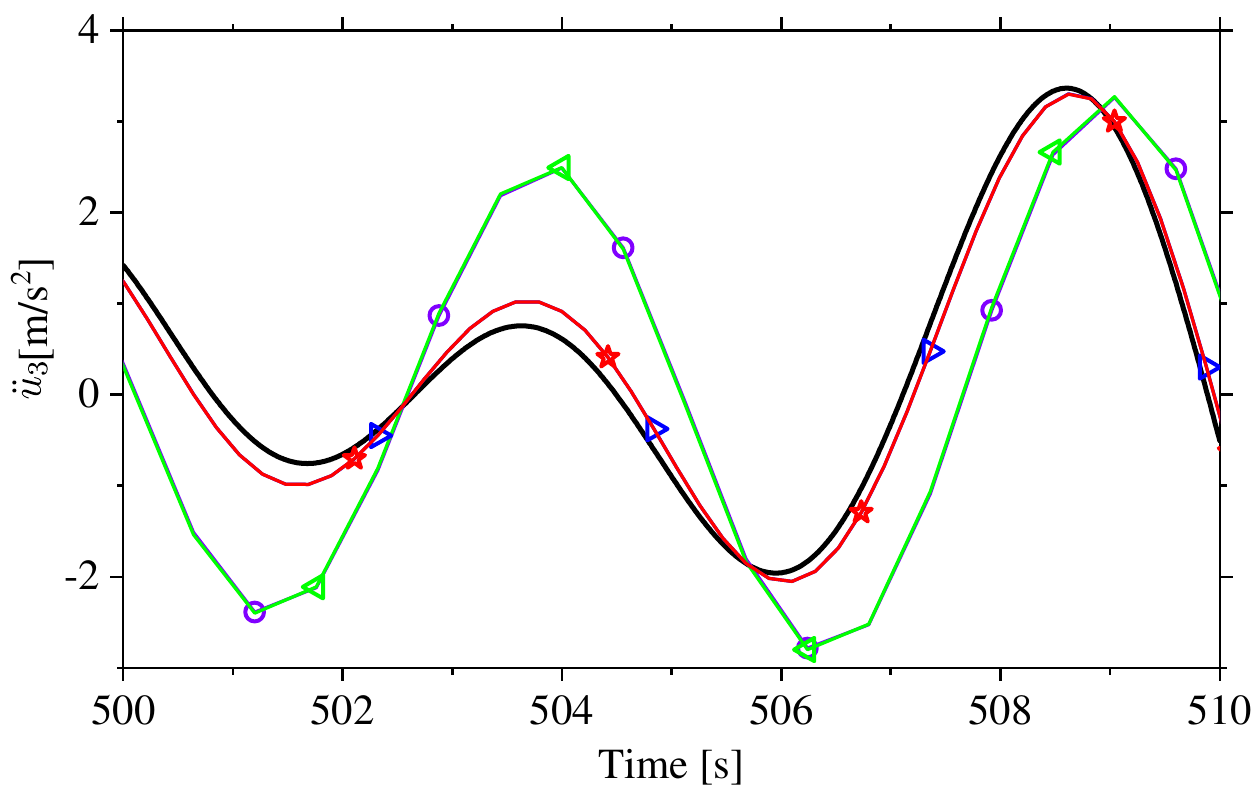}}
	\subfigure[ ]{
		\includegraphics[scale=0.37]{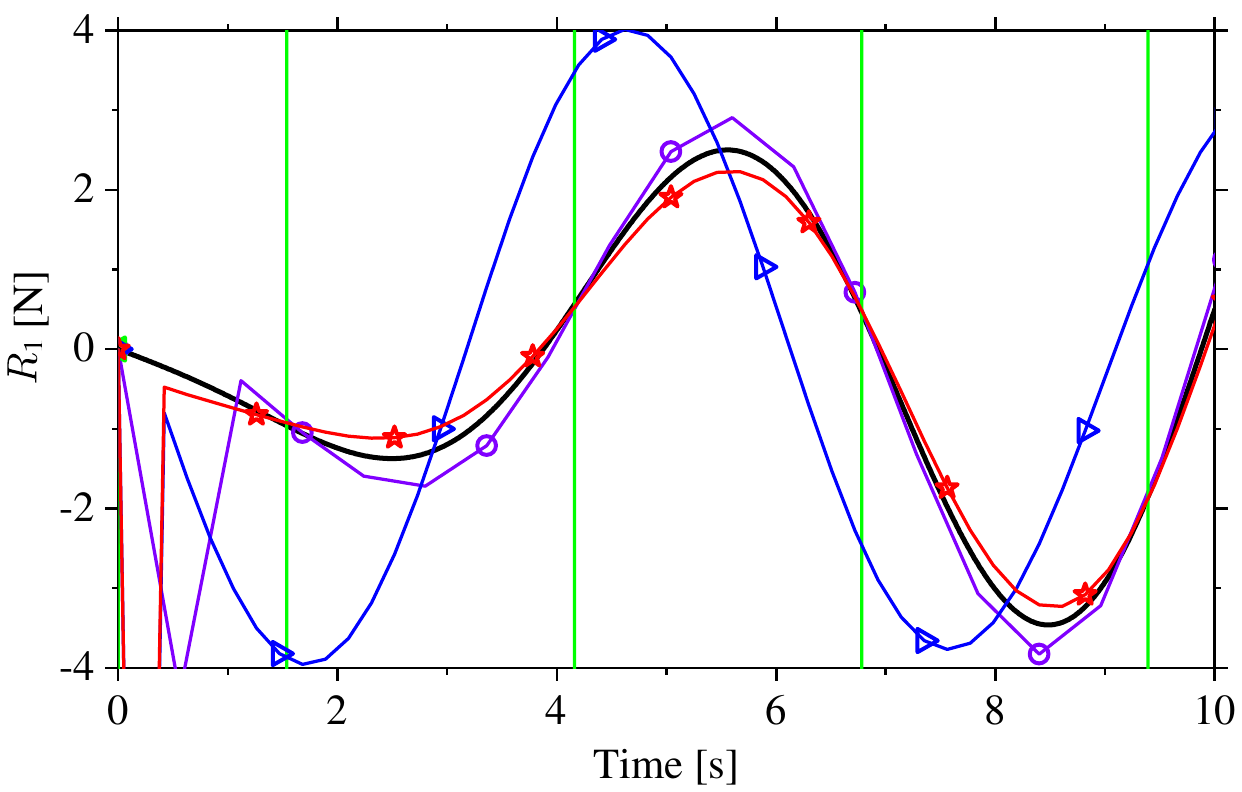}}
	\subfigure[ ]{
		\includegraphics[scale=0.37]{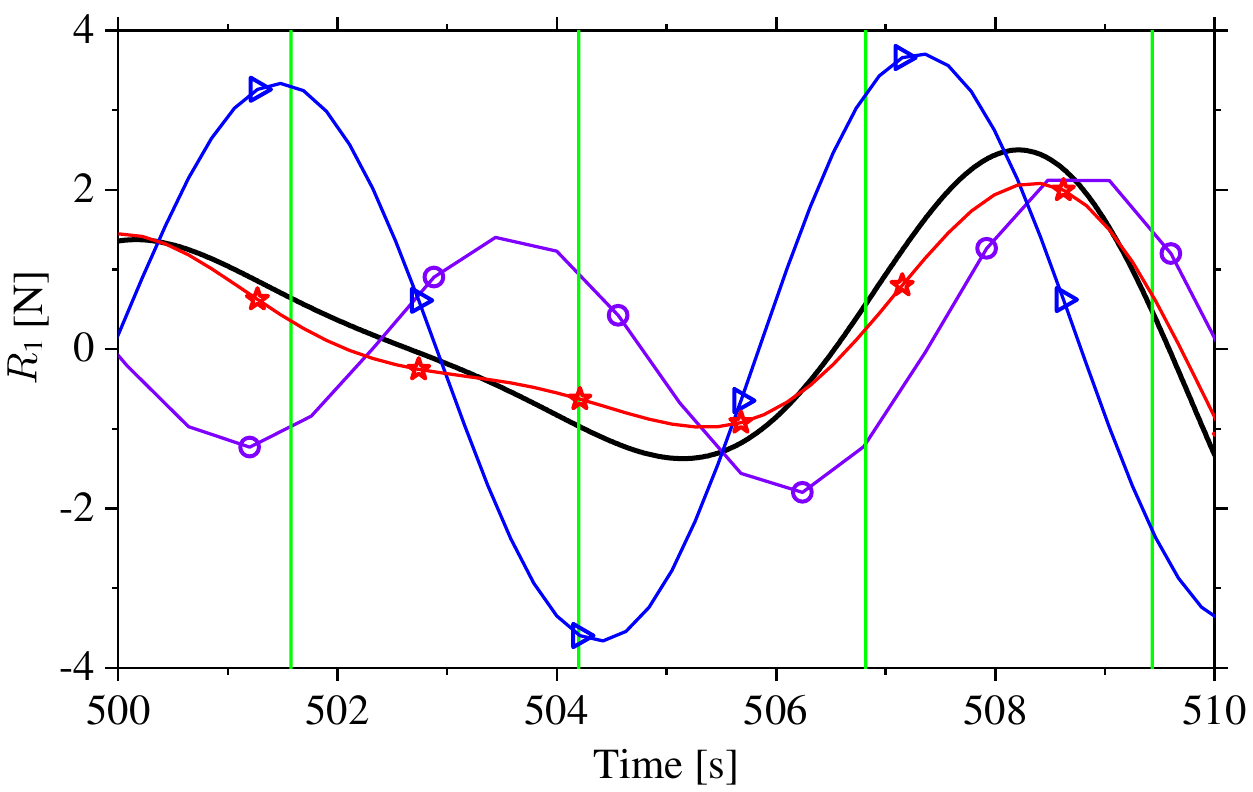}}
	\caption{Numerical displacements, velocities, and accelerations at node 3 and reaction force $R_1$,  predicted by various third-order implicit algorithms with $\rhoinf=0.0$ and $\dt=0.07$s.}
	\label{fig:2dof_3rd_u3_1}
\end{figure}

In the context of the asymptotically annihilating case ($\rhoinf=0.0$), it becomes evident from the findings presented in Figs.~\ref{fig:2dof_3rd_u2_1} and \ref{fig:2dof_3rd_u3_1} that the compared third-order algorithms exhibit enhanced predictive accuracy. It is noteworthy that the third-order $\rhoinf$-Bathe algorithm \cite{choiTimeSplitting2022} necessitates the utilization of complex-valued parameters to achieve $\rhoinf=0.0$. Consequently, its time step is set to $8\dt$, a choice aligned with the approach employed by Choi et al.~\cite{choiTimeSplitting2022} in solving the same model. Figs.~\ref{fig:2dof_3rd_u2_1} and \ref{fig:2dof_3rd_u3_1} elucidate that TR-CS \cite{fungComplextimestepNewmark1998,fungUnconditionallyStable1997} and DSUCI3 \cite{liDirectlySelfstarting2022} manifest impractical accelerations at node 2 and an aberrant reaction force at node 1. Furthermore, the $\rhoinf$-Bathe algorithm with $8\dt$ exhibits substantial amplitude and phase errors when predicting responses at node 3 during $t\in\left[500,~510\right]$s, as shown in Fig.~\ref{fig:2dof_3rd_u3_1}(b,~d,~f). In stark contrast to employing moderate dissipation, SUCI3 emerges as the algorithm with the highest accuracy in solving the model among the compared third-order algorithms when subjected to the most level of numerical high-frequency dissipation.
%When considering the asymptotically annihilating case ($\rhoinf=0.0$), various third-order algorithms can predict more accurate responses, as shown in Figs.~\ref{fig:2dof_3rd_u2_1} and \ref{fig:2dof_3rd_u3_1}. It should be emphasized that the third-order $\rhoinf$-Bathe algorithm \cite{choiTimeSplitting2022} must employ complex-valued parameters to produce $\rhoinf=0.0$, so its time step is set as $8\dt$ that Choi et al. \cite{choiTimeSplitting2022} also adopted to solve this model. Figs.~\ref{fig:2dof_3rd_u2_1} and \ref{fig:2dof_3rd_u3_1} illustrate that TR-CS \cite{fungComplextimestepNewmark1998,fungUnconditionallyStable1997} and DSUCI3 \cite{liDirectlySelfstarting2022} predict unreasonable accelerations at node 2 and the reaction force at node 1. The $\rhoinf$-Bathe algorithm with $8\dt$ produces significant amplitude and phase errors for predicting responses at node 3 during $t\in\left[500,~510\right]$s,  In contrast to using moderate dissipation, SUCI3 predicts the most accurate responses among the compared third-order algorithms when imposing the most level of numerical high-frequency dissipation.

%Hence, as stressed by Choi et al. \cite{choiTimeSplitting2022}, the third-order $\rhoinf$-Bathe method with complex-valued parameters would be useful in cases when the time step size must be small; (iii) SUCI3 predicts competitive solutions with the $\rhoinf$-Bathe algorithm with $3\dt$ and requires less computational cost than the latter.

\begin{figure}[htbp]
	\centering
	\subfigtopskip=2pt %?????????????????
	\subfigbottomskip=-4pt %??????????????????????????????
	\subfigcapskip=-5pt %?????????????
	{
		\includegraphics[scale=0.37]{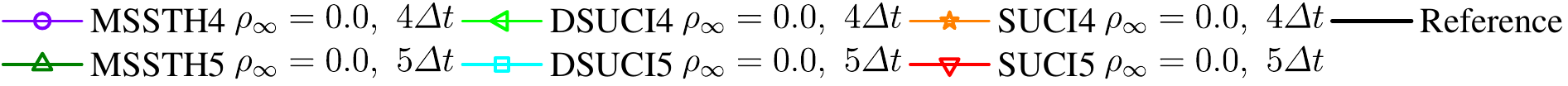}}
	\subfigure[ ]{
		\includegraphics[scale=0.37]{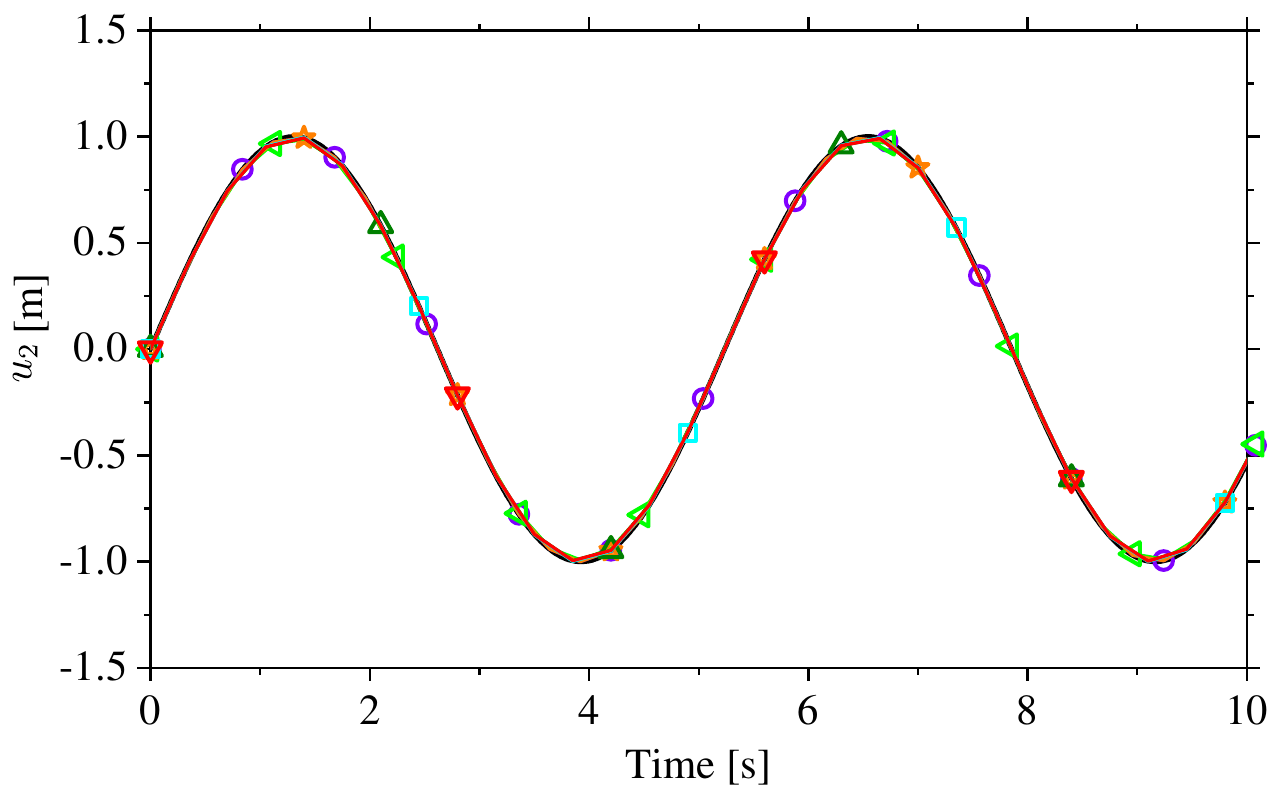}}
	\subfigure[ ]{
		\includegraphics[scale=0.37]{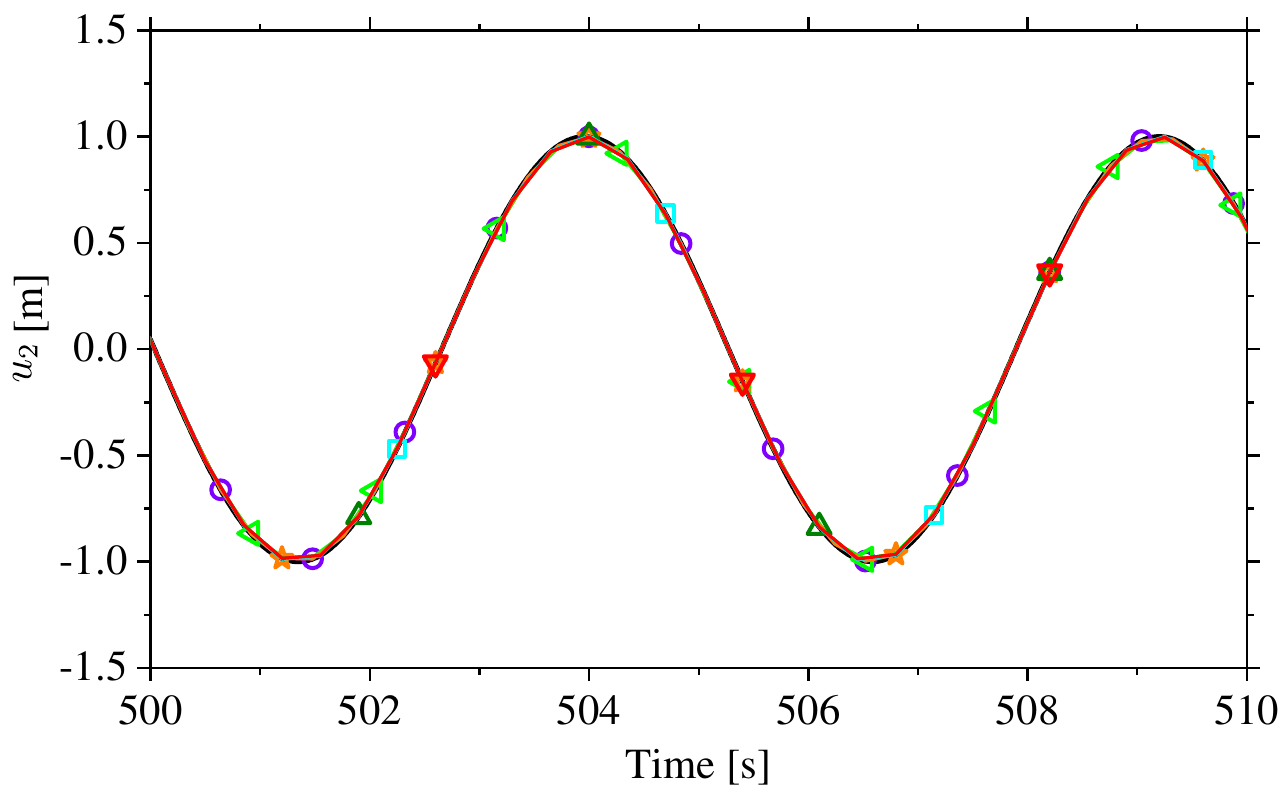}}
	\subfigure[ ]{
		\includegraphics[scale=0.37]{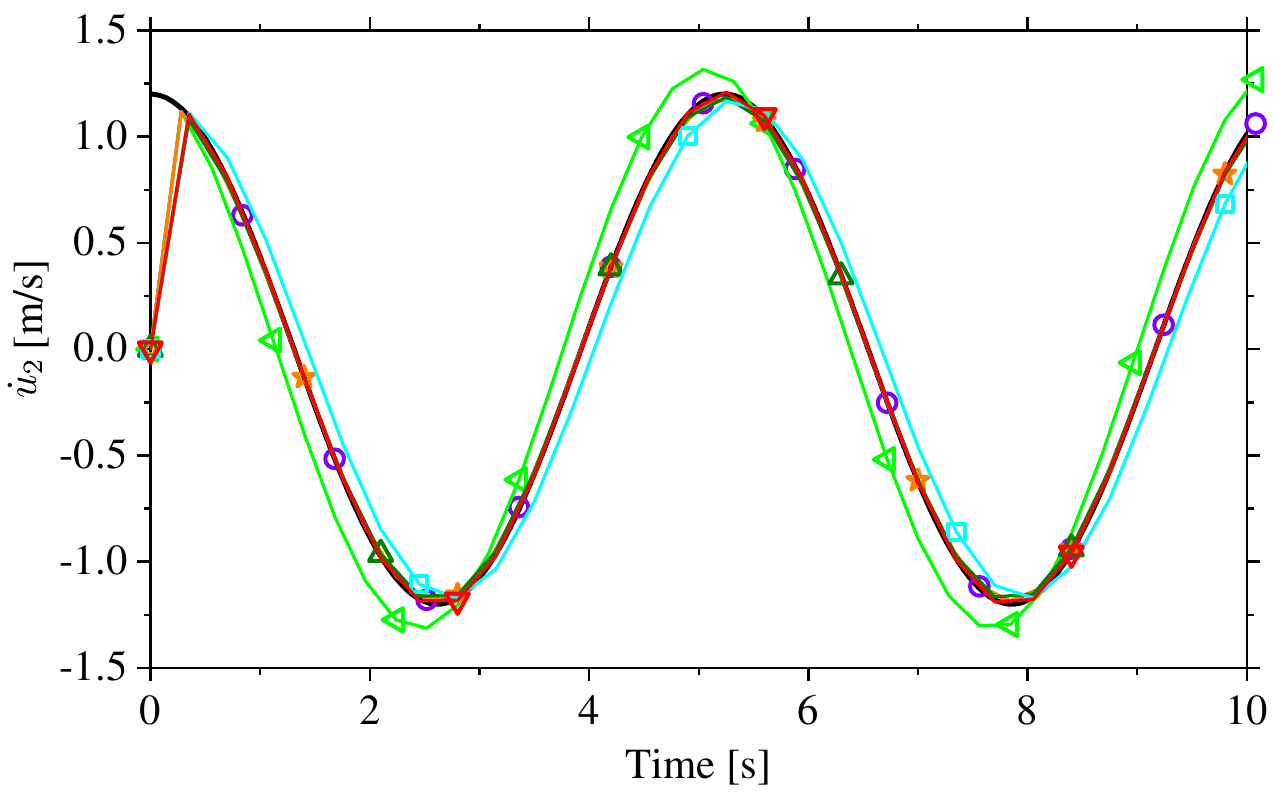}}
	\subfigure[ ]{
		\includegraphics[scale=0.37]{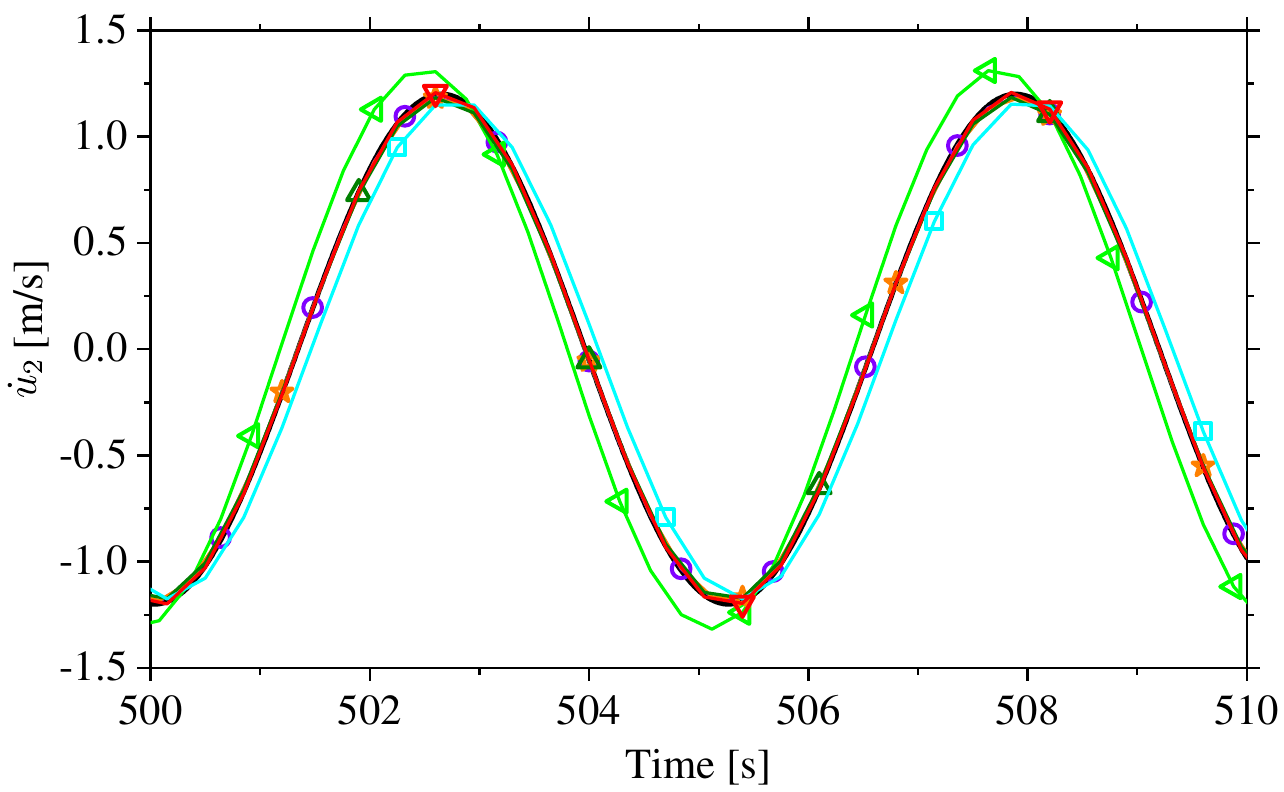}}
	\subfigure[ ]{
		\includegraphics[scale=0.37]{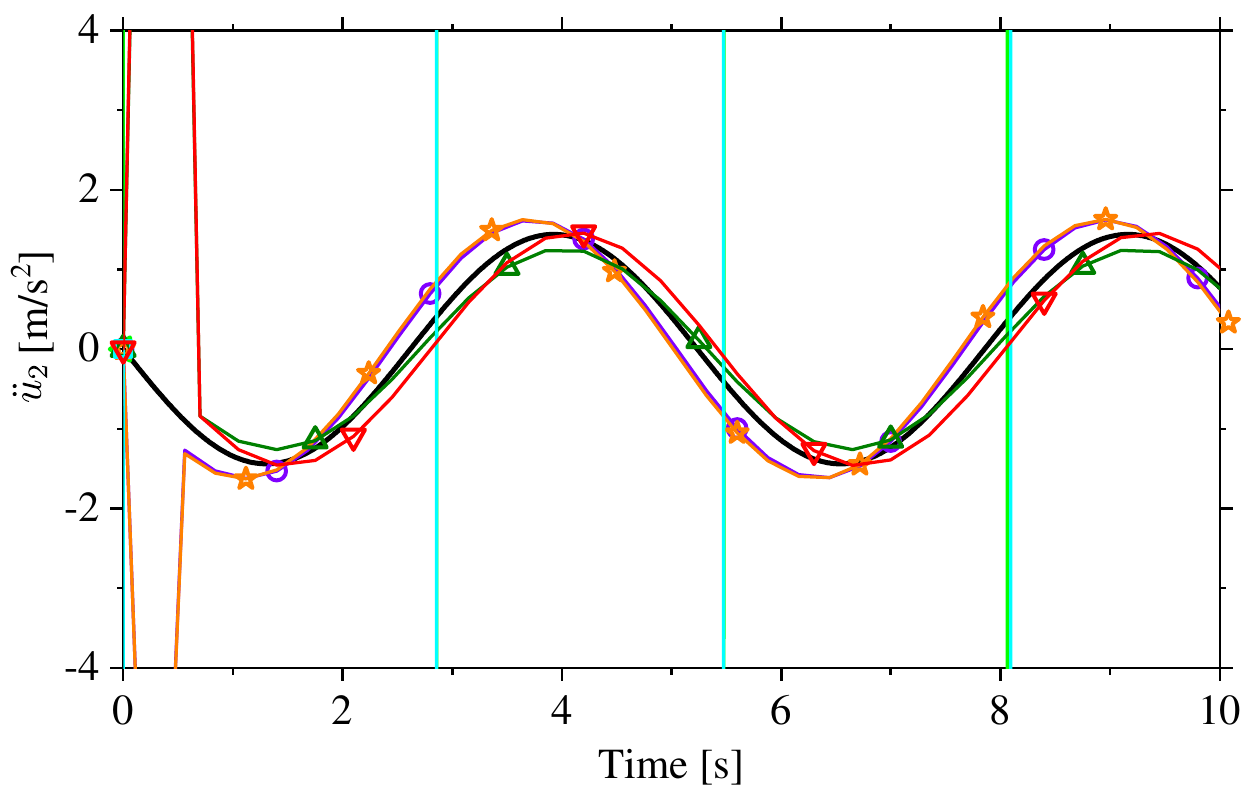}}
	\subfigure[ ]{
		\includegraphics[scale=0.37]{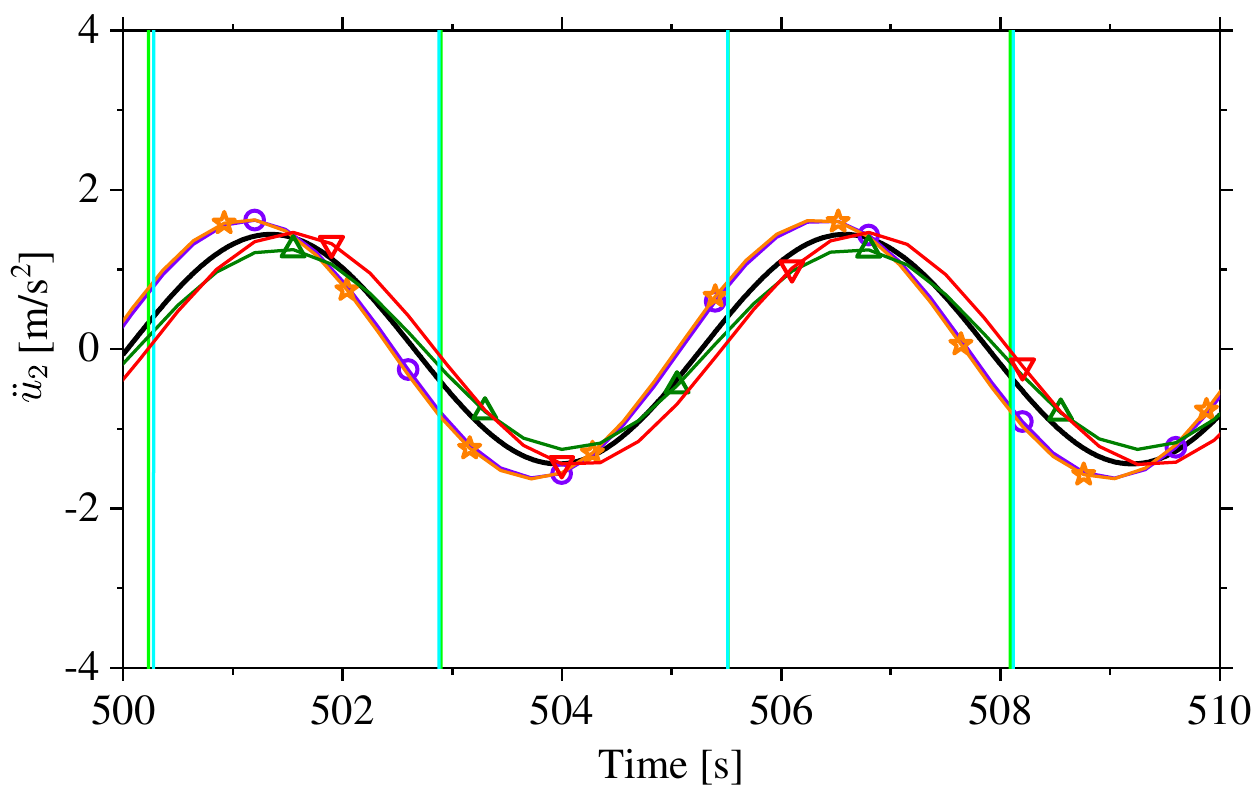}}
	\caption{Numerical displacements, velocities, and accelerations at node 2,  predicted by fourth- and fifth-order implicit algorithms with $\rhoinf=0.0$ and $\dt=0.07$s.}
	\label{fig:2dof_4th_u2}
\end{figure}

\begin{figure}[htbp]
	\centering
	\subfigtopskip=2pt %?????????????????
	\subfigbottomskip=-4pt %??????????????????????????????
	\subfigcapskip=-5pt %?????????????
	{
		\includegraphics[scale=0.37]{2dof_4th_leg_1}}
	\subfigure[ ]{
		\includegraphics[scale=0.37]{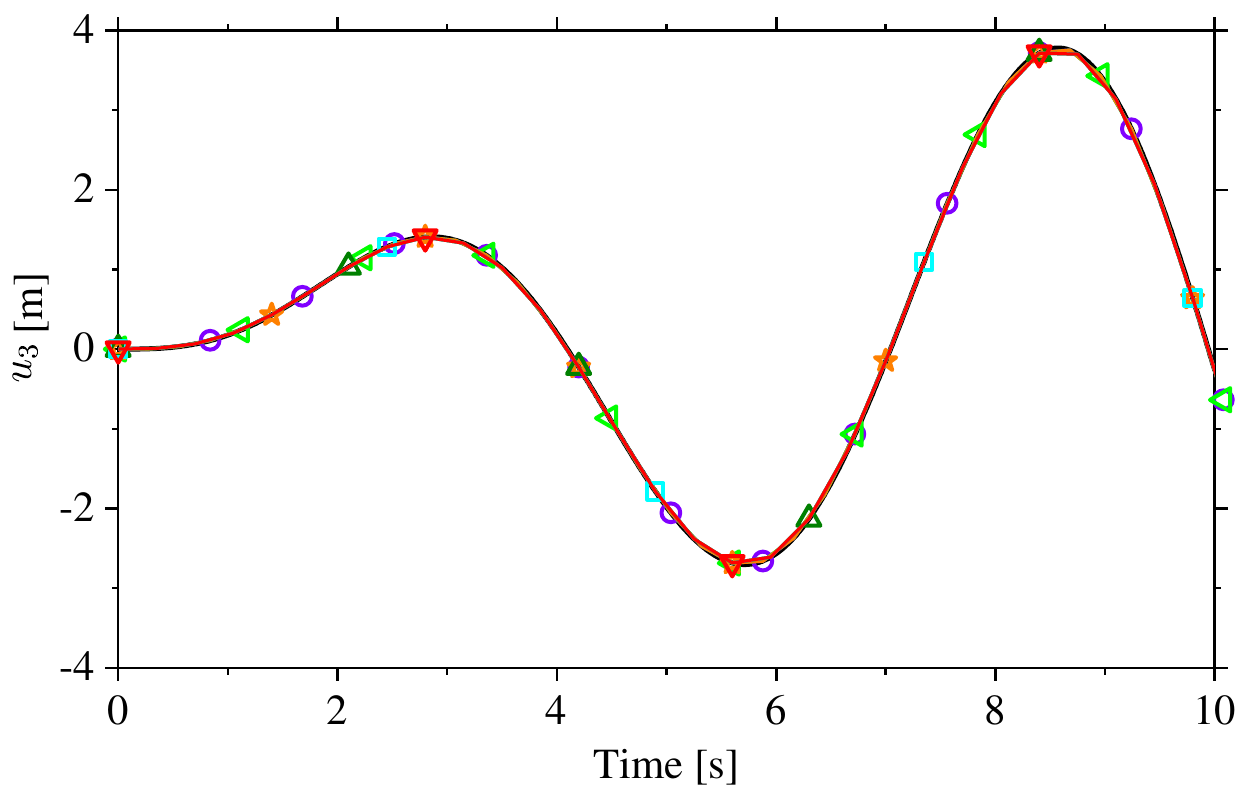}}
	\subfigure[ ]{
		\includegraphics[scale=0.37]{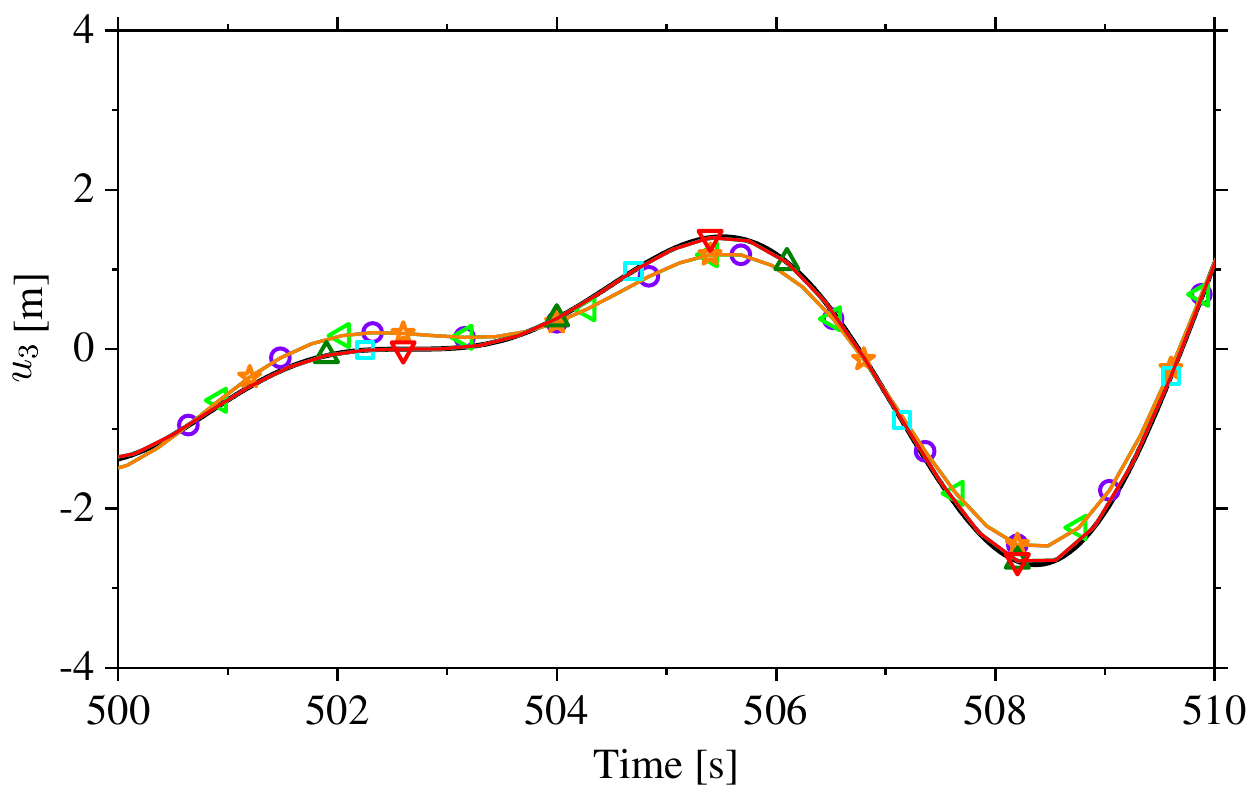}}
	\subfigure[ ]{
		\includegraphics[scale=0.37]{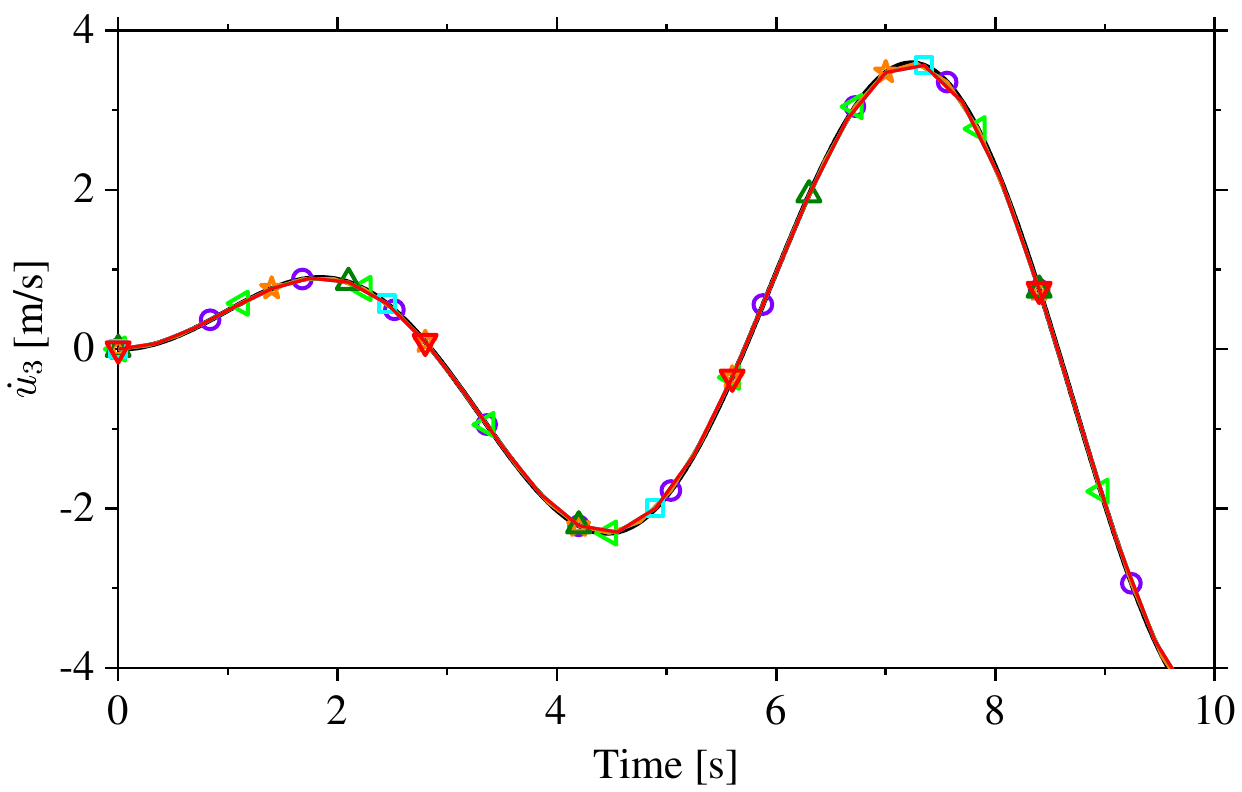}}
	\subfigure[ ]{
		\includegraphics[scale=0.37]{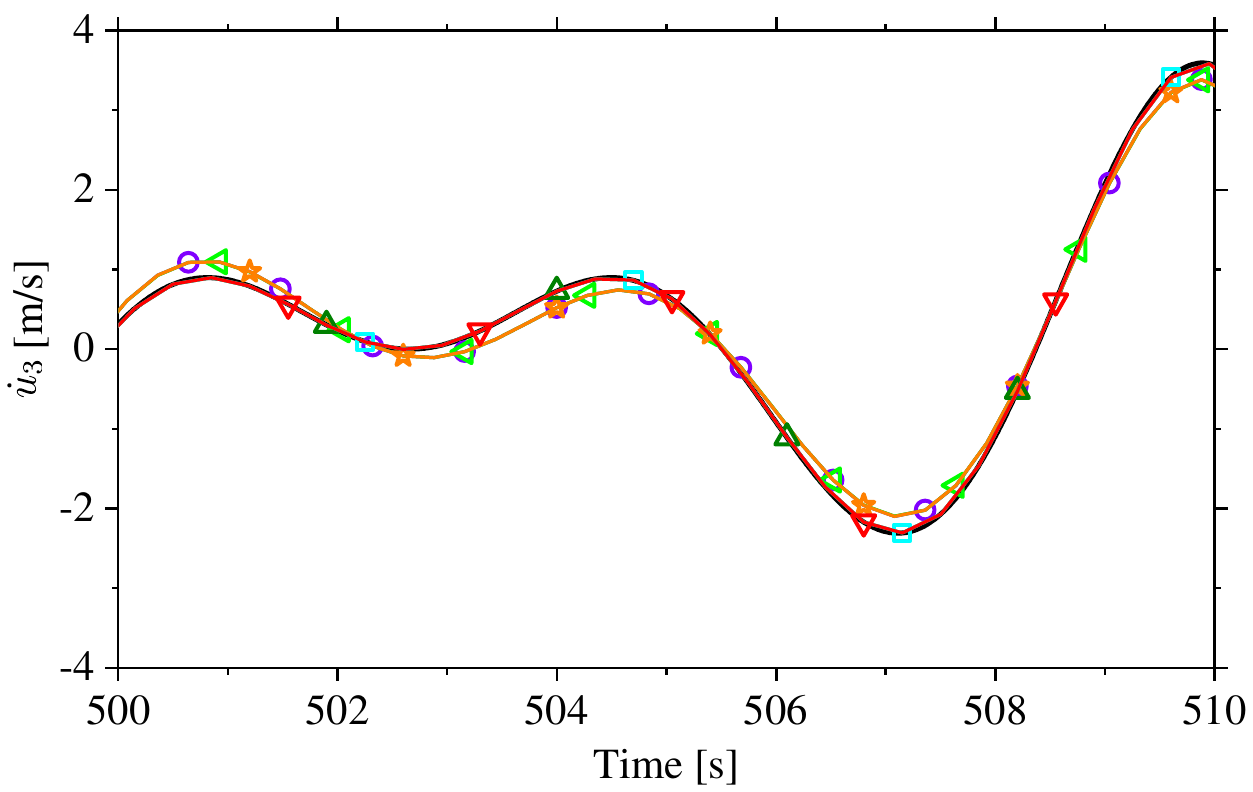}}
	\subfigure[ ]{
		\includegraphics[scale=0.37]{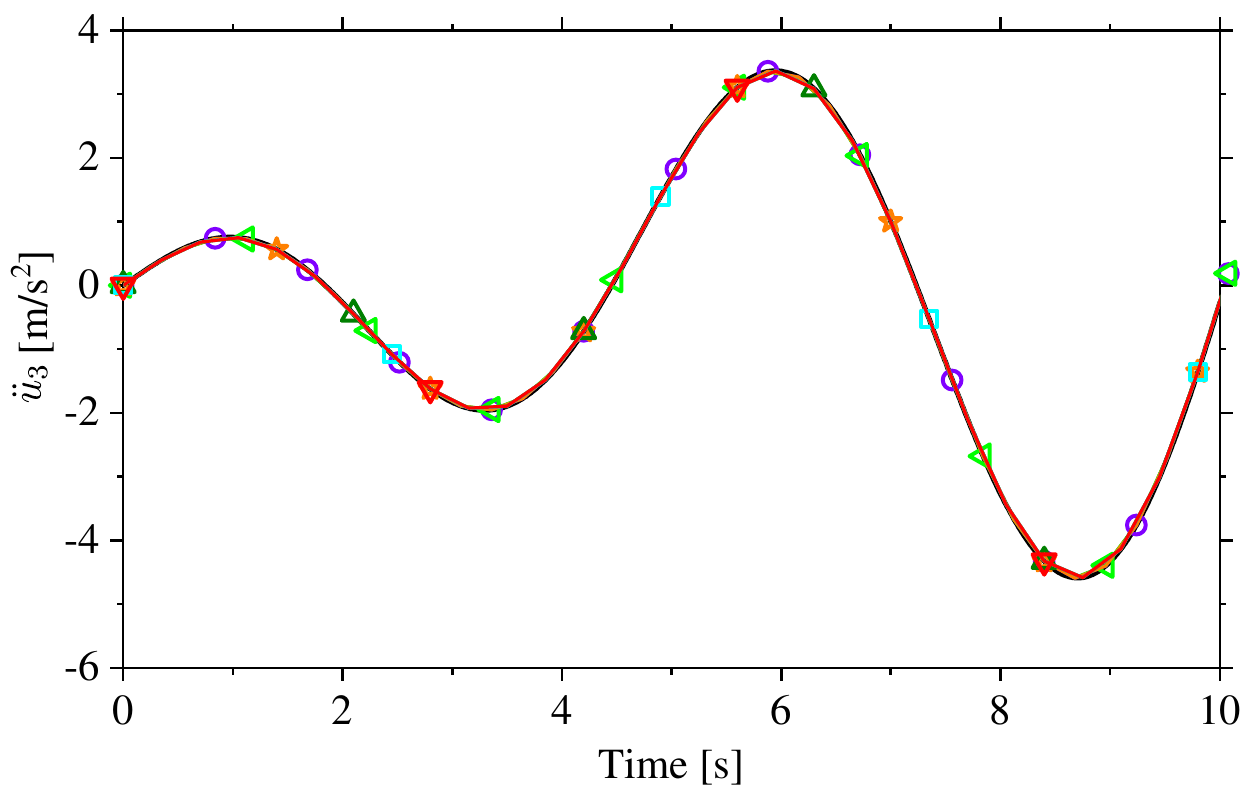}}
	\subfigure[ ]{
		\includegraphics[scale=0.37]{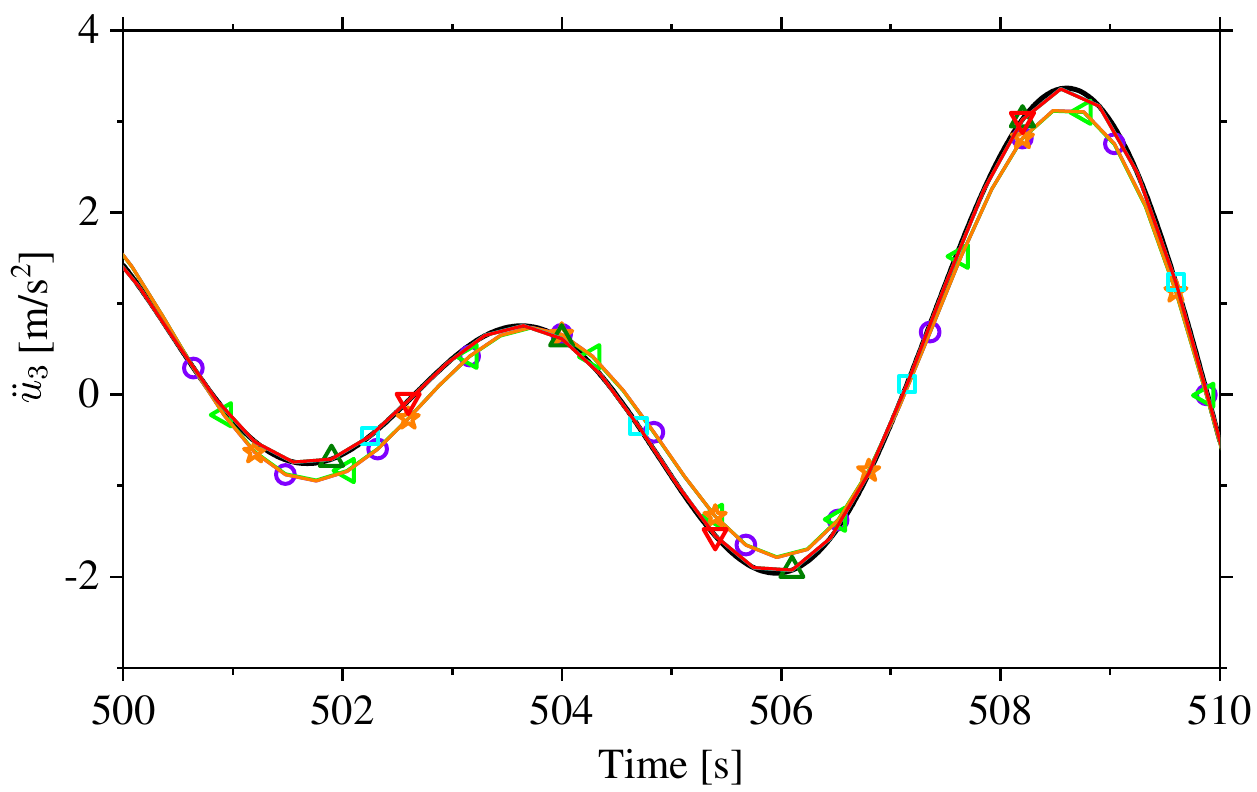}}
	\subfigure[ ]{
		\includegraphics[scale=0.37]{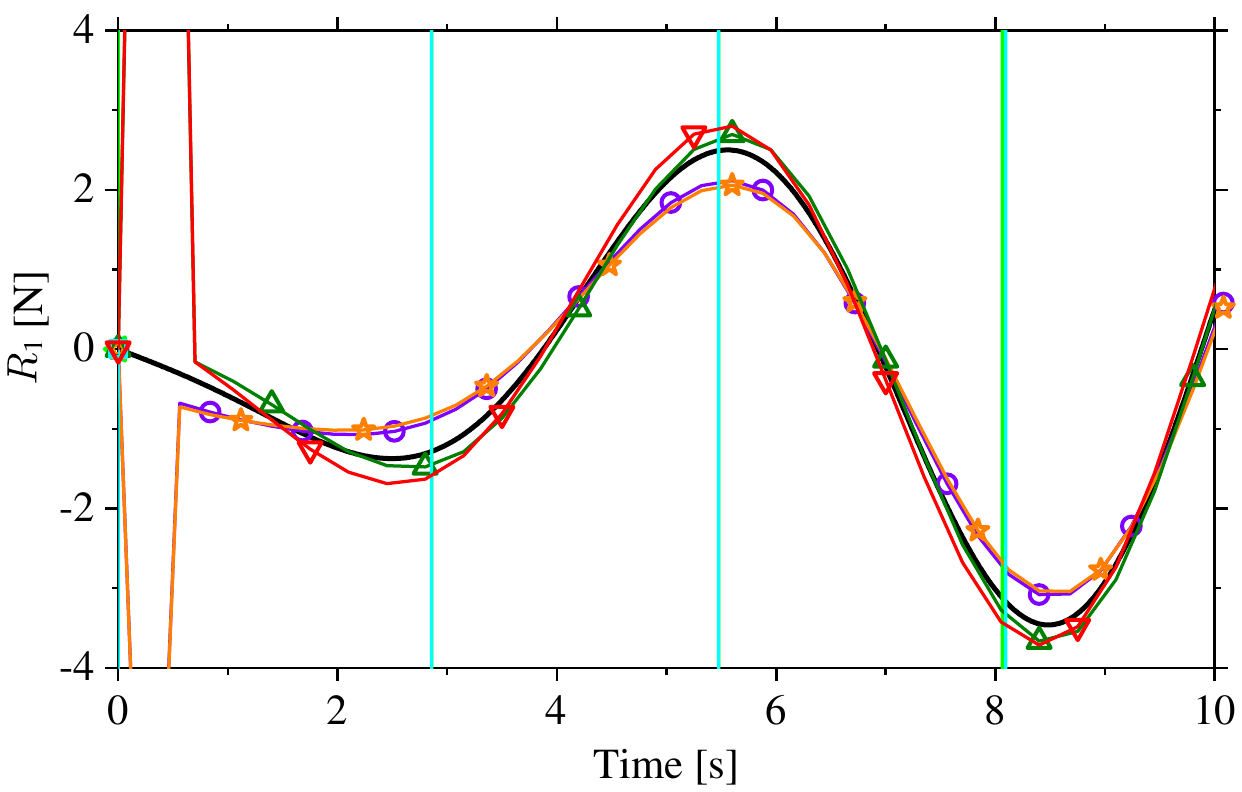}}
	\subfigure[ ]{
		\includegraphics[scale=0.37]{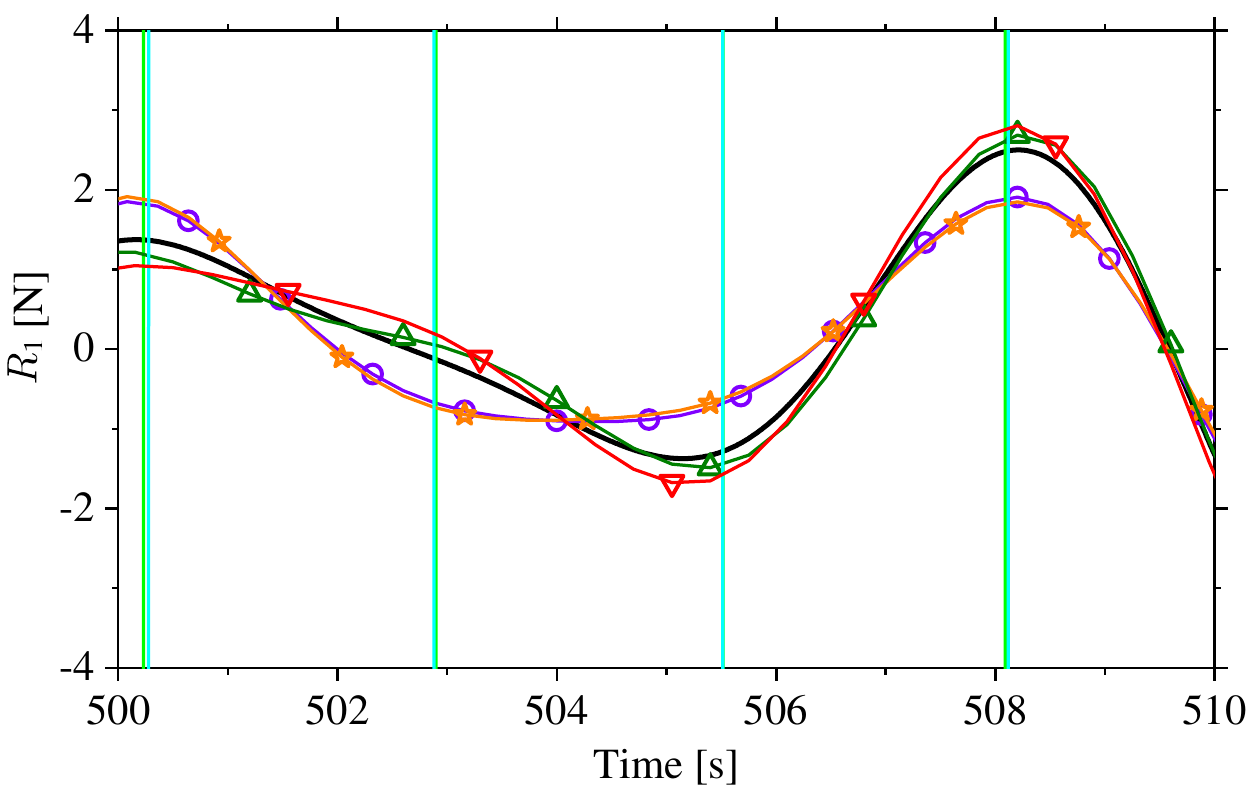}}
	\caption{Numerical displacements, velocities, and accelerations at node 3 and reaction force $R_1$,  predicted by fourth- and fifth-order implicit algorithms with $\rhoinf=0.0$ and $\dt=0.07$s.}
	\label{fig:2dof_4th_u3}
\end{figure}

Figs.~\ref{fig:2dof_4th_u2} and \ref{fig:2dof_4th_u3} depict the numerical responses predicted by fourth- and fifth-order algorithms with $\rhoinf=0.0$ and $\dt=0.07$s. It is discerned that DSUCI$n$ \cite{liDirectlySelfstarting2022} consistently exhibits inferior performance compared to MSSTH$n$ \cite{zhangOptimizationNsubstep2020} and SUCI$n$ when addressing the solution of Eq.~\eqref{eq:standardproblem}. Notably, the fifth-order MSSTH5 \cite{zhangOptimizationNsubstep2020} and SUCI5 algorithms demonstrate slightly more accurate predictions for the responses at node 3 in comparison to their fourth-order counterparts. The disparities in numerical solutions produced by DSUCI6 \cite{liDirectlySelfstarting2022} and SUCI6 closely mirror those observed in Figs.~\ref{fig:2dof_4th_u2} and \ref{fig:2dof_4th_u3}; hence, for the sake of brevity, these results are omitted to save the length of the paper.

%Figs.~\ref{fig:2dof_4th_u2} and \ref{fig:2dof_4th_u3} plot numerical responses predicted by the fourth- and fifth-order algorithms with $\rhoinf=0.0$ and $\dt=0.07$s. It is found that DSUCI$n$ \cite{liDirectlySelfstarting2022} always performs worse than MSSTH$n$ \cite{zhangOptimizationNsubstep2020} and SUCI$n$ for solving Eq.~\eqref{eq:standardproblem}, and the fifth-order MSSTH5 \cite{zhangOptimizationNsubstep2020} and SUCI5 algorithms predict slightly more accurate responses of node 3 than the fourth-order schemes. The numerical differences between DSUCI6 \cite{liDirectlySelfstarting2022} and SUCI6 are similar to those in Figs.~\ref{fig:2dof_4th_u2} and \ref{fig:2dof_4th_u3}, so they are omitted to save the length of the paper. 


\begin{figure}[htbp]
	\centering
	\subfigtopskip=2pt %?????????????????
	\subfigbottomskip=-4pt %??????????????????????????????
	\subfigcapskip=-5pt %?????????????
	\subfigure[ ]{
		\includegraphics[scale=0.45]{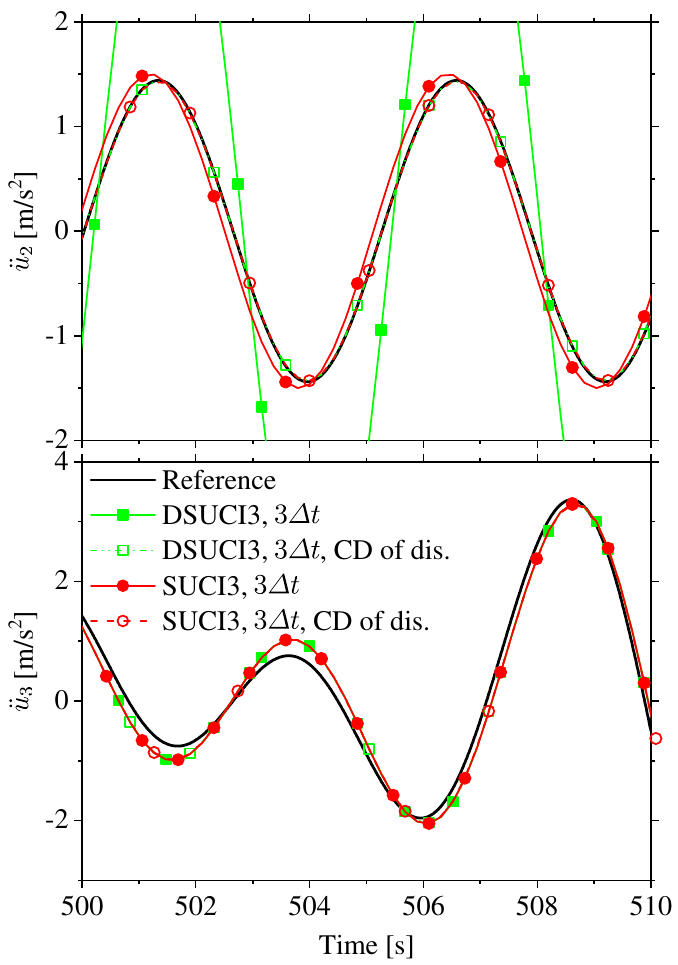}}
	\subfigure[ ]{
		\includegraphics[scale=0.45]{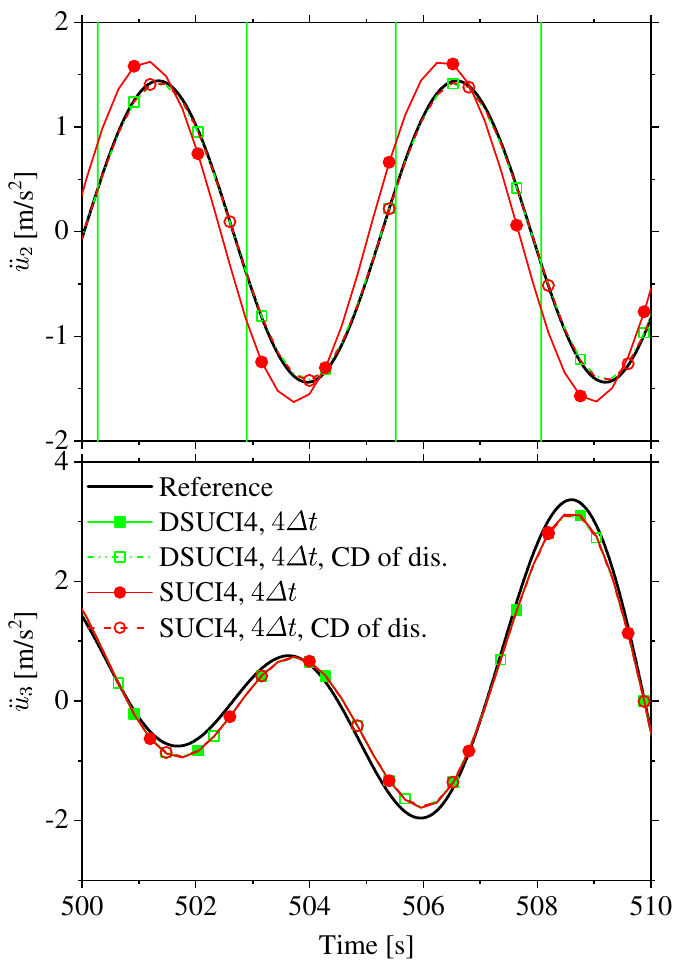}}
	\subfigure[ ]{
		\includegraphics[scale=0.45]{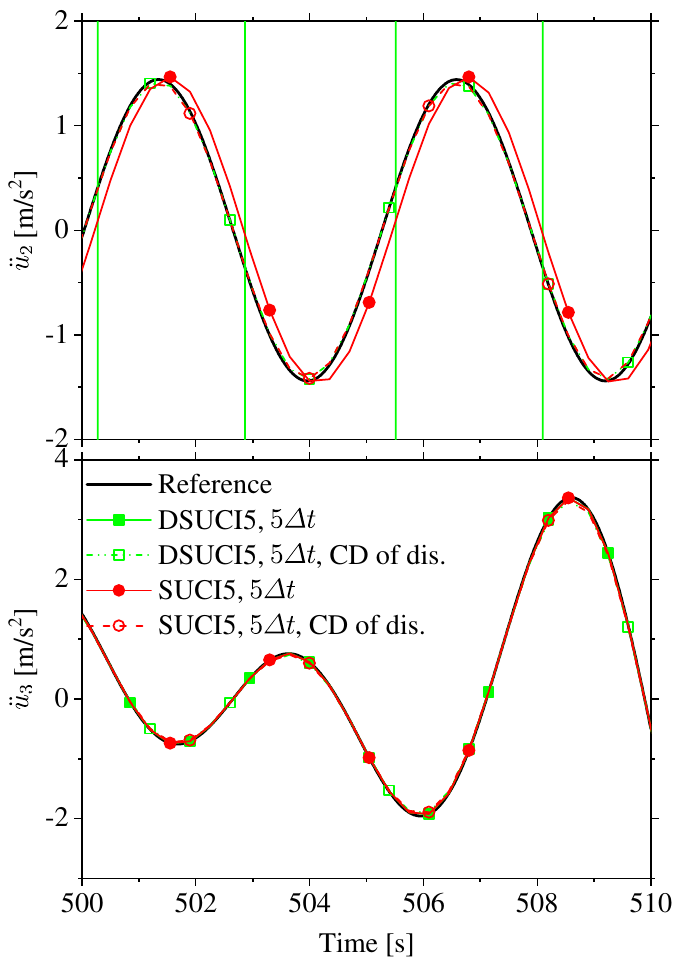}}
	\caption{Numerical accelerations predicted by DSUCI$n$ \cite{liDirectlySelfstarting2022} and SUCI$n$ with $\rhoinf=0.0$ and $\dt=0.07$s.}
	\label{fig:2dof_cd}
\end{figure}

The aforementioned findings highlight that when employing the original strategy for outputting accelerations, DSUCI$n$ \cite{liDirectlySelfstarting2022} can attain the identical high-order accuracy. However, its efficacy diminishes in accurately capturing the steady-state acceleration responses in $\ddot{u}_2$ when addressing the current problem. Conversely, when utilizing the central difference (CD) of displacement, both DSUCI$n$ and SUCI$n$ exhibit the capability to predict accurate acceleration responses ($\ddot{u}_2$), as demonstrated in Fig.~\ref{fig:2dof_cd}. This accuracy stems from the fact that the predicted displacements ($u_2$) are predominantly influenced by the low-frequency steady-state response, thereby facilitating the accurate tracking of low-frequency steady-state acceleration responses. Furthermore, Fig.~\ref{fig:2dof_cd} elucidates that the central difference technique does not yield improvements in the acceleration at node 3. This is attributed to the fact that the acceleration at node 3 is not governed by high-frequency modes, and consequently, the original approach \cite{liDirectlySelfstarting2022} remains proficient in generating accurate accelerations in this context.
%The findings above illustrate that DSUCI$n$ \cite{liDirectlySelfstarting2022} using the original way to output accelerations can achieve identical high-order accuracy, but it cannot follow the steady-state acceleration responses in $\ddot{u}_2$ well for solving the present problem. Alternatively, using the central difference (CD) of displacement, DSUCI$n$ as well as SUCI$n$ provides quite accurate acceleration responses ($\ddot{u}_2$) for users; see Fig.~\ref{fig:2dof_cd}. This is because the predicted displacements ($u_2$) are dominated by the low-frequency steady-state response and thus the resulting accelerations naturally follow the low-frequency steady-state response well. Fig.~\ref{fig:2dof_cd} also reveals that the central difference technique does not improve the acceleration at node 3, as it is not governed by the high-frequency mode and thus the original way \cite{liDirectlySelfstarting2022} can output accurate accelerations.

%This standard double-degree-of-freedom mass-spring system \cite{choiTimeSplitting2022,batheInsightImplicit2012}  demonstrates that when considering the same order of accuracy, SUCI$n$ can perform better than some existing high-order algorithms with respect to the solution accuracy and dissipation control. In particular, using the original parameter settings given in \cite{liDirectlySelfstarting2022}, the directly self-starting DSUCI$n$ algorithms cannot predict accurate numerical responses for solving this benchmark problem (see Appendix C for more discussions), although they share almost the same spectral properties as SUCI$n$. Therefore, the authors recommend the use of SUCI$n$ for dynamic problems with stiff and flexible components.
The standard double-degree-of-freedom mass-spring system \cite{choiTimeSplitting2022,batheInsightImplicit2012} serves as an illustrative example revealing that, when considering the same order of accuracy, SUCI$n$ can outperform certain existing high-order algorithms in terms of both solution accuracy and dissipation control. Specifically, under the original parameter settings outlined in \cite{liDirectlySelfstarting2022}, the directly self-starting DSUCI$n$ algorithms exhibit limitations in accurately predicting numerical responses for solving this benchmark problem. Notably, despite possessing nearly identical spectral properties to SUCI$n$, DSUCI$n$ fails to achieve the desired accuracy in numerical predictions. Consequently, based on these observations, the authors advocate for the adoption of SUCI$n$ when addressing dynamic problems featuring both stiff and flexible components. This recommendation stems from the demonstrated superior performance of SUCI$n$ in achieving accurate and controlled solutions compared to its counterparts, particularly DSUCI$n$, in the context of the discussed mass-spring system.

\subsection{Nonlinear dynamics}
This section will solve three nonlinear examples to show the numerical performance of SUCI$n$ on nonlinear dynamics. For solving nonlinear dynamics of form 
	\begin{equation}\label{eq:non_MCK}
		\mbfM\mbfa(t)+\mbf{R}(\mbfv(t),~\mbfu(t))=\mbfF(t)
	\end{equation}
	where $\mbf{R}$ includes damping and stiffness nonlinearities, the novel method \eqref{eq:nsubstep} requires an iterative scheme, such as the Newton-Raphson scheme used herein, to find satisfactory numerical solutions, as shown in Algorithm \ref{code:suci_nonlinear}. The tangent damping $\mbfC_{n+\gamma_i}^{(k)} $ and stiffness $\mbfK_{n+\gamma_i}^{(k)}$ matrices at the $k$th iteration are calculated, respectively, by
	\begin{equation}\label{key}
		\mbfC_{n+\gamma_i}^{(k)}=\dfrac{\partial\mbf{R}\left(\mbfv_{n+\gamma_i}^{(k)},~\mbfu_{n+\gamma_i}^{(k)}\right)}{\partial \mbfv_{n+\gamma_i}^{(k)}}\quad\text{and}\quad\mbfK_{n+\gamma_i}^{(k)}=\dfrac{\partial\mbf{R}\left(\mbfv_{n+\gamma_i}^{(k)},~\mbfu_{n+\gamma_i}^{(k)}\right)}{\partial \mbfu_{n+\gamma_i}^{(k)}}.
	\end{equation} 
	Note that Algorithm \ref{code:suci_nonlinear} uses the Newton-Raphson iterative scheme, and some mild modifications should be made if other iterative schemes, such as the modified Newton-Raphson, are used. The convergence in equilibrium iterations is reached when one of the following inequalities is satisfied:
\begin{equation}\label{key}
	\big|\big|\mbf{r}^{(k)}\big|\big|_2\le\text{RTOL}\quad\text{or}\quad \big|\big|\otherDelta\ddot{\mbf{U}}^{(k)}\big|\big|_2\le\text{ATOL}
\end{equation}
where $ \mbf{r}^{(k)} $ and $ \otherDelta\ddot{\mbf{U}}^{(k)} $ are the residuals and increment of the acceleration vector in the $ k $th iteration, respectively; RTOL and ATOL are two user-specified convergence tolerances and both of them are given as $ 1.0\times 10^{-8} $ in this paper; $ \big|\big|\cdot\big|\big|_2 $
denotes the Euclidean norm, also known as the 2-norm.

%{\fr 
	\begin{breakablealgorithm}%[h]
		\caption{The novel $s$-sub-step implicit method \eqref{eq:nsubstep} for solving nonlinear problems \eqref{eq:non_MCK}}
		\begin{algorithmic}[1]
			\State Select the number of sub-steps $s\in\{1,~2,~3,~4,~5,~6\}$, set $\rhoinf\in[0,~1]$, and calculate $\gamma_i,~\alpha_{ij}$.
			\State Solve: $\mbfa_0$ by $\mbfM\mbfa_0=\mbfF(t_0)-\mbf{R}(\mbfv_0,~\mbfu_0)$.
			\For{$n=0$ to $n=N$} {\color{dcolor}\qquad\qquad\qquad\ \ \   //Loops for integration steps.}
			\For{$i=1$ to $i=s$}\label{alg:nextsubstep} {\color{dcolor}\qquad\qquad\qquad\  //Loops for $s$ sub-steps.}
			\State Predict: $\mbfa_{n+\gamma_i}^{(1)}=\mbfa_{n+\gamma_{i-1}}$. {\color{dcolor}\qquad\qquad   //Recall that $\gamma_0=0$, i.e., $\mbfa_{n+\gamma_0}=\mbfa_n$.}
			\For{$k=1$ to $\max$-$iter$} {\color{dcolor}\qquad\qquad //Loops for Newton-Raphson iterations.}
			\State Compute: $\mbfv_{n+\gamma_i}^{(k)}=\mbfv_n+\dt\left(\sum_{j=0}^{i-1}\alpha_{ij}\mbfa_{n+\gamma_j}+\alpha_{ii}\mbfa_{n+\gamma_i}^{(k)}\right)$.
			\State Compute: $\mbfu_{n+\gamma_i}^{(k)}=\mbfu_n+\dt\left(\sum_{j=0}^{i-1}\alpha_{ij}\mbfv_{n+\gamma_j}+\alpha_{ii}\mbfv_{n+\gamma_i}^{(k)}\right)$.
			\State Compute: $\mbf{R}\left(\mbfv_{n+\gamma_i}^{(k)},~\mbfu_{n+\gamma_i}^{(k)}\right)$ and $\widetilde{\mbf{K}}_{ii}=\mbfM+\alpha_{ii}\dt\mbfC_{n+\gamma_i}^{(k)}+{\alpha_{ii}}^2\dt^2\mbfK_{n+\gamma_i}^{(k)}$.
			\State Solve: $\varDelta \mbfa^{(k)}$ by $\widetilde{\mbf{K}}_{ii}\varDelta \mbfa^{(k)}=\mbfF(t_n+\gamma_i\dt)-\mbf{R}\left(\mbfv_{n+\gamma_i}^{(k)},~\mbfu_{n+\gamma_i}^{(k)}\right)-\mbfM\mbfa_{n+\gamma_i}^{(k)}$.
			\If{convergence}
			\State Compute: $\mbfa_{n+\gamma_i}=\mbfa_{n+\gamma_i}^{(k)}$,  $\mbfv_{n+\gamma_i}=\mbfv_{n+\gamma_i}^{(k)}$, and $\mbfu_{n+\gamma_i}=\mbfu_{n+\gamma_i}^{(k)}$.
			\State \textbf{Go to} step \ref{alg:nextsubstep}. {\color{dcolor}\qquad\qquad  //Performing solutions in the next sub-step scheme.}
			\Else
			\State Update: $\mbfa_{n+\gamma_i}^{(k)}\leftarrow\mbfa_{n+\gamma_i}^{(k)}+\varDelta \mbfa^{(k)}$.
			\EndIf
			\EndFor {\color{dcolor}\quad // Newton-Raphson iterations.}
			\EndFor {\color{dcolor}\qquad\ \  // Sub-step schemes.}
			\State {\color{dcolor}// Solutions at the discrete instants $t_n$ are provided by the $s$th sub-step scheme due to $\gamma_s=1$.}
			\EndFor
		\end{algorithmic}\label{code:suci_nonlinear}
	\end{breakablealgorithm}

\subsubsection{A simple pendulum}
A classical nonlinear simple pendulum \cite{kimNewFamily2017,rezaiee-pajandModifiedDifferential2017} is considered as the first SDOF system with the strong nonlinearity.
%\begin{figure}[htbp]
%	\centering
%	\includegraphics[scale=0.4]{pendulum}
%	\caption{The simple nonlinear pendulum with $ m=1.0 $.}
%	\label{fig:simplePendulum}
%\end{figure}
The governing equation of motion for this pendulum with length $ L=1 $m subject to a nonzero initial velocity $ \dot{\theta}_0 $ is written as
\begin{equation}\label{eq:kw}
	\ddot{\theta}(t)+\frac{g}{L}\sin(\theta(t))=0
\end{equation}
where $ g $ is the acceleration of gravity and assumed to be unity in this test so that $ g/L=1 $s$^{-2}$ is satisfied in Eq.~(\ref{eq:kw}). It is necessary to point out that $ \dot{\theta}_0=1.999999238456499 $rad/s has been widely used to synthesize the highly nonlinear pendulum \cite{kimNewFamily2017} (\textit{the minimum initial velocity to make the pendulum rotating is calculated as $ \dot{\theta}_{\min}=2 $rad/s}). In the present case, this simple pendulum should oscillate between two peak points ($ \theta=\pm\pi $) instead of making the complete rotation since the initial total energy is slightly less than the minimum total energy to make complete turns. The reference solutions of the pendulum are obtained from the work \cite{kimNewFamily2017}. 

%\begin{figure}[htpb]
%	\centering
%	\includegraphics[scale=0.65]{simPen1}\vskip -2mm
%	\caption{Numerical displacements and absolute errors of the simple pendulum using the TPO/G-$\alpha$ \cite{shaoThreeParameters1988,chungTimeIntegration1993} and SUCI2 \cite{liNovelFamily2020} algorithms with the same $ \dt=0.02 $s.}
%	\label{fig:simPen1}
%\end{figure}
\begin{figure}[htbp]
	\centering
	\includegraphics[scale=0.65]{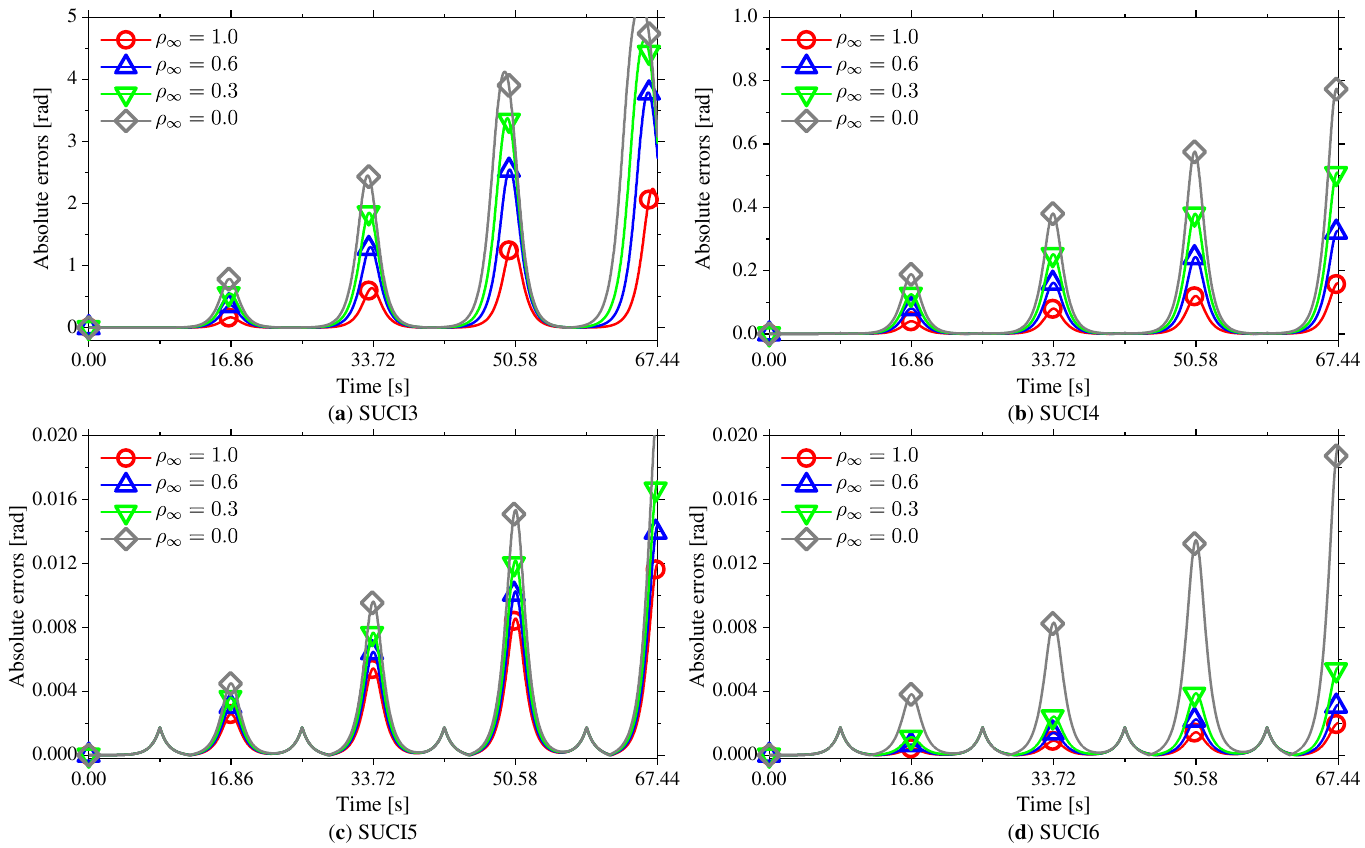}%\vskip -2mm
	\caption{Absolute displacement errors of the simple pendulum using SUCI$n$ with the same $ \dt=0.02 $s.}
	\label{fig:simPen2}
\end{figure}
%The two second-order accurate algorithms, namely single-step TPO/G-$\alpha$ \cite{shaoThreeParameters1988,chungTimeIntegration1993} and two-sub-step SUCI2 \cite{liNovelFamily2020}, are firstly used to solve this pendulum and their numerical results are presented in Fig.~\ref{fig:simPen1}. Obviously, using the current time step $ \dt=0.02 $s, these two algorithms cannot predict reasonable solutions. Particularly, SUCI2 provides the complete turns and TPO/G-$\alpha$ also presents adverse phase errors. To amplify differences among different algorithms or parameter settings, Fig.~\ref{fig:simPen1} also plots absolute displacement errors, and thus the complete turns make errors tend to infinity with the increase of time.

Fig.~\ref{fig:simPen2} presents the absolute displacement errors associated with the novel high-order algorithms, SUCI$n$, for the analysis of the simple pendulum. Two notable observations emerge from the results. Firstly, when employing the same time step size, the higher-order algorithms consistently outperform the lower-order ones for this nonlinear example. It is imperative to emphasize that this enhanced performance of high-order algorithms comes at the expense of heightened computational demands and increased coding complexity. In instances where the time step size is set as $ \dt=0.01\times n $, with $ n $ representing the number of sub-steps, the SUCI$n$ algorithms exhibit larger absolute errors. Nevertheless, they produce a consistent trend with the outcomes depicted in Fig.~\ref{fig:simPen2}. Secondly, it is noteworthy that the non-dissipative scheme, characterized by $ \rhoinf=1 $, yields significantly fewer numerical errors in comparison to the dissipative ones. This discrepancy arises due to the absence of spurious high-frequency components in this simple pendulum.
%The absolute displacement errors of the simple pendulum using the novel high-order algorithms (SUCI$n$) are given in Fig.~\ref{fig:simPen2}, where two expected facts are revealed. One is that when using the same time step size the higher-order algorithms generally perform better than the lower-order ones for solving nonlinear dynamics. It is necessary to point out that such an advantage of high-order algorithms is at the cost of more computational costs and coding complexity. When adopting the time step size as $ \dt=0.01\times n $ where $ n $ denotes the number of sub-steps, SUCI$n$ predicts larger absolute errors but still follow the same tendency as those in Fig.~\ref{fig:simPen2}. The other is that the non-dissipative ($ \rhoinf=1 $) scheme produces significantly fewer numerical errors than the dissipative schemes since there are no spurious high-frequency components in this simple pendulum.

\begin{figure}[htbp]
	\centering
	\subfigtopskip=2pt %?????????????????
	\subfigbottomskip=-4pt %??????????????????????????????
	\subfigcapskip=-5pt %?????????????
	\subfigure[Third-order algorithms]{
		\includegraphics[scale=0.7]{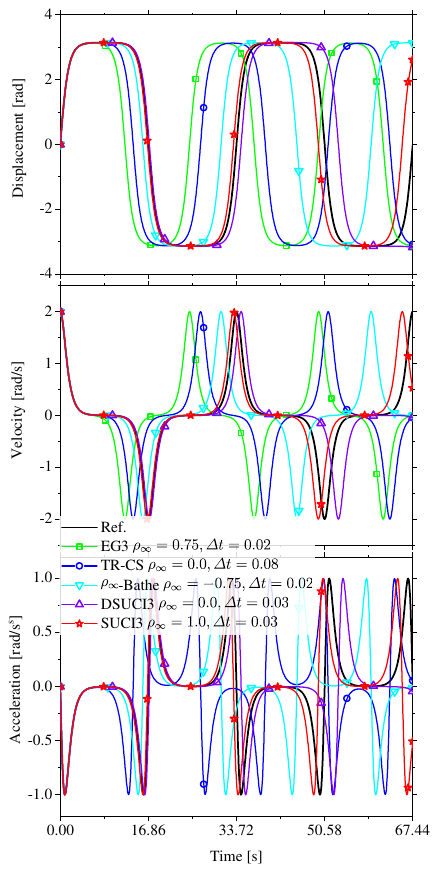}}
	\subfigure[Fourth-order algorithms]{
		\includegraphics[scale=0.7]{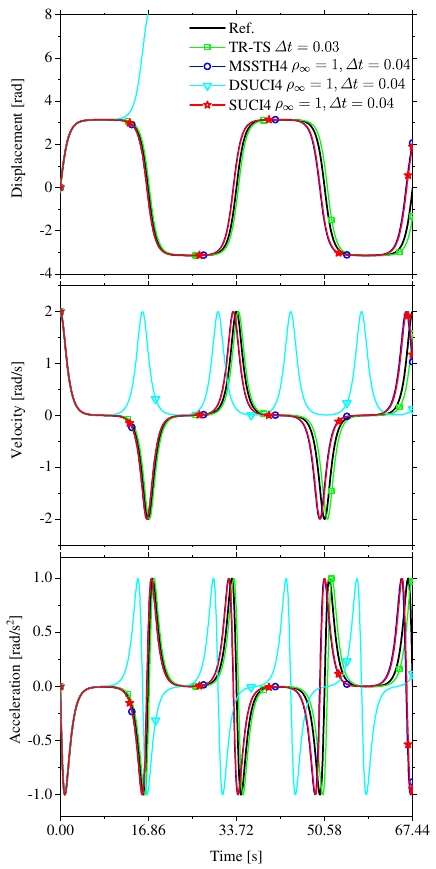}}
	\subfigure[Fifth/sixth-order algorithms]{
		\includegraphics[scale=0.7]{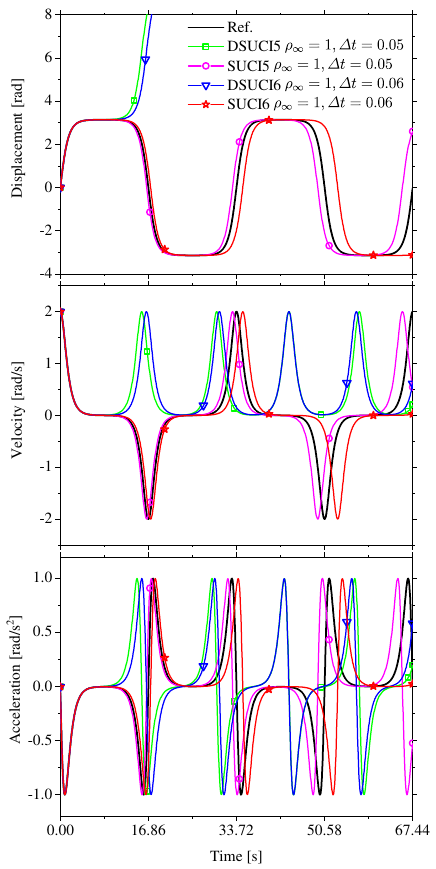}}
	\caption{Numerical solutions of the simple pendulum predicted by various high-order implicit methods.}
	\label{fig:pen_com}
\end{figure}
%\begin{figure}[htbp]
%	\centering
%	\includegraphics[scale=0.77]{simPen3}\vskip -2mm
%	\caption{Absolute displacement errors of the simple pendulum using various high-order algorithms with the same $ \dt=0.02 $s.}
%	\label{fig:simPen3}
%\end{figure}
Fig.~\ref{fig:pen_com} compares numerical solutions among various high-order implicit algorithms. In Fig.~\ref{fig:pen_com}, the time steps $\dt$ for the compared algorithms are typically set as $0.01\times n$, where $n$ represents the number of sub-steps. However, TR-CS \cite{fungComplextimestepNewmark1998,fungUnconditionallyStable1997} adopts $\dt=0.01\times 8$, following the approach used by Choi et al.~\cite{choiTimeSplitting2022}. This decision is made due to TR-CS employing two sub-steps and complex-valued parameters. Using different values of $\rhoinf$ in Fig.~\ref{fig:pen_com}(a) serves the purpose of optimizing the performance of each implicit method. For instance, as demonstrated in \cite{liDirectlySelfstarting2022}, DSUCI3 with $\rhoinf=0$ yields more accurate responses compared to other values of $\rhoinf$.  Notably, the SUCI$n$ algorithms consistently provide more accurate numerical responses than DSUCI$n$ \cite{liDirectlySelfstarting2022}. Among the third-order algorithms, SUCI3 produces the most accurate predictions. It is worth noting that the fourth-order TR-TS method outperforms SUCI4 and MSSTH4 \cite{zhangOptimizationNsubstep2020}, despite its inability to control numerical high-frequency dissipation. Overall, the proposed SUCI$n$ algorithms are superior to the published high-order schemes for solving this simple pendulum.

%It should be emphasized that although SUCI3 is seemly inferior to TR-CS \cite{fungComplextimestepNewmark1998,fungUnconditionallyStable1997}, the latter requires more computational costs due to using the complex-valued parameter. In comparison with DSUCI$n$ \cite{liDirectlySelfstarting2022}, it is found that DSUCI(4-6) perform worse than other high-order methods in the case of $ \dt=0.02 $s, while DSUCI3 predicts fewer absolute errors than SUCI3. Since the displacements predicted by DSUCI(4-6) have deviated from exact quantities, the accelerations predicted by the central difference of displacement naturally deviate from exact solutions and thus are omitted in Fig.~\ref{fig:pen_com}.
% When decreasing the integration step $ \dt $ into $ 0.01 $s, DSUCI$n$ can also predict quite accurate solutions. Therefore, 
%Overall, the present high-order methods are superior to the directly self-starting high-order methods \cite{liDirectlySelfstarting2022} for solving this nonlinear problem since the present methods utilize the acceleration vector $ \mbf{\ddot{U}}_n $ to obtain better robustness.

%----------------------------------------------------------------

\subsubsection{An $ N $-degree-of-freedom mass-spring system}
\begin{figure}[tbhp]
	\centering
	\includegraphics[scale=0.4]{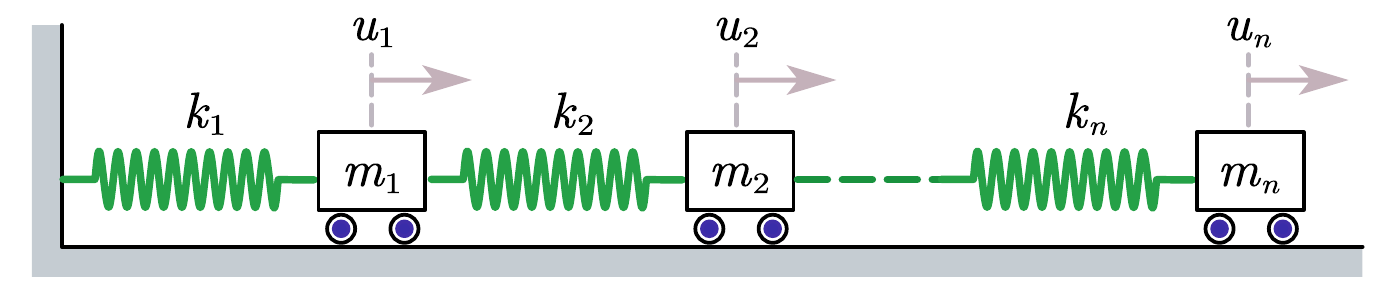}
	\caption{The $ N $-degree-of-freedom mass-spring model.}
	\label{fig:ndof}
\end{figure}

An $ N $-degree-of-freedom mass-spring system \cite{rezaiee-pajandFamilySecondorder2018}, shown in Fig.~\ref{fig:ndof}, is solved to test the computational cost among various integration algorithms. In Fig.~\ref{fig:ndof}, all masses are set as $ m_i=1 $kg, and stiffness coefficients $ k_i $ are assumed as
\begin{equation}\label{key}
	k_i=\begin{cases}
		k                          & i=1         \\
		k\left[1+\alpha(u_i-u_{i-1})^2\right] & 2\le i\le N
	\end{cases}
\end{equation}
where $ k=10^5 $N/m and $ \alpha=2 $m$^{-2}$ is adopted to simulate the nonlinear stiffness hardening system. In addition, the model is excited by a force of $ \sin(t) $ for all masses. The completely non-dissipative trapezoidal rule with $ \dt=1.0\times10^{-7} $s is employed to provide reference solutions.

\begin{table}[htbp]
	%	\centering
	\caption{The elapsed computational time and global displacement errors using the completely non-dissipative ($ \rhoinf=1 $) algorithms.}
	\small\begin{tabular}{lccccccc}\toprule
		\multirow{2}{*}{Algorithms} &\multirow{2}{*}{$ \dt $ [s]} & \multicolumn{2}{l}{ $ N=1000 $} & \multicolumn{2}{l}{$ N=5000 $} & \multicolumn{2}{l}{$ N=10000 $} \\ \cmidrule(lr){3-4}\cmidrule(lr){5-6}\cmidrule(lr){7-8}
		&	& CPU time [s] & Global errors &CPU time [s] & Global errors &CPU time [s] & Global errors \\ 
		\midrule 
		\rowcolor{mygray}
		TPO/G-$ \alpha $ \cite{shaoThreeParameters1988,chungTimeIntegration1993} & 0.02 &  0.856 & $2.3693\times10^{-4}$ &
		3.717 &	$6.9573\times10^{-5}$ &
		11.952 &	$4.7319\times10^{-5}$\\
		\rowcolor{mygray}
		SUCI2 \cite{liNovelFamily2020}  & 0.02 & 1.213 &	$6.1301\times10^{-5}$ &
		5.975 &	$1.9814\times10^{-5}$ &
		21.248 &	$1.3714\times10^{-5}$\\ 
		
		TR-CS \cite{fungComplextimestepNewmark1998,fungUnconditionallyStable1997} & 0.02 & 2.195 &	$2.3159\times10^{-7}$&
		12.153 &	$4.2864\times10^{-8}$&
		34.097 &$	1.3736\times10^{-8}$\\ 
		
		DSUCI3 \cite{liDirectlySelfstarting2022} & 0.02 & 1.569  & $ 7.5232\times 10^{-7} $ &  9.523& $8.6127\times10^{-8}$ & 27.159&	$7.0891\times10^{-8}$ \\
		
		SUCI3 & 0.02 & 1.573 &	$7.6711\times10^{-7}$&
		9.935 &	$8.7909\times10^{-8}$&
		27.927&	$7.0679\times10^{-8}$\\ 
		\rowcolor{mygray}
		TR-TS \cite{tarnowHowRender1994} & 0.02 & 		2.021 &	$3.1656\times10^{-7}$ &
		10.265 &$	5.2774\times10^{-8}$ &
		28.086 &	$1.7341\times10^{-8}$ \\ 
		\rowcolor{mygray}
		MSSTH4 \cite{zhangOptimizationNsubstep2020} & 0.02 & 3.044 & $ 2.3373\times10^{-7} $ & 15.836 & $ 3.9261\times10^{-8} $ & 33.098 & $ 9.7313\times10^{-9} $\\ 
		\rowcolor{mygray}
		DSUCI4 \cite{liDirectlySelfstarting2022} & 0.02 & 2.993  & $ 1.2882\times10^{-7} $ & 14.273 & $ 2.1149\times10^{-8} $ & 33.108 & $ 3.3021\times10^{-9} $ \\
		\rowcolor{mygray}
		SUCI4  & 0.02 &3.121 & $ 1.3047\times10^{-7} $& 14.929 & $ 2.1416\times10^{-8} $& 33.597 & $ 3.3267\times10^{-9} $\\
		
		MSSTH5 \cite{zhangOptimizationNsubstep2020} & 0.02 & 4.232 &  $4.8580\times10^{-7}$& 16.557& $ 1.0683\times10^{-8} $ & 45.706 &$ 3.3453\times10^{-8} $ \\ 
		DSUCI5 \cite{liDirectlySelfstarting2022} &  0.02 & 3.997  & $ 4.3207\times10^{-8} $& 16.893 & $ 7.2260\times10^{-9} $ & 48.891 & $ 4.3192\times10^{-9} $ \\
		SUCI5  & 0.02 & 3.782 & $ 4.0110\times10^{-8} $& 16.499 & $ 6.5481\times10^{-9} $ & 47.538 & $ 2.1996\times10^{-9} $ \\
		
		\rowcolor{mygray}
		DSUCI6 \cite{liDirectlySelfstarting2022} & 0.02 & 4.980  & $ 2.0359\times10^{-8} $ & 22.002 & $ 3.3719\times10^{-9} $ & 70.942 & $ 4.1389\times10^{-10} $ \\
		\rowcolor{mygray}
		SUCI6  & 0.02 & 5.133 & $ 2.0897\times10^{-8} $ & $ 21.155 $ & $ 3.3846\times10^{-9} $ & 70.707 & $ 4.1463\times10^{-10} $ \\ 
		\bottomrule
	\end{tabular}
	\label{tab:rhoinf1}
\end{table}

It should be realized that higher-order accurate integration algorithms often require more computational cost than the common second-order ones when the same time step $ \dt $ is used. Hence, this example focuses mainly on comparing computational efficiency among the same accurate algorithms. Table~\ref{tab:rhoinf1} firstly records the elapsed computational CPU time and global displacement errors of this mass-spring system with $ N=1000,~5000 $, and $ 10000 $ using the non-dissipative implicit algorithms in this paper (\emph{the third-order EG3 \cite{fungExtrapolatedGalerkin1996} and $ \rhoinf $-Bathe \cite{kwonSelectingLoad2021} algorithms with real-valued parameters cannot use $ \rhoinf=\pm 1 $ and they are not thus compared herein}). It is evident that when adopting the same integration step $ \dt=0.02 $s, the higher-order integration algorithms need more computational CPU time than the lower-order ones since more sub-steps are used to achieve higher-order accuracy, but these schemes, in turn, produce smaller global errors. As expected, the computational time of each implicit algorithm sharply increases as the degree-of-freedom increases. Notice also that the third-order complex-sub-step TR-CS scheme \cite{fungComplextimestepNewmark1998,fungUnconditionallyStable1997} requires significantly more computational time than SUCI3, although it also provides slight smaller global errors. Considering the existing four-sub-step (MSSTH4) and five-sub-step (MSSTH5) schemes \cite{zhangOptimizationNsubstep2020}, it has been shown in \cite{liDirectlySelfstarting2022} that they suffer from the order reduction for solving the forced vibrations. Hence, MSSTH4 and MSSTH5 produce significantly larger global errors than SUCI4 and SUCI5, but almost the same computational time is observed since they use the same number of sub-steps. The fourth-order TR-TS \cite{tarnowHowRender1994} consumes less computational time than MSSTH4 and SUCI4 since it is essentially a composite three-sub-step scheme. The directly self-starting high-order algorithms \cite{liDirectlySelfstarting2022} (DSUCI$n$) show almost the same numerical performance as SUCI$n$ since they share similar numerical characteristics. The most significant difference between DSUCI$n$ and SUCI$n$ is that DSUCI$n$ is designed to be directly self-starting, and thus the acceleration output is additionally constructed. It should be pointed out that the computational CPU time of the high-order $ s $-sub-step methods, including the published schemes \cite{zhangOptimizationNsubstep2020,tarnowHowRender1994,liDirectlySelfstarting2022,fungComplextimestepNewmark1998,fungUnconditionallyStable1997}, is $ s $ times less than the single-step TPO/G-$\alpha$'s \cite{shaoThreeParameters1988,chungTimeIntegration1993} computational time. One of the possible reasons for this numerical behavior is that the TPO/G-$\alpha$ method achieves only first-order accurate accelerations, thus slowing down convergence speeds.

\begin{table}[htbp]
	%	\centering
	\caption{The elapsed computational time and global displacement errors using the most dissipative ($ \rhoinf=0 $) algorithms.}
	\small\begin{tabular}{lccccccc}\toprule
		\multirow{2}{*}{Algorithms} &\multirow{2}{*}{$ \dt $ [s]} & \multicolumn{2}{l}{ $ N=1000 $} & \multicolumn{2}{l}{$ N=5000 $} & \multicolumn{2}{l}{$ N=10000 $} \\ \cmidrule(lr){3-4}\cmidrule(lr){5-6}\cmidrule(lr){7-8}
		&	& CPU time [s]  & Global errors &CPU time [s] & Global errors &CPU time [s] & Global errors \\
		\midrule 
		\rowcolor{mygray}
		TPO/G-$ \alpha $ \cite{shaoThreeParameters1988,chungTimeIntegration1993} & 0.02 &  0.810 & $1.2973\times10^{-3}$ &
		8.062 &	$3.8377\times10^{-4}$ &
		12.797 &	$2.4291\times10^{-4}$\\
		\rowcolor{mygray}
		SUCI2 \cite{liNovelFamily2020}  & 0.02 & 1.306 &	$1.1864\times10^{-4}$ &
		12.618 &	$3.8427\times10^{-5}$ &
		20.352  &	$2.6614\times10^{-5}$\\ 
		
		TR-CS \cite{fungComplextimestepNewmark1998,fungUnconditionallyStable1997} & 0.02 & 3.058 &	$2.5343\times10^{-6}$&
		15.455 &	$4.7898\times10^{-7}$&
		40.748 &$	1.9784\times10^{-7}$\\ 
		DSUCI3 \cite{liDirectlySelfstarting2022} & 0.02 & 2.231 &	$4.3721\times10^{-6}$&
		10.231 &	$8.3518\times10^{-7}$&
		30.191 & $	3.4987\times10^{-7}$\\
		SUCI3 & 0.02 & 2.207 &	$4.3874\times10^{-6}$&
		10.645 &	$8.3989\times10^{-7}$&
		30.848 &	$3.5654\times10^{-7}$\\ 
		
		\rowcolor{mygray}
		MSSTH4 \cite{zhangOptimizationNsubstep2020} & 0.02 & 3.132 &	$3.3477\times10^{-6}$&
		17.290 &	$5.2566\times10^{-7}$&
		43.115 &$	9.4042\times10^{-8}$\\ 
		\rowcolor{mygray}
		DSUCI4 \cite{liDirectlySelfstarting2022} & 0.02 & 3.204 &	$1.3111\times10^{-6}$&
		17.541 & $2.2489\times10^{-7}$& 
		43.815&	$3.4101\times10^{-8}$\\
		\rowcolor{mygray}
		SUCI4  & 0.02 & 3.259 &	$1.3474\times10^{-6}$&
		17.603 &	$2.2565\times10^{-7}$&
		43.741&	$3.3932\times10^{-8}$\\ 
		
		MSSTH5 \cite{zhangOptimizationNsubstep2020} & 0.02 & 4.337 &  $7.4887\times10^{-7}$& 21.877& $ 2.4608\times10^{-8} $ & 57.162 &$ 4.6910\times10^{-8} $ \\ 
		DSUCI5 \cite{liDirectlySelfstarting2022} & 0.02 & 4.298 &	$6.8381\times10^{-8}$&
		20.839 & $3.1109\times10^{-9} $ & 57.479 & $ 3.7719\times10^{-9} $\\
		SUCI5  & 0.02 & 4.375 & $ 6.8380\times10^{-8} $& 21.161 & $3.1171\times10^{-9} $ & 57.565 & $ 3.7500\times10^{-9} $ \\ 
		
		\rowcolor{mygray}
		DSUCI6 \cite{liDirectlySelfstarting2022} & 0.02 & 5.012 &	$4.6142\times10^{-8}$&
		26.997 &	$ 7.5503\times10^{-9}$&
		77.212 & $ 7.4311\times10^{-10} $\\
		\rowcolor{mygray}
		SUCI6  & 0.02 & 5.039 & $ 4.5919\times10^{-8} $ & $ 27.439 $ & $ 7.4504\times10^{-9} $ & 77.612 & $ 7.2195\times10^{-10} $ \\ 
		\bottomrule
	\end{tabular}
	\label{tab:rhoinf0}
\end{table}

When considering the most dissipative ($ \rhoinf=0 $) case, the elapsed computational time and global errors among various dissipative algorithms are listed in Table \ref{tab:rhoinf0}. The same conclusions as those in Table \ref{tab:rhoinf1} can be given. Besides, it is also found that when imposing numerical dissipation in the high-frequency range, these integration algorithms produce larger global errors. This observation is in well agreement with the spectral analysis in Section \ref{sec:sp}, where period errors increase as the parameter $ \rhoinf $ decreases from unity into zero, as shown Fig.~\ref{fig:peSUCI2}.

%Tables \ref{tab:rhoinf1} and \ref{tab:rhoinf0} mainly compare the numerical performance of high-order implicit methods with the same order of accuracy and the same time step size. 
It is interesting to compare the numerical performance among various implicit methods with the same sub-step size. Table \ref{tab:rhoinf01} records the elapsed computational time and global errors using various algorithms ($ \rhoinf=0.0 $) with the same sub-step size. Note that each implicit method in Table \ref{tab:rhoinf01} adopts the integration step size $ \dt=0.02\times s $ to guarantee each sub-step size as $ 0.02 $s. With such parameter settings, all implicit methods in Table \ref{tab:rhoinf01} should elapse almost the same computational time for solving linear structures. However, this is not the case for solving nonlinear problems. As one can observe from Table \ref{tab:rhoinf01}, the multi-sub-step methods are significantly superior to the traditional single-step TPO/G-$\alpha$ \cite{shaoThreeParameters1988,chungTimeIntegration1993} method in terms of the computational time and solution accuracy. When the multi-sub-step methods achieve high-order accuracy, such an advantage is further enhanced. As the number of DOFs increases, various implicit methods with the same sub-step size elapse the approximately same CPU time, such as the case of $ N=10000 $. In addition, some useful conclusions from Tables \ref{tab:rhoinf1} and \ref{tab:rhoinf0} are still valid. For instance, the proposed SUCI$n$ method still outperforms MSSTH$n$ \cite{zhangOptimizationNsubstep2020} since the latter produces larger global errors due to the order reduction for solving forced vibrations. 
\begin{table}[htbp]
	%	\centering
	\caption{The elapsed computational time and global displacement errors using various algorithms ($ \rhoinf=0.0 $) with the same sub-step size.}
	\small\begin{tabular}{lccccccc}\toprule 
		\multirow{2}{*}{Algorithms} &\multirow{2}{*}{$ \dt $ [s]} & \multicolumn{2}{l}{ $ N=1000 $} & \multicolumn{2}{l}{$ N=5000 $} & \multicolumn{2}{l}{$ N=10000 $} \\ \cmidrule(lr){3-4}\cmidrule(lr){5-6}\cmidrule(lr){7-8}
		&	& CPU time [s]  & Global errors &CPU time [s] & Global errors &CPU time [s] & Global errors \\
		\midrule 
		\rowcolor{mygray}
		TPO/G-$ \alpha $ \cite{shaoThreeParameters1988,chungTimeIntegration1993} & 0.02 &  0.810 & $1.2973\times10^{-3}$ &
		8.062 &	$3.8377\times10^{-4}$ &
		12.797 &	$2.4291\times10^{-4}$\\
		
		\rowcolor{mygray}
		SUCI2 \cite{liNovelFamily2020}  & 0.04 & 0.894 &	$4.7066\times10^{-4}$ &
		5.112 &	$1.5328\times10^{-4}$ &
		8.251  &	$1.0640\times10^{-4}$\\ 
		
		TR-CS \cite{fungComplextimestepNewmark1998,fungUnconditionallyStable1997} & 0.06 & 0.551 &	$4.7977\times10^{-5}$&
		5.830 &	$9.8774\times10^{-6}$&
		8.856 &$ 4.6496\times10^{-6}$\\ 
		DSUCI3 \cite{liDirectlySelfstarting2022} & 0.06 & 0.635 &	$8.2892\times10^{-5}$&
		5.012 &	$1.7344\times10^{-5}$&
		8.555 & $ 8.4964\times10^{-6}$\\
		SUCI3 & 0.06 & 0.640 &	$8.3278\times10^{-5}$&
		5.131 &	$1.7509\times10^{-5}$&
		8.519 &	$8.4836\times10^{-6}$\\ 
		
		\rowcolor{mygray}
		MSSTH4 \cite{zhangOptimizationNsubstep2020} & 0.08 & 0.691 &	$6.5893\times10^{-5}$&
		5.422 &	$4.1022\times10^{-5}$&
		8.073 &$ 6.2424\times10^{-6}$\\ 
		\rowcolor{mygray}
		DSUCI4 \cite{liDirectlySelfstarting2022} & 0.08 & 0.681 &	$6.3889\times10^{-5}$&
		5.515 & $1.0912\times10^{-5}$& 
		8.478&	$3.9765\times10^{-6}$\\
		\rowcolor{mygray}
		SUCI4  & 0.08 & 0.679 &	$6.3938\times10^{-5}$&
		5.428 &	$1.1061\times10^{-5}$&
		8.113&	$4.1018\times10^{-6}$\\ 
		
		MSSTH5 \cite{zhangOptimizationNsubstep2020} & 0.10 & 0.706 &  $9.2543\times10^{-6}$& 6.331& $ 3.7125\times10^{-6} $ & 8.373 &$ 5.6442\times10^{-7} $ \\ 
		DSUCI5 \cite{liDirectlySelfstarting2022} & 0.10 & 0.740 &	$7.2882\times10^{-6}$&
		6.853 & $1.1921\times10^{-6} $ & 8.607 & $ 3.9446\times10^{-7} $\\
		SUCI5  & 0.10 & 0.719 & $ 7.2879\times10^{-6} $& 6.708 & $1.1921\times10^{-6} $ & 8.653 & $ 3.9450\times10^{-7} $ \\ 
		
		\rowcolor{mygray}
		DSUCI6 \cite{liDirectlySelfstarting2022} & 0.12 & 0.811 &	$6.9486\times10^{-6}$&
		7.304 &	$ 9.4334\times10^{-7}$&
		9.188 & $ 3.2105\times10^{-7} $\\
		\rowcolor{mygray}
		SUCI6  & 0.12 & 0.813 & $ 6.9481\times10^{-6} $ &  7.297  & $ 9.4433\times10^{-7} $ & 9.180 & $ 3.1510\times10^{-7} $ \\ 
		\bottomrule
	\end{tabular}
	\label{tab:rhoinf01}
\end{table}

\begin{figure}[htbp]
	\begin{minipage}[t]{0.5\linewidth}
		\centering
		\includegraphics[scale=0.25]{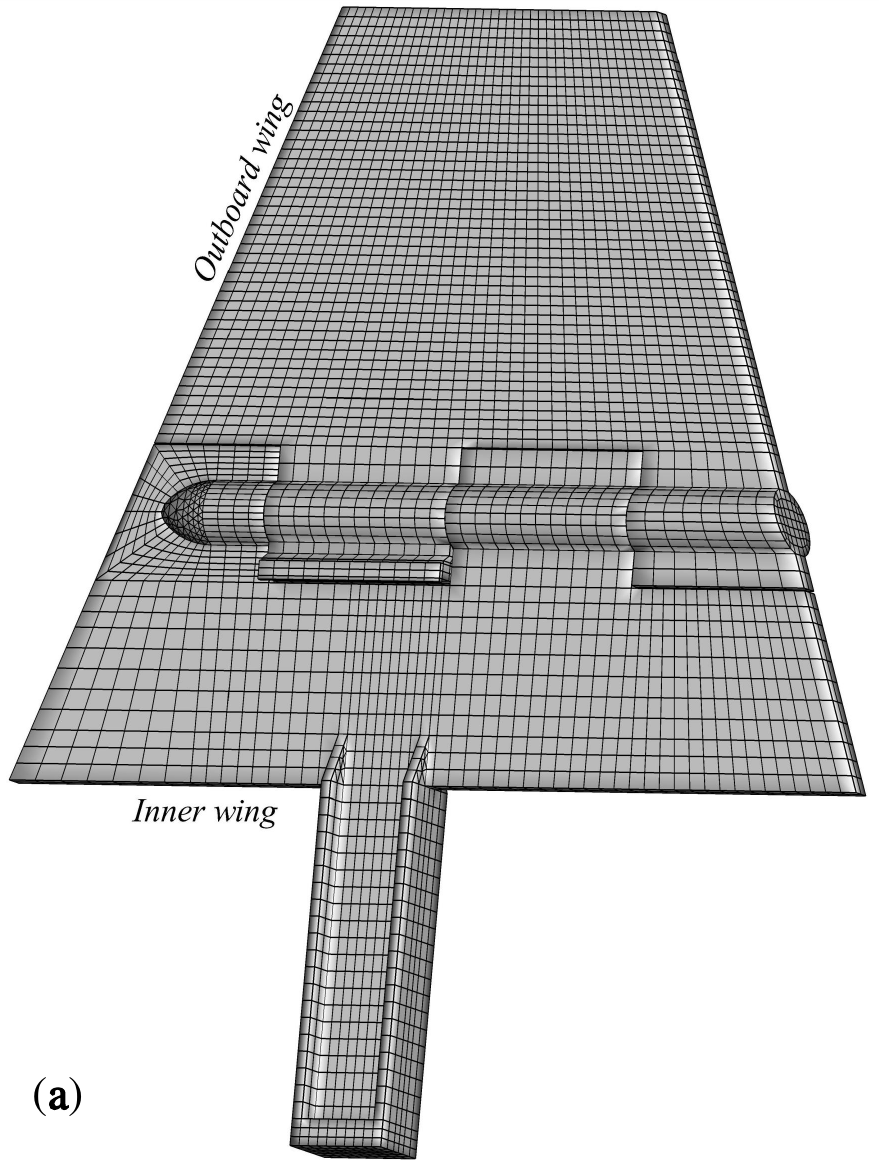}
	\end{minipage}%
	\begin{minipage}[t]{0.5\linewidth}
		\centering
		\includegraphics[scale=0.5]{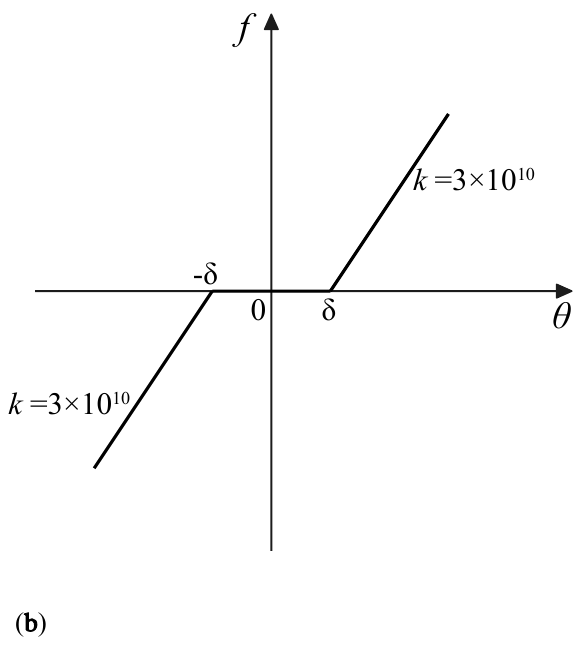}
	\end{minipage}
	\caption{(\textbf{a}) Finite element model of the folding wing with free-play nonlinearity, and (\textbf{b}) model of the nonlinear torsional spring with $ \delta=0.1 $mm.}
	\label{fig:wing}
\end{figure}
\subsubsection{A folding wing with free-play nonlinearity}
To further test the numerical performance of SUCI$n$ on the practical engineering problem, a folding wing \cite{heNonlinearAeroelastic2020,liTrulySelfstarting2020} with free-play nonlinearity consisting of the outboard and inner wings shown in Fig.~\ref{fig:wing}(a) is solved by various high-order implicit methods. The connection of two wings has been modeled by three same nonlinear torsional springs \cite{heNonlinearAeroelastic2020} (see Fig.~\ref{fig:wing}(b)). The detailed mathematical model of this folding wing can refer to the literature \cite{heNonlinearAeroelastic2020}. In this finite element model, there are 9627 solid elements and 80800 DOFs. The initial displacements and velocities of all DOFs are assumed to be zero except that three torsional springs possess the nonzero velocity ($20$mm/s) to make the outboard and inner wings meet with each other.

\begin{figure}[htbp]
	\centering
	\includegraphics[scale=0.7]{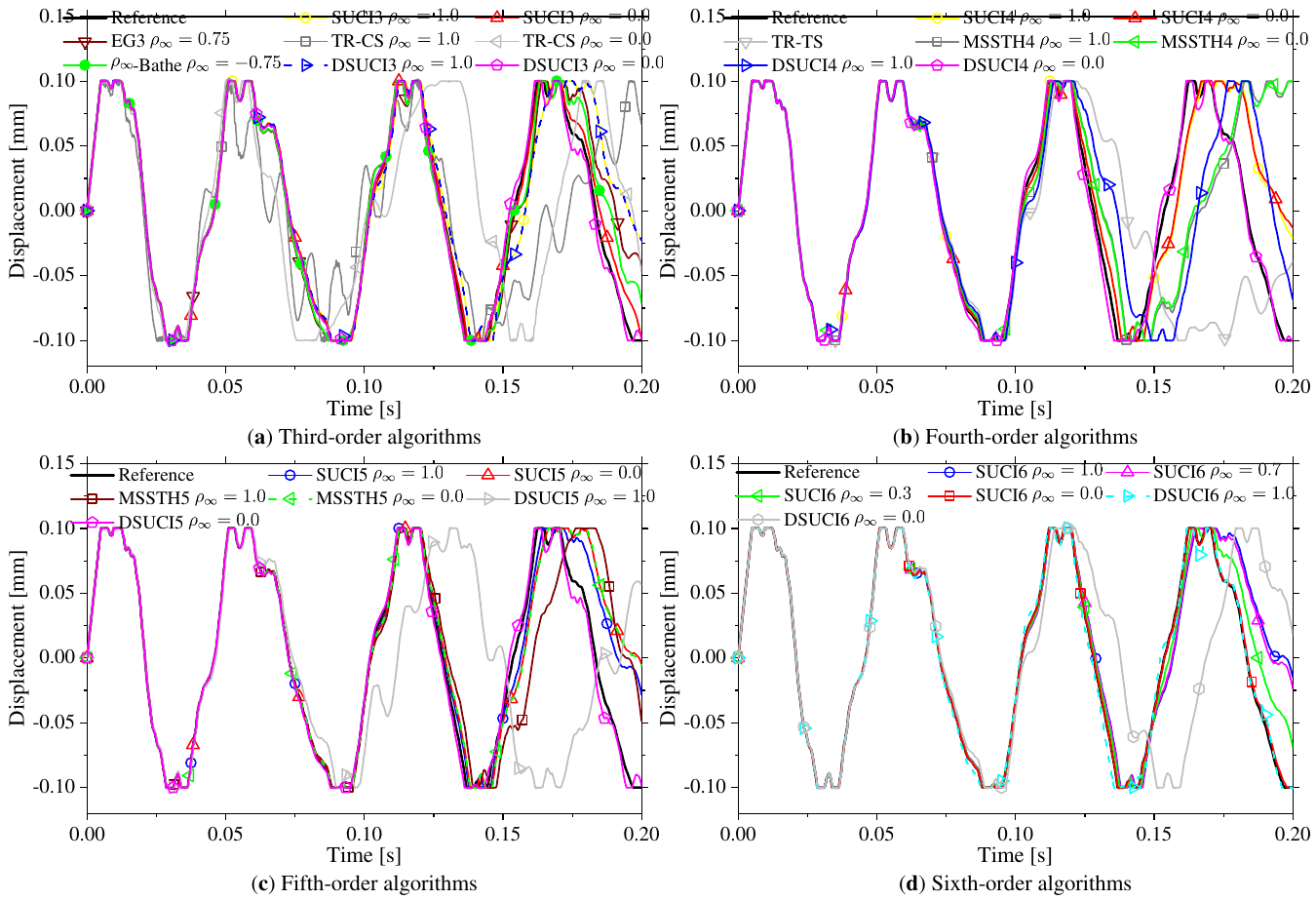}\vskip -2mm
	\caption{Numerical displacements of the third torsional spring using the high-order integration algorithms with the same $ \dt=2.5\times10^{-5} $s.}
	\label{fig:wingsol}
\end{figure}
Numerical displacements of the third torsional spring using the high-order integration algorithms are plotted in Fig.~\ref{fig:wingsol}, where the reference solutions are obtained by using the sixth-order accurate SUCI6 scheme with a smaller time step $ \dt=1.0\times10^{-7} $s. It follows that the third-order complex-sub-step TR-CS scheme \cite{fungComplextimestepNewmark1998,fungUnconditionallyStable1997} produces adverse phase errors compared with SUCI3, EG3 \cite{fungExtrapolatedGalerkin1996}, and $ \rhoinf $-Bathe \cite{kwonSelectingLoad2021}. Moreover, other SUCI$n$ schemes in Fig.~\ref{fig:wingsol}(b-d) also show better robustness than MSSTHn \cite{zhangOptimizationNsubstep2020} and TR-TS \cite{tarnowHowRender1994}. DSUCI4 \cite{liDirectlySelfstarting2022} and DSUCI5 \cite{liDirectlySelfstarting2022} perform better than SUCI4 and SUCI5, respectively, in the present settings. As emphasized previously, except for the acceleration output, the previous DSUCI$n$ \cite{liDirectlySelfstarting2022} method possesses comparative numerical performance with SUCI$n$. On the other hand, this nonlinear example does not strictly follow the previous experience that when using the same integration step size the higher-order schemes can generally predict more accurate numerical responses than the lower-order ones. For instance, the third-order SUCI3($ \rhoinf=0 $) and sixth-order SUCI6($ \rhoinf=0 $) methods perform better than other schemes for solving the present problem.

\begin{figure}[htbp]
	\centering
	\subfigtopskip=2pt %?????????????????
	\subfigbottomskip=-4pt %??????????????????????????????
	\subfigcapskip=-5pt %?????????????
	\subfigure[]{
		\includegraphics[scale=0.3]{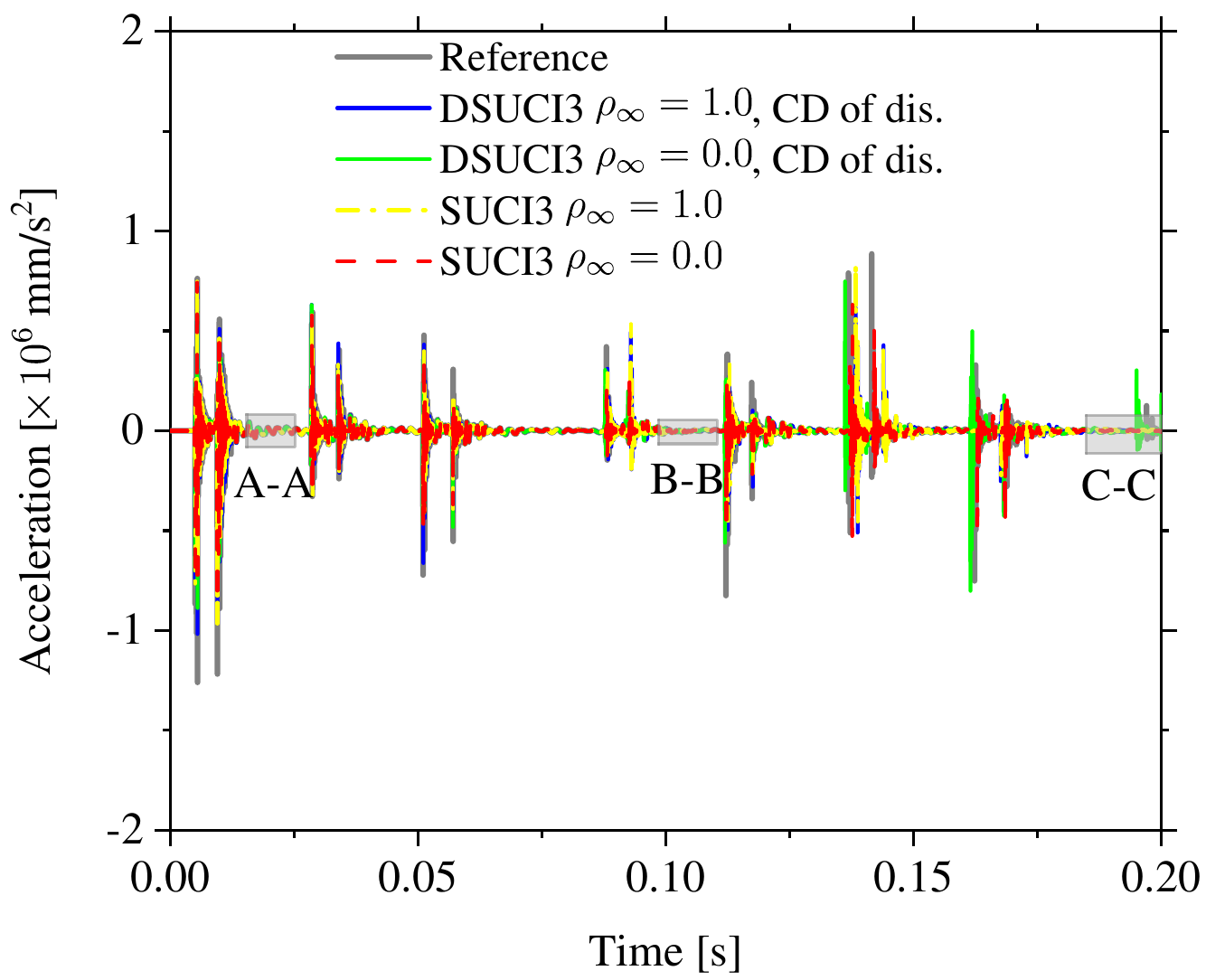}}
	\subfigure[A-A]{
		\includegraphics[scale=0.3]{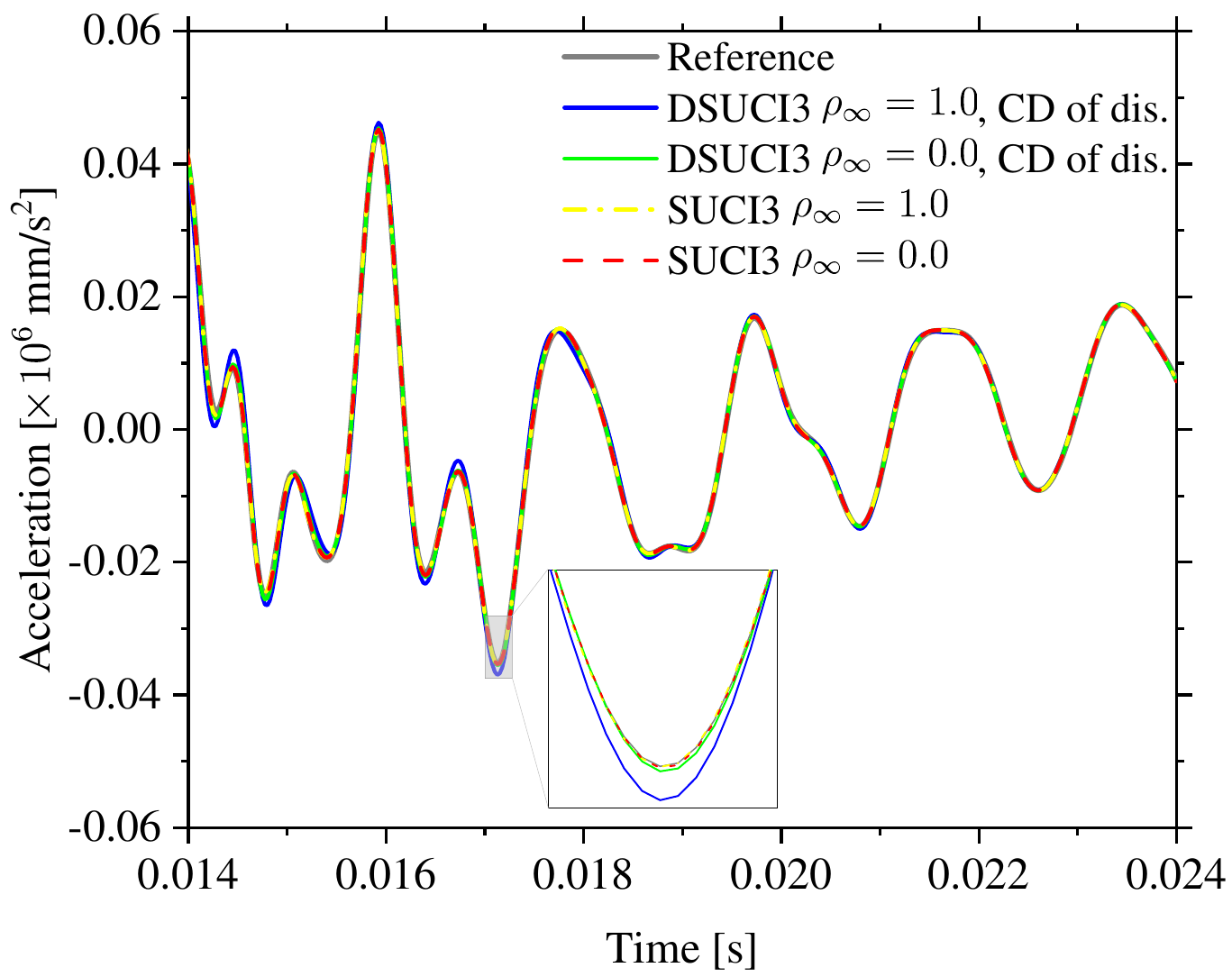}}
	\subfigure[B-B]{
		\includegraphics[scale=0.3]{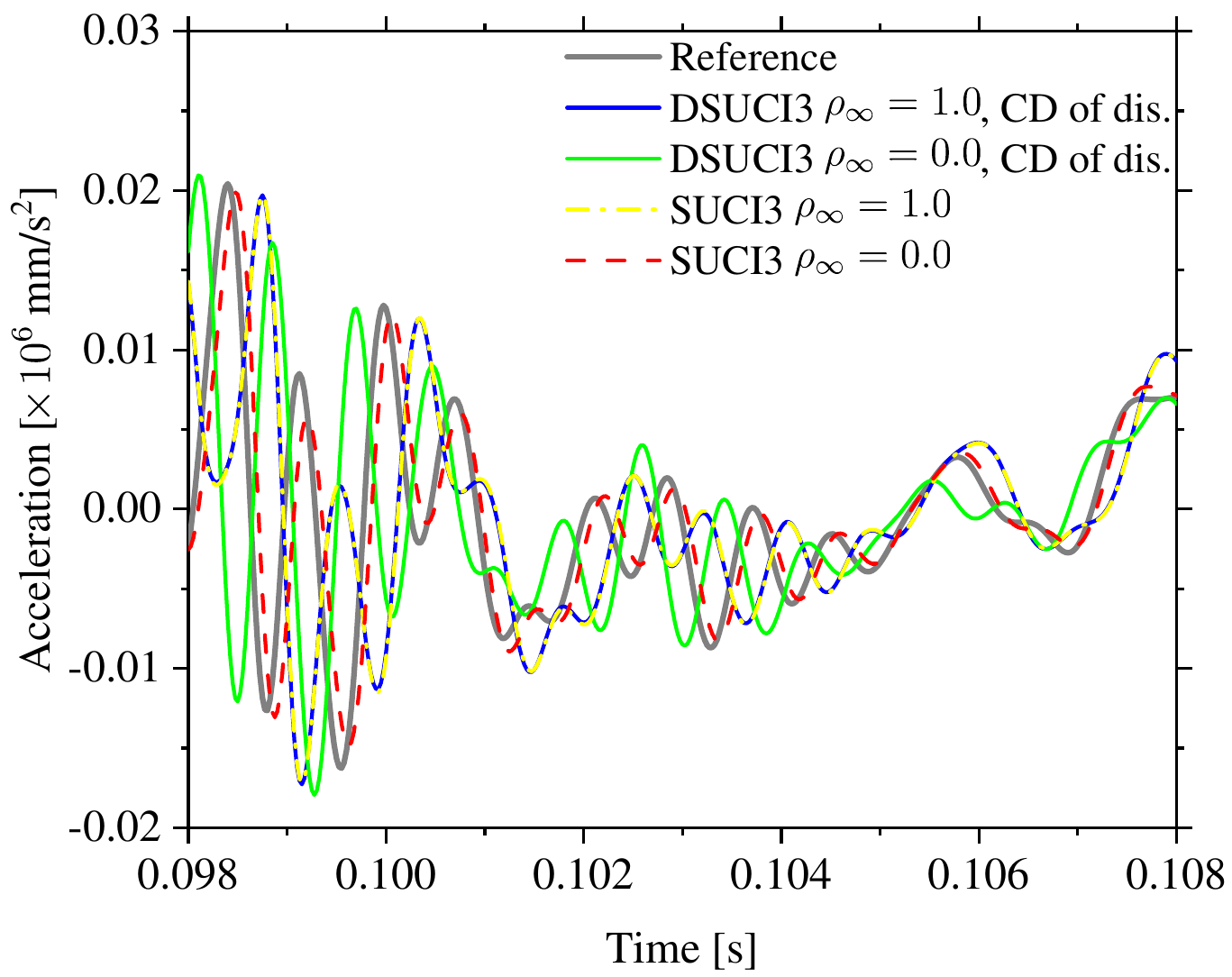}}
	\subfigure[C-C]{
		\includegraphics[scale=0.3]{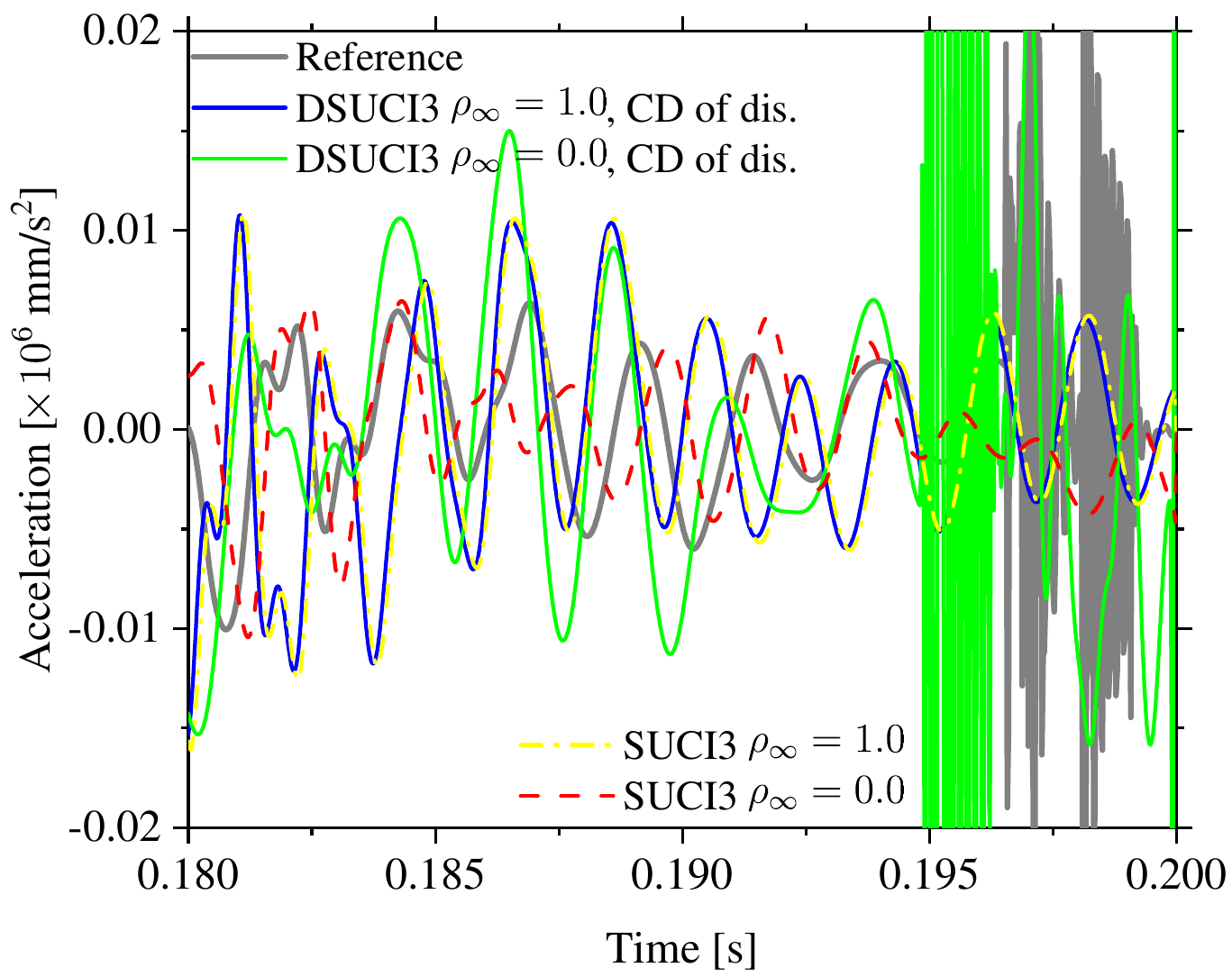}}
	\caption{Numerical accelerations of the third torsional spring using DSUCI3 \cite{liDirectlySelfstarting2022} and SUCI3 with the same $\dt=2.5\times10^{-5}$s. DSUCI3 adopts the central difference of displacement to output accelerations.}
	\label{fig:wing_ds1}
\end{figure}
\begin{figure}[htbp]
	\centering
	\subfigtopskip=2pt %?????????????????
	\subfigbottomskip=-4pt %??????????????????????????????
	\subfigcapskip=-5pt %?????????????
	\subfigure[]{
		\includegraphics[scale=0.23]{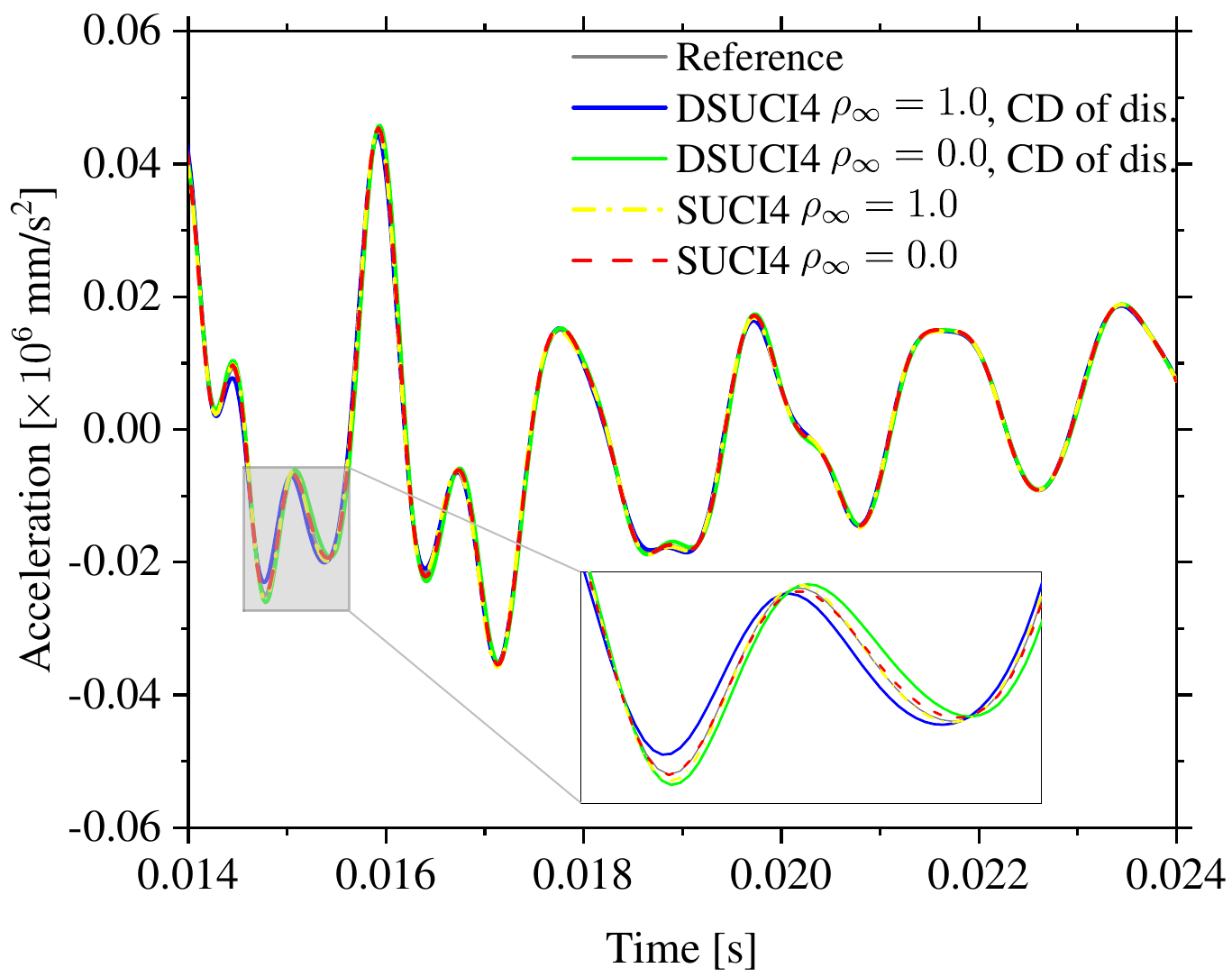}}
	\subfigure[]{
		\includegraphics[scale=0.23]{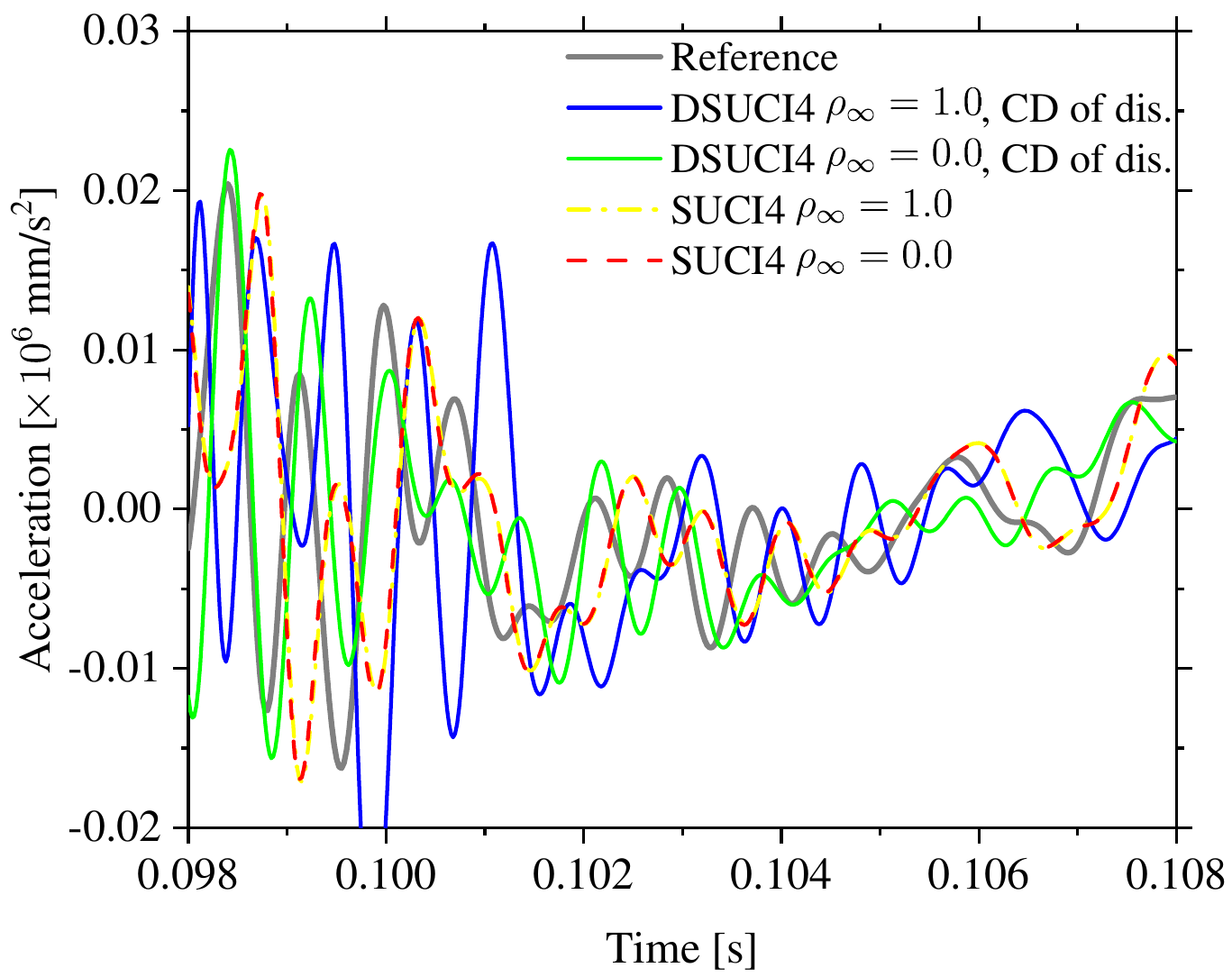}}
	\subfigure[]{
		\includegraphics[scale=0.23]{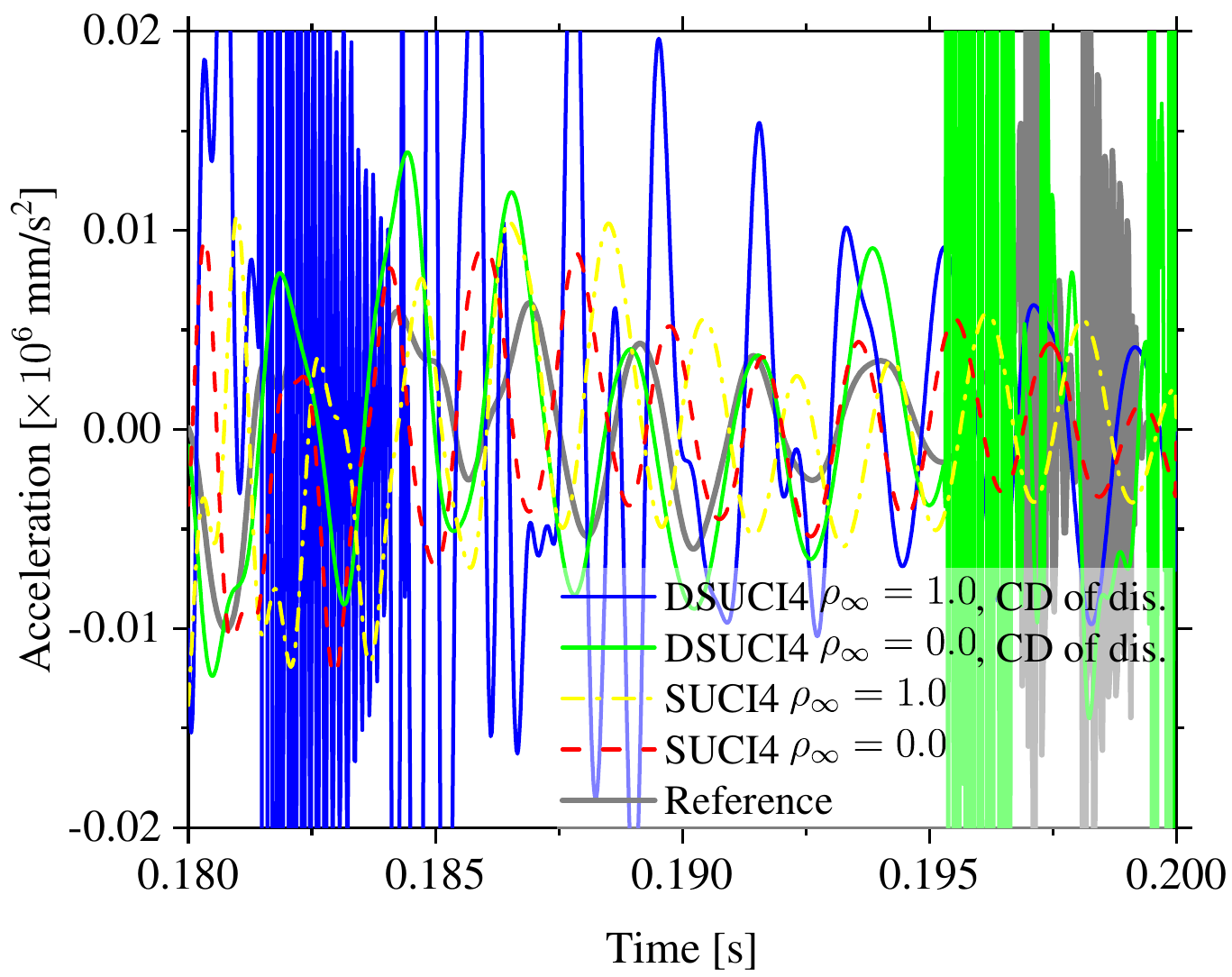}}
	\subfigure[]{
		\includegraphics[scale=0.23]{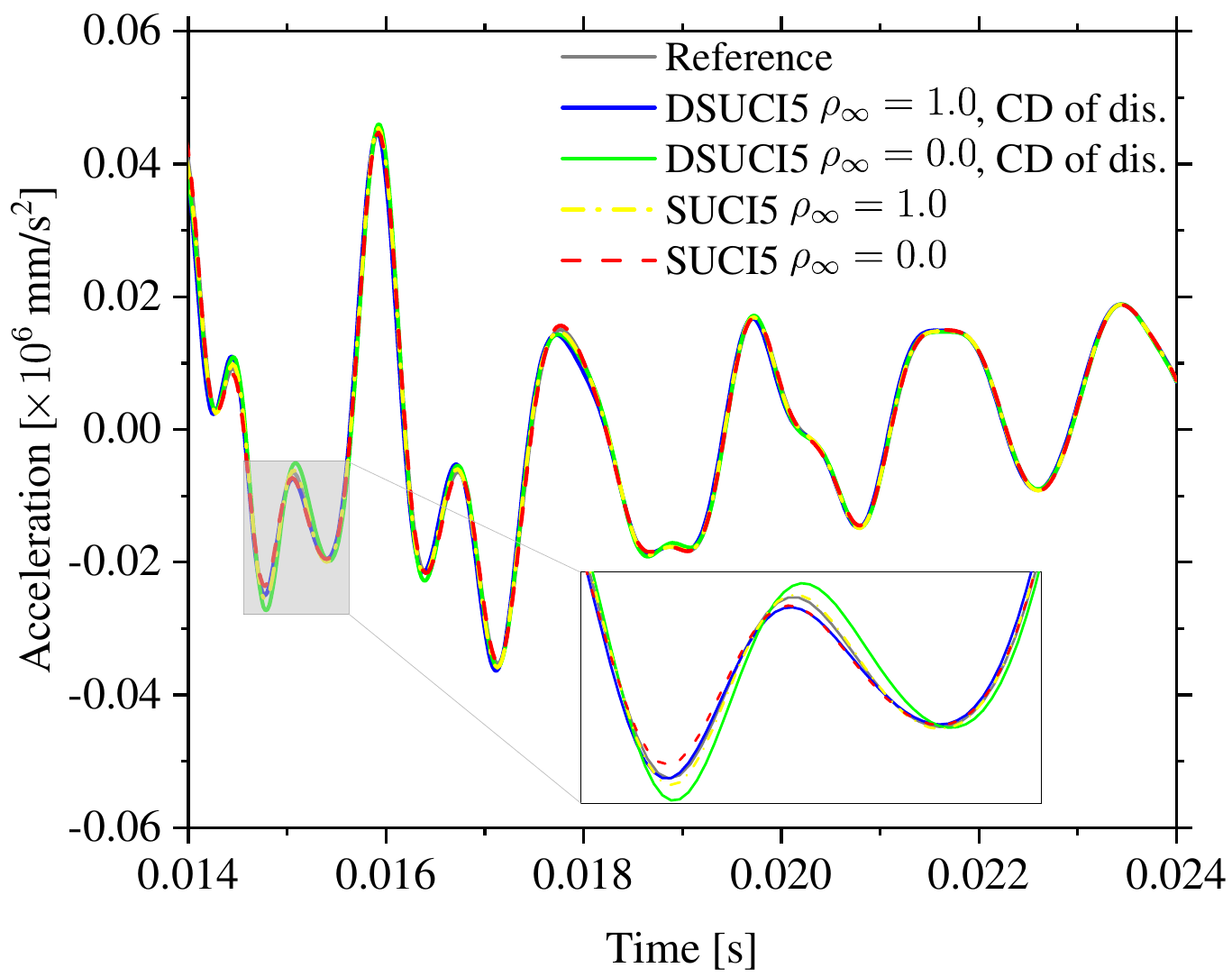}}
	\subfigure[]{
		\includegraphics[scale=0.23]{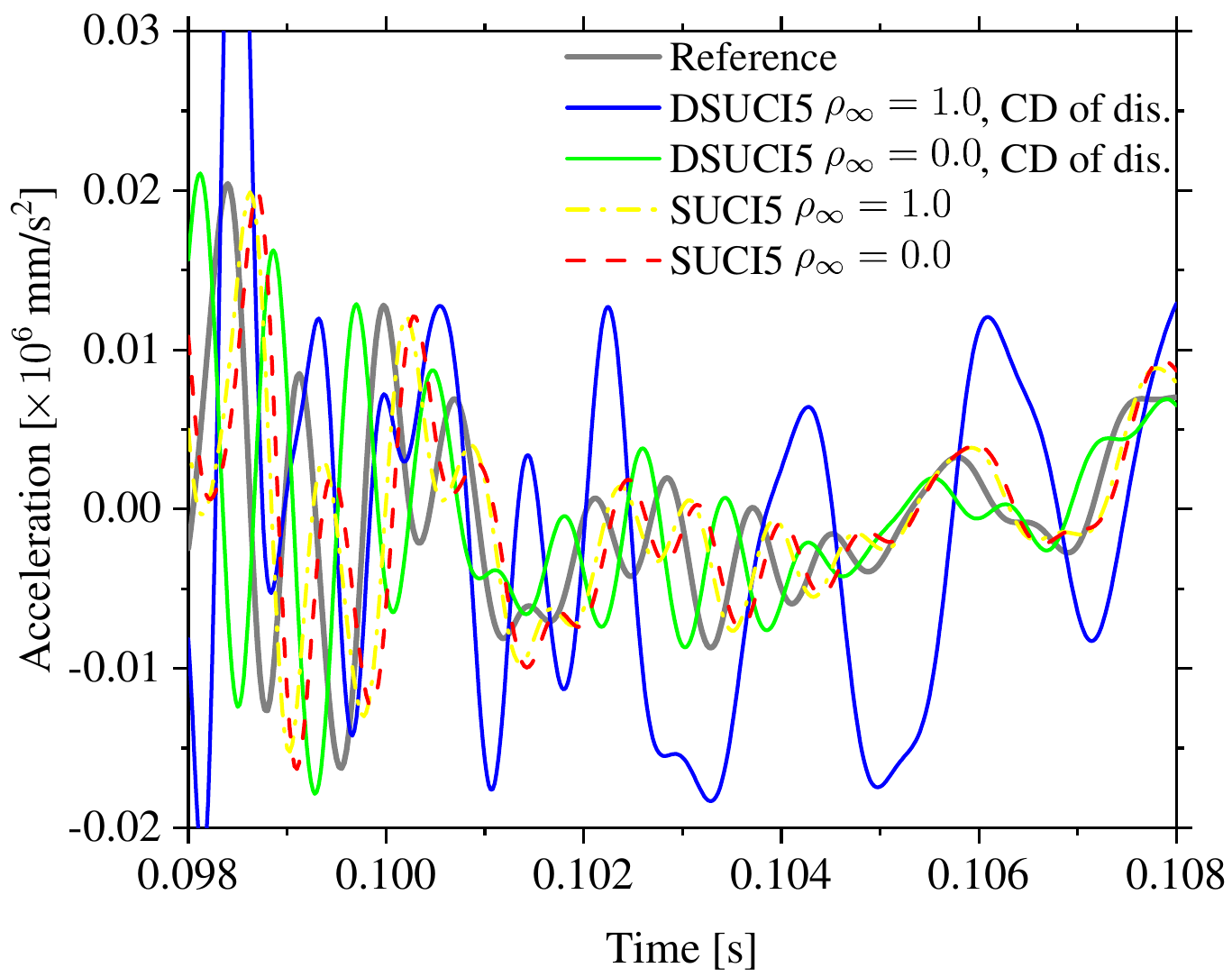}}
	\subfigure[]{
		\includegraphics[scale=0.23]{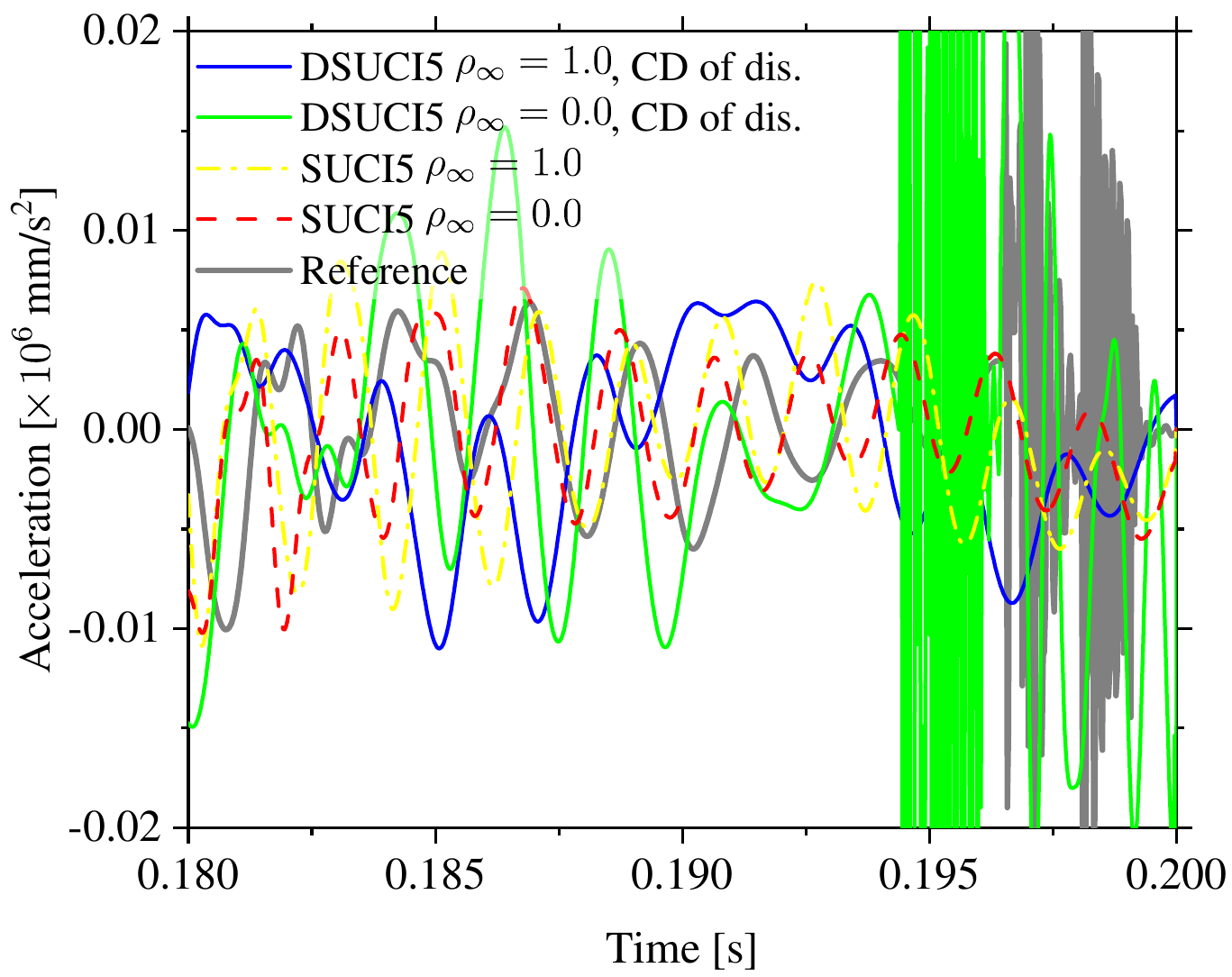}}
	\subfigure[]{
		\includegraphics[scale=0.23]{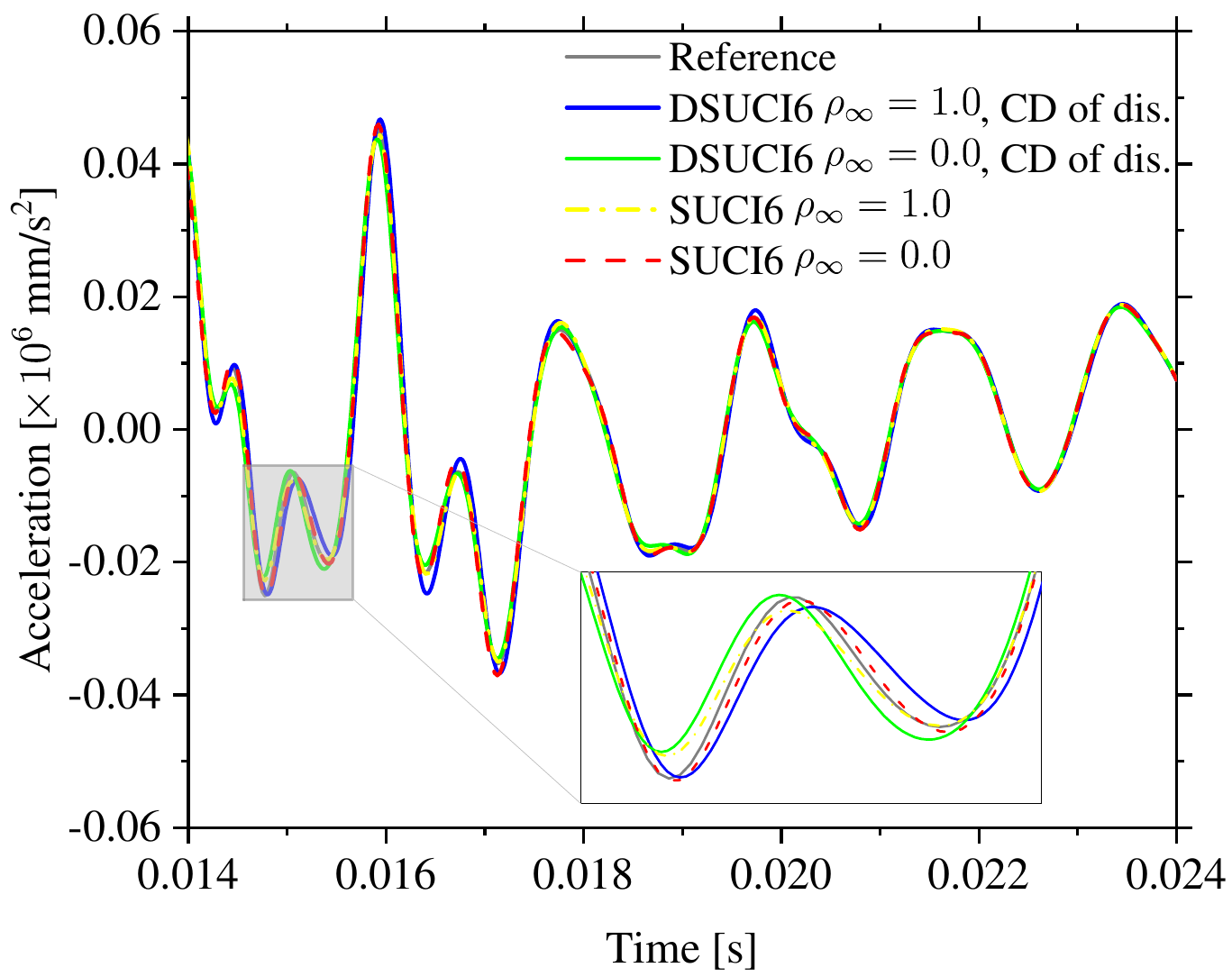}}
	\subfigure[]{
		\includegraphics[scale=0.23]{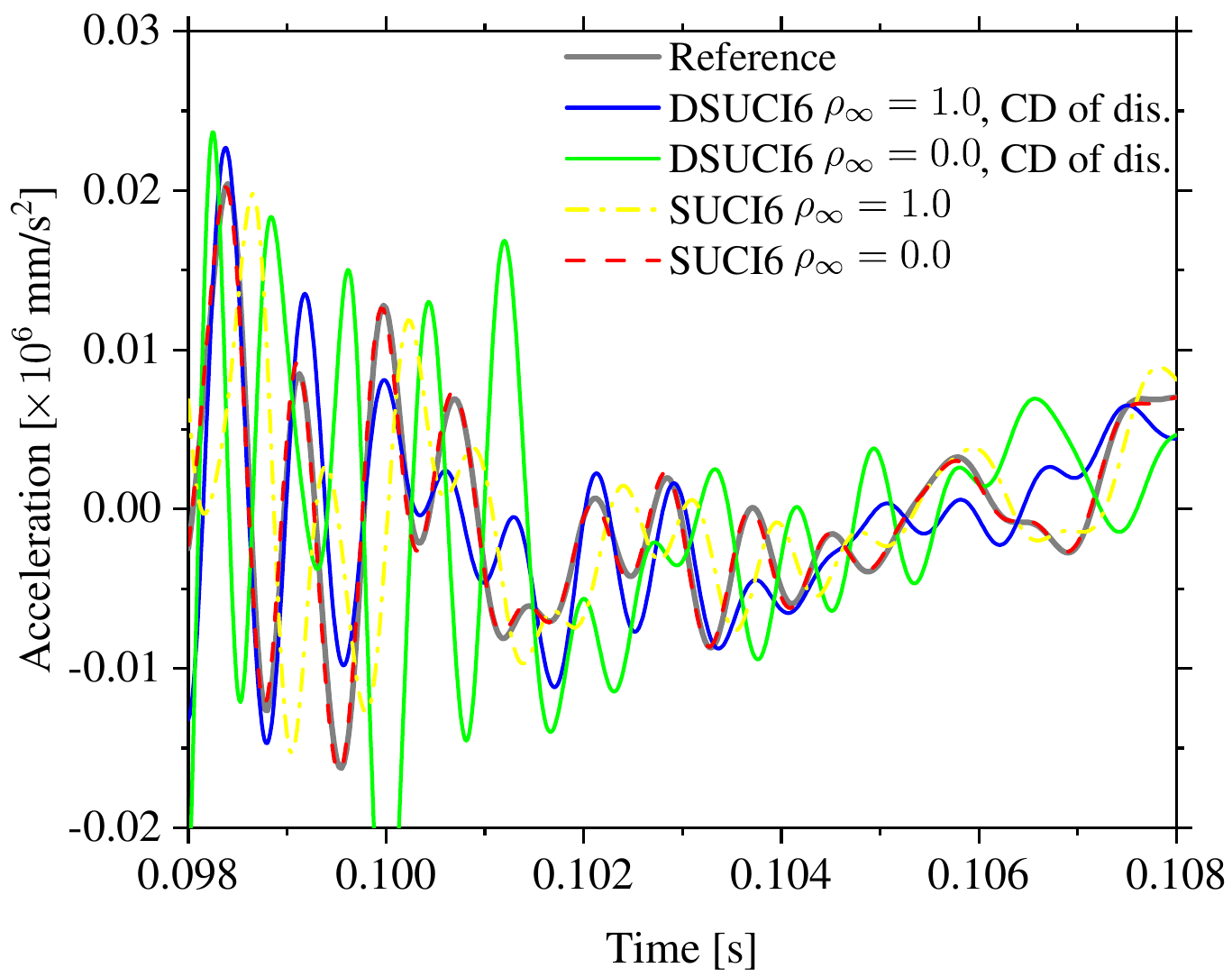}}
	\subfigure[]{
		\includegraphics[scale=0.23]{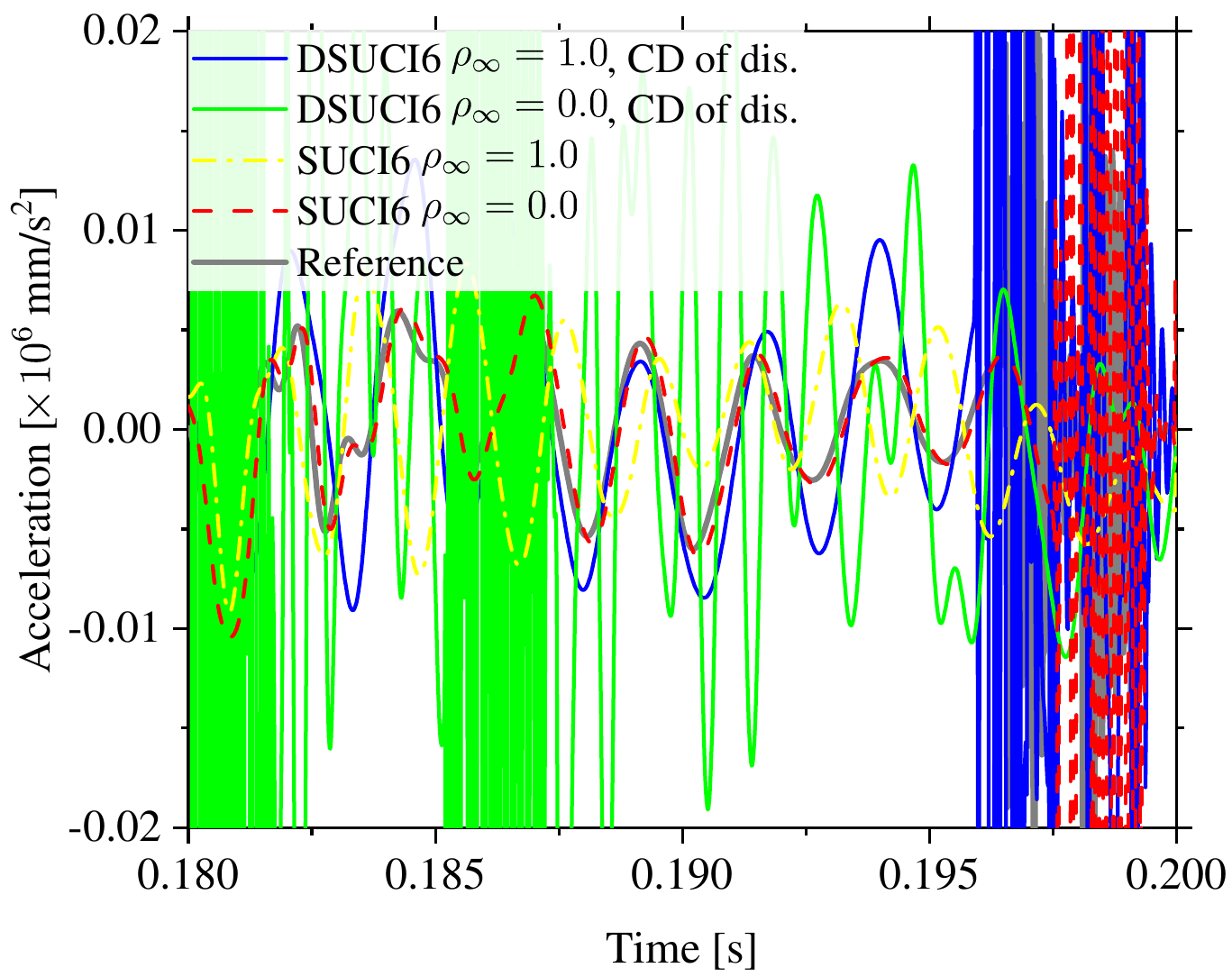}}
	\caption{Numerical accelerations of the third torsional spring using DSUCI(4-6) \cite{liDirectlySelfstarting2022} and SUCI(4-6) with the same $\dt=2.5\times10^{-5}$s. DSUCI(4-6) adopt the central difference of displacement to output accelerations.}
	\label{fig:wing_ds45}
\end{figure}

It is interesting to compare numerical accelerations between DSUCI$n$ \cite{liDirectlySelfstarting2022} and SUCI$n$, and the acceleration responses of DSUCI$n$ are produced by the central difference (CD) of displacement. Using the original post-processing way \cite{liDirectlySelfstarting2022}, DSUCI$n$ does not give more accurate accelerations than using the CD of displacement. This is the main reason for using the latter. Fig.~\ref{fig:wing_ds1} depicts numerical accelerations of the third torsional spring predicted by DSUCI3 and SUCI3. %Note that DSUCI3 provides almost the same acceleration responses when using the central difference of displacement and the original post-processing technique \cite{liDirectlySelfstarting2022}, so acceleration responses from the central difference of displacement are only compared herein. This is the case for SUCI3. 
	Subplots (b-d) of Fig.~\ref{fig:wing_ds1} illustrate that SUCI3 has better acceleration solution accuracy than DSUCI3 at the early stage of the analysis and this advantage gradually fades over the integration time. In particular, although DSUCI3 with $\rhoinf=0$ follows the reference solution well in displacement (see Fig.~\ref{fig:wingsol}(a)), it cannot maintain this advantage in acceleration by using the CD of displacement. Fig.~\ref{fig:wing_ds45} further compares acceleration responses between DSUCI(4-6) and SUCI(4-6), and the same conclusions as those in Fig.~\ref{fig:wing_ds1} can be found. For the present folding wing, SUCI$n$ can predict more accurate acceleration responses than DSUCI$n$ \cite{liDirectlySelfstarting2022} at the early stage of the analysis, and this superiority persists when DSUCI$n$ uses either the original technique \cite{liDirectlySelfstarting2022} or the CD of displacement to output accelerations. Note also that as the integration time increases, these two families of high-order algorithms find it difficult to follow the reference acceleration solution well.

In practical applications, users often need to simulate the analyzed structure more than once via selecting different integration algorithms, algorithmic parameters, and time steps. Based on the spectral analysis in Section \ref{sec:sp} and numerical results in this section, the proposed SUCI$n$ schemes are highly recommended for various dynamic problems in priority.

\section{Conclusions}\label{sec:conclusion}
%This paper constructs a general composite $ s $-sub-step implicit method (\ref{eq:nsubstep}) grounded in the explicit singly diagonally implicit
%Runge-Kutta (ESDIRK) family for second-order hyperbolic problems. Via the accuracy and dissipation analysis, the four novel high-order implicit algorithms denoted by SUCI$n$ are developed to achieve the following numerical characteristics.
This paper formulates a composite $ s $-sub-step implicit method (\ref{eq:nsubstep}), based on the explicit singly diagonally implicit Runge-Kutta (ESDIRK) family, specifically designed for addressing second-order hyperbolic problems. Through a thorough examination of accuracy, dissipation and stability, we develop four innovative high-order implicit algorithms, denoted as SUCI$n$, which have been devised to attain specific numerical characteristics.
\begin{itemize}%[(a)]
	\item The novel methods are identically higher-order accurate and avoid the order reduction for solving forced vibrations. The analysis demonstrates that the $ s $-sub-step implicit method (\ref{eq:nsubstep}) can attain $ s $th-order accuracy within the range $ 2 \leq s \leq 6 $; beyond this range, additional sub-steps are required to achieve higher-order accuracy while incorporating dissipation control and ensuring unconditional stability. Consequently, this paper focuses exclusively on the development of four cost-optimized high-order implicit methods corresponding to three, four, five, and six sub-steps.
	\item The novel methods proficiently manage numerical dissipation in the high-frequency range through the user-specified parameter $ \rhoinf\in\left[0,~1\right] $. This adjustment proves highly effective in eliminating spurious high-frequency components, thereby enabling the precise integration of critical low-frequency modes.
	\item The novel methods achieve identical effective stiffness matrices within each sub-step. Numerous investigations have demonstrated that this characteristic not only serves to significantly economize computational expenses in the solution of linear structures but also yields optimal spectral properties, specifically maximizing high-frequency dissipation while minimizing period errors.
\end{itemize}
%Apart from these, SUCI$n$ also possesses some primary numerical properties. For instance, the unconditional stability is ensured during the process of achieving controllable numerical dissipation, which is also confirmed in the spectral analysis. In addition, the third-order consistency within each sub-step is required and the trapezoidal rule has to thus be employed in the first sub-steps. Of course, four novel methods are naturally self-starting and do not suffer from overshoots.
In addition to these characteristics, the SUCI$n$ algorithms exhibit several primary numerical properties. For instance, they ensure unconditional stability, a feature corroborated through spectral analysis. Furthermore, the third-order consistency within each sub-step necessitates the utilization of the trapezoidal rule in the initial sub-steps. Importantly, the novel methods are inherently self-starting and exhibit immunity to overshoots.

%Numerical examples are solved to confirm the numerical performance and advantages of SUCI$n$. Two typical but straightforward linear systems, namely a standard SDOF damped system subjected to the external force and a double-degree-of-freedom mass-spring model with the spurious high-frequency mode, are solved to test the numerical accuracy and dissipation control of SUCI$n$. SUCI$n$ and the published schemes are also used to solve three nonlinear problems. When adopting the same integration step size, the higher-order methods are generally superior to the lower-order ones with respect to solving nonlinear dynamics. Furthermore, when considering the same order of accuracy and computational costs, the four novel methods perform better than the published high-order schemes.

Numerical examples are employed to validate the numerical performance and advantages of SUCI$n$. Two prototypical yet straightforward linear systems---a standard SDOF damped system subjected to an external force and a 2-DOFs mass-spring model featuring a spurious high-frequency mode---are solved to assess the numerical accuracy and dissipation control of SUCI$n$. Moreover, the SUCI$n$ and published schemes are employed to address three nonlinear problems. Notably, when employing same time step sizes, the higher-order methods consistently outperform lower-order ones in solving nonlinear dynamics. In general, considering the same order of accuracy and computational costs, the four novel methods demonstrate superior performance compared to the existing high-order schemes.

\section*{CRediT authorship contribution statement}
\textbf{Jinze Li:} Writing-Original Draft, Formal analysis, Methodology, Software. \textbf{Hua Li:} Writing-Review \& Editing, Supervision. \textbf{Kaiping Yu:} Supervision, Funding acquisition. \textbf{Rui Zhao:} Supervision, Funding acquisition.

\section*{Acknowledgments}
This work is supported by the National Natural Science Foundation of China (No.~12102103 and 12272105), Fundamental Research Funds for the Central Universities (No.~HIT.NSRIF.2020014) and National Postdoctoral Fellowship Program (No.~GZC20233464). In addition, the first author acknowledges the financial support by the China Scholarship Council (No.~202006120104).

\section*{Conflict of interest}
The authors declare that they have no known competing financial interests or personal relationships that could have appeared to influence the work reported in this paper. 

\section*{Data availability}
Data sharing not applicable to this article as no datasets were generated or analyzed during the current study.

\begin{appendices}
	\renewcommand{\theequation}{A\arabic{equation}}
\section*{Appendix A: The external load in accuracy analysis}\label{app:A}
In the development of high-order algorithms, the effect of external loads on the accuracy must be taken into account; otherwise, the developed high-order algorithms may suffer from the order reduction. For example, the EG3 \cite{fungExtrapolatedGalerkin1996} and MSSTH$n$ \cite{zhangOptimizationNsubstep2020} algorithms do not provide the designed order of accuracy for solving general structures.
Inspired by the work of Fung \cite{fungUnconditionallyStable1996}, the authors briefly explain the reason for using $f(t)=\exp(t)$ in accuracy analysis and give the necessary conditions that the external load $f(t)$ should fulfill in accuracy analysis. 

When advancing numerical solutions at $t_n$ to those at $t_{n+1}=t_n+\dt$ and analyzing local truncation errors during this process, the external load $f(t)$ can be expanded into Taylor series at $t_n$ over a small time increment $\dt$ as
\begin{equation}\label{eq:ft}
	f(\tau)=f(t_n)+\dot{f}(t_n)\tau+\dfrac12\ddot{f}(t_n)\tau^2+\dfrac{1}{3!}\dddot{f}(t_n)\tau^3+\cdots+\dfrac{1}{n!}f^{(n)}(t_n)\tau^n+\cdots\quad\text{for}\quad 0\le\tau\le\dt\quad\text{and}\quad\tau=t-t_n.
\end{equation}
If $f(t)$ is a special function such as a step or delta function, some necessary treatments are made. In general, the accuracy analysis of time integration methods assumes that the external load $f(t)$ satisfies the required smoothness and differentiability. Since the SDOF system \eqref{eq:sdof} is linear, only one term in Eq.~\eqref{eq:ft} needs to be considered without loss of generality. Therefore, assuming that  $f(\tau)=c\cdot \tau^n$, where $c$ is the amplitude coefficient and $n$ is an integer, the considered Eq.~\eqref{eq:sdof} is reduced to 
\begin{equation}\label{eq:tau}
	\ddot{u}(\tau)+2\xi\omega\dot{u}(\tau)+\omega^2u(\tau)=c\cdot \tau^n
\end{equation}
where $\tau=t-t_n$. As well-known, the exact solution of Eq.~\eqref{eq:tau} with initial conditions $u_n$ and $\dot{u}_n$ is divided into two parts: the homogeneous and particular solutions. The homogeneous solution corresponds to the free vibration (i.e., $c=0$), and its expression is omitted herein. One of the particular solutions of Eq.~\eqref{eq:tau} with initial conditions $u_n=0$ and $\dot{u}_n=0$ can be calculated \cite{fungUnconditionallyStable1996} as
\begin{equation}\label{eq:up}
	\overline{u}_p(\tau)=c\times \left\{\dfrac{n!\tau^{n+2}}{(n+2)!}-\frac{2n!\xi\omega\tau^{n+3}}{(n+3)!}-\dfrac{n!(1-4\xi^2)\omega^2\tau^{n+4}}{(n+4)!}+\dfrac{n!4\xi(1-2\xi^2)\omega^3\tau^{n+5}}{(n+5)!}+\cdots \right\}.
\end{equation}
The equation above illustrates that the leading term of displacement for the external load $c\cdot \tau^n$ is proportional to $\tau^{n+2}$. Furthermore, the exact velocity is determined from the time derivative of Eq.~\eqref{eq:up} with respect to $\tau$, and thus the leading term for velocity is proportional to $\tau^{n+1}$. 

On the other hand, Proposition \ref{pos:accuracy} states that the novel method \eqref{eq:nsubstep} achieves $p$th-order accuracy for the standard SDOF system \eqref{eq:sdof} if and only if both $u_{n+1}-u(t_{n+1})=O(\dt^{p+1}) $ and $\dot{u}_{n+1}-\dot{u}(t_{n+1})=O(\dt^{p+1}) $. As a result, the external load $f(\tau)$ proportional to $\tau^p$, $\tau^{p+1}$, etc., does not need to be considered since the terms $\tau^{j}~(j\ge p)$ do not contribute to necessary conditions. In other words, \textbf{the terms from $\tau^0$ to $\tau^{p-1}$ are required and included in the external load $f(\tau)$ to derive complete conditions for $p$th-order accuracy}. 

Therefore, the external load of form  $f(t)=\sum_{n=0}^{2}c_nt^n$ with $c_n\neq 0$ can be used to develop the third-order algorithm, but it cannot be used for more than third-order algorithms. To avoid constructing the specific external load for each high-order algorithm, the exponential function $f(t)=\exp(t)$ can be used uniformly because it contains all items required by each high-order algorithm due to $\exp(t)=\sum_{n=0}^{\infty}{t^n}/{n!}$. This is the main reason why $f(t)=\exp(t)$ is used in this paper to derive complete conditions for achieving $p$th-order accuracy. The above analysis demonstrates that other appropriate external loads, such as $f(t)=\sin(t)+\cos(t)$, can also be used in accuracy analysis, but only using $f(t)=\sin(t)$ or $f(t)=\cos(t)$ cannot achieve identical high-order accuracy due to losing even terms $t^{2k}~(k=0,~1,~2,~\cdots)$ or odd terms $t^{2k+1}~(k=0,~1,~2,~\cdots)$. 

\setcounter{equation}{0}
\renewcommand{\theequation}{B\arabic{equation}}
\section*{Appendix B: Six-sub-step sixth-order scheme: SUCI6}\label{app:B}
In this appendix, the six-sub-step implicit algorithm given by Eq.~\eqref{eq:6substep} is developed. %Numerical schemes in the first four sub-steps are the same as those of SUCI5, and numerical schemes in the last two sub-steps are explicitly written as follows.
\begin{equation}\label{eq:6substep}
	\begin{BMAT}(b){c|c}{c|c}
		\mbf{c} & \mbf{A}\\ %\hline
		& \mbf{b}
	\end{BMAT}=\begin{BMAT}(@,25pt,10pt){c|ccccccc}{ccccccc|c}
		0 & 0 &&&&&&\\
		\gamma_1 & \dfrac{\gamma_1}{2} & \dfrac{\gamma_1}{2}&&&&&\\
		\gamma_2 & \alpha_{20} & \alpha_{21} & \dfrac{\gamma_1}{2}&&&&\\
		\gamma_3 & \alpha_{30} & \alpha_{31} & \alpha_{32} & \dfrac{\gamma_1}{2}&&&\\
		\gamma_4 & \alpha_{40} & \alpha_{41} & \alpha_{42} & \alpha_{43} & \dfrac{\gamma_1}{2}&&\\
		\gamma_5 & \alpha_{50} & \alpha_{51} & \alpha_{52} & \alpha_{53} & \alpha_{54}& \dfrac{\gamma_1}{2}&\\
		1 & \alpha_{60} & \alpha_{61} & \alpha_{62} & \alpha_{63} & \alpha_{64}& \alpha_{65} & \dfrac{\gamma_1}{2}\\
		&  \alpha_{60} & \alpha_{61} & \alpha_{62} & \alpha_{63} & \alpha_{64}& \alpha_{65} & \dfrac{\gamma_1}{2}
	\end{BMAT}
\end{equation}
%\begin{subequations}\label{eq:6substep}
%	\begin{align}
%		\mbf{M}\ddot{\mbf{U}}_{n+\gamma_5}&+\mbf{C}\dot{\mbf{U}}_{n+\gamma_5}+\mbf{K}\mbf{U}_{n+\gamma_5}  =\mbf{F}(t_{n}+\gamma_5\dt)   &\mbf{M}\ddot{\mbf{U}}_{n+1}&+\mbf{C}\dot{\mbf{U}}_{n+1}+\mbf{K}\mbf{U}_{n+1}                       =\mbf{F}(t_{n+1})                                                                                                                                                                                                                              \\
%		\mbf{U}_{n+\gamma_5}                                                                             & =\mbf{U}_n+\dt\left(\sum_{j=0}^4\alpha_{5j}\dot{\mbf{U}}_{n+\gamma_j}+\frac{\gamma_1}{2}\dot{\mbf{U}}_{n+\gamma_5}\right)    & \mbf{U}_{n+1}                                                                                    & =\mbf{U}_n+\dt\left(\sum_{j=0}^5\alpha_{6j}\dot{\mbf{U}}_{n+\gamma_j}+\frac{\gamma_1}{2}\dot{\mbf{U}}_{n+1}\right)                     \\
%		\dot{\mbf{U}}_{n+\gamma_5}                                                                       & =\dot{\mbf{U}}_n+\dt\left(\sum_{j=0}^4\alpha_{5j}\ddot{\mbf{U}}_{n+\gamma_j}+\frac{\gamma_1}{2}\ddot{\mbf{U}}_{n+\gamma_5}\right) &
%		\dot{\mbf{U}}_{n+1}                                                                              & =\dot{\mbf{U}}_n+\dt\left(\sum_{j=0}^5\alpha_{6j}\ddot{\mbf{U}}_{n+\gamma_j}+\frac{\gamma_1}{2}\ddot{\mbf{U}}_{n+1}\right)
%	\end{align}
%\end{subequations}
Notice that the above scheme has achieved identical effective stiffness matrices within each sub-step and the non-dissipative trapezoidal rule has been set in the first sub-step. The conditions (\ref{eq:accuracy}) for achieving sixth-order accuracy, as well as the additional constraints (\ref{eq:2orderi}), are used to determine $ \alpha_{ij} $, which are\vskip -10pt
\begin{subequations}\label{eq:alpha6sub1}
	\begin{align}
		\alpha_{20}&=\frac{-\gamma_1^2+3\gamma_1\gamma_2-\gamma_2^2}{2\gamma_1}                          \hspace{6cm} \alpha_{21}=\frac{\gamma_2(\gamma_2-\gamma_1)}{2\gamma_1}                                            \\
		\alpha_{30}&=\frac{-\gamma_1^2+(3\gamma_3-2\alpha_{32})\gamma_1+2\alpha_{32}\gamma_2-\gamma_3^2}{2\gamma_1}                                    \hspace{3.4cm} \alpha_{31}=\frac{-2\alpha_{32}\gamma_2-\gamma_1\gamma_3+\gamma_3^2}{2\gamma_1}                      \\
		\alpha_{40}&=\frac{-\gamma_1^2+(3\gamma_4-2\alpha_{42}-2\alpha_{43})\gamma_1+2\alpha_{42}\gamma_2+2\alpha_{43}\gamma_3-\gamma_4^2}{2\gamma_1}  \hspace{1.05cm} \alpha_{41}=\frac{-2\alpha_{42}\gamma_2-2\alpha_{43}\gamma_3-\gamma_1\gamma_4+\gamma_4^2}{2\gamma_1}\\
		\alpha_{50} & =\frac{-\gamma_1^2+(3\gamma_5-2\alpha_{52}-2\alpha_{53}-2\alpha_{54})\gamma_1+2\alpha_{52}\gamma_2+2\alpha_{53}\gamma_3+2\alpha_{54}\gamma_4-\gamma_5^2}{2\gamma_1}  \\
		\alpha_{51} & =\frac{-2\alpha_{52}\gamma_2-2\alpha_{53}\gamma_3-2\alpha_{54}\gamma_4-\gamma_1\gamma_5+\gamma_5^2}{2\gamma_1}\\
		\alpha_{60} & =\frac{-\gamma_1^2+(3-2\alpha_{62}-2\alpha_{63}-2\alpha_{64}-2\alpha_{65})\gamma_1+2\alpha_{62}\gamma_2+2\alpha_{63}\gamma_3+2\alpha_{64}\gamma_4+2\alpha_{65}\gamma_5-1}{2\gamma_1}  \\
		\alpha_{61} & =\frac{-2\alpha_{62}\gamma_2-2\alpha_{63}\gamma_3-2\alpha_{64}\gamma_4-2\alpha_{65}\gamma_5-\gamma_1+1}{2\gamma_1}\\
		\alpha_{32} & =\frac{45\gamma_1^5-225\gamma_1^4+300\gamma_1^3-150\gamma_1^2+30\gamma_1-2}{720\alpha_{65}\alpha_{54}\alpha_{43}\gamma_2(\gamma_1-\gamma_2)}\\
		\alpha_{42}&=\frac{\left\{\begin{aligned}
				&-45  (\alpha_{64}\alpha_{43}+\alpha_{65}\alpha_{53})\gamma_1^5+(225\alpha_{64}\alpha_{43}+225\alpha_{65}\alpha_{53}-90\alpha_{65}\alpha_{54}\alpha_{43})\gamma^4+720\alpha_{65}^2\alpha_{54}^2\alpha_{43}^2\gamma_3^2 \\
				& +(360\alpha_{65}\alpha_{54}\alpha_{43}-300\alpha_{64}\alpha_{43}-300\alpha_{65}\alpha_{53})\gamma_1^3    
				+(150\alpha_{64}\alpha_{43}+150\alpha_{65}\alpha_{53}-360\alpha_{65}\alpha_{54}\alpha_{43})\gamma_1^2                                                       \\
				& +(120\alpha_{65}\alpha_{54}\alpha_{43}-30\alpha_{64}\alpha_{43}-30\alpha_{65}\alpha_{53}-720\alpha_{65}^2\alpha_{54}^2\alpha_{43}^2\gamma_3)\gamma_1        
				-12\alpha_{65}\alpha_{54}\alpha_{43}+2\alpha_{64}\alpha_{43}+2\alpha_{65}\alpha_{53}
			\end{aligned}\right\}}{720\alpha_{65}^2\alpha_{54}^2\alpha_{43}\gamma_2(\gamma_1-\gamma_2)}\\
		\alpha_{52} & =\frac{\left\{\begin{aligned}
				&45  (\alpha_{64}^2\alpha_{43}+\alpha_{65}\alpha_{64}\alpha_{53}-\alpha_{65}\alpha_{63}\alpha_{54})\gamma_1^5-12\alpha_{65}\alpha_{54}\alpha_{43}(5\alpha_{65}\alpha_{54}-\alpha_{64})+2\alpha_{65}(\alpha_{63}\alpha_{54}-\alpha_{64}\alpha_{53}) \\
				& +(90\alpha_{65}\alpha_{64}\alpha_{54}\alpha_{43}-225\alpha_{64}^2\alpha_{43}-225\alpha_{65}\alpha_{64}\alpha_{53}+225\alpha_{65}\alpha_{63}\alpha_{54})\gamma_1^4+720\alpha_{65}^3\alpha_{54}^2\alpha_{43}(\alpha_{54}\gamma_4^2+\alpha_{53}\gamma_3^2)  \\
				& +\left[180\alpha_{65}\alpha_{54}\alpha_{43}(\alpha_{65}\alpha_{54}-2\alpha_{64})+300(\alpha_{64}^2\alpha_{43}+\alpha_{65}\alpha_{64}\alpha_{53}-\alpha_{65}\alpha_{63}\alpha_{54}) \right]\gamma_1^3-2\alpha_{43}\alpha_{64}^2   \\
				& +\left[180\alpha_{65}\alpha_{54}\alpha_{43}(2\alpha_{64}-3\alpha_{65}\alpha_{54})-150(\alpha_{64}^2\alpha_{43}+\alpha_{65}\alpha_{64}\alpha_{53}-\alpha_{65}\alpha_{63}\alpha_{54}) \right]\gamma_1^2 \\
				& +\left[360\alpha_{65}^2\alpha_{54}^2\alpha_{43}(1-2\alpha_{65}(\alpha_{54}\gamma_4+\alpha_{53}\gamma_3))-30\alpha_{64}\alpha_{43}(4\alpha_{65}\alpha_{54}-\alpha_{64})+30\alpha_{65}(\alpha_{64}\alpha_{53}-\alpha_{63}\alpha_{54}) \right]\gamma_1
			\end{aligned}\right\}}{720\alpha_{65}^3\alpha_{54}^2\alpha_{43}\gamma_2(\gamma_1-\gamma_2)}                                                                \\
		\alpha_{62} & =\frac{6\alpha_{63}\gamma_3(\gamma_1-\gamma_3)+6\alpha_{64}\gamma_4(\gamma_1-\gamma_4)+6\alpha_{65}\gamma_5(\gamma_1-\gamma_5)+3\gamma_1^2-6\gamma_1+2}{6\gamma_2(\gamma_2-\gamma_1)}\\
		\alpha_{43} & =\frac{
			\left\{\begin{aligned}
				60 & \alpha_{65}\alpha_{53}\gamma_3(\gamma_2-\gamma_3)(\gamma_3-\gamma_1)+60\alpha_{65}\alpha_{54}\gamma_4(\gamma_2-\gamma_4)(\gamma_4-\gamma_1) 
				+5\gamma_1^2(3\gamma_1\gamma_2-3\gamma_1-9\gamma_2+7)\\
				&+30\gamma_1\gamma_2-20\gamma_1-5\gamma_2+3
			\end{aligned} \right\}
		}{60\alpha_{64}\gamma_3(\gamma_2-\gamma_3)(\gamma_1-\gamma_3)} \\
		\alpha_{53} & =\frac{\left\{\begin{aligned}
				240 & \alpha_{65}^2\alpha_{54}^2\gamma_4(\gamma_2-\gamma_4)(\gamma_4-\gamma_1)+30\alpha_{64}\gamma_1^4(\gamma_2-1) 
				+\left[60\alpha_{65}\alpha_{54}(\gamma_2-1)+30\alpha_{64}(3-4\gamma_2)\right]\gamma_1^3                      \\
				& +\left[20\alpha_{65}\alpha_{54}(7-9\gamma_2)+15\alpha_{64}(8\gamma_2-5) \right]\gamma_1^2                    
				+\left[40\alpha_{65}\alpha_{54}(3\gamma_2-2)+2\alpha_{64}(11-20\gamma_2)\right]\gamma_1                      \\
				& +4\alpha_{65}\alpha_{54}(3-5\gamma_2)+2\alpha_{64}(2\gamma_2-1)
			\end{aligned}\right\}}{240\alpha_{65}^2\alpha_{54}\gamma_3(\gamma_2-\gamma_3)(\gamma_1-\gamma_3)}              \\
		\alpha_{63} & =\frac{
			12  \alpha_{64}\gamma_4(\gamma_2-\gamma_4)(\gamma_4-\gamma_1)+12\alpha_{65}\gamma_5(\gamma_2-\gamma_5)(\gamma_5-\gamma_1) 
			-6\gamma_1(\gamma_1(\gamma_2-1)-2\gamma_2)-10\gamma_1-4\gamma_2+3
		}{12\gamma_3(\gamma_2-\gamma_3)(\gamma_1-\gamma_3)}                                       \\
		\alpha_{54} & =\frac{\left\{\begin{aligned}
				15 & (\gamma_3-1)(\gamma_2-1)\gamma_1^3+\left[35\gamma_3-30+5\gamma_2(7-9\gamma_3) \right]\gamma_1^2 
				+\left[15-20\gamma_3+10\gamma_2(3\gamma_3-2)\right]\gamma_1\\
				&+(3-5\gamma_3)\gamma_2+3\gamma_3-2
			\end{aligned}\right\}}{60\alpha_{65}\gamma_4(\gamma_3-\gamma_4)(\gamma_2-\gamma_4)(\gamma_1-\gamma_4)}         \\
		\alpha_{64} & =\frac{\left\{\begin{aligned}
				30 & (1-\gamma_3)(\gamma_2-1)\gamma_1^2+\left[45-50\gamma_3+10\gamma_2(6\gamma_3-5)-60\alpha_{65}\gamma_5(\gamma_3-\gamma_5)(\gamma_2-\gamma_5) \right]\gamma_1 \\
				& +\left[15-20\gamma_3+60\alpha_{65}\gamma_5^2(\gamma_3-\gamma_5) \right]\gamma_2-60\alpha_{65}\gamma_5^3(\gamma_3-\gamma_5)+15\gamma_3-12
			\end{aligned}\right\}}{60\gamma_4(\gamma_3-\gamma_4)(\gamma_2-\gamma_4)(\gamma_1-\gamma_4)}                    \\
		\alpha_{65} & =\frac{\left\{\begin{aligned}
				30 & (1-\gamma_4)(\gamma_3-1)(\gamma_2-1)\gamma_1^2+(45-50\gamma_4+10\gamma_3(6\gamma_4-5))\gamma_1\gamma_2+3\gamma_3(5\gamma_4-4) \\
				& +(45\gamma_4-42+5\gamma_3(9-10\gamma_4))\gamma_1+(15\gamma_4-12+5\gamma_3(3-4\gamma_4))\gamma_2-12\gamma_4+10
			\end{aligned}\right\}}{60\gamma_5(\gamma_4-\gamma_5)(\gamma_3-\gamma_5)(\gamma_2-\gamma_5)(\gamma_1-\gamma_5)}
	\end{align}
\end{subequations}

The remaining algorithmic parameters are only five splitting ratios of sub-step size, namely $ \gamma_i~(i=1,~\cdots,~5) $. Like the previous sub-step methods, the first splitting ratio $ \gamma_1 $ is given to control numerical dissipation in the high-frequency range. In the high-frequency limit ($ \omega\to\infty $), the characteristic polynomial (\ref{eq:cp1}) is simplified using Eq.~\eqref{eq:alpha6sub1} as
\begin{equation}\label{key}
	\left(\zeta_{\infty}-\frac{45\gamma_1^6-540\gamma_1^5+1350\gamma_1^4-1200\gamma_1^3+450\gamma_1^2-72\gamma_1+4}{45\gamma_1^6}\right)^2=0.
\end{equation}
Then, the conditions given by Eq.~(\ref{eq:optimaldissipation}) to achieve controllable numerical dissipation require
\begin{equation}\label{eq:rhoinf6}
	\frac{45\gamma_1^6-540\gamma_1^5+1350\gamma_1^4-1200\gamma_1^3+450\gamma_1^2-72\gamma_1+4}{45\gamma_1^6}=\rhoinf.
\end{equation}
On the other hand, the unconditional stability given by Eq.~\eqref{eq:uc} also imposes constraints on the parameter $\gamma_1$, which is
\begin{equation}\label{eq:suci6_g1}
	\gamma_1\in\left[0.5681292760,~1.081813756\right].
\end{equation}
Hence, Eqs.~\eqref{eq:rhoinf6} and \eqref{eq:suci6_g1} give an proper value of $\gamma_1$ to achieve simultaneously dissipation control and unconditional stability for each given $\rhoinf\in\left[0,~1\right]$, as shown in Fig.~\ref{fig:sucin_g1}(d). These selected values of $\gamma_1$ are given in Table \ref{tab:gamma1}. Other splitting ratios $ \gamma_i~(i=2,~\cdots,~5) $ are taken by default as $ \gamma_i=i\cdot\gamma_1~(i=2,~\cdots,~5) $ in this paper.

\end{appendices}

{\small
\bibliography{Reference}

\begin{thebibliography}{10}
\expandafter\ifx\csname url\endcsname\relax
  \def\url#1{\burl{#1}}\fi
\expandafter\ifx\csname urlprefix\endcsname\relax\def\urlprefix{URL }\fi
\providecommand{\bibinfo}[2]{#2}
\providecommand{\eprint}[2][]{\url{#2}}
\providecommand{\doi}[1]{\url{https://doi.org/#1}}
\bibcommenthead

\bibitem{hughesFiniteElement2000}
\bibinfo{author}{Hughes, T. J.~R.}
\newblock \emph{\bibinfo{title}{The {{Finite Element Method}}: {{Linear
  Static}} and {{Dynamic Finite Element Analysis}}}} Dover {{Civil}} and
  {{Mechanical Engineering}} (\bibinfo{publisher}{{Dover Publications}},
  \bibinfo{year}{2000}).

\bibitem{rezaiee-pajandMoreAccurate2015}
\bibinfo{author}{{Rezaiee-Pajand}, M.} \& \bibinfo{author}{{Karimi-Rad}, M.}
\newblock \bibinfo{title}{More accurate and stable time integration scheme}.
\newblock \emph{\bibinfo{journal}{Engineering with Computers}}
  \textbf{\bibinfo{volume}{31}}, \bibinfo{pages}{791--812}
  (\bibinfo{year}{2015}).

\bibitem{liIdenticalSecondorder2021}
\bibinfo{author}{Li, J.}, \bibinfo{author}{Yu, K.} \& \bibinfo{author}{Li, X.}
\newblock \bibinfo{title}{An identical second-order single step explicit
  integration algorithm with dissipation control for structural dynamics}.
\newblock \emph{\bibinfo{journal}{International Journal for Numerical Methods
  in Engineering}} \textbf{\bibinfo{volume}{122}}, \bibinfo{pages}{1089--1132}
  (\bibinfo{year}{2021}).

\bibitem{nohExplicitTime2013}
\bibinfo{author}{Noh, G.} \& \bibinfo{author}{Bathe, K.-J.}
\newblock \bibinfo{title}{An explicit time integration scheme for the analysis
  of wave propagations}.
\newblock \emph{\bibinfo{journal}{Computers \& Structures}}
  \textbf{\bibinfo{volume}{129}}, \bibinfo{pages}{178--193}
  (\bibinfo{year}{2013}).

\bibitem{liTwoThirdorder2022}
\bibinfo{author}{Li, J.}, \bibinfo{author}{Yu, K.} \& \bibinfo{author}{Zhao,
  R.}
\newblock \bibinfo{title}{Two third-order explicit integration algorithms with
  controllable numerical dissipation for second-order nonlinear dynamics}.
\newblock \emph{\bibinfo{journal}{Computer Methods in Applied Mechanics and
  Engineering}} \textbf{\bibinfo{volume}{395}}, \bibinfo{pages}{114945}
  (\bibinfo{year}{2022}).

\bibitem{rezaiee-pajandNewExplicit2016}
\bibinfo{author}{{Rezaiee-Pajand}, M.} \& \bibinfo{author}{{Karimi-Rad}, M.}
\newblock \bibinfo{title}{A new explicit time integration scheme for nonlinear
  dynamic analysis}.
\newblock \emph{\bibinfo{journal}{International Journal of Structural Stability
  and Dynamics}} \textbf{\bibinfo{volume}{16}}, \bibinfo{pages}{1550054}
  (\bibinfo{year}{2016}).

\bibitem{newmarkMethodComputation1959}
\bibinfo{author}{Newmark, N.~M.}
\newblock \bibinfo{title}{A method of computation for structural dynamics}.
\newblock \emph{\bibinfo{journal}{Journal of Engineering Mechanic Division}}
  \textbf{\bibinfo{volume}{85}}, \bibinfo{pages}{67--94}
  (\bibinfo{year}{1959}).

\bibitem{hilberImprovedNumerical1977}
\bibinfo{author}{Hilber, H.~M.}, \bibinfo{author}{Hughes, T. J.~R.} \&
  \bibinfo{author}{Taylor, R.~L.}
\newblock \bibinfo{title}{Improved numerical dissipation for time integration
  algorithms in structural dynamics}.
\newblock \emph{\bibinfo{journal}{Earthquake Engineering \& Structural
  Dynamics}} \textbf{\bibinfo{volume}{5}}, \bibinfo{pages}{283--292}
  (\bibinfo{year}{1977}).

\bibitem{woodAlphaModification1980}
\bibinfo{author}{Wood, W.}, \bibinfo{author}{Bossak, M.} \&
  \bibinfo{author}{Zienkiewicz, O.}
\newblock \bibinfo{title}{An alpha modification of {{Newmark}}'s method}.
\newblock \emph{\bibinfo{journal}{International Journal for Numerical Methods
  in Engineering}} \textbf{\bibinfo{volume}{15}}, \bibinfo{pages}{1562--1566}
  (\bibinfo{year}{1980}).

\bibitem{shaoThreeParameters1988}
\bibinfo{author}{Shao, H.} \& \bibinfo{author}{Cai, C.}
\newblock \bibinfo{title}{A three parameters algorithm for numerical
  integration of structural dynamic equations}.
\newblock \emph{\bibinfo{journal}{Chinese Journal of Applied Mechanics}}
  \textbf{\bibinfo{volume}{5}}, \bibinfo{pages}{76--81} (\bibinfo{year}{1988}).

\bibitem{chungTimeIntegration1993}
\bibinfo{author}{Chung, J.} \& \bibinfo{author}{Hulbert, G.~M.}
\newblock \bibinfo{title}{A time integration algorithm for structural dynamics
  with improved numerical dissipation: {{The}} generalized-$\alpha$ method}.
\newblock \emph{\bibinfo{journal}{Journal of Applied Mechanics}}
  \textbf{\bibinfo{volume}{60}}, \bibinfo{pages}{371--375}
  (\bibinfo{year}{1993}).

\bibitem{batheCompositeImplicit2005}
\bibinfo{author}{Bathe, K.~J.} \& \bibinfo{author}{Baig, M. M.~I.}
\newblock \bibinfo{title}{On a composite implicit time integration procedure
  for nonlinear dynamics}.
\newblock \emph{\bibinfo{journal}{Computers \& Structures}}
  \textbf{\bibinfo{volume}{83}}, \bibinfo{pages}{2513--2524}
  (\bibinfo{year}{2005}).

\bibitem{LiFurtherAssessment2021}
\bibinfo{author}{Li, J.}, \bibinfo{author}{Yu, K.} \& \bibinfo{author}{Tang,
  H.}
\newblock \bibinfo{title}{Further assessment of three {Bathe} algorithms and
  implementations for wave propagation problems}.
\newblock \emph{\bibinfo{journal}{International Journal of Structural Stability
  and Dynamics}} \textbf{\bibinfo{volume}{21}}, \bibinfo{pages}{2150073}
  (\bibinfo{year}{2021}).

\bibitem{batheConservingEnergy2007}
\bibinfo{author}{Bathe, K.~J.}
\newblock \bibinfo{title}{Conserving energy and momentum in nonlinear dynamics:
  {{A}} simple implicit time integration scheme}.
\newblock \emph{\bibinfo{journal}{Computers \& Structures}}
  \textbf{\bibinfo{volume}{85}}, \bibinfo{pages}{437--445}
  (\bibinfo{year}{2007}).

\bibitem{bankTransientSimulation1985}
\bibinfo{author}{Bank, R.~E.}, \bibinfo{author}{Coughran, W.~M.},
  \bibinfo{author}{Grosse, E.~H.}, \bibinfo{author}{Rose, D.~J.} \&
  \bibinfo{author}{Kentsmith, R.}
\newblock \bibinfo{title}{Transient simulation of silicon devices and
  circuits}.
\newblock \emph{\bibinfo{journal}{IEEE Transaction on Electron Devices}}
  \textbf{\bibinfo{volume}{32}}, \bibinfo{pages}{16} (\bibinfo{year}{1985}).

\bibitem{liThreeOpitmal2023}
\bibinfo{author}{Li, J.}, \bibinfo{author}{Yu, K.}, \bibinfo{author}{Zhao, R.}
  \& \bibinfo{author}{Fang, Y.}
\newblock \bibinfo{title}{Three optimal families of three-sub-step dissipative
  implicit integration algorithms with either second, third, or fourth-order
  accuracy for second-order nonlinear dynamics}.
\newblock \emph{\bibinfo{journal}{International Journal for Numerical Methods
  in Engineering}} \textbf{\bibinfo{volume}{124}}, \bibinfo{pages}{3733--3766}
  (\bibinfo{year}{2023}).

\bibitem{nohBatheTime2019}
\bibinfo{author}{Noh, G.} \& \bibinfo{author}{Bathe, K.-J.}
\newblock \bibinfo{title}{The {{Bathe}} time integration method with
  controllable spectral radius: {{The}} $\rho_\infty$-{{Bathe}} method}.
\newblock \emph{\bibinfo{journal}{Computers \& Structures}}
  \textbf{\bibinfo{volume}{212}}, \bibinfo{pages}{299--310}
  (\bibinfo{year}{2019}).

\bibitem{liNovelFamily2019}
\bibinfo{author}{Li, J.}, \bibinfo{author}{Yu, K.} \& \bibinfo{author}{Li, X.}
\newblock \bibinfo{title}{A novel family of controllably dissipative composite
  integration algorithms for structural dynamic analysis}.
\newblock \emph{\bibinfo{journal}{Nonlinear Dynamics}}
  \textbf{\bibinfo{volume}{96}}, \bibinfo{pages}{2475--2507}
  (\bibinfo{year}{2019}).

\bibitem{rezaiee-pajandMixedMultistep2010}
\bibinfo{author}{{Rezaiee-Pajand}, M.} \& \bibinfo{author}{Sarafrazi, S.~R.}
\newblock \bibinfo{title}{A mixed and multi-step higher-order implicit time
  integration family}.
\newblock \emph{\bibinfo{journal}{Archive Proceedings of the Institution of
  Mechanical Engineers Part C Journal of Mechanical Engineering Science}}
  \textbf{\bibinfo{volume}{224}}, \bibinfo{pages}{2097--2108}
  (\bibinfo{year}{2010}).

\bibitem{liAlternativeBathe2019}
\bibinfo{author}{Li, J.} \& \bibinfo{author}{Yu, K.}
\newblock \bibinfo{title}{An alternative to the {{Bathe}} algorithm}.
\newblock \emph{\bibinfo{journal}{Applied Mathematical Modelling}}
  \textbf{\bibinfo{volume}{69}}, \bibinfo{pages}{255--272}
  (\bibinfo{year}{2019}).

\bibitem{liTrulySelfstarting2020}
\bibinfo{author}{Li, J.} \& \bibinfo{author}{Yu, K.}
\newblock \bibinfo{title}{A truly self-starting implicit family of integration
  algorithms with dissipation control for nonlinear dynamics}.
\newblock \emph{\bibinfo{journal}{Nonlinear Dynamics}}
  \textbf{\bibinfo{volume}{102}}, \bibinfo{pages}{2503--2530}
  (\bibinfo{year}{2020}).

\bibitem{liSimpleTruly2020}
\bibinfo{author}{Li, J.} \& \bibinfo{author}{Yu, K.}
\newblock \bibinfo{title}{A simple truly self-starting and {{L-stable}}
  integration algorithm for structural dynamics}.
\newblock \emph{\bibinfo{journal}{International Journal of Applied Mechanics}}
  \textbf{\bibinfo{volume}{12}}, \bibinfo{pages}{1--29} (\bibinfo{year}{2020}).

\bibitem{malakiyehBatheTime2019}
\bibinfo{author}{Malakiyeh, M.~M.}, \bibinfo{author}{Shojaee, S.} \&
  \bibinfo{author}{Bathe, K.-J.}
\newblock \bibinfo{title}{The {{Bathe}} time integration method revisited for
  prescribing desired numerical dissipation}.
\newblock \emph{\bibinfo{journal}{Computers \& Structures}}
  \textbf{\bibinfo{volume}{212}}, \bibinfo{pages}{289--298}
  (\bibinfo{year}{2019}).

\bibitem{liNovelFamily2020}
\bibinfo{author}{Li, J.} \& \bibinfo{author}{Yu, K.}
\newblock \bibinfo{title}{A novel family of composite sub-step algorithms with
  desired numerical dissipations for structural dynamics}.
\newblock \emph{\bibinfo{journal}{Archive of Applied Mechanics}}
  \textbf{\bibinfo{volume}{90}}, \bibinfo{pages}{737--772}
  (\bibinfo{year}{2020}).

\bibitem{fungExtrapolatedGalerkin1996}
\bibinfo{author}{Fung, T.~C.}, \bibinfo{author}{Fan, S.~C.} \&
  \bibinfo{author}{Sheng, G.}
\newblock \bibinfo{title}{Extrapolated {{Galerkin}} time finite elements}.
\newblock \emph{\bibinfo{journal}{Computational Mechanics}}
  \textbf{\bibinfo{volume}{17}}, \bibinfo{pages}{398--405}
  (\bibinfo{year}{1996}).

\bibitem{tarnowHowRender1994}
\bibinfo{author}{Tarnow, N.} \& \bibinfo{author}{Simo, J.~C.}
\newblock \bibinfo{title}{How to render second order accurate time-stepping
  algorithms fourth order accurate while retaining the stability and
  conservation properties}.
\newblock \emph{\bibinfo{journal}{Computer Methods in Applied Mechanics and
  Engineering}} \textbf{\bibinfo{volume}{115}}, \bibinfo{pages}{233--252}
  (\bibinfo{year}{1994}).

\bibitem{kimEffectiveHigherOrder2017}
\bibinfo{author}{Kim, W.} \& \bibinfo{author}{Reddy, J.~N.}
\newblock \bibinfo{title}{Effective higher-order time integration algorithms
  for the analysis of linear structural dynamics}.
\newblock \emph{\bibinfo{journal}{Journal of Applied Mechanics}}
  \textbf{\bibinfo{volume}{84}}, \bibinfo{pages}{071009}
  (\bibinfo{year}{2017}).

\bibitem{fungUnconditionallyStable1997}
\bibinfo{author}{Fung, T.~C.}
\newblock \bibinfo{title}{Unconditionally stable higher-order {{Newmark}}
  methods by sub-stepping procedure}.
\newblock \emph{\bibinfo{journal}{Computer Methods in Applied Mechanics and
  Engineering}} \textbf{\bibinfo{volume}{147}}, \bibinfo{pages}{61--84}
  (\bibinfo{year}{1997}).

\bibitem{fanComprehensiveUnified1997a}
\bibinfo{author}{Fan, S.~C.}, \bibinfo{author}{Fung, T.~C.} \&
  \bibinfo{author}{Sheng, G.}
\newblock \bibinfo{title}{A comprehensive unified set of single-step algorithms
  with controllable dissipation for dynamics {Part I.} {Formulation}}.
\newblock \emph{\bibinfo{journal}{Computer Methods in Applied Mechanics and
  Engineering}} \textbf{\bibinfo{volume}{145}}, \bibinfo{pages}{87--98}
  (\bibinfo{year}{1997}).

\bibitem{fungComplextimestepNewmark1998}
\bibinfo{author}{Fung, T.~C.}
\newblock \bibinfo{title}{Complex-time-step {{Newmark}} methods with
  controllable numerical dissipation}.
\newblock \emph{\bibinfo{journal}{International Journal for Numerical Methods
  in Engineering}} \textbf{\bibinfo{volume}{41}}, \bibinfo{pages}{65--93}
  (\bibinfo{year}{1998}).

\bibitem{mancusoEfficientIntegration2003}
\bibinfo{author}{Mancuso, M.} \& \bibinfo{author}{Ubertini, F.}
\newblock \bibinfo{title}{An efficient integration procedure for linear
  dynamics based on a time discontinuous {{Galerkin}} formulation}.
\newblock \emph{\bibinfo{journal}{Computational Mechanics}}
  \textbf{\bibinfo{volume}{32}}, \bibinfo{pages}{154--168}
  (\bibinfo{year}{2003}).

\bibitem{krenkConservativeFourthorder2015}
\bibinfo{author}{Krenk, S.}
\newblock \bibinfo{title}{Conservative fourth-order time integration of
  non-linear dynamic systems}.
\newblock \emph{\bibinfo{journal}{Computer Methods in Applied Mechanics and
  Engineering}} \textbf{\bibinfo{volume}{295}}, \bibinfo{pages}{39--55}
  (\bibinfo{year}{2015}).

\bibitem{zhangOptimizationNsubstep2020}
\bibinfo{author}{Zhang, H.}, \bibinfo{author}{Zhang, R.},
  \bibinfo{author}{Xing, Y.} \& \bibinfo{author}{Masarati, P.}
\newblock \bibinfo{title}{On the optimization of $n$-sub-step composite time
  integration methods}.
\newblock \emph{\bibinfo{journal}{Nonlinear Dynamics}}
  \textbf{\bibinfo{volume}{102}}, \bibinfo{pages}{1939--1962}
  (\bibinfo{year}{2020}).

\bibitem{rezaiee-pajandImplicitHigherorder2008}
\bibinfo{author}{{Rezaiee-Pajand}, M.} \& \bibinfo{author}{Alamatian, J.}
\newblock \bibinfo{title}{Implicit higher-order accuracy method for numerical
  integration in dynamic analysis}.
\newblock \emph{\bibinfo{journal}{Journal of Structural Engineering}}
  \textbf{\bibinfo{volume}{134}}, \bibinfo{pages}{973--985}
  (\bibinfo{year}{2008}).

\bibitem{rezaiee-pajandHighlyAccurate2018}
\bibinfo{author}{{Rezaiee-Pajand}, M.}, \bibinfo{author}{Esfehani, S. A.~H.} \&
  \bibinfo{author}{{Karimi-Rad}, M.}
\newblock \bibinfo{title}{Highly accurate family of time integration method}.
\newblock \emph{\bibinfo{journal}{Structural Engineering and Mechanics}}
  \textbf{\bibinfo{volume}{67}}, \bibinfo{pages}{603--616}
  (\bibinfo{year}{2018}).

\bibitem{liDirectlySelfstarting2022}
\bibinfo{author}{Li, J.}, \bibinfo{author}{Zhao, R.}, \bibinfo{author}{Yu, K.}
  \& \bibinfo{author}{Li, X.}
\newblock \bibinfo{title}{Directly self-starting higher-order implicit
  integration algorithms with flexible dissipation control for structural
  dynamics}.
\newblock \emph{\bibinfo{journal}{Computer Methods in Applied Mechanics and
  Engineering}} \textbf{\bibinfo{volume}{389}}, \bibinfo{pages}{114274}
  (\bibinfo{year}{2022}).

\bibitem{kimNewFamily2017}
\bibinfo{author}{Kim, W.} \& \bibinfo{author}{Reddy, J.~N.}
\newblock \bibinfo{title}{A new family of higher-order time integration
  algorithms for the analysis of structural dynamics}.
\newblock \emph{\bibinfo{journal}{Journal of Applied Mechanics}}
  \textbf{\bibinfo{volume}{84}}, \bibinfo{pages}{071008--17}
  (\bibinfo{year}{2017}).

\bibitem{rezaiee-pajandNumericalTime2008}
\bibinfo{author}{{Rezaiee-Pajand}, M.} \& \bibinfo{author}{Alamatian, J.}
\newblock \bibinfo{title}{Numerical time integration for dynamic analysis using
  a new higher order predictor-corrector method}.
\newblock \emph{\bibinfo{journal}{Engineering Computations}}
  \textbf{\bibinfo{volume}{25}}, \bibinfo{pages}{541--568}
  (\bibinfo{year}{2008}).

\bibitem{wangOverviewHighOrder2021}
\bibinfo{author}{Wang, Y.}, \bibinfo{author}{Tamma, K.},
  \bibinfo{author}{Maxam, D.}, \bibinfo{author}{Xue, T.} \&
  \bibinfo{author}{Qin, G.}
\newblock \bibinfo{title}{An overview of high-order implicit algorithms for
  first-/second-order systems and novel explicit algorithm designs for
  first-order system representations}.
\newblock \emph{\bibinfo{journal}{Archives of Computational Methods in
  Engineering}} \textbf{\bibinfo{volume}{28}}, \bibinfo{pages}{3593--3619}
  (\bibinfo{year}{2021}).

\bibitem{fungWeightingParameters1999a}
\bibinfo{author}{Fung, T.~C.}
\newblock \bibinfo{title}{Weighting parameters for unconditionally stable
  higher-order accurate time step integration algorithms. {{Part}} 2 ---
  {Second}-order equations}.
\newblock \emph{\bibinfo{journal}{International Journal for Numerical Methods
  in Engineering}} \textbf{\bibinfo{volume}{45}}, \bibinfo{pages}{971--1006}
  (\bibinfo{year}{1999}).

\bibitem{idesmanNewHighorder2007}
\bibinfo{author}{Idesman, A.~V.}
\newblock \bibinfo{title}{A new high-order accurate continuous {{Galerkin}}
  method for linear elastodynamics problems}.
\newblock \emph{\bibinfo{journal}{Computational Mechanics}}
  \textbf{\bibinfo{volume}{40}}, \bibinfo{pages}{261--279}
  (\bibinfo{year}{2007}).

\bibitem{argyrisDynamicResponse1973}
\bibinfo{author}{Argyris, J.~H.}, \bibinfo{author}{Dunne, P.~C.} \&
  \bibinfo{author}{Angelopoulos, T.}
\newblock \bibinfo{title}{Dynamic response by large step integration}.
\newblock \emph{\bibinfo{journal}{Earthquake Engineering \& Structural
  Dynamics}} \textbf{\bibinfo{volume}{2}}, \bibinfo{pages}{185--203}
  (\bibinfo{year}{1973}).

\bibitem{liStructuralDynamic1996}
\bibinfo{author}{Li, X.} \& \bibinfo{author}{Wiberg, N.}
\newblock \bibinfo{title}{Structural dynamic analysis by a time-discontinuous
  {{Galerkin}} finite element method}.
\newblock \emph{\bibinfo{journal}{International Journal for Numerical Methods
  in Engineering}} \textbf{\bibinfo{volume}{39}}, \bibinfo{pages}{2131--2152}
  (\bibinfo{year}{1996}).

\bibitem{defrutosEasilyImplementable1992}
\bibinfo{author}{{de Frutos}, J.} \& \bibinfo{author}{{Sanz-Serna}, J.~M.}
\newblock \bibinfo{title}{An easily implementable fourth-order method for the
  time integration of wave problems}.
\newblock \emph{\bibinfo{journal}{Journal of Computational Physics}}
  \textbf{\bibinfo{volume}{103}}, \bibinfo{pages}{160--168}
  (\bibinfo{year}{1992}).

\bibitem{songHighorderImplicit2022}
\bibinfo{author}{Song, C.}, \bibinfo{author}{Eisentr{\"a}ger, S.} \&
  \bibinfo{author}{Zhang, X.}
\newblock \bibinfo{title}{High-order implicit time integration scheme based on
  {{Pad\'e}} expansions}.
\newblock \emph{\bibinfo{journal}{Computer Methods in Applied Mechanics and
  Engineering}} \textbf{\bibinfo{volume}{390}}, \bibinfo{pages}{114436}
  (\bibinfo{year}{2022}).

\bibitem{soaresStraightforwardHighorder2020}
\bibinfo{author}{Soares, D.}
\newblock \bibinfo{title}{A straightforward high-order accurate time-marching
  procedure for dynamic analyses}.
\newblock \emph{\bibinfo{journal}{Engineering with Computers}}
  \textbf{\bibinfo{volume}{38}}, \bibinfo{pages}{1659--1677}
  (\bibinfo{year}{2022}).

\bibitem{kwonSelectingLoad2021}
\bibinfo{author}{Kwon, S.-B.}, \bibinfo{author}{Bathe, K.-J.} \&
  \bibinfo{author}{Noh, G.}
\newblock \bibinfo{title}{Selecting the load at the intermediate time point of
  the $\rho_\infty$-{{Bathe}} time integration scheme}.
\newblock \emph{\bibinfo{journal}{Computers \& Structures}}
  \textbf{\bibinfo{volume}{254}}, \bibinfo{pages}{106559}
  (\bibinfo{year}{2021}).

\bibitem{choiTimeSplitting2022}
\bibinfo{author}{Choi, B.}, \bibinfo{author}{Bathe, K.-J.} \&
  \bibinfo{author}{Noh, G.}
\newblock \bibinfo{title}{Time splitting ratio in the $\rho_\infty$-{Bathe}
  time integration method for higher-order accuracy in structural dynamics and
  heat transfer}.
\newblock \emph{\bibinfo{journal}{Computers \& Structures}}
  \textbf{\bibinfo{volume}{270}}, \bibinfo{pages}{106814}
  (\bibinfo{year}{2022}).

\bibitem{hilberCollocationDissipation1978}
\bibinfo{author}{Hilber, H.~M.} \& \bibinfo{author}{Hughes, T. J.~R.}
\newblock \bibinfo{title}{Collocation, dissipation and `overshoot' for time
  integration schemes in structural dynamics}.
\newblock \emph{\bibinfo{journal}{Earthquake Engineering \& Structural
  Dynamics}} \textbf{\bibinfo{volume}{6}}, \bibinfo{pages}{99--117}
  (\bibinfo{year}{1978}).

\bibitem{liSecondorderSsubstep2023}
\bibinfo{author}{Li, J.}, \bibinfo{author}{Li, H.}, \bibinfo{author}{Zhao, R.}
  \& \bibinfo{author}{Yu, K.}
\newblock \bibinfo{title}{On second-order $s$-sub-step explicit algorithms with
  controllable dissipation and adjustable bifurcation point for second-order
  hyperbolic problems}.
\newblock \emph{\bibinfo{journal}{European Journal of Mechanics - A/Solids}}
  \textbf{\bibinfo{volume}{97}}, \bibinfo{pages}{104829}
  (\bibinfo{year}{2023}).

\bibitem{liSuiteSecondorder2023}
\bibinfo{author}{Li, J.}, \bibinfo{author}{Li, H.}, \bibinfo{author}{Lian, Y.},
  \bibinfo{author}{Zhao, R.} \& \bibinfo{author}{Yu, K.}
\newblock \bibinfo{title}{A suite of second-order composite sub-step explicit
  algorithms with controllable numerical dissipation and maximal stability
  bounds}.
\newblock \emph{\bibinfo{journal}{Applied Mathematical Modelling}}
  \textbf{\bibinfo{volume}{114}}, \bibinfo{pages}{601--626}
  (\bibinfo{year}{2023}).

\bibitem{jiUnconditionallyStable2021}
\bibinfo{author}{Ji, Y.}, \bibinfo{author}{Xing, Y.} \&
  \bibinfo{author}{Wiercigroch, M.}
\newblock \bibinfo{title}{An unconditionally stable time integration method
  with controllable dissipation for second-order nonlinear dynamics}.
\newblock \emph{\bibinfo{journal}{Nonlinear Dynamics}}  (\bibinfo{year}{2021}).

\bibitem{butcherNumericalMethods2016}
\bibinfo{author}{Butcher, J.~C.}
\newblock \emph{\bibinfo{title}{Numerical {{Methods}} for {{Ordinary
  Differential Equations}}}} \bibinfo{edition}{Third} edn
  (\bibinfo{publisher}{{Wiley}}, \bibinfo{year}{2016}).

\bibitem{hulbertErrorAnalysis1987}
\bibinfo{author}{Hulbert, G.~M.} \& \bibinfo{author}{Hughes, T. J.~R.}
\newblock \bibinfo{title}{An error analysis of truncated starting conditions in
  step-by-step time integration: {{Consequences}} for structural dynamics}.
\newblock \emph{\bibinfo{journal}{Earthquake Engineering \& Structural
  Dynamics}} \textbf{\bibinfo{volume}{15}}, \bibinfo{pages}{901--910}
  (\bibinfo{year}{1987}).

\bibitem{rezaiee-pajandImprovingStability2011}
\bibinfo{author}{{Rezaiee-Pajand}, M.}, \bibinfo{author}{Sarafrazi, S.~R.} \&
  \bibinfo{author}{Hashemian, M.}
\newblock \bibinfo{title}{Improving stability domains of the implicit higher
  order accuracy method}.
\newblock \emph{\bibinfo{journal}{International Journal for Numerical Methods
  in Engineering}} \textbf{\bibinfo{volume}{88}}, \bibinfo{pages}{880--896}
  (\bibinfo{year}{2011}).

\bibitem{rezaiee-pajandNovelTime2017}
\bibinfo{author}{{Rezaiee-Pajand}, M.}, \bibinfo{author}{Hashemian, M.} \&
  \bibinfo{author}{Bohluly, A.}
\newblock \bibinfo{title}{A novel time integration formulation for nonlinear
  dynamic analysis}.
\newblock \emph{\bibinfo{journal}{Aerospace Science \& Technology}}
  \textbf{\bibinfo{volume}{69}}, \bibinfo{pages}{625--635}
  (\bibinfo{year}{2017}).

\bibitem{rezaiee-pajandEfficientWeighted2021}
\bibinfo{author}{{Rezaiee-Pajand}, M.}, \bibinfo{author}{Esfehani, S. A.~H.} \&
  \bibinfo{author}{Ehsanmanesh, H.}
\newblock \bibinfo{title}{An efficient weighted residual time integration
  family}.
\newblock \emph{\bibinfo{journal}{International Journal of Structural Stability
  and Dynamics}} \textbf{\bibinfo{volume}{21}}, \bibinfo{pages}{2150106}
  (\bibinfo{year}{2021}).

\bibitem{batheInsightImplicit2012}
\bibinfo{author}{Bathe, K.-J.} \& \bibinfo{author}{Noh, G.}
\newblock \bibinfo{title}{Insight into an implicit time integration scheme for
  structural dynamics}.
\newblock \emph{\bibinfo{journal}{Computers \& Structures}}
  \textbf{\bibinfo{volume}{98--99}}, \bibinfo{pages}{1--6}
  (\bibinfo{year}{2012}).

\bibitem{rezaiee-pajandModifiedDifferential2017}
\bibinfo{author}{{Rezaiee-Pajand}, M.} \& \bibinfo{author}{Hashemian, M.}
\newblock \bibinfo{title}{Modified differential transformation method for
  solving nonlinear dynamic problems}.
\newblock \emph{\bibinfo{journal}{Applied Mathematical Modelling}}
  \textbf{\bibinfo{volume}{47}}, \bibinfo{pages}{76--95}
  (\bibinfo{year}{2017}).

\bibitem{rezaiee-pajandFamilySecondorder2018}
\bibinfo{author}{{Rezaiee-Pajand}, M.} \& \bibinfo{author}{{Karimi-Rad}, M.}
\newblock \bibinfo{title}{A family of second-order fully explicit time
  integration schemes}.
\newblock \emph{\bibinfo{journal}{Computational \& Applied Mathematics}}
  \textbf{\bibinfo{volume}{37}}, \bibinfo{pages}{3431--3454}
  (\bibinfo{year}{2018}).

\bibitem{heNonlinearAeroelastic2020}
\bibinfo{author}{He, H.} \emph{et~al.}
\newblock \bibinfo{title}{Nonlinear aeroelastic analysis of the folding fin
  with freeplay under thermal environment}.
\newblock \emph{\bibinfo{journal}{Chinese Journal of Aeronautics}}
  \textbf{\bibinfo{volume}{33}}, \bibinfo{pages}{2357--2371}
  (\bibinfo{year}{2020}).

\bibitem{fungUnconditionallyStable1996}
\bibinfo{author}{Fung, T.~C.}
\newblock \bibinfo{title}{Unconditionally stable higher-order accurate
  {Hermitian} time finite elements}.
\newblock \emph{\bibinfo{journal}{International Journal for Numerical Methods
  in Engineering}} \textbf{\bibinfo{volume}{39}}, \bibinfo{pages}{3475--3495}
  (\bibinfo{year}{1996}).

\end{thebibliography}
}% common bib file

\end{document}